\documentclass[a4paper,11pt,leqno,twoside]{article}
\usepackage{amssymb}
\usepackage[mathscr]{eucal}

\makeatletter
\renewcommand{\subsubsection}{\@startsection{subsubsection}{3}{0mm}{-\baselineskip}{-0.01\baselineskip}{\bf}}
\makeatother

\def\ali{\subsubsection{}}
\def\alin#1{\ali{\sl #1}\ : }

\def\ini{\setcounter{equation}{\value{subsubsection}}\addtocounter{subsubsection}{1}}

\newtheorem{theo}[subsubsection]{Th{\'e}or{\`e}me}
\newtheorem{lemme}[subsubsection]{Lemme}
\newtheorem{prop}[subsubsection]{Proposition}
\newtheorem{coro}[subsubsection]{Corollaire}

\newtheorem{fact}[subsubsection]{Fait}

\newtheorem{conj}[subsubsection]{Conjecture}
\newtheorem{DEf}[subsubsection]{D{\'e}finition}

\newtheorem{rema}[subsubsection]{Remarque}  
\newtheorem{nota}[subsubsection]{Notation}  

\newtheorem{conjintro}{Conjecture}

\newtheorem{theointro}[conjintro]{Th{\'e}or{\`e}me}

\newtheorem{theosec}[subsection]{Th{\'e}or{\`e}me}

\newtheorem{propsec}[subsection]{Proposition}

\newenvironment{boite}[1]
               {\smallskip \noindent {\bf #1 } \rm }
               {\normalsize\smallskip}

\newcommand{\findem}{\hfill$\Box$\par\medskip}
\newcommand{\dem}{\noindent {\sl Preuve :} \rm }

\newenvironment{proof}{\dem}{\findem}

\usepackage[french]{babel}
\usepackage[arrow,matrix]{xy} 
\usepackage[latin1]{inputenc}

\title{Correspondance de Langlands locale et monodromie des espaces de Drinfeld}
\author{J.-F. Dat}

\def\NM{{\mathbb{N}}}
\def\EM{{\mathbb{E}}}
\def\MM{{\mathbb{M}}}
\def\UM{{\mathbb{U}}}

\def\GM{{\mathbb{G}}}
\def\PM{{\mathbb{P}}}
\def\RM{{\mathbb{R}}}

\def\QM{{\mathbb{Q}}}
\def\FM{{\mathbb{F}}}
\def\ZM{{\mathbb{Z}}}
\def\CM{{\mathbb{C}}}
\def\AM{{\mathbb{A}}}
\def\XM{{\mathbb{X}}}
\def\MM{{\mathbb{M}}}
\def\MG{{\mathfrak M}}
\def\BG{{\mathfrak B}}

\def\mG{{\mathfrak m}}

\def\ZG{{\mathfrak Z}}

\def\XG{{\mathfrak X}}

\def\AC{{\mathcal A}}
\def\XC{{\mathcal X}}
\def\JC{{\mathcal J}}

\def\CC{{\mathcal C}}
\def\HC{{\mathcal H}}
\def\NC{{\mathcal N}}
\def\ZC{{\mathcal Z}}

\def\IC{{\mathcal I}}
\def\OC{{\mathcal O}}
\def\MC{{\mathcal M}}

\def\PC{{\mathcal P}}

\def\LC{{\mathcal L}}

\def\FC{{\mathcal F}}
\def\VC{{\mathcal V}}
\def\GC{{\mathcal G}}
\def\UC{{\mathcal U}}
\def\BC{{\mathcal B}}
\def\TC{{\mathcal T}}
\def\DC{{\mathcal D}}

\def\ssi{si et seulement si}
\def\para{sous-groupe parabolique }
\def\paras{sous-groupes paraboliques }
\def\levi{sous-groupe de Levi }
\def\levis{sous-groupes de Levi }
\def\borel{sous-groupe de Borel }
\def\borels{sous-groupes de Borel }

\def\subsetneq{\varsubsetneq}

\def\simto{\buildrel\hbox{$\sim$}\over\longrightarrow}

\def\leq{\leqslant}
\def\geq{\geqslant}
\def\injo{\hookrightarrow}
\def\id{\mathop{\mathrm{Id}}\nolimits}
\def\equ{\Leftrightarrow}

\def\ba{\backslash}
\def\wt{\widetilde}
\def\wh{\widehat}
\def\o#1{\overline{#1}}

\def\application#1#2#3#4#5{\begin{array}{rcl}
                            #1 \;\;\; #2 & \to &  #3 \\
                              #4 & \mapsto & #5 
                            \end{array}} 
\def\cas#1#2#3#4#5{\begin{array}{rcl} #1 \; & = &
    \left\{\begin{array}{rcl} #2 & \hbox{ si } & #3 \\
                             #4 & \hbox{ si } & #5 \end{array}
                         \right. \end{array}}

\def\To#1{\buildrel\hbox{\tiny{$#1$}}\over\longrightarrow}
\def\to{\rightarrow}

\def\ker{\mathop{\hbox{\sl ker}\,}}
\def\coker{\mathop{\hbox{\sl coker}\,}}
\def\im{\mathop{\hbox{\sl im}\,}}

\def\hom#1#2#3{\hbox{\sl Hom}_{#3}\>\!\left(#1,#2\right)}
\def\endo#1#2{\hbox{\sl End}_{#1}\>\!\left(#2\right)}
\def\aut#1#2{\hbox{\sl Aut}_{#1}(#2)}

\def\Mo#1#2{\mathop{\hbox {\sl Mod}_{#1}(#2)}}
\def\Irr#1#2{\mathop{\hbox {\sl Irr}_{#1}(#2)}}
\def\Cu#1#2{\mathop{\hbox {\sl Cusp}_{#1}(#2)}}
\def\Disc#1#2{\mathop{\hbox {\sl Disc}_{#1}(#2)}}
\def\Ell#1#2{\mathop{\hbox {\sl Ell}_{#1}(#2)}}

\def\Ind#1#2{\hbox {\sl Ind}_{#1}^{#2}}
\def\ind#1#2#3{\hbox {\sl Ind}_{#1}^{#2}\>\!\left(#3\right)}  
\def\cInd#1#2{\hbox {\sl ind}_{#1}^{#2}}
\def\cind#1#2#3{\hbox {\sl ind}_{#1}^{#2}\>\!\left(#3\right)} 
\def\ip#1#2#3{\hbox {\sl i}_{#1}^{#2}\>\!(#3)}  
\def\Ip#1#2{\hbox {\sl i}_{#1}^{#2}}

\def\res#1#2#3{\hbox {\sl Res}_{#1}^{#2}\>\!\left(#3\right)}
\def\Res#1#2{\hbox {\sl Res}_{#1}^{#2}}              
\def\rp#1#2#3{\hbox {\sl r}_{#1}^{#2}\>\!\left(#3\right)}    
\def\Rp#1#2{\hbox {\sl r}_{#1}^{#2}}

\def\dim{\mathop{\mbox{\sl dim}}\nolimits}

\def\val{\mathop{\mbox{\sl val}}\nolimits}

\def\gal{\mathop{\mbox{\sl Gal}}\nolimits}

\def \ext#1#2#3#4{\mbox{\sl Ext}^{#1}_{#4}\>\!\left(#2,#3\right)}
\def \tor#1#2#3#4{\mbox{\sl Tor}^{#4}_{#1}\>\!(#2,#3)}

\def \limi#1{\lim\limits_{\displaystyle\longrightarrow\atop {#1}}}
\def \limp#1{\lim\limits_{\displaystyle\longleftarrow\atop {#1}}}
\def \limproj{{\lim\limits_{\longleftarrow}}}

\def\exp{\mathop{\hbox{\sl Exp}\,}}

\begin{document}
\maketitle
\bibliographystyle{plain}

\def\la{\langle}
\def\ra{\rangle}
\def\knr{{\wh{K^{nr}}}}
\def\ka{\wh{K^{ca}}}


\abstract{Les travaux de Harris, Boyer, Harris-Taylor et Hausberger
  ont concr{\'e}tis{\'e} les conjectures de Carayol dans [{\em Non-abelian Lubin-Tate theory.} Automorphic forms, Shimura varieties and $L$-functions, vol II: 15--39, Academic Press,1990] selon lesquelles la cohomologie {\'e}tale $l$-adique de certains espaces
  de modules de groupes formels introduits par Drinfeld permettrait de r{\'e}aliser simultan{\'e}ment les
  correspondances locales de Langlands et Jacquet-Langlands pour les
  repr{\'e}sentations {\em supercuspidales} de $GL_n$ sur un corps local
  non-archim{\'e}dien. 
  Mais selon des r\'esultats annonc\'es r\'ecemment par Boyer et Faltings permettant le calcul des    groupes de cohomologie individuels de ces espaces, la m\^eme strat\'egie ne fonctionne plus  
  pour les repr{\'e}sentations non supercuspidales, m{\^e}me
  pour les s{\'e}ries discr{\`e}tes : la cohomologie ne suffit visiblement plus. Notre but ici est de montrer comment, \`a partir de ce calcul et en  utilisant le formalisme des cat{\'e}gories d{\'e}riv{\'e}es, on peut tout-de-m\^eme obtenir une
  r{\'e}alisation simultan\'ee des deux correspondances pour toute repr{\'e}sentation {\em
    elliptique} de $GL_n$, g\'en\'eralisant ainsi la conjecture de Carayol.}

\def\mdro{\MC_{Dr,0}}
\def\mdrn{\MC_{Dr,n}}
\def\mdr{\MC_{Dr}}
\def\mlto{\MC_{LT,0}}
\def\mltn{\MC_{LT,n}}
\def\mlt{\MC_{LT}}
\def\mltK{\MC_{LT,K}}

\tableofcontents

\section{Introduction}

\def\dd{D_d^\times}

 Soit $K$ un corps local de caract{\'e}ristique r{\'e}siduelle $p$, et $d$ un entier $\geq 1$.
On s'int{\'e}resse aux trois groupes suivants et {\`a} leurs repr{\'e}sentations :
\begin{itemize}
\item le groupe de Weil $W_K$, topologis{\'e} de la mani{\`e}re habituelle,
\item le groupe $G_d:=GL_{d}(K)$, muni de sa topologie localement
  $p$-profinie,
\item le groupe $D_d^\times$  des
{\'e}l{\'e}ments inversibles de l'alg{\`e}bre centrale simple sur $K$ d'invariant
$\frac{1}{d}$, lui aussi localement $p$-profini.
\end{itemize}
Les repr{\'e}sentations de ces trois groupes sont reli{\'e}es par des
{\em correspondances} qui font partie du vaste programme de
g{\'e}n{\'e}ralisation non-ab{\'e}lienne de la th{\'e}orie du corps de classes propos{\'e}
par Langlands. 
Ces correspondances sont en g{\'e}n{\'e}ral exprim{\'e}es en termes de
repr{\'e}sentations {\em complexes}, mais comme on s'int{\'e}resse ici {\`a} leur
r{\'e}alisation cohomologique, on consid{\`e}rera plut{\^o}t leur variante
$l$-adique. Nous fixerons donc dor{\'e}navant un nombre premier $l\neq p$
et ne consid{\`e}rerons dans cette introduction que des $\o\QM_l$-repr{\'e}sentations.
La correspondance de Jacquet-Langlands \cite{DKV}
\cite{Badu1}, que nous rappelons
en \ref{JacLan},
est une bijection entre les classes de repr{\'e}sentations lisses ``de la s{\'e}rie
discr{\`e}te''  de $G_d$ et celles de $\dd$. Nous la noterons $\pi \mapsto JL_d(\pi)$.
La correspondance de
Langlands locale \cite{LRS} \cite{HaTay} \cite{HeLang} \cite{HeSMF}
que nous rappelons en \ref{corLan} induit en particulier une
bijection entre classes de repr{\'e}sentations 
{\em lisses} irr{\'e}ductibles de $G_d$ et classes de repr{\'e}sentations 
{\em continues} de dimension $d$ de $W_K$. Nous la noterons $\pi \mapsto \sigma_d(\pi)$.

Ces deux correspondances ont des vies s{\'e}par{\'e}es et  chacune est
caract{\'e}ris{\'e}e par des propri{\'e}t{\'e}s de pr{\'e}servation d'invariants de
nature arithm{\'e}tico-analytique. Mais un des aspects les plus fascinants
est l'existence d'objets g{\'e}om{\'e}triques -- alg{\'e}briques sur un corps
global ou analytiques sur un corps local -- 
dont la cohomologie $l$-adique permet
de ``r{\'e}aliser'', ou d'``incarner'', certains cas de ces correspondances.
Notons $\ka$ la compl{\'e}tion d'une cl{\^o}ture alg{\'e}brique $K^{ca}$ de $K$
pour l'unique extension de la norme $|.|_K$.
Dans ce texte, on s'int{\'e}resse {\`a} deux tours d'espaces $\ka$-analytiques
d{\'e}finies par Drinfeld.
\begin{itemize}
\item la tour $(\MC^{d/K}_{Dr,n})_{n\in\NM}$ au-dessus du $K$-espace sym{\'e}trique
de Drinfeld de dimension $d-1$, d{\'e}finie dans \cite{Drincov} et
rappel{\'e}e en \ref{tourDr},
\item  la tour $(\MC^{d/K}_{LT,n})_{n\in\NM}$ 
``de Lubin-Tate'' au-dessus de la boule ouverte de dimension $d-1$ sur $K$,
d{\'e}finie dans \cite{Drinell} et rappel{\'e}e en \ref{tourLT}.
\end{itemize}
Chacun des pro-$\ka$-espaces rigides $\MC^{d/K}_{Dr}$ et $\MC^{d/K}_{LT}$ obtenus en
passant {\`a} la limite est muni
d'une action de $G_d\times \dd\times W_K$, donc leur ``cohomologie {\'e}tale
$l$-adique {\`a} supports compacts'' fournit des repr{\'e}sentations de ce produit triple. La
d{\'e}finition habituelle de celle-ci est simplement, pour $?=Dr$ ou $LT$, 
$$H^i_c(\MC^{d/K}_?,\o\QM_l):= \limi{n} H^i_c(\MC^{d/K}_{?,n},\o\QM_l),$$
compte
tenu de la finitude des morphismes de transition.
On sait par des r{\'e}sultats g{\'e}n{\'e}raux de Berkovich que l'action de
$G_d\times \dd$ sur ces espaces est lisse. 
Les quatre articles \cite{HaCusp} \cite{Boyer1} \cite{HaTay}
\cite{Hausb} {\'e}tudient la partie {\em supercuspidale} de ces espaces de cohomologie, pour chacun 
des espaces $\MC^{d/K}_{Dr}$ ou $\MC^{d/K}_{LT}$ en caract{\'e}ristique nulle ou en {\'e}gales
caract{\'e}ristiques. Ils d\'emontrent en particulier la fameuse conjecture 
que Carayol a \'enonc\'ee dans \cite{CarAA} \`a la suite de travaux pionniers de Deligne et Drinfeld pour $d=2$ :
 soit $\pi$ une
repr{\'e}sentation supercuspidale $l$-adique de $G_d$, de caract{\`e}re central
d'ordre fini ; alors on a, en normalisant convenablement les actions, {\em cf} \cite{HaICM}, 
$$\begin{array}{rcl} \hom{H_c^i(\MC^{d/K}_{Dr},\o\QM_l)}{\pi}{G_d} \; &
  \mathop{\simeq}\limits_{\dd\times W_K} &
\hom{H_c^i(\MC^{d/K}_{LT},\o\QM_l)}{\pi}{G_d} \\ &
  \mathop{\simeq}\limits_{\dd\times W_K} &
    \left\{\begin{array}{rcl} JL_d(\pi)\otimes
  \sigma_d(\pi)(?) & \hbox{ si } & i=d-1 \\
                             0 & \hbox{ si } & i\neq d-1 \end{array}
                         \right. \end{array}$$ 
o{\`u} $(?)$ d{\'e}signe une torsion {\`a} la Tate. \footnote{Pour obtenir
  vraiment l'isomorphisme ci-dessus dans le cas de $\MC^{d/K}_{LT}$, il 
convient de conjuguer par l'automorphisme ext{\'e}rieur $g\mapsto
{^t}g^{-1}$  de $G_d$ l'action d{\'e}finie par Carayol dans \cite{CarAA}.}
Les m{\'e}thodes employ{\'e}es dans ces articles reposent sur la propri{\'e}t{\'e}
d'uniformisation $p$-adique de certaines vari{\'e}t{\'e}s globales (de Shimura
ou de Drinfeld) par ces
espaces et les suites spectrales de type Hochschild-Serre ou de cycles
{\'e}vanescents
associ{\'e}es. Elles ne permettent  en g{\'e}n{\'e}ral que d'obtenir des
informations sur la somme altern{\'e}e des $H^i_c$. 

Pendant longtemps,  seul le calcul de Schneider-Stuhler dans \cite{SS1} pour l'espace sym\'etrique de Drinfeld fournissait des renseignements sur la partie non-cuspidale de la cohomologie. Ce calcul a inspir\'e \`a Harris une conjecture (non publi\'ee) sur la forme explicite des groupes de  cohomologie individuels du c\^ot\'e $Dr$ sur un corps $p$-adique. Tout  r\'ecemment, Boyer \cite{Boyer2} a annonc\'e une preuve de cette conjecture... mais du c\^ot\'e $LT$ et en \'egales caract\'eristiques. Un peu auparavant, Faltings avait annonc\'e dans \cite{FaltDrin} que les  cohomologies des deux tours devaient \^etre isomorphes.
Comme ses  arguments sont tr\`es r\'esum\'es et valables seulement sur un corps $p$-adique, ils seront compl\'et\'es et \'etendus en  \'egales caract\'eristiques dans un travail en cours de Fargues, Genestier et Lafforgue.


 Nous
exposerons 
la description pr\'ecise des $H^i_c$ -- qui
reste donc partiellement conjecturale --  en  \ref{conjHarris}.  
Nous nous bornons
ici {\`a} souligner que les repr{\'e}sentations de $G_d$ qui apparaissent sont
toutes {\em elliptiques} au sens o{\`u} la restriction de leur
caract{\`e}re-distribution aux {\'e}l{\'e}ments semi-simples r{\'e}guliers elliptiques
est non-nulle (voir \ref{defell} pour d'autres caract{\'e}risations des
repr{\'e}sentations elliptiques de $G_d$).
Mais par ailleurs et comme le montre d{\'e}ja le calcul de
Schneider-Stuhler pour l'espace
sym{\'e}trique de Drinfeld,  seules {\em certaines}  repr{\'e}sentations
elliptiques apparaissent.
Notons que pour ces repr{\'e}sentations elliptiques -- qui ne sont pas n{\'e}cessairement
temp{\'e}r{\'e}es -- on peut  {\'e}tendre de mani{\`e}re naturelle l'{\'e}nonc{\'e} de la
correspondance de Jacquet-Langlands, {\em cf} \ref{JacLan}.

La cons{\'e}quence majeure de la description des $H^i_c$ est que ceux-ci
r{\'e}alisent bien la correspondance de Jacquet-Langlands ainsi {\'e}tendue,
mais {\em ne r{\'e}alisent
pas} la correspondance de Langlands
pour les repr{\'e}sentations non supercuspidales qui y apparaissent, m{\^e}me les s{\'e}ries
discr{\`e}tes.
Pour {\^e}tre un peu plus pr{\'e}cis, soit $\pi$  une
repr{\'e}sentation elliptique de $G_d$ ; on verra en \ref{corlanell} que la semisimplifi{\'e}e de 
$\sigma_d(\pi)$ est une somme de torsions {\`a} la Tate d'une m{\^e}me
repr{\'e}sentation $l$-adique irr{\'e}ductible de $W_K$ que nous notons ci-dessous
$\sigma'(\pi)$. Si 
 la contribution de $\pi$ {\`a} $H^i_c$ est non-nulle, alors la
 description de Harris-Boyer donne 
$$ \hom{H^i_c(\MC^{d/K},\o\QM_l)}{\pi}{G_d} \mathop{\simeq}\limits_{\dd\times W_K}
JL_d(\pi)\otimes\sigma'(\pi)(?)  $$
o{\`u} $(?)$ d{\'e}signe une certaine torsion {\`a} la Tate et $\MC^{d/K}$ d{\'e}signe
$\MC^{d/K}_{Dr}$ ou $\MC^{d/K}_{LT}$.
On voit en particulier que l'op{\'e}rateur de monodromie (celui que l'on sait associer {\`a} toute
repr{\'e}sentation $l$-adique continue de dimension finie de $W_K$ par le th{\'e}or{\`e}me
de Grothendieck) est toujours nul sur les composantes isotypiques des
$H^i_c$. Mais la description ci-dessus sugg{\`e}re aussi que cet op{\'e}rateur
est essentiellement la seule chose qui manque pour obtenir une r{\'e}alisation
de la correspondance de Langlands, au moins pour les repr{\'e}sentations
qui apparaissent dans ces $H^i_c$. L'id{\'e}e principale de ce texte est
que cet op{\'e}rateur est ``cach{\'e}'' dans ``le'' complexe de cohomologie
$R\Gamma_c(\MC^{d/K}_{Dr},\o\QM_l)$
(resp. $R\Gamma_c(\MC^{d/K}_{LT},\o\QM_l)$) vu comme objet de la cat{\'e}gorie
d{\'e}riv{\'e}e de la cat{\'e}gorie ab{\'e}lienne des $\o\QM_l(G_d\times
\dd)$-repr{\'e}sentations lisses. Oublions un instant les 
difficult{\'e}s de d{\'e}finition et d'{\'e}tude que cela pose et {\'e}non{\c c}ons
nos r{\'e}sultats principaux. Ce complexe $R\Gamma_c$ dans chacun des cas est muni d'une
action de $W_K$. Par cons{\'e}quent pour toute repr{\'e}sentation lisse $\pi$ de
$G_d$, resp. $\rho$ de $\dd$, les complexes de $\o\QM_l$-espaces vectoriels
$$ R\hom{R\Gamma_c}{\pi}{D^b(G_d)} \in D^b(\o\QM_l), \;\;\hbox{ resp.}\;\; 
R\hom{R\Gamma_c}{\pi\otimes\rho}{D^b(G_d\times \dd)} \in D^b(\o\QM_l)$$
sont munis d'une action de $\dd\times W_K$, resp. $W_K$. Rappelons que
la cat{\'e}gorie triangul{\'e}e $D^b(\o\QM_l)$ est {\'e}quivalente {\`a} la cat{\'e}gorie
triangul{\'e}e  des $\o\QM_l$-espaces vectoriels
$\ZM$-gradu{\'e}s ``{\`a} support fini'', une {\'e}quivalence {\'e}tant donn{\'e}e par le foncteur $\HC^*:
\CC^\bullet \in D^b(\o\QM_l) \mapsto \bigoplus_{i\in\ZM}
\HC^i(\CC^\bullet)$. Ainsi le $\o\QM_l$-espace vectoriel $\ZM$-gradu{\'e} 
$$\HC^*(R\hom{R\Gamma_c}{\pi}{D^b(G_d)}) \;\;\hbox{ resp.}\;\; 
\HC^*(R\hom{R\Gamma_c}{\pi\otimes\rho}{D^b(G_d\times \dd)})$$ est muni d'une action
de $\dd\times W_K$, resp. $W_K$. On oubliera g{\'e}n{\'e}ralement la
$\ZM$-graduation par la suite pour se retrouver avec des
repr{\'e}sentations au sens usuel.

Pour {\'e}noncer nos r{\'e}sultats, il convient d'introduire le groupe
$GD:=(G_d\times \dd)/\Delta$ o\`u $\Delta=\{(z,z),z\in K^\times\}$. 
On sait que, quitte {\`a} conjuguer l'action de $G_d$  initialement d{\'e}finie
par Drinfeld par un automorphisme ext{\'e}rieur convenable, l'action
g{\'e}om{\'e}trique de $G_d\times \dd$ sur les espaces 
de modules $\MC^{d/K}_{Dr}$ et $\MC^{d/K}_{LT}$ se factorise par ce groupe. De plus,
toute repr{\'e}sentation irr{\'e}ductible de $GD$ est de la forme
$\pi\otimes \rho$ o{\`u} $\pi$ et $\rho$ sont des repr{\'e}sentations
irr{\'e}ductibles de $G_d$ et $\dd$ dont les caract{\`e}res centraux sont
``inverses l'un de l'autre''.

\begin{conjintro} \label{conj}
Notons $\MC^{d/K}$ pour $\MC^{d/K}_{Dr}$ ou $\MC^{d/K}_{LT}$. Soit $\pi$, resp. $\rho$,
une $\o\QM_l$-repr{\'e}sentation 
  lisse  et irr{\'e}ductible de $G_d$, resp. $\dd$. Alors
  \begin{enumerate}
  \item Lorsque $\pi$ est elliptique, on a
$$\HC^*\left(R\hom{R\Gamma_c(\MC^{d/K},\o\QM_l)}{\pi}{D^b(G_d)}\right)
\mathop{\simeq}\limits_{\dd\times
  W_K}JL_d(\pi)\otimes\sigma_d(\pi)|-|^{\frac{d-1}{2}},$$
le terme de gauche {\'e}tant nul si $\pi$ n'est pas elliptique.

\item Lorsque  $\pi$ est elliptique et $\rho\simeq JL_d(\pi^\vee)$, on a
$$\HC^*\left(R\hom{R\Gamma_c(\MC^{d/K},\o\QM_l)}{\pi\otimes
    \rho}{D^b(G D)}\right)
\mathop{\simeq}\limits_{W_K} \sigma_d(\pi)|-|^{\frac{d-1}{2}},$$
le terme de gauche {\'e}tant nul si $\pi$ n'est pas elliptique ou
si $\rho$ n'est pas isomorphe {\`a} $JL_d(\pi^\vee)$.

   \end{enumerate}
\end{conjintro}
La notation $?^\vee$ d{\'e}signe la contragr{\'e}diente de la repr{\'e}sentation
$?$.  Si on veut r{\'e}sumer cette conjecture en quelques mots
(impr{\'e}cis) : l'introduction du formalisme des cat{\'e}gories d{\'e}riv{\'e}es
permet non seulement de r{\'e}cup{\'e}rer la repr{\'e}sentation Galoisienne en
entier (avec la monodromie), mais aussi fait ``apparaitre'' des
repr{\'e}sentations $\pi$ qui n'apparaissent pas dans les groupes de cohomologie. 
On verra en \ref{reduc1} que
les deux points de la conjecture sont essentiellement {\'e}quivalents
(simplement parce que la cat{\'e}gorie
 des $\o\QM_l(\dd/K^\times)$-modules est semi-simple), mais le
second est le plus propice {\`a} une {\'e}ventuelle g{\'e}n{\'e}ralisation {\`a} d'autres
paires de formes int{\'e}rieures agissant sur des espaces de Rapoport-Zink.
Le th\'eor\`eme principal de ce texte concerne le c\^ot\'e $Dr$ et 
s'\'enonce formellement comme suit :
\begin{theointro} \label{thprin}
 Supposons que 
  \begin{enumerate}
  \item la cohomologie $l$-adique $H^i_c(\MC^{d/K}_{Dr},\o\QM_l)$
    est bien d\'ecrite par la conjecture de Harris \ref{conjHarris}.
\item la cohomologie enti{\`e}re $H^i_c(\MC^{d/K}_{Dr},\ZM_l)$ est
  ``admissible-modulo-le-centre''. 
\item La conjecture ``monodromie-poids'' est v{\'e}rifi{\'e}e pour les
  vari{\'e}t{\'e}s uniformis{\'e}es par $\MC^{d/K}_{Dr,n}$ (c'est-\`a-dire les
  quotients de $\MC^{d/K}_{Dr,n}$ par un sous-groupe discret cocompact de $G_d$ y agissant librement). 
  \end{enumerate}
Alors, la conjecture \ref{conj} est v{\'e}rifi{\'e}e pour $\MC^{d/K}=\MC^{d/K}_{Dr}$.
\end{theointro}
Pr\'ecisons un peu le statut de ces hypoth\`eses.
Comme on l'a d\'eja \'ecrit plus haut, les 
travaux de  Faltings \cite{FaltDrin} compl\'et\'es et g\'en\'eralis\'es par ceux de Fargues, Genestier et Lafforgue montrent que $\MC^{d/K}_{Dr}$ et $\MC^{d/K}_{LT}$ ont des cohomologies \`a coefficients constants isomorphes. Ainsi, il suffit de v\'erifier les deux premi\`eres hypoth\`eses du c\^ot\'e $LT$. Pour l'hypoth\`ese ii), c'est "facile", {\em cf} \ref{admicohoent}. Pour l'hypoth\`ese i) c'est beaucoup plus difficile, mais
lorsque $K$ est d'\'egales caract\'eristiques, Boyer \cite{Boyer2} a annonc\'e une preuve de la conjecture de Harris du c\^ot\'e $LT$. Toujours lorsque $K$ est d'\'egales caract\'eristiques, 
la conjecture "monodromie-poids" est connue pour toutes les vari{\'e}t{\'e}s alg{\'e}briques lisses par les travaux de Deligne sur les conjectures de Weil, voir
\cite{Ito2}, et donc l'hypoth\`ese iii) est aussi v\'erifi\'ee. 

Lorsque $K$ est $p$-adique, les arguments de Boyer devraient s'appliquer formellement de la m\^eme mani\`ere, au moins en admettant la semi-simplicit\'e des gradu\'es du complexe des cycles \'evanescents pour la filtration de monodromie des vari\'et\'es de Shimura \'etudi\'ees par Harris et Taylor dans \cite{HaTay}, donc en particulier en admettant la version locale de la conjecture monodromie-poids pour ces vari\'et\'es.
Quoiqu'il en soit, rappelons que pour le rez-de-chauss\'ee de la tour $\MC^{d/K}_{Dr,0}$, la cohomologie \`a coefficients constants est connue  depuis  Schneider et
Stuhler \cite{SS1}, ind\'ependamment des travaux de Boyer et Faltings. 
De plus, Ito a d\'emontr\'e dans
\cite{Ito} la conjecture monodromie-poids dans
 le cas particulier des vari\'et\'es uniformis\'ees par $\MC^{d/K}_{Dr,0}$.

En r\'esum\'e on peut \'enoncer le

\begin{boite}{Corollaire}
{\em  Supposons $K$ de caract{\'e}ristique $p$. Admettant les r{\'e}sultats de
  Boyer  et Faltings, la conjecture \ref{conj} est v{\'e}rifi{\'e}e pour $\MC^{d/K}_{Dr}$.
Lorsque $K$ est de caract{\'e}ristique nulle, un analogue de la
conjecture est v{\'e}rifi{\'e} 
pour $\MC^{d/K}_{Dr,0}$ : l'{\'e}nonc{\'e} pr{\'e}cis est le th{\'e}or{\`e}me \ref{theodp}.
}\end{boite}

En fait, dans le cas $\MC^{d/K}_{Dr,0}$, on peut se passer du r\'esultat d'Ito et en donner une nouvelle d\'emonstration, par nos techniques (en utilisant aussi une  id\'ee de Mokrane dans \cite{Mokrane}), {\em cf} \ref{MPdp}.
Plus g\'en\'eralement,
nous verrons en \ref{propequMP} comment, sous les hypoth\`eses i) et ii) de \ref{thprin}, 
l'hypoth\`ese iii) se ram\`ene \`a une estimation de l'ordre de nilpotence de certains op\'erateurs de monodromie. C'est cette estimation que l'on ne sait pour l'instant obtenir que pour $\MC^{d/K}_{Dr,0}$ en utilisant la suite spectrale de Rapoport-Zink dans ce cas "semi-stable".

Nous reviendrons \`a la fin de cette introduction sur le cas $LT$ de
la conjecture \ref{conj}.

\bigskip

Revenons maintenant aux difficult{\'e}s inh{\'e}rentes {\`a} l'utilisation des
cat{\'e}gories d{\'e}riv{\'e}es. Il y a principalement deux {\'e}cueils {\`a} franchir :
le premier vient de ce qu'on veut un complexe
$R\Gamma_c$ dans la cat{\'e}gorie d{\'e}riv{\'e}e des $G_d\times
\dd$-modules {\em lisses}\footnote{Une construction {\em ad hoc} d'un
  tel complexe appara{\^\i}t dans \cite{HaCusp} pour $\MC^{d/K}_{Dr,n}$, mais nous voulons ici une
  d{\'e}finition g{\'e}n{\'e}rale et fonctorielle en les paires $(X,H)$ o{\`u} $H$ est
  un groupe topologique  agissant sur un
  espace analytique rigide $X$.}, alors que les techniques {\'e}quivariantes
usuelles fourniraient plut{\^o}t un objet de la cat{\'e}gorie d{\'e}riv{\'e}e de tous les
$G_d\times \dd$-modules. Ce probl{\`e}me a {\'e}t{\'e} joliment r{\'e}solu par
Berkovich dans un manuscrit non publi{\'e}, au moins pour la cohomologie
{\'e}tale {\`a} coefficients finis. Nous consacrerons la partie \ref{berkovich}
de ce texte {\`a} exposer sa th{\'e}orie et {\`a} l'{\'e}tendre aux coefficients
$l$-adiques. Mentionnons seulement ici que lorsque
 $X$ est un $K$-espace analytique au sens de \cite{Bic2} muni d'une action
``continue'' au sens de \cite{Bic3} d'un groupe topologique $H$,
Berkovich introduit la notion de faisceau {\'e}tale $H$-{\'e}quivariant
``lisse'', dont les sections globales \`a supports compacts sont en particulier munies d'une
action lisse de $H$. De tels faisceaux forment un topos not{\'e} $(X,H)_{et}^\sim$. 
 Bien-s{\^u}r, cela nous conduit {\`a} utiliser  sa th{\'e}orie
des $K$-espaces analytiques : nous consid{\`e}rerons donc les $\MC^{d/K}_{Dr,n}$,
resp. $\MC^{d/K}_{LT,n}$, comme des $K$-espaces analytiques au sens de
\cite{Bic2} munis d'actions ``continues'' de $G_d$, resp. $\dd$.

Par contre leurs limites projectives ne sont que des pro-objets dans
la cat{\'e}gorie des espaces de Berkovich, et le deuxi{\`e}me {\'e}cueil est qu'aucun formalisme de
cohomologie {\'e}tale n'est connu de l'auteur pour ces objets. Pour d{\'e}finir les
complexes $R\Gamma_c$, nous utiliserons
les morphismes de p{\'e}riodes  
$$\pi^n_{?} : \MC^{d/K}_{?,n} \To{} \PC^{d/K}_? \;\;\hbox{ pour }\;\;?=LT\hbox{ ou
  } Dr
$$ 
o{\`u} $\PC^{d/K}_{LT} := \PM^{d-1}_K$ et $\PC^{d/K}_{Dr}:=
\Omega^{d-1}_K$ est le compl{\'e}mentaire des
hyperplans $K$-rationnels dans $\PM^{d-1}$.
Pour $?=Dr$, le morphisme de p{\'e}riodes  est donn{\'e} par le th{\'e}or{\`e}me de Drinfeld et est
{\'e}quivariant pour l'action ``naturelle'' de $G_d$ sur $\PC^{d/K}_{Dr}$, et pour
$\MC^{d/K}_{LT}$ il est d{\'e}fini dans \cite{HG} et est {\'e}quivariant pour une
certaine action de $\dd$ sur $\PC^{d/K}_{LT}$. Dans chacun des cas, les
morphismes de p{\'e}riodes sont {\em {\'e}tales}. Leur d\'efinition et
leurs propri\'et\'es sont rappel\'ees en \ref{perLT} et \ref{perDr}.
 Pour donner une
d{\'e}finition uniforme du $R\Gamma_c$, notons temporairement
$J_{LT}:=\dd$ et $J_{Dr}:= G_d$.

Nous d{\'e}finissons en \ref{defRG} le complexe
$R\Gamma_c(\MC^{d/K}_{?},\ZM_l)$ comme l'{\'e}valuation du foncteur
d{\'e}riv{\'e} du foncteur qui {\`a} un syst{\`e}me projectif $(\FC_m)_m$ de
faisceaux de $l^m$-torsion  $J_?$-{\'e}quivariants ``lisses'' sur
$\PC^{d/K}_?$  associe le $G_d\times \dd\times W_K$-module  
$$\limi{n} \Gamma_c\left(\MC^{d/K}_{?,n},{\limp{m}}^\infty \pi^{n,*}_?(\FC_m)\right)$$
en le syst{\`e}me projectif de faisceaux constants $(\ZM/l^m\ZM)_{m\in
  \NM}$ (la notation $\limproj^\infty$ est expliqu\'ee en \ref{deflimpinf}).
On pose aussi 
$R\Gamma_c(\MC^{d/K}_{?},\o\QM_l):=R\Gamma_c(\MC^{d/K}_{?},\ZM_l)\otimes_{\ZM_l}
\o\QM_l.$ De par leur d{\'e}finition, ce sont des complexes de $G_d\times
\dd\times W_K$-modules, mais en l'absence d'une cat{\'e}gorie maniable de
$W_K$-modules, nous nous contenterons de les voir comme des objets de la
cat{\'e}gorie d{\'e}riv{\'e}e des $\o\QM_l(GD)$-modules {\em lisses} munis
d'une action de $W_K$. Nous v{\'e}rifierons bien-s{\^u}r en \ref{defRG} que les groupes de
cohomologie de ce $R\Gamma_c(\MC^{d/K}_?,\ZM_l)$ sont bien les
$H^i_c(\MC^{d/K}_?,\ZM_l)$ d{\'e}finis plus haut.

\bigskip

Une fois le complexe $R\Gamma_c$ bien d{\'e}fini, on est en pr{\'e}sence d'un
objet {\em a priori} compliqu{\'e}. Le miracle est que cet objet est en
fait ``aussi simple qu'il peut l'\^etre'' : il est {\em scindable},
au sens de la section \ref{sectionscin}.

\begin{propsec}\label{prop1} 
 Soit $\MC^{d/K}=\MC^{d/K}_{Dr}$ ou $\MC^{d/K}_{LT}$ et
supposons que la  description conjectur{\'e}e par Harris pour les $H^i_c(\MC^{d/K},\o\QM_l)$ est
  valide. Alors il existe un isomorphisme (pas unique) dans $D^b(GD)$,
 $$R\Gamma_c(\MC^{d/K},\o\QM_l)
    \simto \bigoplus_{i\in \NM} H^i_c(\MC^{d/K},\o\QM_l)[-i].$$
 \end{propsec}

Nous donnerons deux preuves de cette assertion. L'une en \ref{prop1bis}
repose sur un argument relevant purement
de la th{\'e}orie des repr{\'e}sentations de $G_d$ : le calcul complet des
groupes d'extensions $\ext{i}{\pi}{\pi'}{G_d}$ pour les couples de repr{\'e}sentations
elliptiques $\pi$ et $\pi'$ de $G_d$ montre que {\em tout} complexe
ayant la m\^eme cohomologie que $\MC^{d/K}$ est scindable comme dans la
proposition. L'autre utilise la forme
particuli\`ere de l'action d'un rel\`evement de Frobenius sur les $H^i_c$.

\bigskip

Signalons bri\`evement un point technique ;
soit $\varpi$ une uniformisante de $K$ et $\varpi^\ZM$ le sous-groupe de $G_d$ engendr\'e
par $\varpi$. Les espaces quotients $\o\MC^{d/K}:=\MC^{d/K}/\varpi^\ZM$ existent et
pr\'esentent l'avantage d'avoir une cohomologie 
$GD$-admissible, alors que celle de $\MC^{d/K}$ est
``admissible-modulo-le-centre''\footnote{Dans le cas $Dr$, ces propri\'et\'es
  d'admissibilit\'e ne sont connues que via le th\'eor\`eme de Faltings.}. 
Nous expliquerons en \ref{reduc2} comment, pour prouver le th\'eor\`eme
\ref{thprin} on se ram\`ene facilement \`a \'etudier le complexe de cohomologie
$R\Gamma_c(\o\MC^{d/K},\o\QM_l)$, vu comme objet de $D^b(\o\QM_l(\o{GD}))$
o\`u $\o{GD}:=GD/\varpi^\ZM$.

\bigskip

On peut encore simplifier l'\'etude en d\'ecomposant
$R\Gamma_c(\o\MC^{d/K},\o\QM_l)$ selon
la th\'eorie des
repr\'esentations du groupe compact $\dd/\varpi^\ZM$. Nous expliquons en
\ref{pardecomp} comment celle-ci fournit en effet  une d\'ecomposition $W_K$-\'equivariante
$$ R\Gamma_c(\o\MC^{d/K},\o\QM_l) \simeq \bigoplus_{\rho\in
  \Irr{\o\QM_l}{\dd/\varpi^\ZM}}
R\Gamma_c(\o\MC^{d/K},\o\QM_l)_\rho. $$
Par exemple, dans le cas $Dr$, le complexe $R\Gamma_c(\o\MC^{d/K}_{Dr},\o\QM_l)_\rho$ est canoniquement isomorphe au complexe $R\Gamma_c(\Omega_K^{d-1},\LC_\rho)$ o\`u $\LC_\rho$ est un syst\`eme local $l$-adique associ\'e \`a $\rho$. 
Les complexes $R\Gamma_c(\o\MC^{d/K},\o\QM_l)_\rho$, pour
$\rho\in\Irr{\o\QM_l}{\dd/\varpi^\ZM}$ sont eux  aussi
scindables dans  $D^b(\o\QM_l(\o{GD}))$, ce qui permet  d'{\'e}crire leurs
alg{\`e}bres d'endomorphismes sous
la forme 
$$ \endo{D^b(\o{GD})}{R\Gamma_c(\o\MC^{d/K},\o\QM_l)_\rho} \simeq \bigoplus_{i\geq j}
\ext{i-j}{H^i_c(\o\MC^{d/K},\o\QM_l)_\rho}{H^j_c(\o\MC^{d/K},\o\QM_l)_\rho}{\o{GD}}$$
le produit sur le $\o\QM_l$-espace vectoriel de droite {\'e}tant donn{\'e} par
le $\cup$-produit. Cette fois, c'est la th{\'e}orie des
repr{\'e}sentations de $G_d/\varpi^\ZM$ qui
permet de d{\'e}crire compl{\`e}tement ce $\cup$-produit, et par cons{\'e}quent
permet de d{\'e}crire l'alg{\`e}bre des endomorphismes de $R\Gamma_c$.
Pour {\'e}noncer le r{\'e}sultat obtenu, fixons $\rho  \in
\Irr{\o\QM_l}{\dd/\varpi^\ZM}$  et notons
$(\sigma_\rho, V_{\sigma_\rho})$ la repr{\'e}sentation $l$-adique de
$W_K$ associ{\'e}e {\`a} $\rho$ 
par correspondance de Langlands et fonctorialit{\'e} de Jacquet-Langlands
: il s'agit d'une repr{\'e}sentation ind{\'e}composable de dimension $d$ et
l'alg{\`e}bre engendr{\'e}e par
l'image de $W_K$ dans $\endo{\o\QM_l}{V_{\sigma_\rho}}$ est une alg{\`e}bre
triangulaire par blocs que nous noterons $\AC_\rho$.
Le calcul de $\cup$-produits mentionn\'e ci-dessus montre alors ({\em cf} \ref{propnoneq})

\begin{propsec} \label{prop2}
 Sous les hypoth{\`e}ses de la proposition pr\'ec\'edente, il existe un
 isomorphisme de $\o\QM_l$-alg{\`e}bres
$$\endo{D^b(\o{GD})}{R\Gamma_c(\o\MC^{d/K},\o\QM_l)_\rho} \simto 
\AC_\rho.$$
\end{propsec}
On veut ensuite {\'e}lucider l'action de $W_K$ sur
le facteur $R\Gamma_c(\o\MC^{d/K},\o\QM_l)_\rho$. 
Pour des raisons tr\`es g\'en\'erales, {\em cf} \ref{actionphi} i),  puisque
l'action de l'inertie sur les $(H^i_c)_\rho$ 
est potentiellement triviale  (par la description de Boyer-Harris) 
alors  celle sur le $(R\Gamma_c)_\rho$ est
potentiellement unipotente. Par contre  une propri{\'e}t{\'e} de ``continuit{\'e}'' sera
n{\'e}cessaire pour montrer qu'apr\`es restriction \`a un sous-groupe d'indice fini, l'inertie agit par son  quotient
$l$-adique et  ``par l'exponentielle
d'un endomorphisme nilpotent $N_\rho$'' ; c'est l{\`a} que nous aurons
besoin de l'admissibilit{\'e} 
de la cohomologie enti{\`e}re (hypoth\`ese ii) de \ref{thprin}). Pour
d{\'e}crire alors pr{\'e}cis{\'e}ment ce $N_\rho$, nous aurons besoin
de minorer son ordre de nilpotence et pour cela
nous utiliserons l'hypoth\`ese iii) du th\'eor\`eme \ref{thprin}, c'est-\`a-dire
la validit\'e de la conjecture
monodromie-poids pour certaines  vari{\'e}t{\'e}s uniformis{\'e}es par les
$\MC^{d/K}_{Dr,n}/\varpi^\ZM$. 

En notant  $\gamma : W_K \To{}
\endo{D^b(\o{GD})}{R\Gamma_c(\o\MC^{d/K},\o\QM_l)_\rho}^\times$ le
morphisme  de groupes donnant l'action de $W_K$ sur
$R\Gamma_c(\o\MC^{d/K},\o\QM_l)_\rho$, on obtient la description
suivante, {\em cf} \ref{preuveprop3}
\begin{propsec} \label{prop3}
  Supposons les $\ZM_l(\o{GD})$-modules $H^i_c(\MC^{d/K}_{Dr}/\varpi^\ZM,\ZM_l)$ admissibles
  et supposons la conjecture monodromie poids v{\'e}rifi{\'e}e pour (certaines) vari{\'e}t{\'e}s
  uniformis{\'e}es par les $\MC^{d/K}_{Dr,n}/\varpi^\ZM$. Alors il existe un
 isomorphisme d'alg{\`e}bres rendant le diagramme suivant commutatif
$$\xymatrix{
\endo{D^b(\o{GD})}{R\Gamma_c(\MC^{d/K}_{Dr}/\varpi^\ZM,\o\QM_l)_\rho}
\ar[r]^-\sim   & \AC_\rho \\ W_K \ar[u]^\gamma \ar[ur]^{\sigma_\rho} &
}$$
\end{propsec}
Dans le cas $\rho=1$, le complexe \`a \'etudier est simplement $R\Gamma_c(\Omega_K^{d-1},\o\QM_l)$. En utilisant la propri\'et\'e de r\'eduction semi-stable des vari\'et\'es uniformis\'ees par $\Omega_K^{d-1}$ et la suite spectrale de Rapoport-Zink-Steenbrink, on peut minorer l'ordre de $N$ et montrer la proposition ci-dessus sans recours \`a la conjecture Monodromie-Poids. Au contraire, on en d\'eduit ensuite une nouvelle preuve de la conjecture MP dans ce cas, cf \ref{MPdp}.

{\`A} partir de cette proposition, le th{\'e}or{\`e}me \ref{thprin} se prouve
en utilisant {\`a} nouveau le calcul explicite des extensions entre
repr{\'e}sentations elliptiques et de leurs $\cup$-produits. Nous n'en
disons pas plus dans cette 
introduction car cela nous plongerait dans la combinatoire de ces
repr{\'e}sentations. 
N\'eanmoins, le lecteur aura compris qu'une fois admis les r\'esultats
profonds et difficiles de Boyer et Faltings, ce calcul explicite est
la clef de la conjecture \ref{conj}. En fait, au moyen de la th\'eorie
de Bushnell-Kutzko, ce calcul pour les repr\'esentations elliptiques
se ram\`ene \`a celui pour les sous-quotients irr\'eductibles de la
repr\'esentation $\CC^\infty_{\o\QM_l}(G_d/B)$ o\`u $B$ est un
\borel. Ce dernier calcul a un sens sur un groupe r\'eductif
$p$-adique quelconque. Il est effectu\'e dans la partie
\ref{extensions} dans le cas d\'eploy\'e. On utilise la
param\'etrisation des sous-quotients irr\'eductibles de
$\CC^\infty_{\o\QM_l}(G/B)$ par les sous-ensembles de l'ensemble $S$
des racines simples d'un tore maximal $T$ de $B$ dans $\hbox{Lie}(B)$  donn\'ee par 
$$(I\subseteq S)\mapsto \pi_I := \CC^\infty_{\o\QM_l}(G/P_I)/\sum_{J\supset I}
\CC^\infty_{\o\QM_l}(G/P_J) $$
o\`u $P_I$ est le parabolique contenant $B$ associ\'e \`a $I$, ({\em
  cf} \ref{defrep}).
Alors l'\'enonc\'e principal de la partie \ref{extensions}
g\'en\'eralise des r\'esultats ant\'erieurs de Borel, Casselman et
Schneider-Stuhler :
\begin{theosec} \label{theoextintro}
Supposons $G$ semi-simple et posons $\delta(I,J)=|I\cup J|-
|I\cap J|$ pour $I,J\subseteq S$.
  \begin{enumerate}
  \item Soient $I,J$ deux sous ensembles de $S$, alors :
    $$
    \cas{\ext{*}{\pi_I}{\pi_J}{G}}{\o\QM_l}{*=\delta(I,J)}{0}{*\neq
      \delta(I,J)}.$$
\item Soient $I,J,K$ trois sous-ensembles de $S$ tels que
  $\delta(I,J)+\delta(J,K)=\delta(I,K)$, alors le cup-produit
$$  \cup \;:\;\;
\ext{\delta(I,J)}{\pi_I}{\pi_J}{G} \otimes_{\o\QM_l}
\ext{\delta(J,K)}{\pi_J}{\pi_K}{G} \To{} 
\ext{\delta(I,K)}{\pi_I}{\pi_K}{G}
$$ est un isomorphisme.
\end{enumerate}
\end{theosec}
En fait, le th\'eor\`eme \ref{theoext} est \'enonc\'e pour des
coefficients plus g\'en\'eraux que $\o\QM_l$. Signalons que le point
i) a \'et\'e obtenu par Orlik dans \cite{Orlikext}, par une m\'ethode
plus directe et moins explicite que la notre.

\bigskip

Revenons maintenant au cas $LT$ de la conjecture
\ref{conj}. La strat\'egie que nous avons esquiss\'ee \`a travers les
trois propositions ci-dessus fonctionne de la m\^eme mani\`ere
jusqu'\`a la d\'efinition de $N_\rho$ (incluse). Le probl\`eme est
alors de trouver un moyen pour minorer l'ordre de nilpotence de
$N_\rho$ ; on ne peut plus uniformiser ``\`a niveau fini'' comme dans le cas $Dr$.
Il serait int\'eressant de trouver un autre moyen que l'uniformisation
pour \'etudier $N_\rho$. Pour l'instant, je ne vois qu'une 
strat\'egie tr\`es indirecte consistant \`a relever au niveau des
complexes l'\'enonc\'e cohomologique principal de Faltings 
dans \cite{FaltDrin} :

{\em Les $G_d\times \dd\times W_K$-modules $H^i_c(\MC^{d/K}_{Dr},\ZM_l)$ et
  $H^i_c(\MC^{d/K}_{LT},\ZM_l)$ sont isomorphes.}

En fait cet {\'e}nonc{\'e} cohomologique est annonc{\'e} comme cons{\'e}quence d'un
{\'e}nonc{\'e} g{\'e}om{\'e}trique qui s'exprime en termes de sch{\'e}mas formels. {\`A} la
lecture de \cite{FaltDrin} et apr\`es discussions avec Fargues, 
 on peut esp\'erer  : 

\begin{boite}{Conjecture}
{\em  Il existe une {\'e}quivalence de topos $\Phi:\;\;(\PC^{d/K}_{LT},J_{LT})_{et}^\sim
  \simto (\PC^{d/K}_{Dr},J_{Dr})_{et}^\sim$ envoyant les faisceaux constants sur les
  faisceaux constants et un isomorphisme de foncteurs
$$ \limi{n}\Gamma_c(\MC^{d/K}_{Dr,n},\Phi(-)) \simto \limi{n}\Gamma_c(\MC^{d/K}_{LT,n},-).$$}
\end{boite}

Remarquons seulement que si les pro-espaces $\MC^{d/K}_{Dr}$ et $\MC^{d/K}_{LT}$
{\'e}taient de ``vrais'' $K$-espaces analytiques $G_d\times \dd\times
W_K$-isomorphes, alors en notant simplement $\MC^{d/K}$ ces espaces,
les propri{\'e}t{\'e}s des morphismes de p{\'e}riodes permettraient d'identifier
$\MC^{d/K}/J_{Dr} \simto \PC_{LT}^{d/K}$ et $\MC^{d/K}/J_{LT} \simto
\PC_{Dr}^{d/K}$, et 
l'{\'e}nonc{\'e} ci-dessus  serait une
cons{\'e}quence formelle de ces identifications. Mais ce n'est pas si simple...
Quoiqu'il en soit, cette conjecture impliquerait bien-s{\^u}r
 l'existence d'un isomorphisme
  $W_K$-{\'e}quivariant dans $D^b(\o\QM_l(GD))$
$$ R\Gamma_c(\MC^{d/K}_{LT},\o\QM_l) \simeq R\Gamma_c(\MC^{d/K}_{Dr},\o\QM_l),$$
de sorte que la validit{\'e} de la conjecture \ref{conj} pour $\MC^{d/K}_{Dr}$ impliquerait la
m{\^e}me pour $\MC^{d/K}_{LT}$.

\bigskip

D\'ecrivons succinctement le contenu des diverses parties. Le but de la premi\`ere partie est d'exposer les principales id\'ees de ce texte en les d\'ebarassant au maximum de toutes les difficult\'es techniques.
On y traite une variante (moralement un cas particulier)
de la conjecture \ref{conj} concernant l'espace   sym\'etrique $\Omega^{d-1}_K$
de Drinfeld (c'est-\`a-dire, essentiellement $\MC^{d/K}_{Dr,0}$). Le r\'esultat principal \ref{theodp} est inconditionnel
et sa preuve sert de mod\`ele pour le cas g\'en\'eral ; 
on y admet certains r\'esultats techniques d\'emontr\'es ult\`erieurement, comme le calcul d'extensions dans la s\'erie principale, l'existence du formalisme $l$-adique "lisse", et quelques crit\`eres de scindages d'un complexe. On y red\'emontre aussi le th\'eor\`eme d'Ito sur la conjecture monodromie-poids pour les vari\'et\'es uniformis\'ees par $\Omega_K^{d-1}$.

\`A partir de la troisi\`eme partie, l'exposition est plus lin\'eaire. 
La partie \ref{extensions} contient le calcul d'extensions et de
$\cup$-produits entre sous-quotients irr\'eductibles d'une induite de
Borel pour un groupe $p$-adique r\'eductif d\'eploy\'e. Le lecteur d\'esirant acc\'eder au plus vite \`a la preuve du th\'eor\`eme \ref{thprin} devrait se contenter de lire le r\'esultat principal.

Dans la  partie \ref{elliptiques} on d\'efinit et caract\'erise les
repr\'esentations elliptiques et on explicite leur comportement \`a
travers les correspondances de Langlands et Jacquet-Langlands.
Beaucoup de choses y sont bien connues.

La partie \ref{berkovich} parle de cohomologie \'etale d'espaces de
Berkovich. Une grande part consiste \`a exposer des r\'esultats de
Berkovich lui-m\^eme. Les deux seules innovations sont l'introduction
d'un $R\Gamma_c$   ``lisse'' $l$-adique pour un espace muni d'une
action continue d'un groupe, et le rel\`evement de la suite spectrale
de type Hochshild-Serre de Fargues dans \cite{Fargues} au niveau des
complexes de cohomologie. Le lecteur press\'e devrait se contenter de lire les d\'efinitions du paragraphe \ref{deflimpinf} et l'\'enonc\'e de \ref{theoHS}.

Dans la  partie \ref{drinfeld} on rappelle la d\'efinition des espaces de
Drinfeld et Lubin-Tate adapt\'ee \`a la conjecture \ref{conj}.
Puis on d\'efinit le $R\Gamma_c$ comme esquiss\'e dans cette
introduction, en utilisant le formalisme de la partie \ref{berkovich}.

La  partie \ref{real} contient la preuve proprement dite du
th\'eor\`eme \ref{thprin}, apr\`es quelques pr\'eliminaires sur les scindages et les endomorphismes de complexes qui sont certainement bien connus des sp\'ecialistes.

Dans la partie \ref{retour}, on revient sur l'exemple du demi-plan : nous donnons un analogue du th\'eor\`eme principal de \ref{dp} \`a coefficients dans un anneau de torsion "fortement banal".

Dans le premier appendice, nous pr\'esentons un nouveau calcul de la cohomologie \`a supports compacts de $\Omega_K^{d-1}$ valable aussi pour les coefficients $l$-adiques.

Enfin on a renvoy\'e au second appendice quelques preuves qui nous
semblaient alourdir inutilement le texte.

\bigskip

Je remercie chaleureusement M. Harris pour les nombreux {\'e}changes
{\`a} propos du contenu de ce texte et pour m'avoir communiqu{\'e} sa
conjecture ainsi qu'un article en pr{\'e}paration avec Taylor. Je remercie
aussi V. Berkovich pour m'avoir permis d'inclure certains de ses
r{\'e}sultats non publi{\'e}s et P. Boyer qui a essay\'e de m'expliquer son r\'esultat.
J'ai aussi \'et\'e aid\'e au cours de ma r\'edaction par de nombreuses conversations avec 
 L. Fargues, 
  G. Henniart, B.C. Ngo, S. Orlik et M. Strauch que je remercie tous. 
   Les id{\'e}es principales de ce texte sont apparues au cours de longues
discussions avec Alain Genestier qui m'a tout appris sur ce
sujet. Sans lui, pas une ligne n'aurait {\'e}t{\'e} {\'e}crite. Ces discussions
ont {\'e}t{\'e} permises par l'environnement exceptionnel de l'IHES ; je
remercie cet institut pour sa longue hospitalit{\'e} et Laurent Lafforgue pour m'y avoir invit\'e en compagnie de ceux que je viens de remercier.

\section{L'exemple du demi-plan} \label{dp}

Le ``demi-plan sup{\'e}rieur'' ou ``espace sym{\'e}trique'' de Drinfeld  appara{\^\i}t dans \cite{Drinell}
o{\`u} il est d{\'e}fini comme un sous-$K$-espace analytique rigide de l'espace
projectif. Ses points sont donn{\'e}s 
par
$$ \Omega_K^{d-1}(\ka) = \PM^{d-1}(\ka) \setminus \bigcup_{H \in \HC_K} H $$
o{\`u} $\HC_K$ d{\'e}signe l'ensemble des hyperplans $K$-rationnels de
$\PM^{d-1}$. Nous d{\'e}signerons par le m{\^e}me symbole $\Omega_K^{d-1}$ le
$K$-espace analytique de Berkovich correspondant dont on trouve une
description dans \cite{BicCRAS} ; ses points s'identifient {\`a}  
des classes d'{\'e}quivalences de semi-normes multiplicatives sur l'anneau de polyn{\^o}mes
$K[X_0,\cdots,X_{d-1}]$ dont la restriction au sous-espace des
polyn{\^o}mes homog{\`e}nes de degr{\'e} $1$ est une norme.

L'espace analytique $\Omega_K^{d-1}$ est naturellement muni d'une action de
$G_d=GL_d(K)$ triviale sur le centre, qui est {\em continue} au sens de
\cite[par. 6-7]{Bic3}. Le 
changement de base $\Omega_K^{d-1,ca}:=\Omega_K^{d-1} \wh\otimes_K \ka$ est aussi muni
d'une action continue du groupe de Weil $W_K$ de $K$ (et m{\^e}me de son
groupe de Galois). 

Comme on l'a annonc\'e dans l'introduction, le but de cette partie est d'exposer le plus rapidement possible dans un cas particulier simple, les arguments que nous utiliserons dans la partie \ref{real}. Nous utilisons en particulier le th\'eor\`eme principal d'Ito dans \cite{Ito}. Nous montrerons dans la partie \ref{retour} comment se passer de, et red\'emontrer, ce th\'eor\`eme. Nous g\'en\'eraliserons aussi en rempla{\c c}ant les coefficients $\o\QM_l$ par un anneau de torsion "fortement banal" au sens de \ref{bonbanal}.

\medskip

\noindent {\bf Convention h\'er\'etique :} Suivant l'usage, notons
$\ZM_l(1):=\limproj_n (\mu_{l^n}(K^{ca}))$ le $\ZM_l$-torseur des
$l$-racines de l'unit\'e. Pour rendre certaines formules plus
concr\`etes, nous fixons une fois pour toutes un g\'en\'erateur
topologique $\mu$ de $\ZM_l(1)$. Cela nous fixe aussi un isomorphisme
de $\ZM_l$-modules $\mu^*:\;\ZM_l(1)\simto \ZM_l$ qui transporte l'action
naturelle de $W_K$ sur $\ZM_l(1)$ sur le caract\`ere $|-|$ de  $W_K$ trivial
sur l'inertie $I_K$ et qui envoie tout rel\`evement de Frobenius g\'eom\'etrique
$\phi$ sur $q^{-1}$.

\subsection{{\'E}nonc{\'e} du th{\'e}or{\`e}me}

Dans \cite{SS1}, Schneider et Stuhler ont calculé les groupes de cohomologie (sans supports) 
de $\Omega_K^{d-1}$ pour toute théorie cohomologique satisfaisant certains axiomes. La cohomologie étale de $\Omega_K^{d-1,ca}$ à coefficients dans un anneau de torsion première à $p$ définie par Berkovich dans \cite{Bic2} satisfait ces axiomes. 
Cependant, la cohomologie sans supports présente deux inconvénients : elle ne fournit pas des représentations lisses de $G_d$ et de plus il n'existe pas à ce jour de définition de la cohomologie $l$-adique. En revanche, on sait par des résultats généraux de Berkovich (voir partie \ref{berkovich}) que la cohomologie étale à supports compacts est lisse pour l'action de $G_d$, et on sait en donner une (des) version(s) $l$-adique(s), lisse elle aussi.
Par le m{\^e}me principe de calcul (cohomologie d'un arrangement
d'hyperplans projectifs), il est facile de d{\'e}duire de \cite{SS1} ce qu'est la
cohomologie  {\`a} supports compacts $H^i_c(\Omega_K^{d-1,ca},\Lambda)$ lorsque $\Lambda$ est de torsion. C'est aussi un cas particulier du preprint \cite{O2} d'Orlik. Le passage aux coefficients $l$-adiques pose de nouveaux probl\`emes. Bien qu'il ne fasse aucun doute que le r\'esultat de Schneider-Stuhler 
que nous d\'ecrivons ci-dessous reste vrai,
nous donnerons dans l'appendice un calcul par une m\'ethode différente de celle de Schneider et Stuhler et qui fonctionne pour les coefficients $l$-adiques.


Afin de décrire ces espaces de cohomologie, nous introduisons les
``repr{\'e}sentations elliptiques de la s{\'e}rie principale''.  Par
d{\'e}finition, ce sont les sous-quotients irr{\'e}ductibles d'une induite
droite $\ind{B}{G_d}{\o\QM_l}$ de la repr{\'e}sentation triviale d'un \borel $B$
de $G_d$.

\alin{Rappel de classification} \label{rappelclassif}
Fixons $T\subset B$ un tore maximal d'un \borel et $S$ l'ensemble des
 racines simples de $T$ dans $\hbox{Lie}(B)$. On a une bijection naturelle 
$$\application{}{\{\hbox{Sous-ensembles de } S\}}{\{\hbox{Sous-groupes
    paraboliques    contenant } B\}}{I\subseteq S}{P_I:=Z_{G_d}(\bigcap_{\alpha\in I} \ker
  \alpha).B}$$
telle que $P_\emptyset=B$ et $P_S=G_d$.
On sait alors, voir \cite[X 4.6-4.11]{BW} et le rappel dans la partie
\ref{extensions}, que l'induite $\ind{P_I}{G_d}{\o\QM_l}$ a un unique 
quotient irr{\'e}ductible que nous noterons $\pi_I$, que c'est un
sous-quotient irr{\'e}ductible de $\ind{B}{G_d}{\o\QM_l}$ puisqu'on a une
injection canonique $\ind{P_I}{G_d}{\o\QM_l}\subset \ind{B}{G_d}{\o\QM_l}$, et que
l'application
$$\application{}{\{\hbox{Sous-ensembles de} S\}}{\{\hbox{Sous-quotients
    irr{\'e}ductibles de } \ind{B}{G_d}{\o\QM_l}\}}{I\subseteq S}{\pi_I}$$
est une bijection. En particulier, $\pi_\emptyset$ est la
repr{\'e}sentation de Steinberg et $\pi_S$ est la repr{\'e}sentation triviale.

Pour {\'e}noncer le r{\'e}sultat de Schneider et Stuhler, il convient de
num{\'e}roter $S$, {\em i.e.} de fixer une bijection $S\simto
\{1,\cdots, d-1\}$. On le fait de telle sorte que pour
$i\in\{1,\cdots,d-1\}$, le parabolique
$P_{S\setminus\{i\}}$ soit le stabilisateur d'un sous-espace de dimension $i$,
pour l'action naturelle de $G_d=GL_d(K)$ sur $K^d$ ; cela revient \`a num\'eroter la sur-diagonale "de haut en bas". Notons toujours $|.|$
le caract{\`e}re de $W_K$ qui envoie les Frobenius g{\'e}om{\'e}triques sur
$q^{-1}$. 
\begin{theo} (Schneider-Stuhler \cite{SS1} et \ref{cohodp}) \label{ss}
  Pour $i=0,\cdots, d-1$, il existe des isomorphismes $G_d\times W_K$-{\'e}quivariants
$$ H^{d-1+i}_c(\Omega_K^{d-1,ca},\o\QM_l) \simto \pi_{\{1,\cdots,i\}}\otimes
|.|^{-i}$$
\end{theo}
Nous conviendrons que pour $i=0$, la repr{\'e}sentation
$\pi_{\{1,\cdots,i\}}$ est $\pi_\emptyset$.

\alin{Correspondance de Langlands pour les $\pi_I$} \label{deftauI}
Nous allons d{\'e}crire la repr{\'e}sentation continue $l$-adique de $W_K$
associ{\'e}e {\`a} $\pi_I$ par correspondance de Langlands. L'espace de la
repr{\'e}sentation sera $V:=\o\QM_l^d$ dont nous noterons
$e_0,\cdots,e_{d-1}$ la base canonique et $E_{ij}\in \endo{\o\QM_l}{V}$
la base correspondante. Pour chaque $I\subseteq S$, on note $N_I$
l'endomorphisme nilpotent de $V$ d{\'e}fini par 
$$N_I:= \sum_{i\in \o{I^c}} E_{i-1,i}$$
o{\`u} $I^c:=S\setminus I$ d{\'e}signe le compl{\'e}mentaire de $I$ dans $S$ et
$\o{I^c}$ est l'image de $I^c$ par l'application $x\in S\simeq\{1,\cdots,d-1\} \mapsto \o{x}:=d-x$.
 ($S$ a \'et\'e num{\'e}rot{\'e} au paragraphe pr{\'e}c{\'e}dent. Plus conceptuellement l'application $\o{ }$ est $-w_0$ si $w_0$ est l'\'el\'ement de plus grande longueur du groupe de Weyl).
On d{\'e}finit d'abord une repr{\'e}sentation semi-simple $\tau_I^{ss}$ 
(ind{\'e}pendente de $I$, en fait) de $W_K$ sur $V$ en posant
$$\tau_I^{ss}(w)(e_i)= |w|^{-i}e_i.$$
On choisit ensuite un rel{\`e}vement de Frobenius g{\'e}om{\'e}trique $\phi$, ce qui
nous permet d'{\'e}crire tout $w\in W_K$ sous la forme unique $w=\phi^{\nu(w)}
i_\phi(w)$ o{\`u} $\nu(w)\in\ZM$ et $i_\phi(w)\in I_K$ (inertie).
On a alors $|w|=q^{-\nu(w)}$.
Notons maintenant 
$$t_\mu : I_K \To{t_l} \ZM_l(1) \To{\mu^*} \ZM_l$$ o\`u $t_l$ est la 
projection vers le $l$-quotient de $I_K$.
Un calcul rapide montre que pour chaque $I\subseteq S$, la formule
$$ \tau_I^{\phi,\mu}(w) := \tau_I^{ss}(w)\hbox{exp}.(N_It_\mu(i_\phi(w))) $$
d{\'e}finit une repr{\'e}sentation de $W_K$.
Comme la notation le sugg\`ere, celle-ci d\'epend du choix de $\phi$ et $\mu$ mais un autre calcul (un peu moins rapide, {\em cf} \cite[8.4.2]{DelAntwerp}) montre que sa classe n'en d\'epend pas. 
Nous noterons donc simplement $\tau_I$ cette classe.
Nous montrerons en \ref{corell} que la correspondante $\sigma_d(\pi_I)$ de $\pi_I$ par la
correspondance de Langlands est la (classe de la) repr{\'e}sentation $\sigma_d(\pi_I)=
\tau_I\otimes |.|^{\frac{d-1}{2}}$.

\medskip

\noindent{\em Remarque :} on calcule que la contragr\'ediente de $(\tau_I\otimes |.|^\frac{d-1}{2})$ est $\tau_{\o{I}}\otimes |.|^\frac{d-1}{2})$. Par compatibilit\'e de la correspondance \`a la contragr\'ediente, on voit alors que $\pi_I^\vee=\pi_{\o{I}}$. Nous montrerons l'analogue pour un groupe d\'eploy\'e g\'en\'eral dans le point iv) du lemme \ref{serprincorps}.

\alin{Le complexe de cohomologie} 
Notons $PG_d:=PGL_d(K)$. L'action de $G_d$ sur $\Omega_K^{d-1,ca}$ se factorise par
$PG_d$. Notons aussi $D^b(\o\QM_lPG_d)$ la cat{\'e}gorie d{\'e}riv{\'e}e de la cat{\'e}gorie
ab{\'e}lienne des $\o\QM_l$-repr{\'e}sentations lisses de $PG_d$.
Les r{\'e}sultats de la partie
\ref{berkovich} permettent de d{\'e}finir un complexe 
$$ R\Gamma_c(\Omega_K^{d-1,ca},\o\QM_l) \in D^b(\o\QM_l PG_d)$$
muni d'une action de $W_K$, et  dont les espaces de cohomologie
s'identifient canoniquement aux espaces $H^i_c(\Omega_K^{d-1,ca},\o\QM_l)$.
Une propri{\'e}t{\'e} importante de ce complexe est la suivante : supposons
que $\Gamma\subset PG_d$ soit un sous-groupe discret  et sans
torsion de $PG_d$. On sait alors   que  $\Gamma$ agit librement sur $\Omega_{K}^{d-1}$ et
qu'on peut munir $\Omega_{K}^{d-1}/\Gamma$ d'une structure de $K$-espace
analytique telle que le quotient $\Omega_{K}^{d-1} \To{} \Omega_{K}^{d-1}/\Gamma$ soit un
rev{\^e}tement analytique  Galoisien. Dans ces conditions, la proposition
\ref{HSder} montre l'existence d'un isomorphisme $W_K$-{\'e}quivariant dans $D^b(\o\QM_l)$
\ini
\begin{equation}
  \label{HS}
  \o\QM_l \otimes^L_{\o\QM_l[\Gamma]} R\Gamma_c(\Omega_K^{d-1,ca},\o\QM_l) \simto
  R\Gamma_c(\Omega_K^{d-1,ca}/\Gamma,\o\QM_l)
\end{equation}
o{\`u} nous avons
fait l'abus de noter encore $R\Gamma_c(\Omega_K^{d-1,ca},\o\QM_l) \in
D^b(\o\QM_l\Gamma)$  l'image du complexe pr{\'e}c{\'e}demment d{\'e}fini par le
foncteur ``d'oubli'' $D^b(\o\QM_l 
PG_d) \To{} D^b(\o\QM_l\Gamma)$, o{\`u} le produit tensoriel est pris pour
l'augmentation $\o\QM_l[\Gamma] \To{}\o\QM_l$, et 
o{\`u} nous avons not{\'e} $R\Gamma_c(\Omega_K^{d-1,ca}/\Gamma,\o\QM_l) \in D^b(\o\QM_l)$
le complexe de cohomologie de l'espace $\Omega_K^{d-1,ca}/\Gamma$.

Rappelons que nous avons d{\'e}fini dans l'introduction le foncteur
$$\application{\HC^*:\;\;}{D^b(\o\QM_l)}{\{\o\QM_l-\hbox{espaces
  vectoriels}\}}
{\CC^\bullet}{\bigoplus_{i\in\ZM} \HC^i(\CC^\bullet)}.$$
Le but de cette partie est d'expliquer la preuve du th{\'e}or{\`e}me suivant,
qui est une variante dans un cas particulier de la conjecture
\ref{conj}.

\begin{theo} \label{theodp}
  Pour tout $I\subseteq S$, il existe  un isomorphisme $W_K$-{\'e}quivariant
$$ \HC^*(R\hom{R\Gamma_c(\Omega_K^{d-1,ca},\o\QM_l)}{\pi_I}{D^b(\o\QM_lPG_d)})
\simto \sigma_d(\pi_I)\otimes |.|^{\frac{d-1}{2}}.$$
\end{theo}

\subsection{Action de $W_K$ sur le complexe de cohomologie}

Dans cette section, nous noterons $\gamma:\;\; W_K \To{}
\endo{D^b(PG_d)}{R\Gamma_c(\Omega_K^{d-1,ca},\o\QM_l)}$ le morphisme qui d{\'e}finit
l'action de $W_K$ sur $R\Gamma_c(\Omega_K^{d-1,ca},\o\QM_l)$.
La premi{\`e}re {\'e}tape vers le th{\'e}or{\`e}me \ref{theodp} consiste {\`a} montrer que le complexe
$R\Gamma_c(\Omega_K^{d-1,ca},\o\QM_l)$ est scindable.
\begin{prop}
  Il existe un isomorphisme dans $D^b(\o\QM_lPG_d)$
$$ R\Gamma_c(\Omega_K^{d-1,ca},\o\QM_l) \simto \bigoplus_{i=0}^{d-1}
\pi_{\{1,\cdots,i\}}[-d+1-i]$$ 
qui induit en cohomologie les isomorphismes de Schneider-Stuhler \ref{ss}.
\end{prop}
Empressons-nous de pr\'eciser qu'un tel isomorphisme est loin d'\^etre unique, mais nous en exhiberons certains meilleurs que les autres.
Nous allons donner deux preuves ind{\'e}pendantes de ce fait. La
premi{\`e}re repose sur un crit{\`e}re g{\'e}n{\'e}ral de scindage et sur un
calcul explicite des groupes d'extensions entre les diverses repr{\'e}sentations
$\pi_I$, $I\subseteq S$ : on obtient un argument purement alg{\'e}brique qui montre
que {\em tout} complexe de repr{\'e}sentations lisses de $PG_d$ ayant pour
cohomologie les $H^i_c(\Omega_K^{d-1,ca},\o\QM_l)$ est scindable. La deuxi{\`e}me est
d'origine plus g{\'e}om{\'e}trique car on utilise l'action d'un
rel{\`e}vement de Frobenius sur $R\Gamma_c(\Omega_K^{d-1,ca},\o\QM_l)$ ; elle donne plus d'informations qui seront utiles par la suite.

\medskip

{\noindent \em Preuve alg{\'e}brique :}
Nous utiliserons le crit{\`e}re suivant permettant de scinder un complexe
cohomologiquement born{\'e} $\CC^\bullet$ d'une cat{\'e}gorie d{\'e}riv{\'e}e assez
g{\'e}n{\'e}rale. La preuve est facile et est donn{\'e}e en
\ref{sectionscin} : c'est le corollaire \ref{scindage} sp{\'e}cialis{\'e}
dans notre situation.

\medskip

\noindent{\em Si pour tout couple $(i,j) \in \ZM^2$, on a
  $\ext{i-j-1}{\HC^i(\CC^\bullet)}{\HC^j(\CC^\bullet)}{}=0$, alors
  $\CC^\bullet$ est scindable.}

\medskip

Pour appliquer ce crit{\`e}re, nous utilisons le calcul des
extensions entre repr{\'e}sentations du type $\pi_I$,  qui est
l'objet de la partie \ref{extensions}, {\em cf} le th{\'e}or{\`e}me \ref{theoext}.

\medskip

\noindent{\em Pour tous $I,J\subseteq S$, on a 
\ini\begin{equation}\label{ee}
\cas{\ext{*}{\pi_I}{\pi_J}{PG_d}}{\o\QM_l}{*=\delta(I,J)}{0}{*\neq
      \delta(I,J)},
\end{equation}
o{\`u} $\delta(I,J)=|I\cup J|-|I\cap J|$ est le cardinal de la diff{\'e}rence
sym{\'e}trique de $I$ et $J$.}

\medskip

Gr{\^a}ce {\`a} ce calcul et {\`a} la description de Schneider-Stuhler, on
obtient alors pour $i,j\in\{0,\cdots,d-1\}$
$$\cas{\ext{*}{H^{d-1+i}_c(\Omega_K^{d-1,ca},\o\QM_l)}{H^{d-1+j}_c(\Omega_K^{d-1,ca},\o\QM_l)}{PG_d}}
{\o\QM_l}{*=|i-j|}{0}{*\neq
  |i-j|}$$
ce qui implique en particulier que le crit{\`e}re de scindage ci-dessus
est v{\'e}rifi{\'e}.
\findem

\noindent{\em Preuve g{\'e}om{\'e}trique :} Soit $\phi$ un rel{\`e}vement
de Frobenius g{\'e}om{\'e}trique.
La description de
Schneider-Stuhler montre que pour $i=0,\cdots, d-1$, l'endomorphisme
$\HC^{d-1+i}(\gamma(\phi))$ induit par $\gamma(\phi)$ sur
$H^{d-1+i}_c(\Omega_K^{d-1,ca},\o\QM_l)$ est la multiplication par $q^{i}$. On
peut alors encore utiliser un r{\'e}sultat tr{\`e}s g{\'e}n{\'e}ral sur les
complexes donn{\'e} par le lemme \ref{actionphi}. Ce r{\'e}sultat assure
que 
\begin{enumerate}
\item $\gamma(\phi)$ est annul{\'e} par le polyn{\^o}me $\prod_{i=0}^{d-1}
  (X-q^{i})$ et par cons{\'e}quent est un endomorphisme {\em
    semi-simple} de $R\Gamma_c(\Omega_K^{d-1,ca},\o\QM_l)$.
\item Il existe un {\em unique} isomorphisme comme dans la proposition
\ini\begin{equation}\label{scfrob}
 \alpha_\phi:\;\; R\Gamma_c(\Omega_K^{d-1,ca},\o\QM_l) \simto \bigoplus_{i=0}^{d-1}
\pi_{\{1,\cdots,i\}}[-d+1-i]
\end{equation}
 tel que de plus, l'endomorphisme
$\alpha_\phi^{-1}\gamma(\phi)\alpha_\phi$  de l'objet $\bigoplus_{i=0}^{d-1}
\pi_{\{1,\cdots,i\}}[-d+1-i]$ soit 
donn{\'e} par la muliplication par $q^{i}$ sur chaque $\pi_{\{1,\cdots,i\}}[-d+1-i]$.
\end{enumerate}
\findem

Un isomorphisme $\alpha$ comme dans la proposition ci-dessus sera
appel{\'e} un {\em scindage}\footnote{La terminologie 
adopt{\'e}e ici co\"incide bien avec celle de la section
  \ref{sectionscin}, malgr\'e la formulation l\'eg\`erement diff\'erente.} de
$R\Gamma_c(\Omega_K^{d-1,ca},\o\QM_l)$.
   Tout
scindage induit en particulier un isomorphisme
$$ \alpha_* :\; \endo{D^b(PG_d)}{R\Gamma_c(\Omega_K^{d-1,ca},\o\QM_l)} \simto
\endo{D^b(PG_d)}{\bigoplus_{i=0}^{d-1} \pi_{\{1,\cdots, i\}}[-i]}.$$
Il se trouve qu'on peut d{\'e}crire tr{\`e}s explicitement le membre de
droite. On peut tout d'abord le r{\'e}{\'e}crire
$$ \endo{D^b(PG_d)}{\bigoplus_{i=0}^{d-1} \pi_{\{1,\cdots, i\}}[-i]}
 = \bigoplus_{0\leq i\leq j\leq d-1}
 \ext{j-i}{\pi_{\{1,\cdots,j\}}}{\pi_{\{1,\cdots,i\}}}{PG_d} $$
le produit sur le terme de droite {\'e}tant le $\cup$-produit. Or, on peut
d{\'e}crire explicitement le $\cup$-produit pour trois repr{\'e}sentations du type $\pi_I$
de la mani{\`e}re suivante, voir la preuve dans la partie \ref{extensions} :

\medskip

\noindent{\em Soient $I,J,K$ trois sous-ensembles de $S$ tels que
  $\delta(I,J)+\delta(J,K)=\delta(I,K)$, alors le cup-produit
\ini\begin{equation}\label{eec}
\application{\cup \;:\;\;}
{\ext{\delta(I,J)}{\pi_I}{\pi_J}{PG_d} \otimes_{\o\QM_l}
\ext{\delta(J,K)}{\pi_J}{\pi_K}{PG_d}} 
{\ext{\delta(I,K)}{\pi_I}{\pi_K}{PG_d}}{\alpha\otimes \beta}{\beta\cup\alpha}
\end{equation}
 est un isomorphisme.}

\medskip

Choisissons alors pour $i=0,\cdots,d-1$ des g{\'e}n{\'e}rateurs
\ini\begin{equation}\label{generateurs}
\beta_{i,i+1} \in \ext{1}{\pi_{\{1,\cdots, i+1\}}}{\pi_{\{1,\cdots,
    i\}}}{PG_d}
\end{equation}
  et posons pour tous $0\leq i < j\leq d-1$
$$\beta_{i,j}:=  \beta_{i,i+1}\cup\cdots\cup\beta_{j-1,j} \in
\ext{j-i}{\pi_{\{1,\cdots,j\}}}{\pi_{\{1,\cdots,i\}}}{PG_d}. $$
C'est donc un g{\'e}n{\'e}rateur du $\o\QM_l$-espace vectoriel de droite. 
Enfin pour $i=0,\cdots, d-1$, soit $\beta_{ii}$ l'{\'e}l{\'e}ment unit{\'e} de
l'anneau $\ext{0}{\pi_{\{1,\cdots,i\}}}{\pi_{\{1,\cdots,i\}}}{PG_d}$.
Notons maintenant $\TC_d$ la $\o\QM_l$-alg{\`e}bre des matrices $d\times d$
triangulaires sup{\'e}rieures, dont nous notons $(E_{ij})_{0\leq i\leq
  j\leq d-1}$ la base ``canonique''. Nous avons tout fait pour que l'application
$\o\QM_l$-lin{\'e}aire 
\ini\begin{equation}\label{beta}
\application{\beta:\;\;}{\TC_d}{  \bigoplus_{0\leq i\leq j\leq d-1}
 \ext{j-i}{\pi_{\{1,\cdots,j\}}}{\pi_{\{1,\cdots,i\}}}{PG_d}}{E_{ij}}{\beta_{ij}}
 \end{equation}
soit un isomorphisme d'alg{\`e}bres.
Celui-ci d\'epend du choix de g\'en\'erateurs \ref{generateurs} ; un autre choix reviendrait \`a conjuguer par une matrice diagonale.

Tout scindage $\alpha$ de $R\Gamma_c(\Omega_K^{d-1,ca},\o\QM_l)$ et tout choix de g\'en\'erateurs \ref{generateurs} induisent donc un
isomorphisme 
$$ \beta^{-1}\alpha_*:\;\; 
\endo{D^b(PG_d)}{R\Gamma_c(\Omega_K^{d-1,ca},\o\QM_l)} \simto \TC_d. $$

On veut maintenant expliciter l'action  de $W_K$ sur $R\Gamma_c(\Omega_K^{d-1,ca},\o\QM_l)$ au
moyen d'un tel isomorphisme. 
Remarquons  que les repr{\'e}sentations $\tau_I^{\phi,\mu} :W_K \To{}
\endo{\o\QM_l}{\o\QM_l^d}$ que nous avons explicitement d{\'e}finies
en \ref{deftauI} se factorisent par des morphismes $\tau_I^{\phi,\mu} :W_K\To{}
\TC_d$.

La deuxi{\`e}me {\'e}tape importante vers le th{\'e}or{\`e}me \ref{theodp} est 
\begin{prop} \label{propdp2}
Pour tout rel\`evement de Frobenius g\'eom\'etrique $\phi$,
il existe un unique scindage $\alpha_\phi$ de $R\Gamma_c(\Omega_K^{d-1,ca},\o\QM_l)$ et un unique choix de g\'en\'erateurs \ref{generateurs} d\'efinissant un isomorphisme $\beta_{\phi,\mu}$ comme en \ref{beta} tels
que le diagramme suivant soit commutatif 
$$\xymatrix{   \endo{D^b(PG_d)}{R\Gamma_c(\Omega_K^{d-1,ca},\o\QM_l)}
  \ar[r]^-{\beta_{\phi,\mu}^{-1}\alpha_{\phi*}} & \TC_d \\ W_K \ar[u]^{\gamma}
  \ar[ur]_{\tau_\emptyset^{\phi,\mu}} & }
$$
De plus, le choix de g\'en\'erateurs (et donc $\beta_{\phi,\mu}$) est en fait ind\'ependant de $\phi$.
\end{prop}
Rappelons ici que $\tau_\emptyset^{\phi,\mu}$ correspond {\`a} la repr{\'e}sentation
``sp{\'e}ciale'' (ind{\'e}composable) de dimension $d$ de $W_K$. En
particulier cette proposition montre que l'inertie n'agit pas
trivialement sur $R\Gamma_c(\Omega_K^{d-1,ca},\o\QM_l)$ bien qu'elle agisse
trivialement sur les groupes de cohomologie.

La preuve de cette proposition occupe le reste de cette section.

\alin{D{\'e}finition du $N$} \label{definitionN} 
Notons $\NC$ le noyau du  morphisme
d'alg{\`e}bres canonique
$$  \endo{D^b(PG_d)}{R\Gamma_c(\Omega_K^{d-1,ca},\o\QM_l)} \To{can}
\prod_{i=0}^{d-1} \endo{PG_d}{H^{d-1+i}_c(\Omega_K^{d-1,ca},\o\QM_l)} $$
Le point i) du lemme \ref{actionphi} montre dans un contexte
g{\'e}n{\'e}ral que $\NC$ est form{\'e}
d'{\'e}l{\'e}ments nilpotents d'ordre $\leq d$, mais dans le cas
pr{\'e}sent c'est aussi une cons{\'e}quence claire de la description
explicite de $\endo{D^b(PG_d)}{R\Gamma_c(\Omega_K^{d-1,ca},\o\QM_l)}$ donn{\'e}e
ci-dessus.

 Comme l'action de $W_K$ sur
les groupes de cohomologie provient de celle sur le complexe
$R\Gamma_c$, la description de Schneider et Stuhler implique en
particulier que 
$ \gamma(I_K) \subset 1 + \NC$. Un scindage $\alpha$ et un choix de g\'en\'erateurs $\beta$ arbitraires induisent
un isomorphisme $\beta^{-1}\alpha_*$ du
groupe $1+\NC$ sur le groupe multiplicatif des matrices unipotentes
dans $\TC_d$. On en d{\'e}duit une structure de $\o\QM_l$-groupe de Lie
unipotent sur $1+\NC$, structure qui ne d{\'e}pend pas de ces  choix.
On voudrait montrer que le morphisme $\gamma:\;I_K \To{} 1+\NC$ est
{\em continu}. Pour cela remarquons qu'il se factorise par
$$ I_K \To{}  \endo{D^b(PG_d)}{R\Gamma_c(\Omega_K^{d-1,ca},\ZM_l)} \To{can} 
 \endo{D^b(PG_d)}{R\Gamma_c(\Omega_K^{d-1,ca},\o\QM_l)} $$
o{\`u} le deuxi{\`e}me morphisme est obtenu par fonctorialit{\'e} via le foncteur
$-\otimes_{\ZM_l} \o\QM_l$.
D'apr\`es \ref{cohodp} (et moralement d'apr\`es Schneider-Stuhler), on a :
$$ H^{d-1+i}_c(\Omega_K^{d-1,ca},\ZM_l) = \CC^\infty(G_d/P_{\{1,\cdots,
    i\}},\ZM_l)/\sum_{I\supset \{1,\cdots i\}} \CC^\infty(G_d/P_I,\ZM_l),$$
en notant $\CC^\infty(-,\ZM_l)$ les fonctions localement constantes {\`a}
valeurs dans $\ZM_l$. Par \cite[Thm 6.8]{SS1}, on sait que ces repr\'esentations de $PG_d$ admettent des r\'esolutions par des repr\'esentations induites \`a supports compacts de $\ZM_l$-repr\'esentations finies de sous-groupes ouverts compacts. Comme par ailleurs ces repr\'esentations sont $\ZM_l$-admissibles, 
on en d\'eduit ({\em cf} preuve du lemme \ref{finit} pour les d\'etails) que le $\ZM_l$-module
$$ \endo{D^b(PG_d)}{R\Gamma_c(\Omega_K^{d-1,ca},\ZM_l)}$$ est de type fini.

Revenons au morphisme $\gamma :\; I_K\To{} 1+\NC$. On vient de voir
que son image est incluse dans un sous-groupe profini
du groupe de Lie $1+\NC$. Un tel sous-groupe est n\'ecessairement pro-$l$ et
$\gamma$ est  automatiquement continu ; il s'ensuit qu'il se
factorise par le quotient $l$-adique de $I_K$ :
$$\gamma:\;\; I_K \To{t_\mu} \ZM_l \To{\o\gamma_\mu} 1+\NC.$$
Posons alors $N_\mu:=\hbox{log}(\o\gamma_\mu(1))\in \NC$, la densit{\'e} de $\ZM$ dans
$\ZM_l$ montre que
$$\forall i\in I_K,\;\; \gamma(i) = \hbox{exp}(N_\mu t_\mu(i)).$$
Si $w\in W$, en {\'e}crivant la formule ci-dessus pour $wiw^{-1}$, on
constate que $N_\mu$ satisfait l'{\'e}quation
\ini
\begin{equation} \label{eqN}
 \gamma(w)N_\mu\gamma(w)^{-1} = q^{-\nu(w)} N_\mu = |w|N_\mu. 
\end{equation}

\alin{Diagonalisation du Frobenius et choix du scindage} \label{diag}
La proposition que l'on veut d\'emontrer comporte une assertion d'unicit\'e. En ce qui concerne le scindage, nous n'avons pas d'autre choix que de prendre celui de \ref{scfrob} puisque c'est l'unique scindage qui "diagonalise" $\phi$, au sens o\`u pour tout choix de $\beta$,
la matrice $\beta^{-1}\alpha_{\phi*}(\phi)$ est diagonale
d'{\'e}l{\'e}ments diagonaux  $\beta^{-1}\alpha_{\phi*}(\phi)_{ii}=q^{i}$ pour
$i=0,\cdots, d-1$.


Il r{\'e}sulte alors imm{\'e}diatement de l'{\'e}quation \ref{eqN} que l'image de $N_\mu$ est de
la forme $\beta^{-1}\alpha_{\phi*}(N_\mu)=\sum_{i=1}^{d-1} a_i E_{i-1,i}$. En
conjuguant par une matrice diagonale, ce qui  revient {\`a} changer
$\beta$, on peut de plus supposer que les $a_i$
sont nuls ou {\'e}gaux {\`a} $1$. 

En r{\'e}sum{\'e}, on a donc trouv{\'e} le scindage $\alpha_\phi$ et un $\beta_{\phi,\mu}$ 
tels que pour un certain ensemble
$I\subseteq S$, le diagramme
$$\xymatrix{   \endo{D^b(PG_d)}{R\Gamma_c(\Omega_K^{d-1,ca},\o\QM_l)}
  \ar[r]^-{\beta_{\phi,\mu}^{-1}\alpha_{\phi*}} & \TC_d \\ W_K \ar[u]^{\gamma}
  \ar[ur]_{\tau^{\phi,\mu}_I} & }
$$
est commutatif. Nous allons maintenant  prouver que $I=\emptyset$. Il
suffit pour cela de montrer que l'ordre de $N_\mu$ est $d$.

\alin{Uniformisation et puret{\'e}} Soit $\Gamma$ un sous-groupe discret
sans-torsion et {\em cocompact} dans $PG_d$. On sait depuis Mustafin
\cite{Mustafin} que le quotient  
$\Omega_K^{d-1}/\Gamma$ que nous avons d{\'e}ja {\'e}voqu{\'e} plus haut s'identifie dans ce cas
 {\`a} l'analytification d'une vari{\'e}t{\'e} propre et lisse sur
 $K$. Par les th{\'e}or{\`e}mes de type GAGA de Berkovich \cite[7.1]{Bic2}, les
 groupes de cohomologie $H^i(\Omega_K^{d-1,ca}/\Gamma,\o\QM_l)$ sont donc des
repr{\'e}sentations de dimension finie $l$-adiques  continues de $W_K$
et sont par cons{\'e}quent 
munis d'une filtration de monodromie et d'un op{\'e}rateur nilpotent
correspondant (par le th{\'e}or{\`e}me de Grothendieck). 
Ces op{\'e}rateurs sont fonctoriellement induits par l'op{\'e}rateur
$N_\mu$ que nous avons d{\'e}fini sur $R\Gamma_c(\Omega_K^{d-1,ca},\o\QM_l)$, via les
isomorphismes $W_K$-{\'e}quivariants
$$ H^p(\Omega_K^{d-1,ca}/\Gamma,\o\QM_l) \simeq \HC^p\left(\o\QM_l
  \otimes^L_{\o\QM_l[\Gamma]} R\Gamma_c(\Omega_K^{d-1,ca},\o\QM_l) \right)$$
que l'on d{\'e}duit de \ref{HS}. Si $\alpha$ est un scindage de
$R\Gamma_c$, il induit des isomorphismes de $\o\QM_l$-espaces vectoriels
\begin{eqnarray*} 
H^p(\Omega_K^{d-1,ca}/\Gamma,\o\QM_l) & \simto  & \bigoplus_{i=0}^{d-1}
\tor{p-d+1-i}{\o\QM_l}{\pi_{\{1,\cdots,i\}}}{\Gamma} \\
& \simeq & \bigoplus_{i=0}^{d-1} \ext{d-1+i-p}{\pi_{\{1,\cdots,i\}}}{\o\QM_l}{\Gamma}^* \\
& \simeq & \bigoplus_{i=0}^{d-1}
\ext{d-1+i-p}{\pi_{\{1,\cdots,i\}}}{\CC^\infty(G_d/\Gamma,\o\QM_l)}{PG_d}^* 
\end{eqnarray*}
On a bien-s{\^u}r utilis{\'e} et d{\'e}riv{\'e} l'isomorphisme
$$\hom{V}{\o\QM_l}{\Gamma}\simeq (V\otimes_{\o\QM_l[\Gamma]}\o\QM_l)^*$$
puis on l'a dualis{\'e} gr{\^a}ce {\`a} la finitude des dimensions de tous
les espaces en jeu, et on a utilis{\'e} le lemme de Shapiro.
Si on a pris soin d'utiliser le scindage $\alpha_\phi$ de \ref{scfrob}, alors l'action de $\phi$ sur le facteur associ{\'e} {\`a} l'indice $i$ est la
multiplication par $q^{i}$. En particulier, supposons que
$\pi_\emptyset \injo \CC^\infty(PG_d/\Gamma,\o\QM_l)$, ce qui est v{\'e}rifi{\'e}
au moins pour $\Gamma$ ``assez petit'' (peut-{\^e}tre pour tout $\Gamma$
?). Alors la formule  \ref{ee} montre que lorsque $p=d-1$, chacun des termes
$$
\ext{d-1+i-p}{\pi_{\{1,\cdots,i\}}}{\CC^\infty(PG_d/\Gamma,\o\QM_l)}{PG_d}^*
$$ est non-nul. En particulier, on a donc
\ini\begin{equation}\label{gr}
 \hbox{Gr}_W^{2d-2}\left( H^{d-1}(\Omega_K^{d-1,ca},\o\QM_l)\right) \neq 0
\end{equation}
en notant $\hbox{Gr}_W$ le gradu\'e pour la filtration par les poids (dont la d\'efinition est rappel\'ee au paragraphe \ref{monpds}).

 \`A ce stade, nous pouvons conclure la preuve de $N_\mu^{d-1}\neq 0$ en invoquant le th\'eor\`eme suivant
 \begin{theo}
 (Ito, \cite{Ito}) la conjecture monodromie poids ({\em cf} \ref{conmp}) est vraie pour les repr\'esentations $l$-adiques $H^i(\Omega_K^{d-1,ca}/\Gamma,\o\QM_l)$.
 \end{theo}
  En effet, par la forme \'equivalente \ref{equiMP} de cette conjecture et \ref{gr} ci-dessus, on obtient que la puissance $d-1$-\`eme de l'op{\'e}rateur de monodromie de $H^{d-1}(\Omega_K^{d-1,ca}/\Gamma,\o\QM_l)$ est non-nulle, donc {\em a fortiori} $N_\mu^{d-1}\neq 0$.

\bigskip

Mais on peut aussi se passer d'invoquer le th\'eor\`eme d'Ito, car la propri\'et\'e que l'on veut prouver peut se d\'eduire "simplement" de la suite spectrale de Rapoport-Zink et ne n\'ecessite pas l'\'etude plus difficile de cette suite spectrale dans \cite{Ito}. 
L'argument qu'on va utiliser n'est pas nouveau ; il a \'et\'e utilis\'e par Mokrane dans \cite{Mokrane} et \'etait certainement connu de Rapoport-Zink. Il m'a \'et\'e signal\'e par Ngo Bao Chau que je remercie vivement. Commen{\c c}ons par quelques "rappels".

\alin{Rappel sur la suite spectrale de Rapoport-Zink} \label{rappelRZ}
Soit $\XG$ un $\OC_K$-sch\'ema propre et r\'egulier de dimension relative $n$ tel que $X:=\XG\times_{\OC_K} K$ soit lisse et $Y:=\XG\otimes_{\OC_K} k$ soit un diviseur \`a croisements normaux dans $\XG$, somme globale de diviseurs  $Y_i$ lisses sur $k$, pour $i\in I$ ensemble fini.
Suivant Rapoport-Zink \cite{RZss}, notons 
$$Y^{(m)} := \bigsqcup_{J\subseteq I, |J|=m} \left( \bigcap_{i\in J} Y_i\right). $$
$Y^{(m)}$ est une somme de $k$-vari\'et\'es lisses de dimension $n+1-m$.
D'apr\`es \cite[Satz 2.10]{RZss}, il existe une suite spectrale $W_K$-\'equivariante 
$$ E_1^{-r,q+r}= \bigoplus_{k\geq \hbox{sup}(0,-r)} H^{q-r-2k}(Y^{(r+2k+1),ca},\o\QM_l)\otimes |-|^{-r-k} \Rightarrow H^q(X^{ca},\o\QM_l)$$
telle que
\begin{enumerate}
        \item les diff\'erentielles $d_1^{-r,q+r}$ sont des sommes altern\'ees de morphismes de restriction et de morphismes de Gysin, voir ci-dessous un \'enonc\'e plus pr\'ecis dans les cas $r=\pm n$ qui nous int\'eresseront.
        \item le terme $E_1^{-r,q+r}$ est pur de poids $q+r$ et donc la suite spectrale d\'eg\'en\`ere fortement en $E_2$, {\em i.e.} $d_2=0$.
        \item (\cite[Satz 1.11]{RZ}) l'action de $i\in I_K$ sur l'aboutissement est donn\'ee par $(1+\nu)^{t_\mu(i)}$ o\`u $\nu$ est l'endomorphisme de la suite spectrale d\'efini par
 $$ \nu_1^{-r,q+r} =\sum_{k \geq \hbox{sup}(0,-r+1)}  \id_{H^{q-r-2k}(Y^{(r+2k+1),ca},\o\QM_l))}. $$
En particulier, puisque $\nu^{n+1}=0$,
 si $N$ d\'esigne   l'op\'erateur de monodromie de $H^q(X^{ca},\o\QM_l)$, alors 
 on a $N^n=\nu^n$ sur le gradu\'e de la filtration d'aboutissement.

\end{enumerate}

\def\Gr{\hbox{Gr}}

\begin{lemme} \label{Nnonnul} Avec les notations pr\'ec\'edents, les assertions suivantes sont \'equivalentes :
\begin{enumerate}
        \item $N^n \neq 0$ sur $H^n(X^{ca},\o\QM_l)$,
        \item $E_2^{-n,2n}=\ker(d_1^{-n,2n})\neq 0$,
        \item $\Gr^{2n}_W(H^n(X^{ca},\o\QM_l))\neq 0$.
\end{enumerate}
\end{lemme}
\begin{proof}
Comme la suite spectrale d\'eg\'en\`ere en $E_2$, et que $E_1^{-n,2n}$ est le seul terme de poids $2n$, les points ii) et iii) sont \'equivalents. De plus tous les $E_1^{-r,q+r}$ sont de poids compris entre $0$ et $2n$ et par cons\'equent $N^n$ est non nul sur $H^n(X^{ca},\o\QM_l)$ seulement si  $\Gr^{2n}_W(H^q(X^{ca},\o\QM_l))$ est non-nul, c'est-\`a-dire  $i) \Rightarrow iii)$.

Prouvons maintenant que $ii) \Rightarrow i)$. Sous l'hypoth\`ese ii), on a en particulier $E_1^{-n,2n}\neq 0$ et donc
 le $0$-squelette $Y^{(n+1)}$ est non vide.
Par d\'efinition, la diff\'erentielle $d_1^{-n,2n}$ et l'op\'erateur $\nu$ s'ins\`erent dans le diagramme suivant :
$$\xymatrix{E_1^{-n,2n} = H^0(Y^{(n+1)},\o\QM_l) \ar[d]_{d_1^{-n,2n}}  \ar[r]^{\nu^n=\id} &  H^0(Y^{(n+1)},\o\QM_l) = E_1^{n,0} \\ E_1^{-n+1,2n} =H^2(Y^{(n)},\o\QM_l) & H^0(Y^{(n)},\o\QM_l)=E_1^{n-1,0} \ar[u]_{d_1^{n-1,0} }
}.$$
D'apr\`es les rappels ii) et iii) au-dessus du lemme, pour montrer que $N^n \neq 0$ sur l'aboutissement, il suffit de prouver que l'application induite par $\nu^n$ sur le terme $E_2$
$$ E_2^{-n,2n} = \ker (d_1^{-n,2n}) \To{} E_2^{n,0}= \coker(d_1^{n-1,0}) $$
est non nulle, ou encore que $\ker(d_1^{-n,2n})\subsetneq \im(d_1^{n-1,0})$.
Or le morphisme $d_1^{-n,2n}$ est une certaine somme altern\'ee (d\'ependant de la combinatoire d'intersection des $Y_i$) de morphismes de Gysin, tandis que $d_1^{n-1,0}$ est une certaine somme altern\'ee de morphismes de restriction. Par d\'efinition, ces deux morphismes sont {\em adjoints} pour la dualit\'e de Poincar\'e entre $H^2(Y^{(n)},\o\QM_l)$ et $H^0(Y^{(n)},\o\QM_l)$, resp. la forme bilin\'eaire "de Poincar\'e" sur  $H^0(Y^{(n+1)},\o\QM_l)$. Pour montrer que $N^n \neq 0$ sur l'aboutissement, il suffira donc de prouver que $\ker(d_1^{-n,2n})\subsetneq \ker(d_1^{-n,2n})^\perp$. 

Or, les trois espaces $H^0(Y^{(n)},\o\QM_l)$, $H^2(Y^{(n)},\o\QM_l)$ et $H^0(Y^{(n+1)},\o\QM_l)$ admettent des bases "canoniques" index\'ees par les  composantes connexes de $Y^{(n)}$, resp. $Y^{(n+1)}$, et ont donc chacun une $\QM$-structure "canonique"
que nous noterons "$H^i(Y^{(n+j)},\QM)$". De plus, le morphisme $d_1^{-n,2n}$  par d\'efinition respecte ces $\QM$-structures, de sorte que $\ker(d_1^{-n,2n})$ provient d'un $\QM$-sous-espace de "$H^0(Y^{(n+1)},\QM)$".
Enfin, la dualit\'e de Poincar\'e sur $H^0(Y^{(n+1)},\o\QM_l)$ provient par d\'efinition d'un accouplement sur "$H^0(Y^{(m-1)},\QM)$" qui y est visiblement d\'efini positif.

On en d\'eduit que pour montrer que $N^n \neq 0$ sur l'aboutissement, il suffit
de prouver que $\ker(d_1^{-n,2n})\neq 0$. C'est-\`a-dire $ii)\Rightarrow i)$.
\end{proof}

Revenons maintenant \`a la vari\'et\'e $X:=\Omega_K^{d-1}/\Gamma$ pour $\Gamma$ discret cocompact tel que la repr\'esentation de Steinberg apparaisse dans $\CC^\infty(G/\Gamma)$. 
Soit $\wh\Omega_K^{d-1}$ le mod\`ele formel semi-stable (d\'efini par Deligne, {\em cf} \cite{BouCar} par exemple) de $\Omega_K^{d-1}$. On sait que le sch\'ema formel quotient
$\wh\Omega_K^{d-1}/\Gamma$ s'alg\'ebrise en un $\OC_K$-sch\'ema $\XG$ et que
quitte \`a rapetisser un peu $\Gamma$ pour qu'il agisse "tr\`es" librement sur l'immeuble,  
$\XG$ satisfait les conditions de \ref{rappelRZ} avec $n=d-1$.
Mais d'apr\`es \ref{gr}, le lemme pr\'ec\'edent montre que $N_\mu^{d-1}\neq 0$ et cl\^ot la preuve de la proposition \ref{propdp2}.

Il n'est pas difficile alors d'en tirer une nouvelle preuve du th\'eor\`eme d'Ito
\begin{coro} \label{MPdp}
La conjecture monodromie poids ({\em cf} \ref{conmp}) est vraie pour les repr\'esentations $l$-adiques $H^i(\Omega_K^{d-1,ca}/\Gamma,\o\QM_l)$.
\end{coro}
Il suffit pour cela d'appliquer la proposition \ref{propequMP} qui donne les arguments n\'ecessaires dans le cas g\'en\'eral.



\alin{Unicit\'e du choix de g\'en\'erateurs de \ref{generateurs}} \label{indepbeta}
La discussion qui pr\'ec\`ede nous a permis de trouver un $\beta_{\phi,\mu}$, associ\'e \`a un certain choix de g\'en\'erateurs,  rendant commutatif le diagramme de la proposition \ref{propdp2}. Changer ce choix revient \`a conjuguer $\beta_{\phi,\mu}$ par une matrice diagonale non centrale. Or, $\beta_{\phi,\mu}^{-1}\alpha_{\phi*}(N_\mu)$ est un nilpotent r\'egulier donc son centralisateur dans le groupe des matrices diagonales est justement le sous-groupe des matrices
centrales. On en d\'eduit l'unicit\'e de $\beta_{\phi,\mu}$.

Montrons maintenant que $\beta_{\phi,\mu}$ est ind\'ependant de $\phi$. Soit $i\in I_K$, rappelons que par \cite[8.4.2]{DelAntwerp} on a 
$$ \forall w\in W_K,\;\; \tau^{i\phi}_\emptyset(w) = 
\tau^{\phi,\mu}_\emptyset(i)^{\frac{1}{1-q}}\tau^{\phi,\mu}_\emptyset(w)\tau^{\phi,\mu}_\emptyset(i)^{\frac{1}{q-1}}. $$
De m\^eme, en utilisant la propri\'et\'e $\gamma(\phi)\gamma(i)\gamma(\phi)^{-1}=\gamma(i)^q$ et le fait que $\gamma(i)$ est unipotent pour d\'efinir ses puissances rationnelles, on v\'erifie que
$$ \alpha_{i\phi}= \alpha_\phi\circ \gamma(i)^{\frac{1}{1-q}}.$$
\`A partir de ces deux \'equations, un calcul montre que $\beta_{i\phi,\mu}$ fait aussi commuter le diagramme de \ref{propdp2} et par unicit\'e on en conclut que $\beta_{\phi,\mu}
=\beta_{i\phi,\mu}$.

\subsection{Preuve du th{\'e}or{\`e}me} \label{preuvetheodp}
Fixons $I\subseteq S$, et choisissons un rel\`evement de Frobenius $\phi$. 
 Le scindage $\alpha_\phi$ de
la proposition \ref{propdp2}  induit le
premier isomorphisme suivant :
\begin{eqnarray*}
  R\hom{R\Gamma_c(\Omega_K^{d-1,ca},\o\QM_l)}{\pi_I}{D^b(\o\QM_lPG_d)} & \simto &
   \bigoplus_{i=0}^{d-1}
   R\hom{\pi_{\{1,\cdots,i\}}}{\pi_I}{D^b(\o\QM_lPG_d)}[d-1+i] \\
& \simeq & \bigoplus_{i=0}^{d-1}
\ext{\delta(i,I)}{\pi_{\{1,\cdots,i\}}}{\pi_I}{PG_d}[d-1+i-\delta(i,I)]
\end{eqnarray*}
et le second est une cons{\'e}quence de la formule \ref{ee} ;  on y a pos{\'e}
$\delta(i,I):=\delta(\{1,\cdots, i\},I)$. 
Soit $e_i$ un g{\'e}n{\'e}rateur du $\o\QM_l$-espace vectoriel de dimension
$1$ $\ext{\delta(i,I)}{\pi_{\{1,\cdots,i\}}}{\pi_I}{PG_d}$. Les
$(e_i)_{i=0,\cdots,d-1}$ forment une base de l'espace qui nous
int{\'e}resse :
$$ \HC^*_I:=\HC^*\left(R\hom{R\Gamma_c(\Omega_K^{d-1,ca},\o\QM_l)}{\pi_I}{D^b(\o\QM_lPG_d)}\right) \simeq 
   \bigoplus_{i=0}^{d-1} \o\QM_l e_i.$$
Celui-ci a donc au moins la bonne dimension. De plus, l'action de
$\phi$ est donn{\'e}e par
$$\forall i,\;\; \phi.e_i = q^{-i}e_i.$$
L'action de $N_\mu$ ({\`a} droite sur $\HC^*_I$) est donn{\'e}e par
$\cup$-produit par l'{\'e}l{\'e}ment 
$$ \alpha_*(N_\mu)= \sum_{i=0}^{d-1} \beta_{i,i+1} \in \bigoplus_{i=0}^{d-1}
\ext{1}{\pi_{\{1,\cdots,i+1\}}}{\pi_{\{1,\cdots,i\}}}{PG_d} $$
o\`u les $\beta_{i,i+1}$ sont le syst\`eme de g\'en\'erateurs \ref{generateurs} associ\'e \`a l'isomorphisme $\beta_{\phi,\mu}$ de la proposition \ref{propdp2}.
 On a donc pour tout $i=0,\cdots,d-1$
$$  N_\mu.e_i = e_i\cup \alpha_*(N_\mu)=e_i\cup \beta_{i,i+1} \in \o\QM_l
e_{i+1}.$$
La formule \ref{eec} nous montre alors que $N_\mu e_i$ est non nul \ssi\ on
a l'additivit{\'e} des $\delta$ suivante :
$$
\delta(i+1,I)=\delta(i,I)+\delta(\{1,\cdots,i\},\{1,\cdots,i+1\})=
\delta(i,I)+1.$$
Mais cette additivit{\'e} est clairement {\'e}quivalente {\`a} $i+1\in I^c$
(compl{\'e}mentaire de $I$). 
Soit $(E_{ij})_{0\leq i,j\leq d-1}$ la base de $\endo{\o\QM_l}{\HC^*_I}$
associ{\'e}e {\`a} la base $(e_i)_{0\leq i\leq d-1}$. La matrice d{\'e}crivant
l'action de  $N_\mu$ a donc la forme 
$$ \sum_{i\in I^c} a_i E_{i,i-1}$$
o{\`u} les $a_i$ sont non nuls. 
Un changement de base homoth\'etique permet de rendre les $a_i$ \'egaux \`a $1$, puis en renversant l'ordre des vecteurs de la base, on met $N_\mu$ sous la forme $\sum_{i\in\o{I^c}} E_{i-1,i}$. Compte tenu du d\'ecalage par $|-|^{d-1}$, on trouve donc
$$
\HC^*\left(R\hom{R\Gamma_c(\Omega_K^{d-1,ca},\o\QM_l)}{\pi_I}{D^b(\o\QM_lPG_d)}\right)
\simeq \tau_I\otimes |.|^{d-1} \simeq \sigma_d(\pi_I)\otimes |.|^{\frac{d-1}{2}}.
$$

\section{Extensions dans la s\'erie principale}\label{extensions}

\subsection{Notations et premier \'enonc\'e}

\ali 
Soit ${\mathbf G}$ un groupe r{\'e}ductif d{\'e}ploy{\'e} sur un corps local
non-archim{\'e}dien $F$
dont on fixe un tore maximal d{\'e}ploy{\'e} ${\mathbf T}$ et un \borel
${\mathbf B_+}$ de radical unipotent ${\mathbf U_+}$. On notera les groupes de
points $F$-rationnels par les caract{\`e}res non 
{\'e}paissis correspondants $G,T,B_+$. En d\'esignant par $X_*({\mathbf
  T})$ le groupe des cocaract\`eres de ${\mathbf T}$,
on pose $X:=X_*({\mathbf T})\otimes \RM$. 
L'ensemble $S$ des  racines simples de ${\mathbf T}$ dans
$\hbox{Lie}({\mathbf U_+})$ sera parfois consid\'er\'e comme un sous-ensemble du dual $X^*$ de
$X$.  
Le groupe de Weyl correspondant est not{\'e} $W$.
 Comme d'habitude, un \para sera dit standard
s'il contient ${\mathbf B_+}$ et un \levi sera dit standard si c'est la
composante de Levi contenant ${\mathbf T}$ d'un \para standard. Les \levis standards sont en
bijection avec les sous-ensembles de $S$.

\ali \label{nota}
 Pour $J\subseteq S$, on notera
${\mathbf M_J}$ le  \levi standard dont le syst{\`e}me
de racines de ${\mathbf T}$ associ{\'e} 
est engendr{\'e} par $J$. On a donc ${\mathbf M_J}={\mathbf \ZC}_{\mathbf G}(\cap_{\alpha\in J}\ker
\alpha)^0$.  
On notera aussi ${\mathbf A_J}$ le centre connexe de ${\mathbf M_J}$ et
$a_J:=X_*({\mathbf A_J})\otimes \RM$. On remarquera que $a_\emptyset=X$. En
g{\'e}n{\'e}ral, l'injection canonique $X_*({\mathbf A_J})\To{}
X_*({\mathbf T})$ identifie $a_J$ {\`a}
$$\{x\in X, \forall \alpha\in J,\;\; \la x,\alpha\ra =0\}.$$

\alin{Repr{\'e}sentations} \label{defrep}Si $R$ est un anneau commutatif tel
que $p\in R^\times$, on dit qu'un $R$-module $V$ muni d'une action $\pi$ de $G$ est
{\em lisse} si le stabilisateur de tout {\'e}l{\'e}ment est ouvert ; on dit de
plus qu'un tel $RG$-module est {\em
  admissible} si pour tout sous-groupe ouvert compact $H$, le
$R$-module $V^H$ des  $H$-invariants est de type fini. On ommetra
souvent le $V$ de nos notations en d{\'e}signant par $\pi$ {\`a} la fois
la repr{\'e}sentation de $G$ et son $R$-module sous-jacent. La cat{\'e}gorie
ab{\'e}lienne des $RG$-modules lisses sera not{\'e}e $\Mo{R}{G}$. Cette cat\'egorie a suffisamment d'objets projectifs et injectifs. Nous noterons simplement $\ext{*}{-}{-}{RG}$ les groupes d'extensions dans cette cat\'egorie.

Soit $J\subset S$. Notons $P_J$ le groupe parabolique standard engendr{\'e} par
$M_J$ et $B_+$. Pour tout $K\supset J$, on a une injection canonique
$\CC_R^\infty(G/P_K) \injo \CC_R^\infty(G/P_J)$. On pose alors
$$ \pi_J^R:= \CC_R^\infty(G/P_J)/ \left(\sum_{K\supset
    J}\CC_R^\infty(G/P_K)\right).$$  
C'est une $R$-repr{\'e}sentation admissible de $G$. Notre but ici est de
d{\'e}montrer le th{\'e}or{\`e}me suivant, dans lequel on utilise la notation
$\delta(I,J)$ pour le cardinal de la diff{\'e}rence
sym{\'e}trique
$\Delta(I,J):=(I\cup J)\setminus(I\cap J)$ 
 entre deux sous-ensembles $I,J$ de $S$.

\begin{theo} \label{theoext}
Supposons ${\mathbf G}$ semi-simple et soit $R$ un anneau 
  fortement banal pour $G$ au sens de \ref{bonbanal} ci-dessous.
  \begin{enumerate}
  \item Soient $I,J$ deux sous ensembles de $S$, alors :
    $$
    \cas{\ext{*}{\pi^R_I}{\pi^R_J}{RG}}{R}{*=\delta(I,J)}{0}{*\neq
      \delta(I,J)}.$$
\item Soient $I,J,K$ trois sous-ensembles de $S$ tels que
  $\delta(I,J)+\delta(J,K)=\delta(I,K)$, alors le cup-produit
$$  \cup \;:\;\;
\ext{\delta(I,J)}{\pi^R_I}{\pi^R_J}{RG} \otimes_R
\ext{\delta(J,K)}{\pi^R_J}{\pi^R_K}{RG} \To{} 
\ext{\delta(I,K)}{\pi^R_I}{\pi^R_K}{RG}
$$ est un isomorphisme.

\end{enumerate}
\end{theo}


Ce th\'eor\`eme sera pr\'ecis\'e par la proposition \ref{explicite} dans laquelle
 nous construisons de
mani{\`e}re explicite des extensions de Yoneda engendrant les $\hbox{Ext}$
non-triviaux qui sont d\'ecrits ci-dessus. 
Il  g{\'e}n{\'e}ralise des r{\'e}sultats pr\'ec\'edents de Casselman et
Schneider-Stuhler  \cite[Prop 5.9]{SS1}. 
Par ailleurs le point i)  a aussi {\'e}t{\'e} obtenu
ind{\'e}pendamment par S. Orlik via une m{\'e}thode plus directe mais
moins explicite, sous l'hypoth\`ese suppl\'ementaire que $R$ est
auto-injectif et pour des groupes pas n\'ecessairement d\'eploy\'es.

\alin{Bonnes caract{\'e}ristiques} \label{bon}
Suivant une terminologie introduite par Vign{\'e}ras, on dira qu'un nombre premier $l$ est
{\em banal} pour $G$ s'il ne divise pas l'ordre d'un sous-groupe compact de
$G$. On notera aussi
$$N_G:= (\hbox{$p'$-partie du pro-ordre de $G$}) \in \NM^\times. $$
Lorsqu'on consid\`ere des repr\'esentations de $G$ \`a valeurs dans un corps, l'hypoth\`ese "de caract\'eristique banale" nous permettra d'utiliser des r\'esultats de finitude cohomologique et de fid\'elit\'e des foncteurs de Jacquet sur les repr\'esentations de la s\'erie principale.
Mais pour des raisons de calculs d'exposants ({\em cf} surtout le lemme \ref{comb2}), cette hypoth\`ese ne suffira pas pour nos arguments.

 Notons $\rho \in X^*({\mathbf T})$ le
d{\'e}terminant de l'action de ${\mathbf 
  T}$ sur $\hbox{Lie} {\mathbf U_+}$. Vu comme {\'e}l{\'e}ment de $X^*$, c'est
aussi la somme des racines $S$-positives et on peut {\'e}crire
$\rho=\sum_{\alpha\in S} n_\alpha \alpha$ pour des entiers positifs
$n_\alpha$. On dira qu'un nombre premier $l$ est {\em bon} pour $G$ si
le cardinal $q_F$  du corps r{\'e}siduel $k_F$ de $F$ est d'ordre
$\geq \sup_\alpha(n_\alpha)$ modulo $l$, ou ce qui est {\'e}quivalent
si $l$ est 
premier {\`a} l'entier 
$$N_W:=  \prod_{k\leq \sup(n_\alpha)} (q_F^k-1).$$

Pour $GL(n)$, on a $(n-1)(n+1)/4\leq \sup(n_\alpha)\leq n^2/4$, et on voit donc que pour $n\geq 4$, "bon" implique "banal". 
Cependant, en petit rang, il peut exister des $l$ non banals et bons au sens
ci-dessus, par exemple pour $GL(2)$.

\begin{DEf} \label{bonbanal}
 Nous dirons dor{\'e}navant qu'un anneau $R$
est {\em banal}, resp. {\em fortement banal} pour $G$ si l'entier $pN_G$,
resp. l'entier $pN_GN_W$
est inversible dans $R$.
\end{DEf}

Nous n'utilisons dans ce texte que le calcul pour un groupe semi-simple (et m\^eme pour $\mathbf{PGL_n}$).
Mentionnons tout de m\^eme  le calcul pour un groupe r\'eductif (d\'eploy\'e).

\begin{coro}
Supposons ${\mathbf G}$ r{\'e}ductif d\'eploy\'e, de centre de dimension $d$, et $R$ fortement banal 
 pour $G$. Alors pour tous sous-ensembles $I,J$ de $S$, on a
$$ \cas{\ext{*}{\pi^R_I}{\pi^R_J}{G}}{R^{^{({}_r^d)}}}{*=\delta(J,I)+r,\; 0\leq r
   \leq d}{0}{* < \delta(J,I) \hbox{ ou } *> \delta(J,I)+d}$$
\end{coro}
\begin{proof} (\`A partir du cas semi-simple)
Le foncteur d'inflation $\Mo{R}{G/Z}\To{} \Mo{R}{G}$ est exact et admet pour adjoint
{\`a} droite le foncteur des $Z$-invariants
$$\application{}{\Mo{R}{G}}{\Mo{R}{G/Z}}{V}{V^Z \simeq \hom{R}{V}{Z}}.$$
Ce dernier  envoie donc injectifs sur injectifs, et est exact {\`a} gauche. On
a donc la compos{\'e}e de foncteurs d{\'e}riv{\'e}s suivante
$$ R\hom{\pi_I^R}{?}{G}\simeq R\hom{\pi_I^R}{R\hom{R}{?}{Z}}{G/Z}$$
qui fournit une suite spectrale  
$(E_r^{pq})_{r\in \NM}$ de deuxi{\`e}me terme 
$$E_2^{pq}=\ext{p}{\pi_I^R}{\ext{q}{R}{\pi_J^R}{Z}}{G/Z}$$
convergeant vers le gradu{\'e} de
$\ext{p+q}{\pi_I^R}{\pi_J^R}{G}$ pour une certaine filtration.
 
Or, on a un isomorphisme de $R(G/Z)$-modules lisses 
$$ \ext{q}{R}{\pi_J^R}{Z} \simeq  \pi_J^R \otimes \ext{q}{R}{R}{Z},$$
le terme {\`a} la droite du $\otimes$ {\'e}tant muni de l'action triviale de
$G/Z$. 
Par le th{\'e}or{\`e}me pr{\'e}c{\'e}dent, il s'ensuit que seule la colonne
$p=\delta(I,J)$ du terme $E_2^{pq}$ 
est non nulle, et par cons{\'e}quent la suite spectrale d{\'e}g{\'e}n{\`e}re
fortement en un isomorphisme 
$$ E_2^{\delta(I,J)q} = \ext{q}{R}{R}{Z} \simto
\ext{\delta(I,J)+q}{\pi_J^R}{\pi_I^R}{G}. $$
Rappelons maintenant que $Z\simeq \ZM^d\times Z^0$ o{\`u} $Z^0$ est le
sous-groupe compact maximal de $Z$, dont le pro-ordre est suppos{\'e}
inversible dans $R$. Ainsi on a un isomorphisme
$$ \ext{q}{R}{R}{Z} \simto \ext{q}{R}{R}{\ZM^d}. $$
Maintenant un calcul classique (par le complexe de Koszul par exemple, voir
\cite[V.6.4 (ii)]{Brown})
montre que 
$\ext{q}{R}{R}{\ZM^d}$ est un $R$-module libre de dimension $({}_q^d)$
pour $q\leq d$ et est nul pour $q>d$.

\end{proof}

\subsection{Construction explicite d'extensions}

Nous commen{\c c}ons cette section par un peu de combinatoire.

\ali 
Soient $J,K$ deux sous-ensembles de $S$. Comme ci-dessus, on note
  $$\Delta(J,K):=
(J\cup K)\setminus (J\cap K)$$  la 
diff{\'e}rence sym{\'e}trique de $J$ et $K$ dans $S$ et $\delta(J,K)$ son
cardinal. 
Par exemple, $\Delta(J,S)=J^c$, le compl\'ementaire de $J$ dans $S$.
On voit facilement que
$$\Delta(J,\Delta(J,K))=K \;\hbox{ et }\; \Delta(J^c,K)=\Delta(J,K)^c,$$
et non moins facilement que
\begin{lemme} \label{comb1}
  Soient $I,J,K$ trois sous-ensembles de $S$. Les propri\'et\'es suivantes sont \'equivalentes :
  \begin{enumerate}
  \item $\Delta(I,J)\supseteq \Delta(J,K)$
  \item $\Delta(I,J)=\Delta(J,K)\sqcup\Delta(K,I)$
  \item $\delta(I,J)=\delta(J,K)+\delta(K,I)$
  \end{enumerate}
\end{lemme}
\begin{proof}
  Notons $\chi_?:\;S\To{}\ZM$ la fonction caract\'eristique d'un sous-ensemble $?$
  de $S$. On a $\chi_{\Delta(I,J)}=\chi_I+\chi_J-2\chi_{I\cap
    J}$. Deux calculs imm\'ediats montrent que chacun des deux
  premiers points est \'equivalent \`a l'\'egalit\'e
$$ \chi_K-\chi_{K\cap I}-\chi_{K\cap J} +\chi_{I\cap J} =0. $$
Notons maintenant pour une fonction $f:\; S\To{} \ZM$ sa ``somme'' par
$\int_S f:=\sum_{s\in S} f(s)$. Alors $\delta(I,J)=\int_S
\chi_{\Delta(I,J)}$. Un nouveau calcul imm\'ediat montre alors que le
point iii) est \'equivalent \`a l'\'egalit\'e
$$ \int_S \left(\chi_K-\chi_{K\cap I}-\chi_{K\cap J} +\chi_{I\cap
    J}\right) =0. $$
Mais la fonction somm\'ee est \`a valeurs $\geq 0$, donc la nullit\'e
de sa somme \'equivaut \`a sa nullit\'e.
\end{proof}

\ali
Soit $J\subseteq S$. On lui associe une fonction signe 
$$\application{\epsilon_J\;:}{S}{\{\pm 1\}}{\alpha}{-1 \;\hbox{\ssi\ }\;
   \alpha\in J}$$
et un c\^one ouvert de $X$ 
 $$X_J:=\{x\in X,\; \forall \alpha\in S, \;\epsilon_J(x)\langle
 x,\alpha \rangle > 0\}.$$

Les $X_J$ sont les composantes connexes de $X$ priv{\'e} des hyperplans
orthogonaux aux racines {\em simples}.
En particulier,  $X_\emptyset$ est la chambre de Weyl associ{\'e}e {\`a}
$B_+$ et $X_S$ est celle associ{\'e}e au Borel oppos{\'e} {\`a} $B_+$. En
g{\'e}n{\'e}ral, une chambre de Weyl est soit contenue dans soit
disjointe d'un  c{\^o}ne $X_J$.

 \label{1.3}

Si maintenant $J$ et $K$ sont deux
sous-ensembles de $S$, et si l'on convient de surligner
les adh{\'e}rences pour la topologie r{\'e}elle de $X$, on
constate que
$\overline{X_K}\cap \overline{X_J}$ est un c{\^o}ne ferm{\'e}
vectoriellement g{\'e}n{\'e}rateur (ou ce qui est {\'e}quivalent,
contenant un sous-ensemble ouvert) de
$$ a_{\Delta(J,K)} =
  \{x\in X, \forall \alpha
\in \Delta(J,K),\; \la x,\alpha \ra
=0\}.$$
Plus pr{\'e}cis{\'e}ment, c'est une r{\'e}union d'adh{\'e}rences de chambres
paraboliques de $a_{\Delta(J,K)}$ (une chambre parabolique de $a_J$
est par d{\'e}finition une composante connexe du
compl{\'e}mentaire de l'union des hyperplans associ{\'e}s aux racines de
$A_J$ dans le parabolique standard associ{\'e} {\`a} $M_J$.
C'est aussi le $a_J$-int{\'e}rieur de 
l'intersection de l'adh{\'e}rence d'une chambre de Weyl de $X$ avec $a_J$,
lorsque cet int{\'e}rieur est non vide.  Ainsi les
chambres paraboliques sont en bijection avec les sous-groupes
paraboliques de $G$ de Levi $M_J$).

\bigskip

\alin{Repr\'esentations induites d'un Borel}
Par commodit\'e, nous allons supposer que $R$ contient une racine carr\'ee de $p$, fix\'ee une fois pour toutes. Ceci n'est pas n\'ecessaire pour les constructions qui suivent mais simplifie l'exposition. Cela nous permet en particulier de d\'efinir
$$\delta:=\delta_{B_+}^{-\frac{1}{2}}$$
o{\`u} $\delta_{B_+} : A_\emptyset \To{} R^\times$ est le caract{\`e}re-module de $B_+$, et de normaliser les induites paraboliques.
Lorsque $B$ est un sous-groupe de Borel contenant $A_\emptyset$, nous notons $\ip{B}{G}{\delta}$ l'induite parabolique normalis\'ee de $\delta$ le long de $B$. On a $\ip{B}{G}{\delta}=\ind{B}{G}{\delta\delta_B^{\frac{1}{2}}}$ ce qui montre que l'on peut d\'efinir cette induite sans racine de $p$.

Rappelons qu'on d\'efinit une "distance" $d(B,B')$ entre deux \borels contenant $A_\emptyset$ par 
$$d(B,B'):= |\Sigma(A_\emptyset,\hbox{Lie}(B))\setminus \Sigma(A_\emptyset,\hbox{Lie}(B'))|.$$
Enfin nous noterons $C(B) \subset X$ la chambre de Weyl de $X$ associ\'ee au \borel $B$.
Le lemme suivant est bien connu dans le cas $R=\CM$.

\begin{lemme} \label{serprin} On suppose que $R$ est un anneau fortement banal pour $G$, {\em cf} \ref{bonbanal} (pour ${\mathbf G}=GL(n)$, banal suffit). Soient $B, B'$ deux sous-groupes de Borel.
\begin{enumerate}
\item $\hom{\ip{B}{G}{\delta}}{\ip{B'}{G}{\delta}}{G}\simeq
  R$ et contient un g{\'e}n{\'e}rateur canonique not{\'e}
  $J_{B'|B}$, et appel\'e simplement "op\'erateur d'entrelacement".
\item Si $d(B,B'')=d(B,B')+d(B',B'')$, alors $J_{B''|B'} \circ J_{B'|B} = J_{B''|B}$.
\item Si $B$ et $B'$ sont contenus dans un \para $P$ 
  dont on note $M$ la composante de Levi (semi-standard), alors on a $J_{B'|B}=\ip{P}{G}{J_{B'_M|B_M}}$ o{\`u}
  $J_{B'_M|B_M}$ est l'op{\'e}rateur d'entrelacement
  $\ip{B\cap M}{M}{\delta}\To{} \ip{B'\cap M}{M}{\delta}$. 
\item  On a {\'e}quivalence entre les assertions suivantes :
  \begin{enumerate}
  \item $\ip{B}{G}{\delta}$ et $\ip{B'}{G}{\delta}$ sont
    isomorphes (ou, ce qui est \'equivalent par i), $J_{B'|B}$ est un isomorphisme).
  \item     $C(B)$ et $C(B')$ sont contenues dans un m{\^e}me $X_J$, pour un certain $J\subseteq
    S$.
\end{enumerate}
\item Si $C(B)\subset X_J$, alors il existe une factorisation
$$ J_{\o{B}|B}:\;\;\ip{B}{G}{\delta} \twoheadrightarrow \pi_J^R \injo \ip{\o{B}}{G}{\delta}.$$
Autrement dit, on a un isomorphisme $\im J_{\o{B}|B} \simto \pi^R_J$. (On note $\o{B}$ le Borel oppos\'e).

\end{enumerate}
\end{lemme}  
Pour all\'eger cette section, nous reportons 
la preuve de ce lemme {\`a} la section \ref{preuveserprin}. Nous
fixons dor\'enavant l'anneau $R$ fortement banal pour $G$ et l'omettrons
de la plupart des notations.

Choisissons maintenant pour chaque $J\subseteq S$ un Borel $B_J$ tel que $C(B_J)\subset X_J$ et posons 
\ini\begin{equation}\label{defIJ}
 I_J:=\ip{B_J}{G}{\delta}.
\end{equation}
D'apr\`es le lemme pr\'ec\'edent, la classe d'isomorphisme de $I_J$ ne d\'epend pas du choix de $B_J$ et ses endomorphismes sont tous scalaires.
Toujours le lemme pr{\'e}c{\'e}dent
nous fournit aussi un op{\'e}rateur d'entrelacement 
$$J_{J|K}:\; I_K\To{} I_J,$$ 
dont la classe d'homoth{\'e}tie ne d\'epend pas des choix. En particulier, la classe d'isomorphisme de son image ne d\'epend pas non plus de ces choix.

\begin{lemme} \label{lemme1}
Soient $I,J, K \subseteq S$  tels que $\Delta(I,J) \supseteq \Delta(K,J)$.
 Alors on a une factorisation {\em {\`a} homoth{\'e}tie inversible pr{\`e}s}
  $$J_{J|I}\in R^\times .(J_{J|K}\circ J_{K|I})$$ 
\end{lemme}

\begin{proof}
Puisqu'on veut une factorisation \`a homoth\'etie pr\`es, on peut jouer sur les choix des \borels $B_I, B_J$ et $B_K$. Il nous suffit donc de montrer qu'on peut les choisir de sorte que 
$$ J_{B_J|B_I}=J_{B_J|B_K}\circ J_{B_K|B_I}.$$
 D'apr{\`e}s le point iii) du lemme
\ref{serprin}, il nous suffit de trouver des chambres $C_I$ et $C_J$
respectivement dans $X_I$ et $X_J$ reli{\'e}es par  une galerie
tendue dont l'une des chambres interm{\'e}diaires est dans $X_K$.
(Nous renvoyons {\`a} \cite{BT1} pour la notion de galerie tendue et l'usage simple que l'on en fait ci-dessous).


Pour cela, on commence par choisir des chambres $C_I$ et $C_K^1$
telles que  $\o{C_I}\cap \o{C_K^1}$ soit une chambre
parabolique de $\o{X_I}\cap \o{X_K}\subset a_{\Delta(I,K)}$, que nous
noterons $F_{IK}$. 
 De la m{\^e}me mani{\`e}re, on choisit  deux chambres $C_J$ et
 $C_K^2$ dont l'intersection des adh{\'e}rences est une chambre
 parabolique $F_{JK}$ de $\o{X_J}\cap \o{X_K}\subset a_{\Delta(J,K)}$. 
Par d{\'e}finition de nos chambres paraboliques, si $x$ est un point
int{\'e}rieur {\`a} $F_{IK}$, 
$$\begin{array}{rcl} \epsilon_K(\alpha) \la x, \alpha \ra  \; & \hbox{
    est }&
    \left\{\begin{array}{rcl} >0 & \hbox{ si } & \alpha \in \Delta(I,K)^c \\
                             = 0 & \hbox{ si } &  \alpha \in
                             \Delta(I,K) \end{array} 
                         \right. \end{array}
$$
De m{\^e}me, si $y$ est un point int{\'e}rieur {\`a} $F_{JK}$,
$$\begin{array}{rcl} \epsilon_K(\alpha) \la y, \alpha \ra  \; & \hbox{
    est }&
    \left\{\begin{array}{rcl} >0 & \hbox{ si } & \alpha \in \Delta(J,K)^c \\
                             = 0 & \hbox{ si } &  \alpha \in
                             \Delta(J,K) \end{array} 
                         \right. \end{array}
$$
Comme $\Delta(I,K)\cap \Delta(J,K)=\emptyset$ (lemme \ref{comb1}), on en d{\'e}duit que pour
tout point $z$ du segment ouvert $]x,y[$ dans $X$, on a 
$$ \epsilon_K(\alpha) \la z, \alpha \ra >0  \hbox{ pour tout }
\alpha\in S, $$
de sorte que $z\in X_K$. 
En cons{\'e}quence, l'enclos des chambres $C_I$ et $C_J$, qui contient
bien s{\^u}r l'enveloppe convexe de $F_{JK}\cup F_{IK}$ a une
intersection non vide avec $X_K$, et donc avec au moins une chambre
$C_K \subset X_K$. Par d{\'e}finition de l'enclos, il existe une galerie
tendue entre $C_I$ et $C_J$ passant par $C_K$.

\end{proof}

\alin{Interm{\`e}de simplicial}
Soit $E$ un ensemble fini. On lui associe une petite cat{\'e}gorie
$\PC(E)$ dont les objets sont tous les sous-ensembles de $E$ et les
fl{\`e}ches sont les inclusions. On appelle syst{\`e}me de coefficients sur
$\PC(E)$ {\`a} valeurs dans une cat{\'e}gorie ab{\'e}lienne $\CC$ tout foncteur
contravariant $\PC(E)\To{\VC} \CC$.

On veut associer un complexe de cha{\^\i}nes {\`a} un tel syst{\`e}me de
coefficients. Il faut pour cela choisir une orientation de
$\PC(E)$. On peut par exemple choisir un ordre total $\leq$ sur $E$ et
d{\'e}finir pour $I\subset I'\subseteq E$ tels que $I'=I\sqcup \{e\}$ un
signe $\epsilon(I,I'):= (-1)^{|\{i\in I,i\leq e\}|}$.

On d{\'e}finit alors le  complexe de cha{\^\i}nes $\CC_*(\PC(E),\VC)$ associ{\'e} au syst{\`e}me de
coefficients $\VC$ :
$$ 0\To{}\VC(E)\To{d_{|E|}} \cdots \To{} \bigoplus_{|I|=n,I\subseteq E} \VC(I)
\To{d_n}\bigoplus_{|I|=n-1} \VC(I) \To{} \cdots \To{} 
\VC(\emptyset) \To{} 0$$
o{\`u} $$d_n = \bigoplus_{|I|=n} \sum_{I'\subset I} \epsilon(I,I')
\VC(I'\subset I). $$

Il est bien connu que le complexe de cha{\^\i}nes associ{\'e} {\`a} un syst{\`e}me de
coefficients constant est acyclique, sauf si $E=\emptyset$.

\ali
Fixons dor\'enavant deux sous-ensembles $I,J$ de $S$. Nous allons d\'efinir un syst\`eme de coefficients $\VC_{I,J}$ sur la cat\'egorie $\PC(\Delta(I,J))$ \`a valeurs dans $\Mo{R}{G}$. 
Pour cela, posons brutalement pour tout $K\subseteq \Delta(I,J)$
$$ \VC_{I,J}(K) := \im \left(J_{J^c|\Delta(I,K)} :\;\; I_{\Delta(I,K)} \To{} I_{J^c} \right).$$
En \'ecrivant $K=\Delta(I,\Delta(I,K))$ et en utilisant le lemme \ref{comb1}, on obtient 
$\Delta(J,\Delta(I,K))=\Delta(I,J)\setminus K$, puis en passant au compl\'ementaire, 
$\Delta(J^c,\Delta(I,K))= K\sqcup \Delta(I,J)^c$. On a donc
\begin{eqnarray*}
        K\subset K' \subseteq \Delta(I,J) & \Rightarrow & 
        \Delta(J^c,\Delta(I,K)) \subset \Delta(J^c,\Delta(I,K'))
\end{eqnarray*}
donc d'apr\`es le lemme \ref{lemme1}, on a une inclusion canonique
$$ \VC_{I,J}(K\subset K') :\; \VC_{I,J}(K') \injo \VC_{I,J}(K) .$$

On a donc ainsi d\'efini un syst\`eme de coefficients $\VC_{I,J}$, qui d\'epend des choix des \borels $B_J$ mais dont la classe d'isomorphisme n'en d\'epend pas.
En supposant fix\'e un ordre total sur $S$, on obtient un complexe de
cha{\^\i}nes $\CC_*(\PC(\Delta(I,J)),\VC_{I,J})$ de longueur $\delta(I,J)$.
 On peut augmenter ce complexe par
le morphisme surjectif suivant : puisque $\ker(J_{J^c|I})\subseteq \ker(J_{I^c|I})$ (par le lemme \ref{lemme1}), le morphisme $I_I\To{} \im(J_{I^c|I})=\pi_I$ du point v) du lemme \ref{serprin} se factorise  en
$$\VC_{I,J}(\emptyset)= \im(I_I\To{J_{J^c|I}} I_{J^c}) \To{} \pi_I$$
 Nous noterons $\CC_*(\VC_{I,J})^+$ le complexe ainsi augment\'e.
En faisant le changement de variable $K\mapsto \Delta(I,K)$ dans la d\'efinition du complexe de cha\^ines, 
le complexe augment\'e $\CC_*(\VC_{I,J})^+$ s'\'ecrit aussi :



\ini\begin{equation} \label{compl}
\begin{array}{rcl}  0\To{} \pi_{J}\simeq\im(J_{J^c|J}) \To{} \cdots  & \To{d_{n+1}} &
  \bigoplus\limits_{\tiny\begin{array}{c} \delta(I,K)=n \\ \delta(J,K)=\delta(I,J)-n\end{array}} \im(J_{J^c|K})  \\
&  \To{d_n} &  \cdots  \To{} \im(J_{J^c|I}) \To{} \pi_I \To{} 0,\end{array}
\end{equation}
le terme $\pi_I$ \'etant en degr\'e $1$ et le terme $\pi_J\simeq \im(J_{J^c|J})$ (cet isomorphisme vient de \ref{serprin} v)) en degr\'e $-\delta(I,J)$.

\bigskip

On peut maintenant \'enoncer une version explicite du th\'eor\`eme \ref{theoext} :

\begin{theo} \label{explicite} Supposons $R$ fortement banal pour $G$. Soient $I,J\subseteq S$.
\begin{enumerate}
        \item Le complexe augment\'e $\CC_*(\VC_{I,J})^+$ de \ref{compl} est acyclique.
        \item Soit $\alpha_{I,J} \in \ext{\delta(I,J)}{\pi_I}{\pi_J}{RG}$ la classe de $\CC_*(\VC_{I,J})^+$. Alors pour tout $H\subseteq S$ tel que $\delta(H,J)=\delta(H,I)+\delta(I,J)$, le $\cup$-produit
        $$  \alpha_{I,J} \cup - :\;\; \ext{*}{\pi_H}{\pi_I}{RG} \To{} \ext{*+\delta(I,J)}{\pi_H}{\pi_J}{RG} $$
        est un isomorphisme pour $*\in \ZM$.
\end{enumerate}
\end{theo}

La preuve occupera les sections suivantes. Pr\'ecisons tout de m\^eme comment cet \'enonc\'e implique le th\'eor\`eme \ref{theoext}. Cela utilise le lemme suivant :

\begin{lemme} \label{dimcoh} Supposons $R$ banal pour $G$ et
soit $I\subseteq S$.
\begin{enumerate}
        \item On a $R\simto \endo{RG}{\pi_I}$.
        \item Si $G$ est semi-simple, il existe une r\'esolution projective de $\pi_I$ dans $\Mo{R}{G}$ de longueur $|S|+1$.
\end{enumerate}
\end{lemme}
\begin{proof}
Pour le point i), choisissons un \borel $B$ tel $C(B)\subset X_I$. D'apr\`es \ref{serprin} v) on a un monomorphisme de $R$-modules
$$ \endo{RG}{\pi_I} \injo \hom{\ip{B}{G}{\delta}}{\ip{\o{B}}{G}{\delta}}{RG}. $$
Mais d'apr\`es \ref{serprin} i), c'est un isomorphisme et le terme de droite est libre de rang $1$ sur $R$.

Pour le point ii), on peut invoquer deux r\'ef\'erences. D'une part, le cas ${\mathbf G}=PGL(n)$ est trait\'e dans \cite{SS1}. D'autre part, le cas  $R$  corps alg\'ebriquement clos de caract\'eristique banale est trait\'e dans \cite{Vigsheaves} (g\'en\'eralisant le cas $R=\CM$ trait\'e dans \cite[II.3.3]{SS2}). Il se trouve que les deux r\'ef\'erences se g\'en\'eralisent \`a notre situation. Pour ce qui est de \cite{SS1}, on remarque que les arguments reposent purement sur la combinatoire des immeubles de Tits et Bruhat-Tits : les r\'esolutions y sont construites pour un anneau quelconque $R$ et deviennent bien projectives d\`es  que $R$ est banal.
Il suffit donc de v\'erifier que les propri\'et\'es  combinatoires de ces immeubles qui sont utilis\'ees sont vraies pour tout groupe d\'eploy\'e. Ceci ne fait aucun doute (et doit m\^eme \^etre implicite dans \cite{SS2}), mais nous ne le ferons pas ici.

En revanche nous g\'en\'eralisons les arguments de \cite{Vigsheaves} dans l'appendice \ref{decomposition}, {\em cf} le corollaire \ref{corofinit}.
\end{proof}

\alin{Le th\'eor\`eme \ref{explicite} implique le th\'eor\`eme \ref{theoext}}
Commen{\c c}ons par remarquer le cas particulier $H=I$ de \ref{explicite} ii) : on obtient un isomorphisme
 $$ \alpha_{I,J}\cup - :\;\;\hom{\pi_I}{\pi_I}{RG} \simto \ext{\delta(I,J)}{\pi_I}{\pi_J}{RG}, $$
 ce qui, compte tenu de l'isomorphisme de \ref{dimcoh} i), montre que $\ext{\delta(I,J)}{\pi_I}{\pi_J}{RG}$ est libre de rang $1$ sur $R$ et engendr\'e par $\alpha_{I,J}$.

Il faut maintenant montrer l'annulation des autres $\hbox{Ext}$. Toujours par le cas particulier ci-dessus du point ii) de \ref{explicite}, il suffit de prouver que pour tout $I$, 
$$ \ext{*}{\pi_I}{\pi_I}{RG} = 0 \;\; \hbox{ d\`es que } \; *>0. $$
Appliquons pour cela ce point ii) au cas $H=I=J^c$. On obtient les isomorphismes
$$ \ext{*}{\pi_I}{\pi_I}{RG}  \simto \ext{|S|+*}{\pi_I}{\pi_{I^c}}{RG}. $$
Or par \ref{dimcoh} ii), le terme de droite est nul d\`es que $*>0$.
On a donc prouv\'e le point i) de \ref{theoext}. Le point ii) de ce m\^eme th\'eor\`eme est alors une cons\'equence imm\'ediate du point ii) de \ref{explicite}.

\alin{L'image de $J_{J|K}$} \label{image} 
Autant pour expliciter un peu les termes du complexe \ref{compl} que pour la preuve du th\'eor\`eme \ref{explicite} dans les sections suivantes,
 nous voulons d\'ecrire assez pr\'ecis\'ement l'image de l'op\'erateur d'entrelacement
 $J_{J|K}:\; I_K\To{} I_J$ associ\'e \`a deux sous-ensembles $J,K$ de $S$.

Pour cela, nous devons malheureusement compliquer encore un peu les notations. Si $T\subseteq S$, le groupe de Levi $M_T$ est un groupe r\'eductif, muni  d'un tore maximal, $A_\emptyset$, et d'une  base de son syst\`eme de racines, $T$, associ\'ee au \borel $B_+\cap M_T$. La construction \ref{defrep} permet donc d'associer \`a tout sous-ensemble $I$ de $T$ une repr\'esentation de $M_T$ que nous noterons $\pi_{I,T}$. 
Par exemple, $\pi_{T,T}$ est la repr\'esentation triviale de $M_T$, $\pi_{\emptyset,T}$ est sa repr\'esentation de Steinberg, et
pour $T=S$ on a $M_S=G$ et on retrouve $\pi_{I,S}=\pi_I$. 
 
\begin{lemme} \label{lemmeimage}
Fixons deux sous-ensembles $J,K$ de $S$ et notons simplement $\Delta:=\Delta(J,K)$.
Pour tout \para $P_{K,J}$ de $G$ associ\'e \`a une chambre parabolique ({\em cf} \ref{1.3}) de $a_\Delta$ contenue dans $\o{X_J}\cap\o{X_K}$, on a 
$$ \im(J_{J|K}) \simeq \ip{P_{K,J}}{G}{\delta_{P_\Delta}^{-\frac{1}{2}}\pi_{_{K\cap\Delta,\Delta}}}. $$
(Rappelons que $P_\Delta$ est le sous-groupe parabolique standard de Levi $M_\Delta$).
\end{lemme}
 \begin{proof}
 Comme la classe d'isomorphisme de l'image de $J_{J|K}$ ne d\'epend  que
de la classe d'homoth{\'e}tie de cet op\'erateur,  nous pouvons  jouer sur les choix de $B_J$ et $B_K$ qui interviennent dans la d\'efinition \ref{defIJ}.
 Choisissons donc ces \borels de sorte que chacun soit contenu dans $P_{K,J}$. Il suffit pour cela que les chambres de Weyl associ\'ees v\'erifient 
 $C(B_J)\subset X_J$, $C(B_K)\subset X_K$ et $\overline{C(B_K)}\cap
\overline{C(B_J)} = \overline{C(P_{K,J})}$ en notant $C(P_{K,J})$ la chambre parabolique de $a_\Delta$ associ\'ee \`a $P_{K,J}$.

On d{\'e}duit alors du point iii)  du lemme \ref{serprin} l'existence
d'un diagramme commutatif 
$$\xymatrix{ I_K=\ip{B_K}{G}{\delta} \ar[r]^{J_{J|K}}
  \ar@{=}[d] & I_J=\ip{B_J}{G}{\delta} \ar@{=}[d] \\
  \ip{P_{K,J}}{G}{\ip{M_\Delta\cap B_K}{M_\Delta}{\delta}} \ar[r]_{\ip{P_{K,J}}{G}{J}} &
 \ip{P_{K,J}}{G}{\ip{M_\Delta\cap B_J}{M_\Delta}{\delta}}
}$$
o{\`u} $$J=J_{M_\Delta\cap B_K|M_\Delta\cap B_J} :\;\ip{M_\Delta\cap
  B_K}{M_\Delta}{\delta}\To{}\ip{M_\Delta\cap
  B_J}{M_\Delta}{\delta} $$ est l'op{\'e}rateur d'entrelacement
canonique. 
On a une factorisation $\delta=\delta_{P_\Delta}^{-\frac{1}{2}} \delta_{M_\Delta\cap B_+}^{-\frac{1}{2}}$, et on est donc ramen\'e \`a prouver que l'image de l'op\'erateur d'entrelacement 
$$ J=J_{M_\Delta\cap B_K|M_\Delta\cap B_J} :\;\ip{M_\Delta\cap
  B_K}{M_\Delta}{\delta_{M_\Delta\cap B_+}^{-\frac{1}{2}}}\To{}\ip{M_\Delta\cap
  B_J}{M_\Delta}{\delta_{M_\Delta\cap B_+}^{-\frac{1}{2}}} $$ 
est isomorphe \`a $\pi_{K\cap\Delta,\Delta}$. 
Mais puisque $J\cap\Delta= \Delta\setminus(K\cap\Delta)$, ceci est
le point v) du lemme \ref{serprin}.
\end{proof}

\begin{coro} \label{chgtbase}
Le syst\`eme de coefficients $\VC_{I,J}$ et le complexe augment\'e $\CC_*(\VC_{I,J})^+$ sont compatibles au changement de scalaires en le sens suivant : si $R\To{\psi} R'$ est un morphisme d'anneaux fortement banals, alors les morphismes canoniques 
$$ \VC_{I,J}^R \otimes_{R,\psi} R' \To{} \VC_{I,J}^{R'}, \hbox{ resp. }  \; \CC_*(\VC_{I,J}^R)^+ \otimes_{R,\psi} R' \To{} \CC^*(\VC_{I,J}^{R'})^+ $$
sont des isomorphismes de syst\`emes de coefficients sur $\Delta(I,J)$ \`a coefficients dans $\Mo{R'}{G}$, resp. de complexes de $R'G$-repr\'esentations lisses.
\end{coro}

\begin{proof}
Il suffit bien-s\^ur de v\'erifier l'assertion pour le syst\`eme de coefficients, puisque le complexe de cha\^ines lui est associ\'e fonctoriellement (seule l'augmentation demande un argument suppl\'ementaire mais qui ne pose pas de probl\`eme). Par le lemme pr\'ec\'edent, on est ramen\'e \`a prouver que pour tout $I\subseteq T\subseteq S$, l'application $\pi_{I,T}^R\otimes_{R,\psi} R' \To{} \pi_{I,T}^{R'} $ est un isomorphisme de $R'G$-repr\'esentations. Mais cela r\'esulte de l'exactitude \`a droite du produit tensoriel et  de la d\'efinition de $\pi_{I,T}$ comme conoyau de la fl\`eche 
$$ \bigoplus_{T\supseteq J\supset I} \ind{P_J \cap M_T}{M_T}{1} \To{} \ind{P_I\cap M_T}{M_T}{1} $$
puisque cette fl\`eche est clairement compatible \`a l'extension des scalaires.
\end{proof}

\subsection{Le cas $R$ corps alg\'ebriquement clos}

Dans cette section on suppose que
$R$ est un corps alg{\'e}briquement clos de caract\'eristique $\neq p$, dans lequel on choisit une racine de $p$ pour normaliser les foncteurs paraboliques. Pour un \para $P$ de Levi $M$, on notera  $\Ip{P}{G}$ et $\Rp{G}{P}$ les foncteurs d'induction et restriction paraboliques normalis\'es.

La caract\'eristique de $R$ sera toujours suppos\'ee {\em banale},
{\em cf} \ref{bonbanal}. La raison essentielle vient du r\'esultat
suivant :
\begin{fact} \label{faitbanal}
  Appelons ``bloc principal'' de $\Mo{R}{G}$ la sous-cat\'egorie pleine
  des objets dont tous les sous-quotients irr\'eductibles sont de la
  s\'erie principale, c'est-\`a-dire apparaissent comme sous-quotients
  d'induites $\ip{B}{G}{\chi}$ pour $\chi$ caract\`ere non ramifi\'e
  de $T$. Lorsque $R$ est de caract\'eristique banale pour $G$, alors 
 cette sous-cat\'egorie est "facteur direct"
  de $\Mo{R}{G}$, et les restrictions des foncteurs $\Rp{G}{B}$, pour $B$
  un \borel,  y  sont {\em fid\`eles}.
\end{fact}
\begin{proof}
Lorsque $R=\CM$, ceci est une partie du th\'eor\`eme principal de Bernstein dans \cite{bernstein}. La preuve est une \'elaboration du fait que si une repr\'esentation irr\'eductible complexe est cuspidale, c'est-\`a-dire annul\'ee par tous les foncteurs de restriction parabolique, alors elle est projective dans la cat\'egorie des repr\'esentations complexes \`a caract\`ere central fix\'e. Cette propri\'et\'e des cuspidales n'est pas vraie en toutes caract\'eristiques, mais d'apr\`es \cite{Vigsheaves}, elle l'est en caract\'eristiques banales. \`A partir de l\`a, il est commun\'ement admis (voir
par exemple \cite{Vigicm}) que les arguments de Bernstein s'appliquent de la m\^eme mani\`ere que dans le cas complexe.
\end{proof}

\alin{S\'eries principales elliptiques}
Nous allons tout d'abord pr\'eciser et renforcer le lemme \ref{serprin} dans le cas d'un corps. 
Appelons "s\'eries principales elliptiques"  les sous-quotients irr{\'e}ductibles de l'induite
 $\ind{B_+}{G}{1}=\ip{B_+}{G}{\delta}$ (bien que cette terminologie soit peut-\^etre un peu usurp\'ee lorsque ${\mathbf G}\neq GL(n)$). 
 Dans le cas $R=\CM$, elles sont classifi{\'e}es par les
sous-ensembles de $S$, comme on peut le d{{\'e}duire} de plusieurs travaux
plus ou moins ind{\'e}pendents  dont ceux de Rodier dans \cite{rodier},
Langlands (``quotient de Langlands'') ou
Bernstein-Zelevinski.
Pour $GL(n)$, Vign\'eras les a classifi\'ees en toutes caract\'eristiques ($\neq p$) et dans le cas (fortement) banal, on obtient encore une param\'etrisation par les sous-ensembles de $S$.
Toutes ces approches montrent en particulier 
que ces repr{\'e}sentations sont toutes de la forme 
$\pi^R_I$ d{\'e}finie en \ref{defrep}.
 Le lemme suivant montre entre autres qu'il en est bien de m{\^e}me sur $R$ de  caract{\'e}ristique {\em fortement banale} pour $G$.

\begin{lemme} \label{serprincorps} $R$ est un corps alg\'ebriquement clos de caract\'eristique fortement banale (pour ${\mathbf G}=GL(n)$, "banale" suffit). 
\begin{enumerate}
\item Pour tout \borel $B$, l'induite $\ip{B}{G}{\delta}$ a un unique quotient
  irr{\'e}ductible et toute s{\'e}rie principale elliptique est isomorphe
  {\`a} un tel quotient.
\item  On a {\'e}quivalence entre les assertions suivantes :
  \begin{enumerate}
  \item $\ip{B}{G}{\delta}$ et $\ip{B'}{G}{\delta}$ sont
    isomorphes.
  \item leurs quotients irr{\'e}ductibles sont isomorphes
  \item     $C(B)$ et $C(B')$ sont contenues dans un m{\^e}me $X_J$, pour $J\subset
    S$.
\end{enumerate}

\item Pour tout $J\subseteq S$, $\pi^R_J$ est irr\'eductible. Si $C(B)\subset X_J$ alors le quotient irr{\'e}ductible de $\ip{B}{G}{\delta}$ est
isomorphe \`a $\pi_J^R$. 
\item Soit $w_0$ l'\'el\'ement de plus grande longueur de $W_G$. Alors la contragr\'ediente de $\pi_J^R$ est $\pi^R_{-w_0(J)}$.
\end{enumerate}
\end{lemme}  
La preuve de ce lemme est donn\'ee avec celle de \ref{serprin} dans la
section \ref{preuveserprin}.
Remarquons que tous les points peuvent \^etre mis en d\'efaut en caract\'eristique non banale, m\^eme l'irr\'eductibilit\'e des $\pi_I^R$.

\begin{rema} \label{parLanpiI} (Param\`etres de Langlands de $\pi_J^\CM$). On suppose ici que $R=\CM$ et ${\mathbf G}=GL(n)$. D'apr\`es le lemme pr\'ec\'edent, $\pi_J$ est l'unique quotient irr\'eductible de $\im(J_{K|J})$ pour tout $K$ et donc en particulier pour $K=S$. Or d'apr\`es le lemme \ref{lemmeimage} on a $\im(J_{S|J})\simeq \ip{\o{P_{J^c}}}{G}{\delta_{P_{J^c}}^{-\frac{1}{2}} \pi_{J\cap J^c, J^c}}$ et on peut expliciter $\pi_{J\cap J^c,J^c}=\hbox{St}_{M_{J^c}}$ (repr\'esentation de Steinberg de $M_{J^c}$). Comme $\o{P_{J^c}}$ est conjugu\'e \`a $P_{-w_0(J^c)}$ par l'\'el\'ement de plus grande longueur $w_0$ de $W_G$, on obtient que $\pi_J$ est l'unique quotient irr\'eductible de l'induite $\ip{P_{-w_0(J^c)}}{G}{\delta_{P_{-w_0(J^c)}}^{\frac{1}{2}}\hbox{St}_{M_{-w_0(J^c)}}}$. Ceci donne les param\`etres de Langlands de $\pi_J$. 
\end{rema}

\alin{Exposants}
 Puisque $R$ est un corps alg{\'e}briquement clos de
caract{\'e}ristique $\neq p$, on sait  que toute
$R$-repr{\'e}sentation lisse irr{\'e}ductible a un caract{\`e}re central. En
particulier,
lorsque $\pi$ est une $R$-repr{\'e}sentation de longueur finie de $M_J$, on
note $\exp(A_J,\pi)$ l'ensemble des caract{\`e}res centraux des
sous-quotients ir{\'e}ductibles de $\pi$. 
Rappelons que dans ces circonstances on a une d{\'e}composition canonique,
dite ``d{\'e}composition isotypique'',
$$\pi\simeq \bigoplus_{\chi\in\exp(A_J,\pi)} \pi_{\chi},\;\hbox{ o{\`u}
  }\; \pi_\chi=\{v\in \pi, \exists n \in \NM, \forall z\in A_J, \;
(\pi(z)-\chi(z))^nv=0 \}.$$
On a aussi $$\exp(A_J,\pi)=\{\chi:A_J\To{} R^\times,\;\; \hom{\chi}{\pi}{A_J}\neq 0\}.$$

\begin{lemme} \label{lemme3}
Soient $K,J\subseteq S$. Les propri\'et\'es suivantes pour $I\subseteq S$ sont \'equivalentes :
\begin{enumerate}
\item $\hom{I_I}{\im(J_{J|K})}{G} \neq 0$,
\item  $\pi_I$ est un sous-quotient irr\'eductible de $\im(J_{J|K})$,
\item $\Delta(I,J)\supseteq \Delta(J,K)$.  
\end{enumerate}
\end{lemme}
\begin{proof}
L'implication $iii)\Rightarrow i)$ est une cons\'equence de la propri\'et\'e de factorisation \`a homoth\'etie pr\`es du lemme \ref{lemme1}.
L'implication $i)\Rightarrow ii)$ est imm\'ediate puisque $\pi_I$ est l'unique quotient irr\'eductible de $I_I$, d'apr\`es le lemme \ref{serprincorps}.

Il nous reste \`a prouver $ii)\Rightarrow iii)$.
Pour cela, nous allons  d'abord calculer $\exp(T,\Rp{G}{B_J}(\pi_I))$.
Par d\'efinition de $\pi_I$, cet ensemble de caract\`eres lisses de $T$ est contenu dans l'orbite $W.\delta$ de $\delta$.
En utilisant la r\'eciprocit\'e de Frobenius-Casselman
\cite[II.3.8-2]{Vig} pour la premi\`ere ligne ci-dessous, on a 
\begin{eqnarray*}
\hom{w(\delta)}{ \rp{G}{{B_J}}{\pi_I}}{RT} \neq 0 & \equ & 
\hom{\ip{\o{B_J}}{G}{\delta}}{\pi_I}{RG}\neq 0 \\
 & \equ &  \hom{\ip{w^{-1}(\o{B_J})}{G}{\delta}}{\pi_I}{RG}\neq 0 \\
& \equ & w^{-1}(-C_J)\subset X_I  \\ 
\end{eqnarray*}
gr\^ace \`a \ref{serprincorps} iv).
On a donc 
\begin{eqnarray*}
\exp(T,\Rp{G}{B_J}(\pi_I)) & = & \{w(\delta), w^{-1}(-C_J)\subset X_I\}
\\
 & = & \{w(\delta), \forall x\in -C_J, \forall \alpha\in S,\;
 \epsilon_I(\alpha) \la x,w(\alpha) \ra >0\}  
\end{eqnarray*}
On peut paraphraser la derni{\`e}re {\'e}galit{\'e} en introduisant le
c{\^o}ne polaire 
$D_J$ de $C_J$ d{\'e}fini par $D_J=\{y\in X^*, \forall x\in C_J,\;
\la x,y\ra >0\}$. C'est aussi le c{\^o}ne engendr{\'e} par les racines simples
correspondant {\`a} la chambre de Weyl $C_J$ (on appelle parfois $D_J$ la
chambre obtuse). On peut donc r{\'e}{\'e}crire

\ini
\begin{equation}\label{exp1}
 \exp(T,\Rp{G}{B_J}(\pi_I)) = \{w(\delta),\; \forall \alpha\in S,
\epsilon_I(\alpha)w(\alpha) \in -D_J\}.\label{eq:4}
\end{equation}

Dans la suite, on utilise les notations du lemme \ref{lemmeimage}. Pour all\'eger un peu ces notations, on pose  $\sigma_J^K:= \delta_{P_\Delta}^{-\frac{1}{2}} \pi_{_{K\cap\Delta,\Delta}}$.
Soit $P_{K,J}$ un \para comme dans \ref{lemmeimage} et $B_J,B_K\subset P_{K,J}$ deux \borels tels que $C(B_J)\subset X_J$ et $C(B_K)\subset X_K$.
On va maintenant calculer $\exp(\Rp{G}{B_J}\circ
\ip{P_{K,J}}{G}{\sigma_J^K})$.
 Pour utiliser le lemme g{\'e}om{\'e}trique,
introduisons le sous-ensemble 
$$W^J=\{w\in W,\; w(D_J^\Delta)\subset D_J\}$$
o{\`u} $D_J^\Delta$ est le c{\^o}ne de $(a_{\Delta(J,K)})^\perp$ engendr{\'e} par l'ensemble $J\cap \Delta$ des racines
simples correspondant au sous-groupe de Borel $B_J\cap M_\Delta$ de
$M_\Delta$ (rappelons que $\Delta=\Delta(J,K)$). On sait
\cite[2.11]{BZ2} que $W^J$ est un  ensemble de 
repr{\'e}sentants privil{\'e}gi{\'e}s des classes {\`a} droite de $W$ modulo
$W(\Delta)$, et le lemme g{\'e}om{\'e}trique \cite[2.12]{BZ2} nous assure que

\ini\begin{equation}\label{exp2}
\exp(T,\Rp{G}{B_J}\circ
\ip{P_{K,J}}{G}{\sigma_J^K}) = \bigsqcup_{v\in W^J}
v.\exp(T,\Rp{M_\Delta}{M_\Delta\cap B_J}{\sigma_J^K})
\end{equation}
Par ailleurs, la d{\'e}finition de $\sigma_J^K$ nous permet de calculer
comme \ref{exp1}

\ini\begin{equation} \label{exp3}
 \exp(T,\Rp{M_\Delta}{M_\Delta\cap B_J}{\sigma_J^K})  =  \{w_\Delta(\delta)
 , w_\Delta\in W(\Delta) \hbox{ et } \forall \alpha\in \Delta,
 \epsilon_K(\alpha)w_\Delta(\alpha)\in -D_J^\Delta\}
\end{equation}

Supposons maintenant que $\pi_I$ soit un sous-quotient irr{\'e}ductible de
$\im(J_{J|K})$. En particulier, on doit avoir 
$$\exp(T,\Rp{G}{B_J}(\pi_I))\subset\exp(T,\Rp{G}{B_J}\circ
\ip{P_{K,J}}{G}{\sigma_J^K}). $$
Fixons alors $\chi$ dans le terme de gauche. Par \ref{exp2} et
\ref{exp3}, on peut {\'e}crire $\chi=vw_\Delta(\delta)$ avec $v\in W^J$ et
$w_\Delta\in W(\Delta)$.

 Soit $\alpha\in \Delta$, par \ref{exp3}, on a
 $\epsilon_K(\alpha)w_\Delta(\alpha) \in -D_J^\Delta$. Par la
 d{\'e}finition de $W^J$, on a donc aussi
 $\epsilon_K(\alpha)vw_\Delta(\alpha) \in -D_J$. Mais par \ref{exp1},
 ceci entraine que $\epsilon_K(\alpha)=\epsilon_I(\alpha)$. En
 d'autres termes : la restriction de  $\epsilon_I$ {\`a} $\Delta(J,K)$ est
 l'oppos{\'e}e de celle de $\epsilon_J$. Ou autrement dit,
 $\Delta(I,J)\supseteq \Delta(J,K)$.

\end{proof}

\alin{Acyclicit\'e des complexes $\CC_*(\VC_{I,J})^+$ de \ref{compl}}
Fixons deux sous-ensembles $I,J$ de $S$. Nous voulons ici montrer, sous l'hypoth\`ese que $R$ est un corps alg\'ebriquement clos de caract\'eristique fortement banale, l'acyclicit\'e du complexe de cha\^ines augment\'e $\CC_*(\VC_{I,J})^+$.
 Comme la caract{\'e}ristique de $R$ est  banale, le foncteur
$\Rp{G}{B_+}$ est exact et fid{\`e}le sur 
le bloc principal, {\em cf} \ref{faitbanal}.
 Il suffit  donc de
 v{\'e}rifier que $\Rp{G}{B_+}(\CC_*(\VC_{I,J}^+))$ est exact. Ce dernier complexe
 est un complexe de repr{\'e}sentations non-ramifi{\'e}es de longueur finie du tore maximal $T$,
 et se d{\'e}compose donc en une somme  directe de ses composantes
 $\chi$-isotypiques, pour $\chi$ caract{\`e}re non ramifi{\'e} de $T$.
 $$\Rp{G}{B_+}(\CC_*(\VC_{I,J})^+)=\bigoplus_{\chi\in W.\delta}\left(\Rp{G}{B_+}{\CC_*(\VC_{I,J})^+}\right)_\chi$$
 Pour calculer ces composantes isotypiques, on remarque d'abord que
 les repr{\'e}sentations $\rp{G}{B_+}{\im(J_{J|K})}$ sont toutes de multiplicit{\'e}
 $1$, et donc semisimples. Il s'ensuit simplement que 
$$\left(\Rp{G}{B_+}{\CC_*(\VC_{I,J})^+}\right)_\chi \simeq \hom{\chi}{\Rp{G}{B_+}{\CC_*(\VC_{I,J})^+}}{T}.$$
Posons maintenant pour tout $K\subseteq \Delta(I,J)$
$$\VC_{I,J,\chi}(K):=\hom{\chi}{\Rp{G}{B_+}\VC_{I,J}(K)}{T}.$$
Les fl{\`e}ches de transition $\VC_{I,J}(K\subset K')$
induisent par fonctorialit{\'e} des fl{\`e}ches
 $$\VC_{I,J,\chi}(K\subset K'):\;\;\VC_{I,J,\chi}(K')\To{} \VC_{I,J,\chi}(K)$$ et
on obtient ainsi un syst{\`e}me de coefficients
$\VC_{I,J,\chi}$ sur la cat\'egorie $\PC(\Delta(I,J))$ \`a valeurs dans la cat\'egorie des $R$-espaces vectoriels.
On peut augmenter le complexe de cha\^ines $\CC_*(\VC_{I,J,\chi})$ associ\'e par l'application
$$ \VC_{I,J,\chi}(\emptyset)=\hom{\chi}{\Rp{G}{B_+}(I_I)}{T} \To{} \hom{\chi}{\Rp{G}{B_+}(\pi_I)}{T} $$ et en notant $\CC_*(\VC_{I,J,\chi})^+$ le complexe ainsi augment\'e, on a bien-s\^ur
$$ \CC_*(\VC_{I,J,\chi})^+ = \left(\Rp{G}{B_+}\CC_*(\VC_{I,J})^+\right)_\chi.$$
Il nous suffit donc de montrer que pour tout $\chi$ le complexe de gauche est acyclique.

Commen{\c c}ons par remarquer que les fl\`eches de transition $\VC_{I,J,\chi}(K\subset K')$ sont injectives (par exactitude de $\Rp{G}{B_+}$ et exactitude {\`a} gauche de $\hom{\chi}{.}{T}$, et parce que les fl\`eches $\VC_{I,J}(K\subset K')$ le sont).
Comme les espaces $\VC_{I,J,\chi}(K)$ sont de dimension $0$ ou $1$, il nous suffira donc de pr{\'e}ciser cette dimension pour d\'eterminer le syst\`eme de coefficients $\VC_{I,J,\chi}$. Pour cela, notons $L$ l'unique sous-ensemble de $S$ tel que
$I_L\simeq
\ip{\o{B_+}}{G}{\chi}$ (ou ce qui est {\'e}quivalent, tel que $\pi_L$ soit
un quotient de $\ip{\o{B_+}}{G}{\chi}$). 
Par  r\'eciprocit\'e de Frobenius-Casselman \cite[II.3.8-2]{Vig}, on a  $\VC_{I,J,\chi}(K)\neq 0$ \ssi\ $\hom{I_L}{\im(J_{J^c|\Delta(I,K)})}{G}
\neq 0$.
D'apr{\`e}s le lemme \ref{lemme3}, ceci est encore {\'e}quivalent {\`a}
$\Delta(L,J^c)\supseteq  \Delta(J^c,\Delta(I,K))=K \sqcup \Delta(I,J^c)$ (la derni\`ere \'egalit\'e a \'et\'e expliqu\'ee lors de la d\'efinition des $\VC_{I,J}$).

En particulier, le syst\`eme de coefficients $\VC_{I,J,\chi}$ et son complexe augment\'e ne sont non nuls que si $\Delta(L,J^c)\supseteq \Delta(I,J^c)$.
Supposons tout d'abord que $\Delta(L,J^c)\supsetneq \Delta(I,J^c)$.
Dans ce cas-l\`a, $L\neq I$ donc le complexe augment\'e $\CC_*(\VC_{I,J,\chi})^+$ co\"incide avec le complexe non-augment\'e  car $\hom{\chi}{\Rp{G}{B_+}(\pi_I)}{T}=0$.
Or, la discussion pr{\'e}c{\'e}dente montre que le syst{\`e}me de coefficients
$\VC_{I,J,\chi}$
est support{\'e} par le sous-simplexe des sous-ensembles de $\Delta(L,J^c)\setminus \Delta(I,J^c) = \Delta(I,L)$ (lemme \ref{comb1})
et y est constant. En d'autres termes
$$ \CC_*(\PC(\Delta(I,J)),\VC_{I,J,\chi}) \simeq \CC_*(\PC(\Delta(I,L)),R)$$
et comme on l'a d\'eja rappel{\'e}, le complexe de droite est acyclique, puisque $I\neq L$.

Supposons maintenant $\Delta(L,J^c)=\Delta(I,J^c)$, c'est-\`a-dire $L=I$. Alors le syst{\`e}me de coefficients 
$\VC_{I,J,\chi}$ est support{\'e} par $\{\emptyset\}$
mais l'augmentation est ici un isomorphisme, de sorte que
$$ \CC_*(\VC_{I,J,\chi})^+ \simeq \left( R\To{id} R \right)$$
 dont l'homologie est bien s{\^u}r nulle.

\bigskip

C'est essentiellement dans le lemme suivant qu'intervient vraiment
l'hypoth\`ese de {\em bonne} caract\'eristique de \ref{bon}. 

\begin{lemme}
  \label{comb2} 
Supposons que la caract{\'e}ristique de $R$ soit bonne pour $G$,
c'est-{\`a}-dire ne divise pas l'entier $N_W$ de \ref{bon}.
Soit $J\subseteq S$ et $w\in W$ tel que $w(\delta)_{|A_J}=\delta_{|A_J}$.
Alors $w\in W(J)$ o{\`u} $W(J)$ est le groupe de Weyl de $M_J$.
\end{lemme}
\begin{proof}
Remarquons tout-d'abord que puisque $\delta$ est trivial sur le centre
de $G$, il suffit de montrer l'assertion pour $G^{ad}:={\mathbf G}/{\mathbf
  Z({\mathbf G})}(K)$. Nous supposerons donc ${\mathbf G}$ adjoint.

Pour un {\'e}l{\'e}ment $x$ quelconque de $X^*$ nous noterons $x_J\in a_J^*$
son image par la projection canonique $X^*\To{} a_J^*$.

\noindent {\em Premi{\`e}re {\'e}tape : sous l'hypoth{\`e}se de caract{\'e}ristique, on
  a pour tout $J$ et tout $w$ :
  $$w(\delta)_{|A_J}=\delta_{|A_J} \equ w(\rho)_J =\rho_J. $$}
Rappelons pour cela que $\delta$ est d{\'e}fini par $\forall t\in T, \delta(t):=
q_F^{-{1\over 2}\val_F(\rho(t))}$, ou ce qui est {\'e}quivalent : 
 pour tout
cocaract{\`e}re $\tau:\GM_m\To{} {\mathbf T}$, 
$$ \delta(\tau(\varpi_F)):= q_F^{-{1\over 2}\la \rho,\tau \ra} $$
o{\`u} $\val_F$, $\varpi_F$ et $q_F$ sont respectivement la valuation,
une uniformisante et le cardinal du corps r{\'e}siduel de $F$, et $\la
.,.\ra : X^*({\mathbf T}) \times X_*({\mathbf T})$ est 
l'accouplement canonique.

Comme on s'est ramen\'e au cas adjoint, on peut consid\'erer la
famille  $\{\omega_\alpha\}_{\alpha\in S}$ 
 des co-poids fondamentaux. On a 
\begin{eqnarray*} \left(w(\delta)_{|A_J}=\delta_{|A_J}\right) & \hbox{
    \ssi\ }  & \left(\forall \alpha \in
J^c,\;
w(\delta)(\omega_\alpha(\varpi_F)) = \delta(\omega_\alpha(\varpi_F))\right)
\\ & \hbox{
    \ssi\ }  & (\forall \alpha \in
J^c,\; q_F^{\frac{1}{2}\la \rho-w(\rho),\omega_\alpha \ra} =1)
\end{eqnarray*}
D'autre part, identifiant $X_*({\mathbf T})$ avec un r{\'e}seau de $X$,
 on  a
 \begin{eqnarray*} w(\rho)_J =\rho_J & \hbox{
    \ssi\ }  & \forall \alpha \in
J^c,\;
 \la \rho-w(\rho),\omega_\alpha \ra =0
\end{eqnarray*}

Il nous suffira donc de montrer que pour tout $\alpha\in S$, on a 
$$ q_F^{\frac{1}{2}\la \rho-w(\rho),\omega_\alpha \ra} =1  \Rightarrow \la
\rho-w(\rho),\omega_\alpha \ra =0. $$
Pour cela, on remarque que si  $w_0$ est l'{\'e}l{\'e}ment de plus grande
longueur de $W$,  alors $\frac{1}{2}{(\rho-w(\rho))}$ est le
d{\'e}terminant de l'action de ${\mathbf T}$ sur $\hbox{Lie}({\mathbf U}_+\cap
{\mathbf U}_+^{ww_0})$, de sorte que 
$$0\leq \frac{1}{2}\la
\rho-w(\rho),\omega_\alpha \ra \leq \la \rho,\omega_\alpha \ra = n_\alpha .$$

\bigskip

\noindent{\em Deuxi{\`e}me {\'e}tape : $w(\rho)_J=\rho_J \Rightarrow w\in W(J)$.} 

 Notons
momentan{\'e}ment $a^{J*}$ le noyau de $X^* \To{} a_J^*$ ;  c'est aussi le
sous-espace vectoriel de $X^*$ engendr{\'e} par les racines $\alpha\in J$.
La projection $X^* \To{} a_J^*$ admet une section canonique gr\^ace au diagramme suivant
$$ a_J^* \buildrel\sim\over\longleftarrow X^*({\mathbf M_J})\otimes \RM \To{res} X^*({\mathbf M_\emptyset}) \otimes \RM = X^*, $$
et on notera encore $a_J^*$ l'image de cette section. Si $[-,-]$ est
un produit scalaire sur $X$ invariant sous le groupe de Weyl $W$, alors
$a_J^*$ est l'orthogonal de $a^{J*}$.
On a alors les d\'ecompositions orthogonales 
$\rho=\rho^J+\rho_J \in a^{J*}\oplus a_J^*$ o\`u $\rho^J$ 
est la demi-somme des racines positives de $M_J$, et 
$w(\rho)=w(\rho)_J+w(\rho)^J$ 

Soit $v$ un \'el\'ement de $W(J)$ tel que $v(w(\rho)^J)$ soit dans la chambre de
Weyl positive de $a^{J*}$, c'est {\`a} dire : $\forall \alpha\in J, [
\alpha, v(w(\rho)^J) ] \geq 0$.
Comme le produit scalaire est $W$-invariant, on a $v(w(\rho)^J)=v(w(\rho))^J$ et $v(w(\rho))_J =w(\rho)_J$ donc $v(w(\rho))_J=\rho_J$ par notre
hypoth{\`e}se. 
Comme par ailleurs $vw(\rho)\in \rho + \sum_{\alpha>0} \RM_- \alpha$,
on obtient en soustrayant $\rho_J$ et compte tenu de $a^J\cap
\sum_{\alpha>0} \RM_- \alpha = \sum_{\alpha\in J} \RM_-\alpha$ :
$$ vw(\rho)^J \in \rho^J +\sum_{\alpha \in J} \RM_- \alpha$$
En prenant le produit scalaire avec $vw(\rho)^J$ on obtient donc :
$$ ||vw(\rho)^J||^2 \leq [ vw(\rho)^J,\rho^J ] $$
Mais puisque $||vw(\rho)^J||=||\rho^J||$, il s'ensuit que
$vw(\rho)^J=\rho^J$. 

On a donc obtenu $vw(\rho)=\rho$, ce qui {\'e}quivaut {\`a} $vw=1$.

\end{proof}

\subsection{Preuve du th\'eor\`eme \ref{explicite}}

\alin{Acyclicit\'e des complexes $\CC_*(\VC_{I,J})^+$} Nous avons d\'eja trait\'e le cas o\`u $R$ est un corps alg\'ebriquement clos dans la section pr\'ec\'edente. Nous allons nous ramener \`a ce cas-l\`a. Pour cela, rappelons le 
\begin{fact} \label{faitlibre}
Les $R$-modules sous-jacents aux repr\'esentations $\pi_I^R$ sont libres.
\end{fact}
\begin{proof}
Dans le cas ${\mathbf G}=GL(n)$ ceci est prouv\'e dans  \cite[Cor. 4.5]{SS1}
(et m\^eme pour $R=\ZM$). 
L'argument de {\em loc. cit. } repose sur la d{\'e}composition de Bruhat et 
fonctionne de la m{\^e}me mani{\`e}re pour n'importe quel groupe
d{\'e}ploy{\'e}.
\end{proof}
Soit maintenant $H$ un pro-$p$-sous-groupe ouvert de $G$, distingu\'e dans un compact sp\'ecial $H_0$ tel que $G=H_0B_+$, et ayant des d\'ecompositions d'Iwahori relativement \`a chaque \para semi-standard. Ces hypoth\`eses permettent de calculer facilement les $H$-invariants d'une induite parabolique. On sait que la famille de ces sous-groupes engendre la topologie de $G$ et il suffit donc de prouver l'acyclicit\'e des complexes $(\CC_*(\VC_{I,J}^R)^+)^H$ pour tout tel $H$.
Mais gr\^ace au fait ci-dessus et au lemme \ref{lemmeimage},  on s'aper{\c c}oit que le complexe $(\CC_*(\VC_{I,J}^R)^+)^H$ est un complexe de $R$-modules projectifs de type fini.
Comme on sait, \ref{chgtbase}, que le complexe est compatible aux changements de scalaires, il nous suffit donc de traiter le cas universel $R_u:=\ZM[\frac{1}{\sqrt{p}N_G}]$. 
Dans ce dernier cas, on a pour tout nombre premier $l$ banal une suite exacte de complexes
$$ 0\To{} \CC_*(\VC_{I,J}^{R_u})^+ \To{\times l} \CC_*(\VC_{I,J}^{R_u})^+ \To{} \CC^*(\VC_{I,J}^{R_u/l})^+ \To{} 0.$$
Comme $R_u/l$ est isomorphe \`a $\FM_l$ ou $\FM_l\times \FM_l$ et qu'on a d\'eja trait\'e le cas d'un corps de coefficients fortement banal, le complexe de droite est acyclique.
On en d\'eduit que la multiplication par $l$ est un isomorphisme sur les $R_u$-modules de cohomologie $\HC_*(\VC_{I,J}^{R_u})$. Comme ces derniers sont de type fini, il r\'esulte du lemme de Nakayama qu'ils sont nuls.


\alin{Cup-produits et suites spectrales} \label{yonss}
Afin de prouver le point ii) du th\'eor\`eme \ref{explicite}, nous rappelons quelques sorites sur les groupes
d'extensions, valables dans un contexte beaucoup plus g{\'e}n{\'e}ral que
le notre.

Soit $\CC$ une cat{\'e}gorie ab{\'e}lienne ayant assez 
d'injectifs et soient $\pi,\sigma$ et $\rho$ trois objets de $\CC$. On
se donne aussi un complexe acyclique
$$ C_n=\pi \To{} C_{n-1} \To{} \cdots \To{} C_0 \To{} C_{-1}=\sigma. $$
Ce complexe d{\'e}finit un {\'e}l{\'e}ment
$\alpha\in\ext{n}{\sigma}{\pi}{\CC}$ et on s'int{\'e}resse au
cup-produit par $\alpha$ :
$$ \alpha\cup - \; :\;\; \ext{*}{\rho}{\sigma}{\CC} \To{}
\ext{n+*}{\rho}{\pi}{\CC} $$
(Lorsqu'on pense le $\cup$-produit comme une simple composition de morphismes, l'application ci-dessus est bien le $\cup$-produit {\em \`a gauche} par $\alpha$).

Nous voulons ici souligner comment ce cup-produit
se lit sur la suite spectrale
$$ {E}_1^{pq} = \ext{q}{\rho}{C_p}{\CC} \Rightarrow
\ext{q-p}{\rho}{\sigma}{\CC}$$
obtenue en appliquant le foncteur $\hom{\rho}{-}{\CC}$ \`a une r\'esolution injective du complexe $C_\bullet$ et en filtrant le bicomplexe obtenu par les colonnes.
Pour cela, il faut se rappeler que puisque les diff\'erentielles $d_r$ sont de degr\'es $(-1-r,-r)$, il y a sur le bord $p=n$ de la suite spectrale des inclusions 
$ E_\infty^{n*} \injo E_1^{n*}$.
D'autre part, le terme $E_\infty^{n*}$ est le dernier quotient de la filtration du terme $*-n$ de l'aboutissement 
et on a donc une projection canonique 
$\ext{*-n}{\rho}{\sigma}{\CC} \twoheadrightarrow E_\infty^{n*}$.
La compos{\'e}e de ces deux applications
$$  \ext{*}{\rho}{\sigma}{\CC} \To{} {E}_\infty^{n,n+*} \To{} {E}_1^{n,n+*}= \ext{n+*}{\rho}{\pi}{\CC}$$
est justement le cup-produit \`a gauche par $\alpha$.

\alin{Preuve de \ref{explicite} ii)}
Nous utilisons les remarques du paragraphe \ref{yonss} pr\'ec\'edent : le complexe acyclique $\CC_*(\VC_{I,J})^+$ explicit\'e en \ref{compl} fournit la suite spectrale 
$$
  {E}_1^{pq} =
  \ext{q}{\pi_H}{\bigoplus_{\tiny\begin{array}{c} \delta(I,K)=p \\ \delta(J,K)=\delta(I,J)-p\end{array}} \im(J_{J^c|K})}{G} \Rightarrow
  \ext{q-p}{\pi_H}{\pi_I}{G} $$
  Le cup-produit {\`a} gauche par $\alpha_{I,J}$ consid\'er\'e dans \ref{explicite} ii) s'identifie alors {\`a} la
  compos{\'e}e
  $$
  \ext{*}{\pi_H}{\pi_I}{G} \To{} {E}_\infty^{\delta,\delta+*} \To{}
  {E}_1^{\delta,\delta+*} = \ext{\delta+*}{\pi_H}{\pi_{J}}{G}$$
o{\`u} on a pos{\'e} $\delta:=\delta(I,J)$ (on ne confondra pas avec le caract\`ere $\delta$ de $B_+$ !).
Nous allons montrer que pour tout $p\neq \delta$ et tout $q$, on a
${E_1}^{pq}=0$. Ceci impliquera en particulier la d\'eg\'en\'erescence en $E_1$ de la suite spectrale et donc le fait que la compos{\'e}e ci-dessus
est un isomorphisme pour tout $*\in\ZM$, ce qui est bien ce qu'on cherche \`a prouver.

Soit donc $K\subseteq S$ tel que 
$$ \delta(I,K) \neq \delta \;\; \hbox{ et } \;\; \delta(I,J)=\delta(I,K)+\delta(K,J),$$
ou ce qui est \'equivalent par \ref{comb1}, tel que
$$ K\neq J \;\; \hbox{ et }\;\; \Delta(K,J)\subseteq \Delta(I,J). $$
Choisissons un \para $P_{K,J^c}$ comme dans le lemme \ref{lemmeimage} dont nous reprenons les notations (notamment $\Delta:=\Delta(K,J^c)$).
Par r{\'e}ciprocit{\'e} de Frobenius-Shapiro (c'est une cons\'equence formelle de la r\'eciprocit\'e de Frobenius et de l'exactitude des foncteurs paraboliques)
on a 
$$ \ext{*}{\pi_H}{\im(J_{J^c|K})}{G} =
\ext{*}{\rp{G}{{P_{K,J^c}}}{\pi_H}}{\delta_{P_{\Delta}}^{-\frac{1}{2}}\pi_{_{K\cap\Delta,\Delta}}}{M_{\Delta(J^c,K)}}.$$
Or d'apr{\`e}s le lemme \ref{coro2} ci-dessous, le terme de
droite est non nul seulement si 
$\Delta(H^c,J^c)=\Delta(H,J) \subseteq \Delta(J^c,K)$.
Il est donc non nul seulement si
$\Delta(J^c,K)\supseteq \Delta(H,J)\cup \Delta(J^c,I)$. Or, l'hypoth\`ese $\delta(H,I)+\delta(I,J)=\delta(H,J)$ \'equivaut par \ref{comb1} \`a $\Delta(I,J) \subseteq \Delta(H,J)$. Puisque $\Delta(J^c,I)=\Delta(J,I)^c$, il s'ensuit que le terme de droite est non-nul seulement si $\Delta(J^c,K)=S$, c'est-\`a-dire
si $K=J$, ce que nous avons exclu.

\begin{lemme} \label{coro2}
Fixons $K,J$ et $I$ des sous-ensembles de $S$. Soit $P_{K,J}$ le \para
associ\'e \`a une chambre parabolique de $a_{\Delta(J,K)}$ contenue
  dans $\o{X_J}\cap \o{X_K}$ (comme  dans le lemme \ref{lemmeimage}).
Supposons enfin que $R$ est un anneau fortement banal pour $G$.
\begin{center} Si $\ext{*}{\rp{G}{{P_{K,J}}}{\pi_I}}{\delta_{P_{\Delta}}^{-\frac{1}{2}}\pi_{_{K\cap\Delta,\Delta}}}{M_{\Delta}} \neq 0$, alors  $\Delta(I^c,J)\subseteq \Delta(J,K)$.
\end{center}

\end{lemme}
\begin{proof}
Commen{\c c}ons par 
l'observation g{\'e}n{\'e}rale suivante : si $\CC$
est une cat{\'e}gorie ab{\'e}lienne avec assez d'injectifs et de
projectifs, et $\ZG$ d{\'e}signe le centre de la cat{\'e}gorie $\CC$,
alors les groupes $\ext{*}{V}{W}{\CC}$ sont naturellement et
canoniquement des $\ZG$-modules. De plus si $z\in \ZG$ agit par un
scalaire sur $V$, resp. $W$, il agit par ce m{\^e}me scalaire sur
$\ext{*}{V}{W}{\CC}$. 

\medskip

{\em Strat\'egie} :
Dans notre cas on prend $\CC=\Mo{R}{M_\Delta}$ et on
va produire, sous l'hypoth{\`e}se $\Delta(I^c,J)\setminus \Delta(J,K)\neq \emptyset$, 
 un {\'e}l{\'e}ment de $R[A_\Delta]$ agissant par $0$ sur $\rp{G}{{P_{K,J}}}{\pi_I}$ et par
l'identit{\'e} sur $\delta_{P_{\Delta}}^{-\frac{1}{2}}\pi_{_{K\cap\Delta,\Delta}}$.
La nullit\'e des ${Ext}^*$ entre ces deux objets en r\'esultera imm\'ediatement.


\medskip

{\em Premi\`ere \'etape :}
Tout
caract\`ere  $A_\Delta\To{} R^\times$ induit un morphisme de
$R$-alg\`ebres $R[A_\Delta]\To{} R$. C'est le cas par exemple
pour les restrictions des caract\`eres $w(\delta)$ \`a $A_\Delta$.
D'apr\`es le lemme \ref{comb2}, on a
\ini
\begin{equation} \label{R}
\hbox{Si}\;\;  \ker\nolimits_{R[A_\Delta]}(\delta_{|A_\Delta}) +
  \ker\nolimits_{R[A_\Delta]}(w(\delta)_{|A_\Delta}) \subsetneq R[A_\Delta],\;\hbox{
    alors }\;\; w\in W(\Delta).
\end{equation}
En effet, supposons que la somme des deux noyaux ci-dessus est un id\'eal propre de $R[A_\Delta]$ et 
choisissons un id\'eal maximal $\MC$ de $R$ contenant l'image de $\ker(w(\delta)_{|R[A_\Delta]})$ par $\delta_{|A_\Delta}$. Alors on a
$$\ker\nolimits_{R/\MC[A_\Delta]}(\delta_{|A_\Delta}) +
  \ker\nolimits_{R/\MC[A_\Delta]}(w(\delta)_{|A_\Delta}) \subsetneq R/\MC[A_\Delta]$$
o\`u l'on note encore $\delta$ et $w(\delta)$ pour leurs compos\'ees avec la projection $R\To{} R/\MC$.
Puisque les deux noyaux ci-dessus sont des id\'eaux maximaux, ils co\"incident et on a $$\delta_{|A_\Delta}=w(\delta)_{|A_\Delta} \; \hbox{mod.} \; \MC.$$ On peut donc appliquer \ref{comb2} au corps $R/\MC$.

\medskip

{\em Deuxi\`eme \'etape :} Comme dans la preuve du lemme \ref{lemmeimage}, on choisit un \borel
$B_J$ contenu dans $P_{K,J}$ et tel que $C_J:=C(B_J)\subset X_J$. 
Nous allons montrer que
$$
\hbox{Si }\; 
\{w\in W,\; w(C_J)\subset X_{I^c}\} \cap W(\Delta) \neq \emptyset, \hbox{ alors }\; \Delta(I^c,J)\subseteq\Delta(J,K).
$$
En effet, soit $w$ dans l'intersection ci-dessus. Puisque $w\in W(\Delta)$, il fixe $a_{\Delta(J,K)}$ points par points.
On a donc $\o{X_{I^c}}\cap \o{X_J} \supset w(\o{C_J)}\cap \o{C_J} \supset \o{C_J} \cap
a_{\Delta(J,K)}$. Autrement dit, $\o{X_{I^c}}\cap\o{X_J}$ contient un
c{\^o}ne vectoriellement g{\'e}n{\'e}rateur de $a_{\Delta(J,K)}$ (par notre choix de $B_J$). Comme on sait par ailleurs, \ref{1.3}, que 
$\o{X_{I^c}}\cap\o{X_J}$ est un c{\^o}ne ferm{\'e} g{\'e}n{\'e}rateur de
$a_{\Delta(J,I^c)}$, il s'ensuit que $a_{\Delta(J,K)} \subseteq
a_{\Delta(I^c,J)}$ et donc que $\Delta(I^c,J) \subseteq \Delta(J,K)$.

\medskip

{\em Troisi\`eme \'etape :} Supposons dor\'enavant que $\Delta(I^c,J)\setminus \Delta(J,K)\neq \emptyset$, de sorte que, par l'\'etape pr\'ec\'edente,
$\{w\in W,\; w(C_J)\subset X_{I^c}\} \cap W(\Delta) = \emptyset$.
Appliquons alors \ref{R} dans le cas "universel" de l'anneau $R_u=\ZM[\frac{1}{\sqrt{p}N_WN_G}]$. On obtient une d\'ecomposition 
$$ 1_{R_u[A_\Delta]} = z_\delta + z^\delta \in \ker\nolimits_{R_u[A_\Delta]} (\delta_{|A_\Delta}) + \prod_{w(C_J)\subset X_{I^c}} \ker\nolimits_{R_u[A_\Delta]}(w(\delta)_{|A_\Delta}).
$$

Comme tout anneau fortement banal $R$ re{\c c}oit $R_u$, les \'el\'ements $z_\delta$ et $z^\delta$ induisent des \'el\'ements correspondants dans $R[A_\Delta]$. Il est clair que l'action de $z_\delta$ sur $\delta_{P_{\Delta}}^{-\frac{1}{2}}\pi_{_{K\cap\Delta,\Delta}}$ est nulle.
Montrons que celle de $z^\delta$ sur $\rp{G}{{P_{K,J}}}{\pi_I}$ est nulle aussi. 
Par extension des scalaires, il suffit de traiter le cas "universel" $R_u$. Soit $H$ un pro-$p$-sous-groupe ouvert de $M_\Delta$, on sait que $\rp{G}{{P_{K,J}}}{\pi_I}^H$ est un $R_u$-module  de type fini, puisque c'est le conoyau de la fl\`eche
$$  \bigoplus_{L\supset I} \left(\Rp{G}{P_{K,J}}\circ \ip{P_L}{G}{\delta_{P_L}^{-\frac{1}{2}}}\right)^H \To{} \left(\Rp{G}{P_{K,J}}\circ \ip{P_I}{G}{\delta_{P_I}^{-\frac{1}{2}}}\right)^H,$$
et qu'il est sans torsion, puisque par \ref{serprin} v) on a l'inclusion 
$$ \rp{G}{{P_{K,J}}}{\pi_I}^H \injo \rp{G}{{P_{K,J}}}{I_{I^c}}^H $$ et 
que le terme de droite est sans torsion par la formule de Mackey.
 Pour voir que l'action de $z^\delta$ est nulle, il suffit donc de le v\'erifier apr\`es extension des scalaires $R_u\injo \CM$. Par fid\'elit\'e du foncteur $\Rp{M_\Delta}{B_J\cap M_\Delta}$, il suffit encore de v\'erifier que $z^\delta$ annule $\rp{G}{{P_{K,J}}}{\pi^\CM_I}$.
Mais ceci r\'esulte du calcul des exposants de  $\rp{G}{B_J}{\pi_I}$ sur un corps, effectu\'e au-dessus de \ref{exp1}. 

\end{proof}

\subsection{Preuve des lemmes \ref{serprin} et \ref{serprincorps}} \label{preuveserprin}
Nous utiliserons ici les bases de la th\'eorie des op\'erateurs d'entrelacements sur un anneau de coefficients g\'en\'eral qui sont d\'ecrites dans \cite{nutemp}.
Soit $B$ un \borel et $\o{B}$ son oppos\'e.
Suivant la terminologie de {\em loc. cit.} 2.10, un caract\`ere $\chi$ de $A_\emptyset$ est dit $(B,\o{B})$-r\'egulier si 
$$ \hbox{Ann}_{R[A_\emptyset]}(\chi) + \bigcap_{1\neq w\in W} \hbox{Ann}_{R[A_\emptyset]}(\chi^w) = R[A_\emptyset] . $$
Ici le caract\`ere $\chi$ est vu comme un caract\`ere de la $R$-alg\`ebre $R[A_\emptyset]$ et la notation $\hbox{Ann}$ d\'esigne l'id\'eal annulateur d'un $R[A_\emptyset]$-module.
Cette notion ne d\'epend pas du \borel $B$ et on dira donc simplement que $\chi$ est "r\'egulier". Elle est encore \'equivalente \`a 
$$ \forall w\in W,\;\; \hbox{Ann}_{R[A_\emptyset]}(\chi) + 
 \hbox{Ann}_{R[A_\emptyset]}(\chi^w) = R[A_\emptyset] . $$

\begin{lemme}
Lorsque $R$ est fortement banal pour $G$, le caract\`ere $\delta$ est r\'egulier au sens ci-dessus. Pour $G=GL(n)$, "banal" suffit.
\end{lemme}
\begin{proof}
Lorsque $R$ est un corps, ceci est le cas particulier $J=\emptyset$ du lemme \ref{comb2}. Nous donnons ici une preuve pour $R$ un anneau, qui permet d'am\'eliorer la restriction sur la caract\'eristique. 
Fixons $w\neq 1\in W$ et soit $\alpha\in S$ telle que $w(\alpha)\neq \alpha$. Notons 
$\alpha^\vee :\;\; F^\times \To{} A_\emptyset$ la coracine associ\'ee \`a $\alpha$. Alors 
$\delta(\alpha^\vee(\varpi)) = q $  et $\delta^w(\alpha^\vee(\varpi))= q^l$ o\`u $l$ est la somme des coefficients de $w^{-1}(\alpha)$ dans la base $S$. Soit $h$ le nombre de Coxeter du syst\`eme de racines de $G$, on a donc $|1-l|\leq h$. 
Supposons que l'entier $\prod_{0<r\leq h}(1-q^r)$ soit inversible dans $R$.
Alors
l'entier $q-q^{l}$ est aussi inversible dans $R$ de sorte que l'\'egalit\'e
$$ 1 = \frac{\alpha^\vee(\varpi) - q}{q^{l} -q} + \frac{q^{l}- \alpha^\vee(\varpi)}{q^{l}-q} $$
dans $R[A_\emptyset]$ montre bien que $\hbox{Ann}_{R[A_\emptyset]}(\delta) + 
 \hbox{Ann}_{R[A_\emptyset]}(\delta^w) = R[A_\emptyset] .$
L'hypoth\`ese sur $R$ ci-dessus est en g\'en\'eral plus faible que l'hypoth\`ese "bon et banal". Pour $GL(n)$ elle est \'equivalente \`a "banal", par la formule donnant le cardinal de $GL_n(\FM_q)$.
\end{proof}

\alin{Preuve de \ref{serprin} i), ii)  et iii)}
Fixons deux \borels $B$ et $B'$.
Par r\'eciprocit\'e de Frobenius, on a
$ \hom{\ip{B}{G}{\delta}}{\ip{B'}{G}{\delta}}{G} \simeq \hom{\Rp{G}{B'}\circ \ip{B}{G}{\delta}}{\delta}{A_\emptyset}$. Mais d'apr\`es le lemme g\'eom\'etrique \cite[2.12]{BZ2} et la r\'egularit\'e de $\delta$ du lemme pr\'ec\'edent, on a
$$ \Rp{G}{B'}\circ \ip{B}{G}{\delta} \simeq \bigoplus_{w\in W} \delta^w. $$
On en d\'eduit la premi\`ere assertion de \ref{serprin} i).
Pour l'existence et la d\'efinition d'un op\'erateur d'entrelacement "canonique" dans ces circonstances, nous renvoyons \`a \cite[2.11]{nutemp}.  Pour les propri\'et\'es \ref{serprin} ii) et iii) que ces op\'erateurs satisfont, nous renvoyons \`a la proposition 7.8. de {\em loc. cit}.

\alin{Preuve de \ref{serprincorps} i)} Dans ce paragraphe, $R$ est donc un corps alg\'ebriquement clos de caract\'eristique {\em bonne et banale} (ou simplement banale pour ${GL(n)}$). 
D'apr\`es \ref{faitbanal}, on a  une partition
$$ W.\delta = \bigsqcup_{\pi \in JH(\ip{B}{G}{\delta})} \exp(A_\emptyset,
\rp{G}{\o{B}}{\pi}).$$ 
En particulier, $JH(\ip{B}{G}{\delta})$ est sans multiplicit\'es.
Soit $\pi$ un sous-quotient irr\'eductible de $\ip{B_+}{G}{\delta}$,
et soit $w$ tel que $\delta^w\in
\exp(A_\emptyset,\rp{G}{\o{B_+}}{\pi})$. Par r\'eciprocit\'e de
Frobenius-Casselman, on a un morphisme non nul
$\ip{B_+}{G}{\delta^w}=\ip{B_+^{w^{-1}}}{G}{\delta} \To{} \pi$, ce qui
montre que $\pi$ est un quotient d'une induite du type
$\ip{B}{G}{\delta}$.  

Fixons maintenant $B$, alors par les m\^emes arguments, une repr\'esentation irr\'eductible  $\pi$ de $G$ est un quotient de $\ip{B}{G}{\delta}$ \ssi\ $\delta\in \exp(A_\emptyset,\rp{G}{\o{B}}{\pi})$. Par la partition ci-dessus, une telle repr\'esentation est unique (\`a isomorphisme pr\`es). Puisque $\hom{\ip{B}{G}{\delta}}{\ip{B_+}{G}{\delta}}{G}\neq 0$, c'est bien un sous-quotient de $\ip{B_+}{G}{\delta}$.

\bigskip

\noindent{\em Remarque :} Choisissons un ordre total $\leq$ sur l'ensemble
 $\PC(S)$ des sous-ensembles de $S$, raffinant l'ordre induit par la
relation d'inclusion. On obtient une filtration d\'ecroissante 
 $\hbox{Fil}_I:=\sum_{I\leq K} \ind{P_K}{G}{1}$ de  $\ind{B_+}{G}{1}$ dont le gradu\'e
est $\bigoplus_{I\subseteq S} \pi_I^R$ (on peut par exemple le v\'erifier en calculant les exposants de $\rp{G}{B_+}{\hbox{Fil}_I}$).
Il s'ensuit que dans le groupe de Grothendieck on a l'\'egalit\'e 
$$ \ind{B_+}{G}{1} = \sum_{I\subseteq S} \pi_I^R $$
et que les $\pi_I^R$ sont \`a supports disjoints. En particulier, on a
$\hbox{long}(\ind{B_+}{G}{1})\geq |\PC(S)|$.

\alin{Preuve de \ref{serprincorps} ii)}
$R$ est toujours un corps alg\'ebriquement clos (fortement) banal.
L'implication $(a)\Rightarrow (b)$ est tautologique.
 Supposons que $\ip{B}{G}{\delta}$ et $\ip{B'}{G}{\delta}$ ont
leurs uniques quotients irr{\'e}ductibles isomorphes et appelons $\pi$ (la classe d'isomorphisme de) ce quotient. Montrons que l'op{\'e}rateur
d'entrelacement $J_{B'|B}$, qui est non nul, doit {\^e}tre
surjectif. En effet, son image est non-nulle, et a pour (unique) quotient irr\'eductible $\pi$. Si son conoyau \'etait non-nul, il aurait aussi $\pi$ comme (unique) quotient irr\'eductible, contredisant la multiplicit\'e $1$ dans $JH(\ip{B'}{G}{\delta})$. Donc $J_{B'|B}$ est surjectif.
On en d\'eduit maintenant qu'il est injectif, puisque pour tout sous-groupe ouvert compact $H$ de $G$, on a 
$\dim_{R}(\ip{B}{G}{\delta}^H) =
\dim_{R}(\ip{B'}{G}{\delta}^H) = |B\ba G/H|$.
On a donc prouv\'e $(b) \Rightarrow (a)$.

Il nous reste maintenant {\`a} montrer que $J_{B'|B}$ est un isomorphisme
\ssi\ $C':=C(B')$ et $C:=C(B)$ sont dans un m{\^e}me $X_J$ pour un certain $J\subseteq S$.
{\'E}tudions d'abord le cas particulier o{\`u} les chambres $C$ et $C'$ sont
adjacentes, le mur {\'e}tant associ{\'e} {\`a} une racine $r$. 
Le \levi $M_r:=\ZC_G(\ker r)$ est un groupe r\'eductif d\'eploy\'e de rang $1$ (dont le groupe adjoint n'est autre que $PGL(2)$). 
Le caract\`ere $\delta$ est $(B\cap M_r,B'\cap M_r)$-r\'egulier ({\em cf} \cite[7.8.ii)]{nutemp}) et on dispose donc d'un op\'erateur d'entrelacement non nul $J_{B'_M|B_M}:\; \ip{B\cap M_r}{G}{\delta}\To{} \ip{B'\cap M_r}{G}{\delta}$. 
Apr{\`e}s avoir
  {\'e}tudi{\'e} les s{\'e}ries principales (en caract{\'e}ristique
  banale) de ces groupes
  d{\'e}ploy{\'e}s de rang $1$, voir par exemple \cite{Viggl2}, ou \cite[8.4.i)]{nutemp} coupl\'e \`a la formule pour la mesure de Plancherel des s\'eries principales de $SL(2)$ et $PGL(2)$, on sait que  
les assertions suivantes sont {\'e}quivalentes :
\begin{enumerate}
\item  La repr{\'e}sentation
  $\ip{B\cap M_r}{M_r}{\delta}$ est r{\'e}ductible.
\item
  L'op{\'e}rateur $J_{B'_M|B_M}$  n'est pas inversible.
\item $\delta \circ r^\vee = q_F^{\pm\val_F} $, {\'e}galit{\'e} de
  caract{\`e}res lisses $\GM_m(F)=F^\times \To{} R^\times$.
\end{enumerate}
Si $l(r)$ d\'esigne la somme des coefficients de $r$ dans la base $S$, alors
la derni\`ere assertion est encore \'equivalente \`a : $q^{l(r)}=q^{\pm 1}$. Comme on l'a d\'eja vu un peu plus haut, sous notre hypoth\`ese de caract\'eristique banale, ceci \'equivaut \`a $l(r)=\pm 1$ ou encore $r\in \pm S$.

En utilisant maintenant \ref{serprin} ii) et iii), on d\'eduit de la discussion pr\'ec\'edente  que $J_{B'|B}$ est inversible \ssi\ il existe une galerie
tendue entre $C'$ et $C$ dont tous les murs successifs sont distincts
des murs associ{\'e}s aux racines simples, autrement dit, une galerie
incluse dans une composante connexe de $X$ priv{\'e} des orthogonaux des
racines simples. Une telle composante connexe est un $X_J$ et on a
obtenu $(a)\Rightarrow (c)$. Pour
la r{\'e}ciproque, il faut encore v{\'e}rifier que deux chambres
quelconques dans $X_J$ sont reli{\'e}es par une galerie tendue dont tous les
membres sont dans $X_J$, mais ceci r{\'e}sulte de la convexit{\'e} de $X_J$.

\bigskip

\noindent{\em Remarque :} Par ce que l'on vient de prouver
et l'absence de multiplicit\'e dans $JH(\ind{B_+}{G}{1})$, on obtient
$\hbox{long}(\ind{B_+}{G}{1}) = |\PC(S)|$. On d\'eduit donc de la
remarque pr\'ec\'edente que les repr\'esentations $\pi_I^R$ sont
irr\'eductibles.

\alin{Preuve de \ref{serprin} iv)} Dans ce paragraphe, $R$ est un anneau fortement banal pour $G$ ou simplement banal pour $GL(n)$. Remarquons que l'assertion (a) de \ref{serprin} iv) est \'equivalente \`a "$J_{B'|B}$ est un isomorphisme". Comme la formation de $\ip{B}{G}{\delta}$ et la d\'efinition de $J_{B'|B}$ sont compatibles \`a l'extension des scalaires par un morphisme $R\To{} R'$, il suffit de  consid\'erer  le cas "universel" $R=\ZM[\frac{1}{\sqrt{p}N_GN_W}]$. 

Or, pour tout sous-groupe ouvert compact $H$, le sous-$R$-module $\ip{B}{G}{\delta}^H$ est libre de type fini. Il s'ensuit que $J_{B'|B}$ est un isomorphisme \ssi\ il l'est apr\`es r\'eduction \`a tout corps r\'esiduel de $R$ et donc \ssi\ il l'est apr\`es changement des  scalaires par un morphisme $R\To{} C$ o\`u $C$ est alg\'ebriquement clos (et n\'ecessairement de caract\'eristique banale). 

Ainsi l'\'equivalence que l'on doit montrer est une cons\'equence de celle de \ref{serprincorps} ii) montr\'ee ci-dessus.

\alin{Preuve de \ref{serprincorps} iii) et \ref{serprin} v)}
Dans ce paragraphe, $R$ est un anneau (fortement) banal qui sera parfois suppos\'e \^etre un corps.
Nous allons commencer par traiter le cas $J=S$ des \'enonc\'es que l'on veut  prouver. Dans ce cas, on peut expliciter l'op\'erateur $J_{B_+|\o{B_+}}$ par
$$\application{J_{B_+|\o{B_+}}:\;\;}
{\ip{\o{B_+}}{G}{\delta} = \ind{\o{B_+}}{G}{\delta_{\o{B_+}}}}
{\ind{B_+}{G}{1}=\ip{B_+}{G}{\delta}}
{f}{\int_{G/\o{B_+}} f(g)dg}.
$$
Sous cette forme, on remarque la factorisation
$$ J_{B_+|\o{B_+}}:\;\;\ip{\o{B_+}}{G}{\delta} \twoheadrightarrow 1=\pi_S^R \injo \ip{B_+}{G}{\delta}$$
qui est ce que l'on cherchait.

Fixons maintenant $J\subseteq S$ et choisissons un \borel $B_J$ tel que $C(B_J)\subset X_J$ et $\o{C(B_J)} \cap \o{X_\emptyset} = a_J\cap \o{X_\emptyset}$ (rappelons que $X_\emptyset$ est aussi la chambre de Weyl associ\'ee \`a $B_+$). Alors $B_J$ et $B_+$ sont contenus dans le \para standard $P_J$ de Levi $M_J$ et
par \ref{serprin} iii), on a un diagramme commutatif
$$\xymatrix{ \ip{B_J}{G}{\delta} \ar[r]^{J_{B_+|B_J}} \ar@{=}[d] & \ip{B_+}{G}{\delta} \ar@{=}[d] \\ \ip{P_J}{G}{\ip{\o{B_+}\cap M_J}{M_J}{\delta}} \ar[r]_{\ip{P_J}{G}{J'}} &
\ip{P_J}{G}{\ip{B_+\cap M_J}{M_J}{\delta}}
}
$$
o\`u $J'$ est l'op\'erateur d'entrelacement $J_{\o{B_+}\cap M_J|B_+\cap M_J} : \ip{\o{B_+}\cap M_J}{M_J}{\delta}\To{} \ip{{B_+}\cap M_J}{M_J}{\delta}$. Par la factorisation $\delta=\delta_{P_J}^{-\frac{1}{2}}\delta_{B_+\cap M_J}^{-\frac{1}{2}}$, et le cas $J=S$ trait\'e pr\'ec\'edemment, on obtient donc une factorisation
$$ J_{B_+|B_J}:\;\;\ip{B_J}{G}{\delta} \twoheadrightarrow \ip{P_J}{G}{\delta_{P_J}^{-\frac{1}{2}}}=\ind{P_J}{G}{1} \injo \ind{B_+}{G}{1}= \ip{B_+}{G}{\delta}.$$

Maintenant soit $K\supset J$. Comme dans le lemme \ref{lemme1} on peut trouver un \borel $B_K$ tel que $C(B_K)\subset X_K$ et  $\o{C(B_K)}\cap \o{X_\emptyset}=a_K\cap X_\emptyset$, et de plus $d(B_K,B_+)=d(B_K,B_J)+d(B_J,B_+)$. On obtient alors le diagramme commutatif suivant :
$$\xymatrix{
 \ip{B_J}{G}{\delta} \ar@{->>}[r] & \ind{P_J}{G}{1} \ar@{^(->}[r] & \ip{B_+}{G}{\delta} \ar[r]^{J_{\o{B_J}|B_+}} &  \ip{\o{B_J}}{G}{\delta} \\
 \ip{B_K}{G}{\delta} \ar[u]^{J_{B_J|B_K}} \ar@{->>}[r] & \ind{P_K}{G}{1} \ar@{^(->}[u]  &
& }$$

\begin{lemme}
Pour $K\neq J$, on a $J_{\o{B_J}|B_J}\circ J_{B_J|B_K} = 0$.
\end{lemme}
\begin{proof}
Soient $B$ et $B'$ deux \borels adjacents  tels que $d(B_J,B_K)=d(B_J,B)+1+d(B',B_K)$ et tels que $C(B)\subset X_J$ et $C(B') \cap X_J=\emptyset$.
Comme on a aussi $d(B_J,\o{B_J})=d(B_J,B)+1 +d(B',\o{B_J})$, la compos\'ee envisag\'ee 
s'\'ecrit encore
$$ J_{\o{B_J}|B_J}\circ J_{B_J|B_K} = J_{\o{B_J}|B'}\circ J_{B'|B}\circ J_{B|B_J}\circ J_{B_J|B}\circ J_{B|B'}\circ J_{B'|B_K}.$$
Par \ref{serprin} i), l'endomorphisme $J_{B|B_J}\circ J_{B_J|B}$ de $\ip{B}{G}{\delta}$ est un scalaire et commute donc au reste. Ainsi la compos\'ee envisag\'ee est de la forme
$$ J_{\o{B_J}|B_J}\circ J_{B_J|B_K}  = ? \circ J_{B'|B}\circ J_{B|B'} \circ ? $$
et il nous suffira de prouver que $J_{B'|B}\circ J_{B|B'}=0$. Par hypoth\`ese, la racine $\alpha$ dont le mur s\'epare $C(B)$ et $C(B')$ est dans $S$. Si $P$ est le \para contenant $B$ et $B'$, on a donc par \ref{serprin} iii)
$$ J_{B'|B}\circ J_{B|B'}= \ip{P}{G}{J_{\o{B}\cap M_\alpha|B\cap M_\alpha}J_{B\cap M_\alpha|\o{B}\cap M_\alpha}}.$$
On est donc ramen\'e \`a montrer que lorsque $G$ est de rang semi-simple $1$, on a $J_{\o{B_+}|B_+}\circ J_{B_+|\o{B_+}}=0$. Lorsque $R$ est un corps alg\'ebriquement clos de caract\'eristique banale, ceci est bien connu. Pour passer au cas g\'en\'eral, on se ram\`ene au cas universel $R=\ZM[\frac{1}{\sqrt{p}N_GN_W}]$. Dans ce dernier cas, la compos\'ee envisag\'ee est nulle puisque elle l'est sur tous les points de $\hbox{Spec} R$.

\end{proof}

Reprenons  la preuve de \ref{serprin} v) et \ref{serprincorps} iii) l\`a o\`u on l'a interrompue. Rappelons que $\pi_J$ est d\'efini par la suite exacte
$$ \bigoplus_{K\supset J} \ind{P_K}{G}{1} \To{} \ind{P_J}{G}{1} \To{} \pi_J \To{} 0.$$
Ainsi par le lemme et le diagramme qui le pr\'ec\`ede, on obtient une factorisation
$$ \xymatrix{ J_{\o{B_J}|B_J}:\;\;
\ip{B_J}{G}{\delta} \ar@{->>}[r] & \ind{P_J}{G}{1} \ar@{->>}[d] \ar[r]  & \ip{\o{B_J}}{G}{\delta} \\ & \pi_J^R \ar[ru]^{\theta_R} &
}.$$
Il reste maintenant \`a v\'erifier que $\theta_R$ est
injective. Lorsque $R$ est un corps, cela r\'esulte de
l'irr\'eductibilit\'e de $\pi^R_J$ ({\em cf} les deux remarques
ci-dessus) et de la non-nullit\'e de $J_{\o{B_J}|B_J}$. Pour passer au cas g\'en\'eral, remarquons que la formation de $\pi_J^R$ et $\theta_R$ est compatible au changement de base. De plus on sait par \ref{faitlibre} que pour tout sous-groupe ouvert compact, $(\pi_J^R)^H$ est un $R$-module libre, de m\^eme que $\ip{\o{B_J}}{G}{\delta}^H$.
La propri\'et\'e d'injectivit\'e n'est pas compatible au changement de base, mais la propri\'et\'e d'\^etre un plongement localement scindable l'est. 
Ainsi pour montrer que $\theta_R$ induit un plongement localement scindable sur les $H$-invariants, on peut se ramener au cas universel $R=\ZM[\frac{1}{\sqrt{p}N_G}]$, puis aux localis\'es de ce cas universel. Dans ce dernier cas o\`u $R$ est local principal, la condition de scindabilit\'e est \'equivalente \`a l'existence d'un mineur inversible pour la matrice de $\theta_R$ dans des bases quelconques de $(\pi_J^R)^H$ et $\ip{\o{B_J}}{G}{\delta}^H$. Cette existence est assur\'ee par r\'eduction modulo l'id\'eal maximal, en utilisant le cas des corps, d\'eja trait\'e.

\alin{Preuve de \ref{serprincorps} iv)}
Soit $B$ tel que $C(B)\subset X_J$, de sorte que $\pi_J$ est l'unique quotient irr\'eductible de $\ip{B}{G}{\delta}$. Alors sa contragr\'ediente $\pi_J^\vee$ est l'unique sous-repr\'esentation irr\'eductible de $\ip{B}{G}{\delta^{-1}}=\ip{w_0(B)}{G}{\delta}$. Mais par \ref{serprin} v), celle-ci est aussi l'image de $J_{w_0(B)|\o{w_0(B)}}$. Donc $\pi_J^\vee$ est l'unique quotient irr\'eductible de $\ip{\o{w_0(B)}}{G}{\delta}$ dont la chambre associ\'ee est $C(\o{w_0(B)})=-w_0(C(B))$ et est donc contenue dans $X_{-w_0(J)}$. Donc $\pi_J^\vee\simeq \pi_{-w_0(J)}$.

\section{Repr{\'e}sentations elliptiques et correspondances}\label{elliptiques}

Dans cette partie, nous rappelons des faits bien connus sur les
correspondances de Langlands et Jacquet-Langlands, d'autres un peu moins
connus mais qui le sont certainement des sp{\'e}cialistes, puis nous nous
consacrons {\`a} une description explicite dans le cas des
repr{\'e}sentations elliptiques.

\subsection{Correspondance de Jacquet-Langlands} \label{JacLan}

\alin{Notations}
Il sera commode de consid{\'e}rer  des repr{\'e}sentations {\`a}
coefficients dans un corps $C$ abstraitement isomorphe au corps des
complexes $\CM$. Ce corps pourra parfois {\^e}tre $\CM$ ou un $\o\QM_l$,
selon les besoins topologiques qu'on aura. Pour tout groupe $H$
localement profini, nous notons $\Irr{C}{H}$ l'ensemble des classes 
d'isomorphisme de $C$-repr{\'e}sentations {\em lisses } irr{\'e}ductibles de
$H$. Pour le groupe $G_d$, on isole certains sous-ensembles
remarquables :
\begin{enumerate}
\item On d{\'e}signe par  $\Cu{C}{G_d}\subseteq \Irr{C}{G_d}$ le sous-ensemble
  des repr{\'e}sentations {\em
    cuspidales}, {\em i.e.} dont les coefficients sont {\`a} support
  compact-modulo-le-centre.  
\item On note $\Disc{\CM}{G_d} \subseteq
  \Irr{\CM}{G_d}$ le sous-ensemble des repr{\'e}sentations ``de la s{\'e}rie
  discr{\`e}te'', {\em i.e.} dont les coefficients sont de carr{\'e} int{\'e}grable (au sens
  complexe) modulo-le-centre. 
Malgr{\'e} cette d{\'e}finition de nature analytique, il se trouve que la
notion de ``s{\'e}rie discr{\`e}te'' de $G_d$  est invariante par
automorphismes du corps $\CM$ ;  c'est une cons{\'e}quence des th{\'e}or{\`e}mes
de multiplicit{\'e}s limites de \cite{Rog}, ou plus prosa{\"\i}quement une
cons{\'e}quence de la classification de Bernstein-Zelevinski \cite[Thm
9.3]{Zel}. Cela nous permet d'isoler sans ambigu{\"\i}t{\'e} un sous-ensemble 
$\Disc{C}{G_d} \subseteq
  \Irr{C}{G_d}$ dont nous appellerons les membres ``s{\'e}ries discr{\`e}tes''
  par abus de langage.
\end{enumerate}

On peut consid{\'e}rer la correspondance de Jacquet-Langlands comme le
reflet spectral de la correspondance ``g{\'e}om{\'e}trique''
bien connue 
entre classes de conjugaison elliptiques semi-simples r{\'e}guli{\`e}res de
$G_d$ et de $\dd$, donn{\'e}e par l'{\'e}galit{\'e} des  polyn{\^o}mes caract\'eristiques. L'{\'e}nonc{\'e} spectral
classique concerne les repr{\'e}sentations complexes  : 

\begin{theo} (Correspondance de Jacquet-Langlands, \cite{DKV},\cite{Badu1}) \label{JL1}
  Il existe une bijection $$JL_d:\; \Irr{\CM}{\dd}\simto \Disc{\CM}{G_d}$$ caract{\'e}ris{\'e}e par
  l'{\'e}galit{\'e} de caract{\`e}res $\chi_{JL_d(\rho)}(g) =
  (-1)^{d-1}\chi_{\rho}(x)$ pour toutes classes elliptiques $g\in G_d,
  x\in\dd$ se  correspondant ({\em i.e.} ayant m{\^e}me polyn{\^o}me minimal
  de degr{\'e} $d$).
\end{theo}

Remarquons que l'{\'e}galit{\'e} de caract{\`e}res de  cet {\'e}nonc{\'e}
 peut se tester sur les fonctions localement constantes {\`a} support compact dans
l'ensemble (ouvert) des elliptiques r{\'e}guliers, et ne fait donc
intervenir que le caract{\`e}re-distribution des repr{\'e}sentations de
$\Disc{\CM}{G_d}$. Ce caract{\`e}re-distribution est d{\'e}fini sur n'importe
quel corps de coefficients, en particulier sur $C$. Compte tenu de ce
que l'ensemble $\Disc{\CM}{G_d}$ et la condition d'{\'e}galit{\'e} des
caract{\`e}res sont stables par l'action
des automorphismes du corps $\CM$, l'{\'e}nonc{\'e} ci-dessus se transf{\`e}re
sans ambigu{\"\i}t{\'e} {\`a} un {\'e}nonc{\'e} formellement analogue sur le corps $C$. 
Nous noterons encore $JL_d:\; \Irr{C}{\dd}\simto \Disc{C}{G_d}$ la
bijection obtenue.


\alin{Groupes de Grothendieck}
Notons maintenant $R(G_d)$ et $R(\dd)$ les groupes de Grothendieck des
$C$-repr{\'e}sentations de longueur finies. La bijection de Jacquet-Langlands 
induit une injection $R(\dd)\injo R(G_d)$, v{\'e}rifiant
l'{\'e}galit{\'e} de caract{\`e}res du th{\'e}or{\`e}me \ref{JL1} sur les classes de conjugaison
elliptiques se correspondant. On veut d{\'e}finir une r{\'e}traction pour cette injection,
v{\'e}rifiant la m{\^e}me {\'e}galit{\'e} de caract{\`e}res. 
Soit $\chi$ un caract{\`e}re lisse de $K^\times$, notons $R(G_d,\chi)$,
resp. $R(\dd,\chi)$,  le
sous-groupe de $R(G_d)$, resp. $R(\dd)$, engendr{\'e} par les irr{\'e}ductibles de caract{\`e}re
central $\chi$. 
Les d{\'e}compositions selon le caract{\`e}re central $R(G_d)=\bigoplus_\chi
R(G_d,\chi)$, resp. $R(\dd)=\bigoplus_\chi R(\dd,\chi)$, sont respect{\'e}es par
l'application $JL$. De plus, les groupes $R(G_d,\chi)$ et $R(\dd,\chi)$
 sont munis de formes bilin{\'e}aires enti{\`e}res, voir \cite{SS2},
$$\application{\la\;.\;\ra:\;}{R(G_d,\chi)\times
  R(G_d,\chi)}{\ZM}{(\pi,\pi')}{\sum_i (-1)^i\dim(\ext{i}{\pi}{\pi'}{G_d,\chi})}$$
et respectivement pour $\dd$. 
Ces formes bilin{\'e}aires en induisent une sur la somme
directe $R(G_d)$, resp. $R(\dd)$, que l'on note de la m{\^e}me mani{\`e}re.
Dans le cas de $\dd$, cette forme est non
d{\'e}g{\'e}n{\'e}r{\'e}e et on a simplement $\la\rho,\rho'\ra=\dim\hom{\rho}{\rho'}{\dd}$.
Dans le cas de $G_d$, la situation est sensiblement plus compliqu{\'e}e.
Notons 
$R_I(G_d)$ le sous-groupe de $R(G_d)$ engendr{\'e} par les induite
paraboliques de repr{\'e}sentations de sous-groupes de Levi propres et 
$\o{R}(G_d)$ le quotient $R(G_d)/R_I(G_d)$.

\begin{lemme}
 $\o{R}(G_d)$ est libre sur $\ZM$ et  la forme bilin{\'e}aire $\la\;,\;\ra$  s'y descend 
 et y est non-d{\'e}g{\'e}n{\'e}r{\'e}e, une base orthonormale {\'e}tant donn{\'e}e par les
 images de s{\'e}ries discr{\`e}tes. 
\end{lemme}
\begin{proof} 
Pour cette preuve, nous pouvons identifier $C$ et $\CM$.
L'{\'e}nonc{\'e} est alors une cons{\'e}quence des trois propri{\'e}t{\'e}s suivantes :
  \begin{enumerate}
  \item La classification de Zelevinski \cite{Zel}(ou celle dite du quotient de
    Langlands) montre que
$$R(G_d) =\bigoplus_{(M,\sigma,\psi)} \ZM [\ip{M}{G_d}{\sigma\psi}] $$
o{\`u} $[?]$ est l'{\'e}l{\'e}ment de $R(G_d)$ associ{\'e} {\`a} la repr{\'e}sentation $?$, les
triplets $(M,\sigma,\psi)$ sont form{\'e}s d'un \levi standard de $G_d$
(c'est-{\`a}-dire un produit de $GL_{d_i}$ diagonaux), d'une
repr{\'e}sentation temp{\'e}r{\'e}e de $M$ et d'un caract{\`e}re non-ramifi{\'e} de $M$ ``dans la
chambre de Weil positive'', et sont pris {\`a} $G_d$-conjugaison pr{\`e}s. Le
signe $\Ip{M}{G_d}$ est l'induction parabolique normalis{\'e}e le long du
parabolique triangulaire sup{\'e}rieur dont le Levi est $M$.
Comme on sait de plus que les repr{\'e}sentations temp{\'e}r{\'e}es de $G_d$ sont
soit induites, soit des s{\'e}ries discr{\`e}tes, on en d{\'e}duit
$$R(G_d)= \left(\bigoplus_{\pi\in \Disc{C}{G_d}} \ZM[\pi]\right) \bigoplus
R_I(G_d),$$
ce qui montre que $\o{R}(G_d)$ est libre sur $\ZM$ et qu'une base en est
donn{\'e}e par les images des s{\'e}ries discr{\`e}tes.
\item  Soit $M$ un \levi standard de $G_d$ et $\pi$ une repr{\'e}sentation
  admissible de  $G_d$. Alors pour toute repr{\'e}sentation admissible de
  $G_d$, un argument de d{\'e}formation attribu{\'e} {\`a} Kazhdan dans \cite{SS2} montre que
$$\la \pi,\ip{M}{G_d}{\sigma} \ra = \la \ip{M}{G_d}{\sigma},\pi \ra =0.$$
\item Soit $\pi,\pi'$ deux s{\'e}ries discr{\`e}tes de $G_d$, on a par \cite{Vigext}
$$ \la \pi,\pi'\ra =\delta_{\pi\pi'}.$$
  \end{enumerate}
\end{proof}

D'autre part, la formule des caract{\`e}res induits de Van Dijk
\cite{vdijk} montre  que l'application qui {\`a} un {\'e}l{\'e}ment $x\in R(G_d)$
associe la restriction de son caract{\`e}re-distribution aux {\'e}l{\'e}ments
elliptiques r{\'e}guliers se factorise {\`a} travers le quotient
$\o{R}(G_d)$. On d{\'e}duit alors du th{\'e}or{\`e}me \ref{JL1} le

\begin{coro} \label{JL2}
  L'application $R(\dd)\To{JL_d} R(G_d)$ induit un isomorphisme isom{\'e}trique
  $$R(\dd)\simto \o{R}(G_d)$$ caract{\'e}ris{\'e} par l'{\'e}galit{\'e} de caract{\`e}res du
  th{\'e}or{\`e}me \ref{JL1}.
\end{coro}

En cons{\'e}quence, on a un morphisme dans l'autre sens $R(G_d)\To{}
\o{R}(G_d) \simto R(\dd)$ que nous noterons encore $JL_d$ ! Ceci permet de
d{\'e}finir $JL_d(\pi)$ pour toute repr{\'e}sentation de $G_d$. C'est {\em a
  priori} seulement un {\'e}l{\'e}ment de $R(\dd)$.

\begin{lemme} \label{defell}
Soit $\pi$ une $C$-repr{\'e}sentation irr{\'e}ductible de $G_d$. Les
  propri{\'e}t{\'e}s suivantes sont {\'e}quivalentes :
  \begin{enumerate}
  \item $\pi$ a le m{\^e}me support cuspidal qu'une s{\'e}rie discr{\`e}te.
  \item $\pi$ a un caract{\`e}re non nul sur les {\'e}l{\'e}ments elliptiques
    semi-simples r{\'e}guliers.
  \item L'image de $\pi$ dans $\o{R}(G_d)$ est non nulle.
  \end{enumerate}
 \end{lemme}

Rappelons que le support cuspidal d'une repr{\'e}sentation
irr{\'e}ductible $\pi$ est l'unique classe de conjugaison de couples
$(M,\tau)$ form{\'e}s d'un \levi $M$ et d'une repr{\'e}sentation cuspidale
irr{\'e}ductible $\tau$ de $M$ qui appara{\^\i}t dans le module de Jacquet
normalis{\'e} de  $\pi$ le long d'un parabolique de Levi $M$.

\medskip

\begin{proof}
L'implication $ii) \Rightarrow iii)$ est une cons{\'e}quence de la formule
de Van Dijk \cite{vdijk} qui montre que le caract{\`e}re d'une induite
parabolique en un {\'e}l{\'e}ment semi-simple r{\'e}gulier elliptique est nul.

Pour voir l'implication $iii) \Rightarrow i)$,
{\'e}crivons  $[\pi]=\sum_{(M,\sigma)} [\ip{M}{G_d}{\sigma}]$ (somme
d'induites de repr{\'e}sentations essentiellement temp{\'e}r{\'e}es) comme
nous le permet la classification par le quotient de Langlands.
On peut supposer que les supports cuspidaux de chaque $\sigma$ sont
contenus dans celui de $\pi$.
Comme les temp{\'e}r{\'e}es sont soit induites, soit des s{\'e}ries
discr{\`e}tes, on voit que si $\pi\notin R_I(G_d)$,
alors il y a une s{\'e}rie 
discr{\`e}te de m{\^e}me support cuspidal que $\pi$. Donc $iii)
\Rightarrow i)$.

Pour l'implication $i)\Rightarrow ii)$, on peut utiliser la
combinatoire de la classification de Zelevinski. Soit $\pi^{disc}$
l'unique s{\'e}rie discr{\`e}te de m{\^e}me support cuspidal que
$\pi$. Avec les notations de \cite{Zel}, on sait qu'il existe 
un (unique) segment 
$\Delta=[\tau,\tau']$ o{\`u} $\tau$ est une cuspidale irr{\'e}ductible de
$G_{d'}$ avec $d'$ diviseur de $d$, de sorte que 
\begin{enumerate}
\item L'ensemble partiellement ordonn{\'e} des multisegments de support
  $[\tau,\tau']$ a pour plus petit {\'e}l{\'e}ment $a_{min}=\{\Delta\}$ et
  plus grand {\'e}l{\'e}ment $a_{max}=\{\{\tau'\},\cdots, \{\tau\}\}$.
\item $\pi^{disc}=\la a_{max} \ra$ et toutes les repr{\'e}sentations de
  m{\^e}me support cuspidal sont de la forme $\la a \ra$ pour un
  multisegment $a$ de support $[\tau,\tau']$.
\item La repr{\'e}sentation irr{\'e}ductible $\la b \ra$ associ{\'e}e {\`a} $b$ appara{\^\i}t comme
  sous-quotient de la repr{\'e}sentation induite $\pi(a)$ associ{\'e}e {\`a} $a$
  \ssi\ $b\leq a$, et dans ce cas sa multiplicit{\'e} est $1$. 
\end{enumerate}
Supposons maintenant que $\pi=\la a \ra$.
Pour $b< a$ notons $d(b,a)$ la longueur
d'une cha{\^\i}ne maximale $b<m_1<\cdots< a$ ({\em i.e} le nombre
d'op{\'e}rations ``{\'e}l{\'e}mentaires'' au sens de \cite[7]{Zel} pour passer de
$a$ {\`a} $b$), 
on obtient donc l'{\'e}galit{\'e} dans $R(G_d)$
$$ [\la a \ra] = \sum_{b\leq a} (-1)^{d(b,a)} [\pi(b)] $$
Or, pour $b=a_{min}$, on a $\la a_{min}\ra=\pi(a_{min})$ (c'est la
repr{\'e}sentation de Speh associ{\'e}e {\`a} $\Delta$ et c'est l'image de la
s{\'e}rie discr{\`e}te $\pi^{disc}$ par l'involution de 
Zelevinski) et  pour $a_{min}\neq b$, la repr{\'e}sentation $\pi(b)$ est une
induite ``propre''. 
On obtient donc la congruence $[\pi]= \pm [\la  a_{min}\ra]
\hbox{ mod } R_I(G_d)$ dans $R(G_d)$. En l'appliquant aussi {\`a}
$\pi^{disc}$, on obtient $$[\pi]=\pm [\pi^{disc}] \hbox{ mod }
R_I(G_d).$$ 
La formule
de Van Dijk montre alors que les caract{\`e}res de $\pi$ et $\pi^{disc}$
 co{\"\i}ncident au signe pr{\`e}s sur les {\'e}l{\'e}ments elliptiques. Or, les
formules d'orthogonalit{\'e} pour les s{\'e}ries discr{\`e}tes montrent que le
caract{\`e}re d'une s{\'e}rie discr{\`e}te sur ces {\'e}l{\'e}ments est non nul.



\end{proof}

Une repr{\'e}sentation satisfaisant les propri{\'e}t{\'e}s {\'e}quivalentes ci-dessus
sera dite {\em elliptique} et nous noterons $\Ell{C}{G_d} \subset
\Irr{C}{G_d}$ le sous-ensemble des classes de repr{\'e}sentations
elliptiques. On a bien-s{\^u}r
$$ \Cu{C}{G_d} \subset \Disc{C}{G_d} \subset \Ell{C}{G_d}. $$
Nous allons donner une classification de ces repr{\'e}sentations adapt{\'e}e
aux besoins de ce texte.

\begin{nota} \label{notationsrho}
  Soit $\rho\in \Irr{C}{\dd}$. Nous noterons
\begin{itemize}
\item $d_\rho\in \NM^*$ l'unique diviseur de $d$ et $\tau_\rho^0$
  l'unique repr\'esentation cuspidale irr\'eductible de $G_{d/d_\rho}$
  tels que $JL_d(\rho)$ apparaisse dans l'induite parabolique standard
  normalis\'ee
  $\ip{(G_{d/d_\rho})^{d_\rho}}{G_d}{|det|^{\frac{1-d_\rho}{2}}\tau_\rho^0 \times\cdots\times |det|^{\frac{d_\rho-1}{2}} \tau_\rho^0}$. L'existence de
  $d_\rho$ et $\tau_\rho^0$ est assur{\'e}e par la classification des
  s{\'e}ries discr{\`e}tes de $G_d$ par Bernstein-Zelevinski,
  \cite[Thm 9.3]{Zel}.
\item $M_\rho$ le \levi\ standard $(G_{d/d_\rho})^{d_\rho}$ et
  $\tau_\rho$ la repr\'esentation supercuspidale irr\'eductible $
  \tau_\rho^0 |det|^{(1-d_\rho)/2} \otimes \cdots \otimes \tau_\rho^0
  |det|^{(d_\rho-1)/2}$ de $M_\rho$.  Ainsi la paire $(M_\rho,\tau_\rho)$
  est un repr\'esentant du support cuspidal de $JL_d(\rho)$
\end{itemize}
\end{nota}

\alin{Rappels sur la th{\'e}orie de Bernstein \cite{bernstein}} \label{decbernstein} 
Si $G$ est un groupe r\'eductif $p$-adique, rappelons que $\Mo{C}{G}$ d\'esigne la
cat{\'e}gorie ab{\'e}lienne de toutes les 
$C$-repr{\'e}sentations lisses de $G$.
Soit $(M,\tau)$ une paire Levi-cuspidale, d\'efinissons
$\BG_{M,\tau}^{G}$
la sous-cat{\'e}gorie pleine de $\Mo{C}{G}$ form{\'e}e des objets dont tous
les sous-quotients irr{\'e}ductibles contiennent $(M,\tau\psi)$ dans leur
support cuspidal, pour un certain caract{\`e}re non-ramifi{\'e} $\psi$ de $M$. On
sait que la cat{\'e}gorie $\BG_{M,\tau}^{G}$ est un ``facteur direct
ind{\'e}composable'' (que nous appellerons bloc de Bernstein associ{\'e} {\`a}
$(M,\tau)$) de
$\Mo{C}{G}$ et qu'on a une d{\'e}composition
$$ \Mo{C}{G} \simeq \bigoplus_{(M,\tau)/\sim} \BG_{M,\tau}^{G} $$ 
o{\`u} $(M,\tau)\sim (M',\tau')$ \ssi\ il existe $g\in G$ et $\psi$
caract{\`e}re non ramifi{\'e} de $M$ tels que $M'=M^g$ et $\tau'=(\tau\psi)^g$
({\'e}quivalence ``inertielle''). En particulier, les idempotents centraux
primitifs de $\Mo{C}{G}$ sont en bijection avec les classes
inertielles de paires $(M,\tau)$.

Pour un produit de groupes lin\'eaires, lorsque $M = T$ est un tore
maximal et $\tau$ est un caract{\`e}re 
non-ramifi{\'e} de $T$, 
on appelle parfois le bloc associ\'e $\BG^{G}_{T,1}$
le ``bloc unipotent'' de $G$.

\`A l'oppos\'e, lorsque $M=G$, le bloc est dit "cuspidal". Dans le cas  $G_n=GL_n(K)$, les blocs cuspidaux ont une structure simple : 
\begin{fact} \label{bloccusp}
Si $\tau \in \Cu{C}{GL_n(K)}$, et $f_\tau$ est le nombre de caract\`eres non ramifi\'es de $K^\times$ tels que $(\chi\circ \hbox{det})\otimes \pi=\pi$ alors  pour toute extension finie $K'|K$ de degr\'e r\'esiduel $f_\tau$, il existe une unique \'equivalence de cat\'egories $\BG^{G}_{G,\tau}\To{\alpha_\tau} \BG^{{K'}^\times}_{{K'}^\times,1}=\Mo{C}{{K'}^\times/\OC_{K'}^\times}$
 telle que 
 \begin{enumerate}
 \item $\alpha_{\tau,K'}(\tau)$ est la repr\'esentation triviale de ${K'}^\times$ et
 \item pour tout caract\`ere non ramifi\'e $\chi : K^\times \To{}C^\times$, on a $$\alpha_\tau\left((\chi\circ\hbox{det})\otimes \pi\right)=(\chi\circ N_{K'|K})\otimes \alpha_\tau(\pi).$$
 \end{enumerate}
De plus si $\hbox{deg}(K'|K)=n$, alors pour toute repr\'esentation irr\'eductible $\tau'$ dans $\BG^G_{G,\tau}$, on a $\omega_{\tau'}(\varpi) = \omega_\tau(\varpi) \omega_{\alpha_\tau(\tau')}(\varpi)$
pour toute uniformisante $\varpi$ de $K$.
\end{fact} 

 \begin{proof}
 Cela est implicite dans \cite{BK}, voir notamment la discussion  suivant (7.6.18) et le lemme (6.2.5). Mais c'est en fait beaucoup plus \'el\'ementaire et on n'a pas besoin de la th\'eorie des types simples ! Esquissons les arguments : par Bernstein \cite{bernstein} et un peu de th\'eorie de Clifford (voir \ref{defun} pour quelques d\'etails) on montre que $\BG^G_{G,\tau}$ est \'equivalente \`a  $\Mo{C}{G_\tau/G^0}$ o\`u $G^0$ est le noyau commun de tous les caract\`eres non ramifi\'es de $G$ et $G_\tau$ est le noyau commun de ceux qui stabilisent la classe d'isomorphisme de $\tau$. Si on envoie un g\'en\'erateur de $G_\tau/G^0$ sur un g\'en\'erateur de ${K'}^\times/\OC_{K'}^\times$, on obtient une \'equivalence ayant les propri\'et\'es requises (noter que la derni\`ere assertion est \'equivalente \`a $|\hbox{det}(\varpi)|=|N_{K'|K}(\varpi)|$).
 
Pour l'unicit\'e, on est amen\'e \`a montrer qu'une auto-\'equivalence de $\Mo{C}{\ZM}$ qui "stabilise" les caract\`eres est isomorphe \`a l'identit\'e. Par Morita, une telle auto-\'equivalence est donn\'ee par un $C[X^{\pm 1},Y^{\pm 1}]$-module $P$ qui est projectif de type fini s\'epar\'ement sur $C[X^{\pm 1}]$ et sur $C[Y^{\pm 1}]$. En tant que $C[X^{\pm 1}]$, $P$ est n\'ecessairement de rang $1$ et donc isomorphe \`a $C[X^{\pm 1}]$. L'action de $Y$ sur $P$ est donn\'ee par un polyn\^ome $Q(X)$, et l'hypoth\`ese dit que pour tout $c\in C$, on $Q(c)=c$. Donc $Q(X)=X$.

 \end{proof}

Pour ce qui est des blocs non-cuspidaux, toujours dans le cas $GL_n$, Bushnell et Kutzko en ont donn\'e une belle expression en termes d'alg\`ebres de Hecke dans \cite{BK}. 
Le cas qui nous int\'eresse est celui o\`u $M$ est un produit de $GL$ de m\^eme taille et $\tau$ est un produit tensoriel d'une m\^eme repr\'esentation cuspidale. Dans ce cas la th\'eorie de Bushnell-Kutzko g\'en\'eralise l'\'enonc\'e ci-dessus en montrant qu'un tel bloc est \'equivalent au bloc unipotent d'un groupe lin\'eaire plus petit sur un corps plus gros.

Si $\rho\in\Irr{C}{\dd}$, la description des blocs cuspidaux  s'applique \`a la repr\'esentation $\tau_\rho^0$ de $G_{d_\rho}$. En fait le nombre $f_{\tau_\rho^0}$ est aussi le nombre $f_\rho$ de caract\`eres non ramifi\'es de $\dd$ stabilisant la classe de $\rho$. Par produit, l'\'equivalence $\alpha_{\tau_\rho^0}$ induit une \'equivalence $\alpha_{\tau_\rho} : \BG^{M_\rho}_{M_\rho,\tau_\rho} \simto \BG^{({K'}^\times)^{d_\rho}}_{({K'}^\times)^{d_\rho},1}$.

\begin{fact}\label{equivcat} (Bushnell-Kutzko)
Soit $\rho\in\Irr{C}{\dd}$. 
Il existe une extension $K_\rho$ de $K$ de degr\'e $d/d_\rho$ et degr\'e r\'esiduel $f_\rho$ et 
  une \'equivalence de cat\'egories $\alpha_\rho$ s'inscrivant dans le diagramme commutatif suivant 
$$ \xymatrix{
  \BG_{M_\rho,\tau_\rho}^{G_d} \ar@<.5ex>[d]^{\Rp{G_d}{M_\rho}} \ar[r]^{\alpha_\rho} & 
  \BG^{G'_{d_\rho}}_{T'_{d_\rho},1} \ar@<.5ex>[d]^{\Rp{G'_{d_\rho}}{T'_{d_\rho}}}  \\
  \BG_{M_\rho,\tau_\rho}^{M_\rho}
  \ar@<.5ex>[u]^{\Ip{M_\rho}{G_d}} \ar[r]^{\alpha_{\tau_\rho}} & 
  \BG^{T'_{d_\rho}}_{T'_{d_\rho},1} \ar@<.5ex>[u]^{\Ip{T_{d_\rho}}{G_{d_\rho}}}
} $$ 
o\`u  $G'_{d_\rho}:=GL_{d_\rho}(K')$ et $T'_{d_\rho}$ est son tore diagonal, et les fl\`eches
  verticales d\'esignent les foncteurs paraboliques normalis\'es standard.
De plus, on a
\begin{enumerate}
\item  $\alpha_\rho(JL_d(\rho))$ est isomorphe {\`a} la repr{\'e}sentation de
    Steinberg de $G'_{d_\rho}=GL_{d_\rho}(K_\rho)$.
\item Pout tout caract\`ere non ramifi\'e $\chi$ de $K^\times$, on a $$\alpha_\rho\left((\chi\circ\hbox{det})\otimes \pi\right) \simeq (\chi\circ N_{K_\rho|K}\circ\hbox{det})\otimes\alpha_\rho(\pi). $$
\item Pour toute repr\'esentation $\pi\in \BG^G_{M_\rho,\tau_\rho}$ telle que $\pi(\varpi)$ soit un scalaire, $\alpha_\rho(\pi)(\varpi)$ est aussi un scalaire et on a $\alpha_\rho(\pi)(\varpi)=\pi(\varpi)\omega_\rho(\varpi)^{-1}$. 
\end{enumerate}

\end{fact}
\begin{proof}
  L'existence de $\alpha_\rho$ rendant le diagramme commutatif est 
  donn\'ee par \cite[Thm (7.6.20)+ Cor (7.6.21)]{BK}. En fait, pour \^etre un peu plus pr\'ecis,
  l'\'equivalence de cat\'egories donn\'ee par {\em loc.cit} ne concerne que les
  repr\'esentations admissibles de longueur finie. Pour obtenir
  l'\'equivalence des blocs ``en entier'', il faut utiliser aussi
  \cite{BK1}. L'existence du diagramme commutatif est donn\'ee par
  \cite[7.6.21]{BK} (m\^eme renvoi \`a \cite{BK1} pour les blocs ``en entier'').
  Le fait que l'image de $JL_d(\rho)$ est la Steinberg est une cons\'equence de 
  \cite[Thm 7.7.1]{BK} et du fait que le support cuspidal de cette image contient $(T'_{d_\rho},1)$ par d\'efinition de $\alpha_{\tau_\rho}$ et compatibilit\'e aux foncteurs paraboliques. 
  La compatibilit\'e \`a la torsion par les caract\`eres se d\'eduit de \cite[(7.5.12)]{BK}. La derni\`ere propri\'et\'e est mentionn\'ee en \cite[(7.7.6)]{BK} pour les repr\'esentations admissibles. Elle se g\'en\'eralise aux autres repr\'esentations, en utilisant la discussion de \cite[(7.5.9)]{BK} par exemple.
\end{proof}

D'apr\`es la d\'efinition des repr\'esentations elliptiques,
l'\'equivalence de cat\'egories $\alpha_\rho$ induit une bijection entre l'ensemble
$\Ell{C}{\BG_{M_\rho,\tau_\rho}^{G_d}}$ des repr{\'e}sentations
elliptiques de $G_d$ dans le bloc  $\BG_{M_\rho,\tau_\rho}^{G_d}$ et l'ensemble
$\Ell{C}{\BG^{G'_{d_\rho}}_{T_{d_\rho},1}}$ des repr{\'e}sentations elliptiques de
$G'_{d_\rho}$ dans son bloc unipotent. Notons que cette bijection ne d\'epend pas de $\alpha_\rho$, seulement de l'existence de $\alpha_\rho$ rendant commutatif le diagramme ci-dessus.

\alin{Classification des repr{\'e}sentations elliptiques} \label{classrepell}
L'assertion \ref{equivcat} ram{\`e}ne la classification des repr{\'e}sentations
elliptiques de $G_d$ {\`a} celles des blocs unipotents des $GL_{d'}(K')$ pour
$d'$ divisant $d$ et $K'$ extension de $K$ de degr\'e $\leq d$.  Or,
les repr{\'e}sentations elliptiques du bloc 
unipotent de $GL_{d'}(K')$ sont par d{\'e}finition les sous-quotients de l'induite droite
$\ind{B}{GL_{d'}(K')}{1}$ {\`a} partir d'un \borel de $GL_{d'}(K')$. Nous avons d{\'e}ja
rappel{\'e} la classification de ces repr{\'e}sentations dans la
partie \ref{dp} ; elle ne fait intervenir que le syst\`eme de racines
et pas le corps de d\'efinition $K'$.
Soit $S_{d'}$ l'ensemble des
racines simples de $G_{d'}$ que nous identifierons \`a celui de
$GL_{d'}(K')$.  Pour tout $I\subseteq S_{d'}$, 
nous noterons ici $\pi_{d',I}$, resp. $\pi^{K'}_{d',I}$  la
repr\'esentation de $G_{d'}$, resp $GL_{d'}(K')$, 
d\'efinie en \ref{rappelclassif},  o\`u elle  \'etait not\'ee simplement
$\pi_I$. 
L'application $(I\subseteq S_{d'}) \mapsto
\pi^{K'}_{d',I}$ d{\'e}finit donc une bijection entre l'ensemble $\PC(S_{d'})$ des
sous-ensembles de $S_{d'}$ 
et l'ensemble 
des $C$-repr{\'e}sentations  elliptiques du bloc unipotent $\BG^{GL_{d'}(K')}_{T,1}$
Cette bijection envoie $\emptyset$ sur la repr{\'e}sentation de Steinberg et
$S_{d'}$ sur la repr{\'e}sentation triviale.


\begin{nota} \label{notaell}
Soit $\rho\in \Irr{C}{\dd}$ et $\alpha_\rho$ l'\'equivalence de
cat\'egories de \ref{equivcat}. Nous noterons
$\pi_\rho^I:=\alpha_\rho^{-1}(\pi^{K_\rho}_{d_\rho,I})$, pour $I\subseteq
S_{d_\rho}$. Ce sont les repr{\'e}sentations
elliptiques de $G_d$ dont le support cuspidal est le m{\^e}me que celui de
$JL_d(\rho)$. On a donc en particulier $JL_d(\rho)=\pi_\rho^\emptyset$. 
\end{nota}



L'application $(\rho \in \Irr{C}{\dd},I\subseteq S_{d_\rho}) \mapsto
\pi_\rho^I\in \hbox{Ell}_C(G_d)$ est une bijection. On obtient ainsi une classification de
toutes les repr{\'e}sentations elliptiques de $G_d$.

Comme on l'a mentionn{\'e} {\`a} la fin de la
preuve de \ref{defell},
si $\pi^{disc}$ est la s{\'e}rie discr{\`e}te de m{\^e}me 
support cuspidal que la repr{\'e}sentation $\pi\in \Ell{C}{G_d}$, alors on
a $[\pi]=\pm [\pi^{disc}]$ dans 
$\o{R}(G_d)$ et par cons{\'e}quent : $JL_d(\pi)=\pm JL_d(\pi^{disc})$. On peut
pr{\'e}ciser ce signe en termes de la classification ci-dessus :

\begin{rema}
   Supposons que $\pi^{disc} = JL_d(\rho)=\pi_\rho^\emptyset$ pour $\rho\in \Irr{C}{\dd}$
  et que $\pi = \pi_\rho^I$ pour $I\subseteq S_{d_\rho}$.
Alors on a dans $\o{R}(G_d)$  l'{\'e}galit{\'e}
$[\pi_\rho^I]=(-1)^{|I|}[\pi_\rho^\emptyset]$. Par cons{\'e}quent on a
aussi
$$JL_d[\pi^I_\rho] = (-1)^{|I|} 
[\rho]$$ dans $R(\dd)$.
\end{rema}
\begin{proof}
Gr{\^a}ce {\`a} l'{\'e}quivalence de cat{\'e}gories $\alpha_\rho$ et {\`a} ses propri{\'e}t{\'e}s
de commutation {\`a} l'induction parabolique, on est ramen{\'e} {\`a} prouver
l'{\'e}galit{\'e} $
[\pi_{d_\rho,I}]=(-1)^{|I|}[\pi_{d_\rho,\emptyset}]=[\hbox{St}_{G_{d_\rho}}]$
dans $R(G_{d_\rho})$. Pour prouver celle-ci, notons $P_I$ le
parabolique standard de $G_{d_\rho}$ associ{\'e} {\`a} $I\subseteq
S_{d_\rho}$. Le lemme X.4.6 de \cite{BW} montre que
dans $R(G_{d_\rho})$, on a
$$ [\pi_{d_\rho,I}] = \sum_{I\subseteq J \subseteq S_{d_\rho}} (-1)^{|J\setminus I|}  
   [\ind{P_I}{G_{d_\rho}}{1}]. $$
En appliquant ceci {\`a} $I$ et $\emptyset$, on obtient $[\pi_{d_\rho,I}]
= (-1)^{d_\rho-1-|I|}[1_{G_{d_\rho}}] = 
(-1)^{|I|}[\hbox{St}_{G_{d_\rho}}]$.

\end{proof}

\begin{rema} \label{dualiteell}
La duale de $\pi_\rho^I$ est $\pi_{\rho^\vee}^{\o{I}}$, o\`u l'on identifie $S_{d_\rho}=S_{d_{\rho^\vee}}$ et $\o{I}$ est l'image de $I$ par l'involution $-w_0$ de $S_{d_\rho}$ dans lui-m\^eme en notant $w_0$ l'\'el\'ement de plus grande longueur du groupe de Weyl de $G_{d_\rho}$.
\end{rema}

\begin{proof}
Remarquons tout d'abord que par compatibilit\'e de la correspondance de Jacquet-Langlands \`a la contragr\'ediente, on a $JL_d(\rho^\vee)=JL_d(\rho)^\vee$ et par suite
$M_{\rho^\vee}=M_\rho$ et $\tau_{\rho^\vee}\simeq (\tau_\rho)^\vee$. 
Soit maintenant $P_\rho$ le \para standard de Levi $M_\rho$. L'ensemble des racines simples de $\ZC(M_\rho)$ dans $P_\rho$ s'identifie \`a $S_{d_\rho}$. Ce dernier est donc en bijection $I\mapsto P_{\rho,I}$ avec 
l'ensemble des \paras standards contenant $M_\rho$.

On peut expliciter la d\'efinition de $\pi_\rho^I$ \`a travers l'\'equivalence $\alpha_\rho$ de la mani\`ere suivante :
Soit $M_{\rho,I}$ la composante de Levi standard de $P_{\rho,I}$ et $Q_{\rho,I}$ le \para standard de $M_{\rho,I}$ dont la composante de Levi est $M_\rho$. Alors on a
$$ \pi^I_\rho = \hbox{Cosoc}\left(\Ip{P_{\rho,I}}{G_d}\left(\hbox{Soc}\left(\Ip{Q_{\rho,I}}{M_{\rho,I}}\left(\tau_\rho\right)\right)\right) \right)
$$
o\`u le socle et le cosocle se trouvent \^etre irr\'eductibles.
Par dualit\'e il vient
$$ (\pi^I_\rho)^\vee = \hbox{Soc}\left(\Ip{P_{\rho,I}}{G_d}\left(\hbox{Cosoc}\left(\Ip{Q_{\rho,I}}{M_{\rho,I}}\left((\tau_\rho)^\vee\right)\right)\right) \right).
$$
En appliquant l'\'equivalence $\alpha_{\rho^\vee}$, on obtient
$$\alpha_{\rho^\vee}\left((\pi^I_\rho)^\vee\right) = (\pi_{d_\rho,I})^\vee $$
qui d'apr\`es \ref{serprincorps} iv) n'est autre que $\pi_{d_\rho,\o{I}}$.


\end{proof}

\alin{Extensions entre repr{\'e}sentations elliptiques} On a vu que
l'application ``de Jacquet-Langlands'' $JL_d : R(\dd)\To{} R(G_d)$
pr{\'e}serve les caract{\'e}ristiques 
d'Euler-Poincar{\'e} ({\`a} caract{\`e}re central fix{\'e}). Mais l'ingr{\'e}dient
essentiel de ce texte est le calcul de tous les  groupes d'extensions
 entre repr{\'e}sentations elliptiques. Celui-ci se ram{\`e}ne par
les {\'e}quivalences $\alpha_\rho$ de \ref{equivcat} au calcul de la partie \ref{extensions}
sur les s{\'e}ries principales 
elliptiques de $G_{d_\rho}$.

Fixons une uniformisante $\varpi$ de $K$. Nous noterons par
la m\^eme lettre son image dans le centre de $G_d$ ou celui de $\dd$
par les inclusions $K^\times \subset G_d,\,\dd$ et nous noterons
$\varpi^\ZM$ le sous-groupe qu'elle engendre dans $G_d$ ou $\dd$.
Remarquons que par la caract\'erisation de la correspondance de
Jacquet-Langlands, si
$\rho\in\Irr{\o\QM_l}{\dd/\varpi^\ZM}$ alors $JL_d(\rho)\in
\Irr{\o\QM_l}{G_d/\varpi^\ZM}$ et par cons\'equent toutes les
repr\'esentations $\pi_\rho^I$ ont aussi un caract\`ere central
trivial sur $\varpi$.

\begin{prop} \label{theoextell}
Soient $\rho, \rho' \in \Irr{\o\QM_l}{\dd/\varpi^\ZM}$ et $I\subseteq S_{d_\rho}, I'\subseteq S_{d_{\rho'}}$.
\begin{enumerate}
  \item Pour tout $i\in\NM$, on a
    $$  \cas{\ext{i}{\pi^I_{\rho}}{\pi^{I'}_{\rho'}}{G_d/\varpi^\ZM}}{\o\QM_l}{\rho\simeq \rho' \hbox{ et }i=\delta(I,I')}{0}{\rho\neq \rho' \hbox{ ou } i\neq
      \delta(I,I')}. $$
\item Soient $I,J,K$ trois sous-ensembles de $S_{d_\rho}$ tels que
  $\delta(I,J)+\delta(J,K)=\delta(I,K)$, alors le cup-produit
$$  \cup \;:\;\;
\ext{\delta(I,J)}{\pi^I_{\rho}}{\pi^J_{\rho}}{G_d/\varpi^\ZM} \otimes_{\o\QM_l}
\ext{\delta(J,K)}{\pi^J_{\rho}}{\pi^K_\rho}{G_d/\varpi^\ZM} \To{} 
\ext{\delta(I,K)}{\pi^I_{\rho}}{\pi^K_\rho}{G_d/\varpi^\ZM}
$$ est un isomorphisme.
\end{enumerate}
\end{prop}

\begin{proof}
Commen{\c c}ons par montrer la nullit\'e des Ext lorsque $\rho\neq
\rho'$. Tout d'abord, si ${\omega_\rho}_{|\OC_K^\times}\neq {\omega_{\rho'}}_{|\OC_K^\times}$, alors les $Ext$ sont \'evidemment nuls. Supposant alors $\omega_\rho=\omega_{\rho'}$, on a 
$$ \ext{*}{\pi^I_\rho}{\pi^J_{\rho'}}{G_d,\omega_\rho} 
=\ext{*}{\pi^I_\rho}{\pi^J_{\rho'}}{G_d/\varpi^\ZM} $$
car les objets projectifs de $\Mo{\omega_\rho}{G_d}$ sont encore projectifs dans $\Mo{\o\QM_l}{G_d/\varpi^\ZM}$.
Rappelons maintenant que si $\pi,\pi'$ sont deux repr\'esentations
irr\'eductibles de m\^eme caract\`ere central $\omega$ mais de supports
cuspidaux distincts, alors $\ext{i}{\pi}{\pi'}{G,\omega}=0$ pour tout
$i\in\NM$, voir par exemple \cite[6.1]{Vigext}. Or, le support cuspidal de $\pi^I_\rho$ est
le m\^eme que celui de $\pi^\emptyset_\rho=JL_d(\rho)$. Comme il y a
au plus une s\'erie discr\`ete de support cuspidal donn\'e, on en
d\'eduit l'assertion.

Pour le reste du th\'eor\`eme, nous allons nous ramener au
th\'eor\`eme \ref{theoext} gr\^ace \`a l'\'equivalence
$\alpha_\rho$. Pour cela,
remarquons que par la propri\'et\'e iii) de \ref{equivcat},
 $\alpha_\rho$ induit une \'equivalence entre la sous-cat\'egorie pleine des objets de $\BG^{G_d}_{M_\rho,\tau_\rho}$ triviaux sur $\varpi$ et la sous-cat\'egorie pleine des objets
de $\BG^{G'_{d_\rho}}_{T,1}$ triviaux sur $\varpi$ (rappelons que $G'_{d_\rho}=GL_{d_\rho}(K_\rho)$ avec les notations de \ref{equivcat}). Il s'ensuit que 
$$\ext{*}{\pi^I_\rho}{\pi^J_\rho}{G_d/\varpi^\ZM} = \ext{*}{\pi^{K_\rho}_{d_\rho,I}}{\pi^{K_\rho}_{d_\rho,J}}{G'_{d_\rho}/\varpi^\ZM}.$$
Mais puisque $K_\rho^\times/\varpi^\ZM$ est compact, les objets projectifs de $\Mo{\o\QM_l}{PG'_{d_\rho}}$ sont encore projectifs lorsqu'ils sont vus comme objets de $\Mo{\o\QM_l}{G'_{d_\rho}/\varpi^\ZM}$ de sorte que 
$$\ext{*}{\pi^{K_\rho}_{d_\rho,I}}{\pi^{K_\rho}_{d_\rho,J}}{G'_{d_\rho}/\varpi^\ZM}
=\ext{*}{\pi^{K_\rho}_{d_\rho,I}}{\pi^{K_\rho}_{d_\rho,J}}{PG'_{d_\rho}}.$$
On est donc ramen\'e au th\'eor\`eme \ref{theoext}.
 \end{proof}

\subsection{Correspondance (locale) de Langlands} \label{corLan}
\def\rep{\hbox{Rep}}

Soit $K^{nr}$ la sous-extension non ramifi{\'e}e maximale de $K$ dans
$K^{ca}$ et $\wh{K^{nr}}$ sa compl{\'e}tion. On note toujours $I_K :=
\gal(K^{ca}/K^{nr})\subset W_K$ le groupe d'inertie de $K$.

\alin{Formulation {\`a} la Weil-Deligne}
Rappelons qu'une  "représentation de Weil-Deligne" $\sigma$ de $W_K$ à valeurs dans le corps $C$ est un triplet $(\sigma^{ss},N_\sigma,V_\sigma)$ où 
\begin{itemize}
\item $V_\sigma$ est un espace vectoriel de dimension finie sur $C$,  
\item $\sigma^{ss}: W_K\To{} GL(V_\sigma)$ est une représentation continue (pour la topologie discrète de $C$) et semi-simple de $W_K$.
\item $N_\sigma\in \endo{C}{V_\sigma}$ est un endomorphisme nilpotent de $V$ tel que pour tout $w\in W_K$, on a 
\ini\begin{equation}\label{monod}
\sigma^{ss}(w)N_\sigma\sigma^{ss}(w)^{-1} = |w|N_\sigma.
\end{equation}
\end{itemize}
Notons $\rep^d_C(WD_K)$ l'ensemble des classes d'{\'e}quivalences des
 $C$-repr{\'e}sentations de Weil-Deligne de dimension $d$.
La correspondance de Langlands sur $K$ est une famille de bijections
$(\sigma_d)_{d\in\NM^*}$ 
$$\sigma_d :\;\;\Irr{\CM}{GL_d(K)} \simto \rep^d_\CM(WD_K)$$ 
 qui v{\'e}rifient un certain nombre de propri{\'e}t{\'e}s
suffisantes pour les rendre uniques
(compatibilit{\'e} avec la th{\'e}orie du
corps de classe qui donne aussi le cas $d=1$,  pr{\'e}servation de certains
invariants de nature arithm{\'e}tico-analytique (facteurs $L$ et
$\epsilon$ de paires), etc...). Nous renvoyons {\`a} \cite{HeLang},\cite{HeSMF} pour l'{\'e}nonc{\'e}
pr{\'e}cis de ces propri{\'e}t{\'e}s caract{\'e}ristiques que nous n'utiliserons pas
en totalit{\'e}. Nous utiliserons sans commentaires particuliers la compatibilité à la contragrédiente et la compatibilité à la torsion par les caractères qui dans le cas des caractéres non ramifiés s'exprime formellement ainsi :
$$\forall \pi\in\Irr{\CM}{G_d},\forall a\in\CM^\times,\;\; \sigma_d((|\hbox{det}|^a\pi)=|-|^a \sigma_d(\pi) $$
(en particulier, on a normalisé le corps de classes de manière à ce que les uniformisantes correspondent aux Frobenius géométriques).
Rappelons par ailleurs que 
$$ \sigma_d(\Cu{\CM}{G_d}) = \{\sigma \in \rep^d_\CM(WD_K) \; \hbox{
  irr{\'e}ductibles}\}. $$
En fait, la correspondance est  d{\'e}termin{\'e}e par sa restriction aux
repr{\'e}sentations cuspidales de $G_d$ et irr{\'e}ductibles de $W_K$, la
classification de Zelevinski \cite{Zel} et la classification des
repr{\'e}sentations de Weil-Deligne en fonction des irr{\'e}ductibles de
$W_K$.

La g{\'e}om{\'e}trie alg\'ebrique r{\'e}alise plut\^ot une variante de la
correspondance de Langlands, appel{\'e}e parfois correspondance de Hecke
et qui se d{\'e}duit de celle de Langlands par une simple torsion ``{\`a} la
Tate'' par le caract{\`e}re $w \mapsto |w|^{{1-d}\over 2}$. Il est sous-entendu ici que l'on choisit toujours (si besoin) la racine carrée positive de $q$ dans $\CM$ pour définir une puissance demi-entière de $|-|$.
Pour les représentations cuspidales, cette variante se trouve {\^e}tre compatible 
aux actions des automorphismes de $\CM$ des deux c{\^o}t{\'e}s, cf \cite[(7.4)]{HeTNB}. 
L'extension aux  représentations non cuspidales respecte cette compatibilit\'e, 
comme nous allons le 
préciser maintenant en suivant \cite[(7.4)]{HeTNB}.
Si $\tau$ est un automorphisme de $\CM$ et $\sigma=(\sigma^{ss},N,V)$ une représentation de Weil-Deligne, notons $\sigma^\tau:=(\sigma^{ss}\otimes 1, N\otimes 1, V\otimes_{\CM,\tau}\CM)$.
Cela définit une action de $\aut{}{\CM}$ sur 
$\rep^d_\CM(WD_K)$ et de même on en définit une autre sur $\Irr{\CM}{G_d}$.

\begin{fact}
Pour tout automorphisme $\tau$ du corps $\CM$ et toute représentation $\pi\in\Irr{\CM}{G_d}$, on a
$$ |-|^{(1-d)/2}\sigma_d(\pi^\tau) = (|-|^{(1-d)/2}\sigma_d(\pi))^\tau .$$
\end{fact}
\begin{proof}
Le cas le plus difficile est celui o\`u $\pi$ est cuspidale ; il est prouvé dans \cite[(7.4)]{HeTNB}
(et cela repose sur des compatibilités avec des cas de correspondance globale).

Pour le cas g\'en\'eral, modifions un peu l'\'enonc\'e en introduisant
la notation $\varepsilon_\tau = \tau(\sqrt q)/\sqrt q \in \{\pm 1\}$ et les caract\`eres $$\varepsilon_\tau:\;\; g\in G_d \mapsto \varepsilon_\tau^{\hbox{val}\circ\hbox{det}(g)}  
\;\;\hbox{et}\;\; \varepsilon_\tau:\;\; w\in W_K \mapsto \varepsilon_\tau^{\hbox{log}_q|w|}.$$
L'\'enonc\'e que l'on veut prouver se reformule en
\ini\begin{equation}\label{refor}
 \sigma_d(\pi^\tau)=\varepsilon_\tau^{d-1}\sigma_d(\pi)^\tau. 
 \end{equation}
\'Etant donn\'ees $\pi_1,\cdots,\pi_r$ des repr\'esentations de $G_{d_1},\cdots,G_{d_r}$, notons $\pi_1\times \cdots \times \pi_r$ la repr\'esentation de $G_{d_1+\cdots +d_r}$ induite parabolique standard normalis\'ee. Alors si $d=d_1+\cdots +d_r$, on calcule 
$$ \left(\pi_1\times\cdots \times \pi_r\right)^\tau = \varepsilon_\tau^{d-d_1}\pi_1^\tau\times \cdots \times \varepsilon_\tau^{d-d_r}\pi_r^\tau. $$
Supposons maintenant que cette induite est irr\'eductible et que (\ref{refor}) est connu pour chaque $\pi_i$. On sait alors que
$$ \sigma_d(\pi_1\times\cdots \times \pi_r)=\sigma_{d_1}(\pi_1)\oplus\cdots \oplus\sigma_{d_r}(\pi_r) $$
et on en tire imm\'ediatement l'\'enonc\'e pour l'induite. D'apr\`es la classification de Zelevinski, cela nous ram\`ene au cas o\`u $\pi$ est l'unique quotient irr\'eductible de 
$$\hbox{St}_{k_1}(|det|^{n_1}\pi_g)\times \cdots \times \hbox{St}_{k_r}(|det|^{n_r}\pi_g)$$
avec $n_1\geq\cdots\geq n_r \in\ZM$ et $\pi_g\in\Cu{\CM}{G_g}$. On a not\'e ici 
$\hbox{St}_{k}{\pi_g}$ l'unique quotient irr\'eductible de $\pi_g|det|^{(1-k)/2}\times \cdots \times \pi_g|det|^{(k-1)/2}$ ($k$ facteurs, on passe au suivant en multipliant par $|det|$).
Remarquons que $(\hbox{St}_{k}{\pi_g})^\tau=\hbox{St}_{k}{(\varepsilon_\tau^{kg-g-k+1} \pi_g^\tau)}$.
Par ailleurs, d'apr\`es \cite[2.7]{HeSMF}, on a $\sigma_{kg}(\hbox{St}_{k}{\pi_g})=\sigma_g(\pi_g)\otimes\sigma_k(\hbox{St}_{G_k})$.
On calcule donc 
\begin{eqnarray*}
\sigma_{kg}((\hbox{St}_{k}{\pi_g})^\tau)=\sigma_{kg}(\hbox{St}_{k}{(\varepsilon_\tau^{kg-g-k+1} \pi_g^\tau)}) & = & \varepsilon_\tau^{kg-g-k+1}\sigma_g(\pi_g^\tau)\otimes \sigma_k(\hbox{St}_{G_k}) \\
& = & \varepsilon_\tau^{kg-g-k+1}\varepsilon_\tau^{1-g}\varepsilon_\tau^{1-k}
(\sigma_g(\pi_g)\otimes\sigma_k(\hbox{St}_{G_k}))^\tau \\
& = & \varepsilon_\tau^{1-kg}\sigma_{kg}(\hbox{St}_{k}{\pi_g})^\tau
\end{eqnarray*}
compte tenu de ce que l'on sait d\'eja pour $\pi_g$ et $\hbox{St}_{G_k}$. Ainsi (\ref{refor}) est v\'erifi\'e pour les repr\'esentations du type $\hbox{St}_k(\pi_g)$. 
Revenons \`a notre repr\'esentation irr\'eductible $\pi$ et rappelons que selon \cite[2.9]{HeSMF}, il lui est associ\'e la repr\'esentation galoisienne
$$ \sigma_d(\pi)=\sigma_{k_1g}(\hbox{St}_{k_1}(|det|^{n_1}\pi_g))\oplus \cdots \oplus \sigma_{k_rg}\hbox{St}_{k_r}(|det|^{n_r}\pi_g)). $$
Remarquons par ailleurs que $\pi^\tau$ est l'unique quotient irr\'eductible de la repr\'esentation
$$\varepsilon_\tau^{d-k_1g}(\hbox{St}_{k_1}(|det|^{n_1}\pi_g))^\tau\times \cdots \times \varepsilon_\tau^{d-k_rg}(\hbox{St}_{k_r}(|det|^{n_r}\pi_g))^\tau$$  
qui n'est autre que la repr\'esentation
$$\hbox{St}_{k_1}\left( \varepsilon_\tau^{d-k_1g}\varepsilon_\tau^{k_1g-g-k_1-1} |det|^{n_1}\pi_g^\tau\right)\times\cdots\times
\hbox{St}_{k_1}\left( \varepsilon_\tau^{d-k_1g}\varepsilon_\tau^{k_1g-g-k_1-1} |det|^{n_r}\pi_g^\tau\right).
$$
Il lui est donc associ\'e la repr\'esentation
$$\sigma_d(\pi^\tau)= \sigma_{k_1g}\left(  \hbox{St}_{k_1}(\varepsilon_\tau^{d-g-k_1-1}|det|^{n_1}\pi_g^\tau)\right)\oplus\cdots\oplus \sigma_{k_rg}\left(  \hbox{St}_{k_r}(\varepsilon_\tau^{d-g-k_r-1}|det|^{n_r}\pi_g^\tau)\right)
$$
qui n'est autre que
$$\varepsilon_\tau^{d-k_1g}\sigma_{k_1g}\left((\hbox{St}_{k_1}(|det|^{n_1}\pi_g))^\tau\right)\oplus \cdots \oplus \varepsilon_\tau^{d-k_rg}\sigma_{k_rg}\left((\hbox{St}_{k_r}(|det|^{n_r}\pi_g))^\tau\right).$$  
Ainsi, par le cas d\'eja trait\'e des repr\'esentations $\hbox{St}$, on en d\'eduit (\ref{refor}) pour $\pi$.

\end{proof}

La compatibilit{\'e} aux automorphismes de $\CM$  permet de transposer sans ambigu{\"\i}t{\'e}
la correspondance de Langlands tordue \`a la Hecke au corps abstrait $C$. 
 Pour obtenir une correspondance de Langlands sur $C$,
il faut alors choisir 
(si besoin) une racine du cardinal du corps r{\'e}siduel de $K$.

\alin{Formulation continue $l$-adique} Lorsque $C=\o\QM_l$,
nous avons besoin d'une formulation en termes de repr{\'e}sentations
continues $l$-adiques de $W_K$ plut{\^o}t que de repr{\'e}sentations de Weil-Deligne.
Voici  bri\`evement le lien entre les deux formulations expliqu\'e par Deligne
dans \cite[8]{DelAntwerp}.
Rappelons  que $W_K$ est muni de la topologie d{\'e}finie par la
topologie profinie de $I_K$ et la topologie discr{\`e}te de $W_K/I_K\simeq
\ZM$. On note alors  $\rep^d_l(W_K)$ l'ensemble des classes
d'{\'e}quivalences de $\o\QM_l$-repr{\'e}sentations Frobenius-semisimples qui
sont obtenues par 
extension des scalaires de repr{\'e}sentations {\em continues} de
$W_K$ {\`a} valeurs dans une extension finie de $\QM_l$.

Rappelons que nous avons fix\'e un prog\'en\'erateur $\mu$ de $\ZM_l(1)$ et qu'il lui est associ\'e le morphisme surjectif
$$ t_\mu:\;\; I_K\To{t_l} \ZM_l(1)\To{\mu^*} \ZM_l.$$

Soit $(\sigma^{ss},N_\sigma,V_\sigma)$ une $\o\QM_l$-repr{\'e}sentation de Weil-Deligne
comme dans le paragraphe pr{\'e}c{\'e}dent et soit $\phi$
un rel{\`e}vement de Frobenius g{\'e}om{\'e}trique dans $W_K$. L'équation \ref{monod} montre 
 que la formule
$$\sigma(w):= \sigma^{ss}(w) \hbox{exp}(N_\sigma t_\mu(i_\phi(w))), \;\;\hbox{ o{\`u}
  } w=\phi^{\nu(w)} i_\phi(w) \in \phi^\ZM\ltimes I_K$$
d{\'e}finit une repr{\'e}sentation continue et Frobenius-semisimple de $W_K$
sur l'espace  $V_{\sigma}$.
 Le th{\'e}or{\`e}me ``de la
monodromie $l$-adique''  de Grothendieck montre que l'application 
$\rep^d_{\o\QM_l}(WD_K) \To{} \rep^d_l(W_K)$ ainsi obtenue est une bijection.
Deligne a montr\'e qu'elle ne d\'epend pas des choix de $\phi$ et $\mu$.
L'op{\'e}rateur nilpotent $N_\sigma \in \endo{\o\QM_l}{V_{\sigma}}$ est
appel{\'e} ``monodromie'' et contr{\^o}le le d{\'e}faut de semi-simplicit{\'e} de la
repr{\'e}sentation $\sigma$. En particulier la semisimplifi{\'e}e de
$\sigma$ n'est autre que $\sigma^{ss}$.


La correspondance de Langlands locale, tordue \`a la Hecke, fournit donc une correspondance entre repr\'esentations irr\'eductibles de $G_d$ et repr\'esentations $l$-adiques et  c'est cette derni{\`e}re que la
g{\'e}om{\'e}trie ({\em i.e.} la cohomologie $l$-adique) peut pr{\'e}tendre r{\'e}aliser.
Moyennant le
choix d'une racine du cardinal du corps r{\'e}siduel dans $\o\QM_l$, on obtient la "vraie"
correspondance, que nous noterons toujours
$$\pi \in \Irr{\o\QM_l}{G_d} \mapsto
\sigma_d(\pi) \in \rep^d_l(W_K).$$

\alin{Repr{\'e}sentations elliptiques} \label{corlanell}
Avant d'expliquer ce qu'il advient des repr{\'e}sentations elliptiques de
$G_d$ {\`a} travers la correspondance de Langlands, rappelons  que
$$ \sigma_d(\Disc{\o\QM_l}{G_d}) = \{\sigma \in \rep^d_l(W_K) \; \hbox{
  ind{\'e}composables}\}. $$
Ainsi, 
la composition 
$$ \sigma_d\circ JL_d :\; \Irr{\o\QM_l}{\dd}\To{} {\rep^d_l(W_K)},$$
qui est un exemple  simple de fonctorialit{\'e} de Langlands, 
induit une bijection de $\Irr{\o\QM_l}{\dd}$ sur 
l'ensemble des ind{\'e}composables de $\rep^d_l(W_K)$.
\def\I#1{\{1,\cdots,#1\}}

Maintenant, comme les repr{\'e}sentations elliptiques ont par d{\'e}finition
le m{\^e}me support cuspidal qu'une repr{\'e}sentation de la s{\'e}rie discr{\`e}te,
la compatibilit{\'e} de la correspondance de Langlands {\`a} l'induction
parabolique (qui est un exemple encore plus simple de fonctorialit{\'e} de
Langlands) donne la caract{\'e}risation
$$ \sigma_d(\Ell{\o\QM_l}{G_d}) = \{\sigma \in \rep^d_l(W_K), \exists \sigma'  \hbox{
  ind{\'e}composable telle que }  \sigma^{ss}={\sigma'}^{ss}\}. $$

Donnons maintenant une description explicite {\`a} partir de la classification
des repr{\'e}sentations
elliptiques donn{\'e}e plus haut. 

Pour $d'$ divisant $d$, nous fixons une num\'erotation de $S_{d'}$,
{\em i.e.} une bijection $S_{d'} \simto \{1,\cdots,d'-1\}$,
comme en \ref{rappelclassif}. Ceci nous permet de d\'efinir comme en
\ref{deftauI} une (classe de) $\o\QM_l$-repr\'esentation continue $\tau_I$ de
$W_K$ que nous noterons ici 
plus pr\'ecis\'ement  $\tau_{d',I}$.


\begin{lemme} \label{corell}
  Soit $\rho \in \Irr{\o\QM_l}{D}$ et 
$I\subseteq S_{d_\rho}$. Utilisant les notations \ref{notationsrho} on a 
$$ \sigma_d(\pi^I_\rho) \simeq \sigma_{d/d_\rho}(\tau_\rho^0)\otimes
\sigma_{d_\rho}(\pi_{d_\rho,I}) \simeq 
\sigma_{d/d_\rho}(\tau_\rho^0)\otimes\tau_{d_\rho,I} \otimes |.|^{\frac{d_\rho-1}{2}}.$$
\end{lemme}
\begin{proof}
Pour calculer la repr\'esentation $\sigma_{d_\rho}(\pi_{d_\rho,I})$, il nous faut exprimer les param\`etres de Langlands de $\pi_{d_\rho,I}$, c'est-\`a-dire pr\'esenter celle-ci comme unique quotient d'une repr\'esentation du type
$$ |det|^{n_1}\hbox{St}_{G_{d_1}}\times \cdots \times |det|^{n_r}\hbox{St}_{G_{d_r}} $$
avec $n_1\geq \cdots \geq n_r$. Supposons ceci fait, alors par la d\'efinition qu'on a donn\'ee de la repr\'esentation $\pi_\rho^I$ au moyen de l'\'equivalence de cat\'egories \ref{equivcat} (qui est compatible aux foncteurs d'inductions paraboliques et qui envoie temp\'er\'ees sur temp\'er\'ees), celle-ci sera l'unique quotient de la repr\'esentation
$$ |det|^{n_1}\hbox{St}_{d_1}(\tau_\rho^0)\times\cdots\times |det|^{n_r}\hbox{St}_{d_r}(\tau_\rho^0) .$$
Ainsi on aura, par \cite[2.9]{HeSMF},
\begin{eqnarray*}
\sigma_d(\pi_\rho^I) & =  & |-|^{n_1}\sigma_{dd_1/d_\rho}(\hbox{St}_{d_1}(\tau_\rho^0))\oplus\cdots\oplus |-|^{n_r}\sigma_{dd_r/d_\rho}(\hbox{St}_{d_r}(\tau_\rho^0)) \\
 & = & |-|^{n_1} \sigma_{d/d_\rho}(\tau_\rho^0)\otimes \sigma_{d_1}(\hbox{St}_{G_{d_1}}) \oplus\cdots\oplus |-|^{n_r} \sigma_{d/d_\rho}(\tau_\rho^0)\otimes \sigma_{d_r}(\hbox{St}_{G_{d_r}}) \\
 & = & \sigma_{d/d_\rho}(\tau_\rho^0)\otimes \sigma_{d_\rho}(\pi_{d_\rho,I})
 \end{eqnarray*}
d'o\`u la premi\`ere \'egalit\'e de l'\'enonc\'e.

Les param\`etres de Langlands de la repr\'esentation $\pi_{d_\rho,I}$ sont donn\'es par la remarque \ref{parLanpiI}.
Avec les notations courantes que l'on all\`ege en abr\'egeant
$d=d_\rho$ et $\pi_I=\pi_{d_\rho,I}$ et en posant $\o{I^c}$ l'image de $I^c=S\setminus I$ par l'application $x\in S\mapsto (d-x)\in S$, on a donc :
\begin{center}
$\pi_I$ est l'unique quotient irr\'eductible de l'induite 
$\ip{M_{\o{I^c}}}{G_d}{\delta_{P_{\o{I^c}}}^\frac{1}{2} \hbox{St}_{M_{\o{I^c}}}}$.
\end{center}
Pour traduire ceci en termes de produits \`a la Zelevinski, introduisons un peu de combinatoire : sur l'ensemble $\{0,\cdots,d-1\}$, on met une structure de graphe orient\'e en d\'eclarant que $x\rightarrow y$ \ssi\ ($x\leq y$ et $\{x+1,\cdots, y\}\subset \o{I^c}$). En particulier, on a $x\rightarrow x$ pour tout $x$. La partition de $\{0,\cdots, d-1\}$ en composantes connexes s'\'ecrit $\{0,\dots, d_1-1\}\sqcup\{d_1,\cdots, d_1+d_2-1\} \sqcup\cdots $ pour une unique suite $d_1, \cdots d_r$ d'entiers $\geq 1$ et de somme $d$. Alors le \levi $M_{\o{I^c}}$ n'est autre que $G_{d_1}\times\cdots\times G_{d_r}$ et un calcul \'el\'ementaire du module de $P_{\o{I^c}}$ montre que 
$$\ip{M_{\o{I^c}}}{G_d}{\delta_{P_{\o{I^c}}}^\frac{1}{2} \hbox{St}_{M_{\o{I^c}}}}
= |det|^{n_1} \hbox{St}_{G_{d_1}} \times\cdots\times|det|^{n_r}\hbox{St}_{G_{d_r}}$$
en posant pour tout $i=1,\cdots, r$,
$$ n_i= (-d_1-\cdots -d_{i-1} +d_i+\cdots +d_r)/2. $$
D'un autre c\^ot\'e, revenant \`a la d\'efinition de la repr\'esentation $\tau_{I}$ en \ref{deftauI}, on constate que
$$ |-|^{(d-1)/2}\tau_I= |-|^{m_1} \sigma_{d_1}(\hbox{St}_{G_{d_1}}) \oplus \cdots \oplus 
|-|^{m_r} \sigma_{d_r}(\hbox{St}_{G_{d_r}}) $$
avec pour tout $i=1,\cdots,r$,
$$ m_i=(d-d_i)/2-(d_1+\cdots d_{i-1}) .$$
On v\'erifie imm\'ediatement que $n_i=m_i$ pour tout $i$. Comme $\pi$ est l'unique quotient de $\ip{M_{\o{I^c}}}{G_d}{\delta_{P_{\o{I^c}}}^\frac{1}{2} \hbox{St}_{M_{\o{I^c}}}}$ et que les $n_i$ sont d\'ecroissants, on en conclut que $\sigma_d(\pi_I)=|-|^{(d-1)/2}\tau_I$.

\end{proof}

\subsection{Remarques de th{\'e}orie des repr{\'e}sentations}
Ce paragraphe contient des consid{\'e}rations techniques ayant peu
d'int{\'e}r{\^e}t en soi, pour r{\'e}f{\'e}rence ult{\'e}rieure. Il est conseill{\'e} de le
sauter dans une premi{\`e}re lecture.

\alin{Le contexte} \label{notationsnu}
On consid\`ere un groupe localement profini $H$ muni d'un morphisme
surjectif $H\To{\nu_H} \ZM$ dont on note $H^0$ le noyau. 
En notant $\nu_K$ la valuation normalis\'ee de $K^\times\To{} \ZM$, les exemples qui nous int\'eressent sont :
\begin{itemize}
        \item $G_d$ muni de $\nu_{G}:= \nu_K\circ \hbox{det}$,
        \item $\dd$ muni de $\nu_{D} := \nu_K\circ \hbox{Nr}_{D_d/K}$ (norme r\'eduite),
        \item $W_K$ muni de $\nu_{W} :\;\; W_K\To{} \gal(k^{ca}/k)\simto \ZM$ o\`u l'on choisit un Frobenius arithm\'etique comme g\'en\'erateur de $\gal(k^{ca}/k)$.
        \item $G_d\times \dd$ et $(G_d\times \dd)/K^\times_{diag}$ munis de $\nu_{GD} := -\nu_{G}+\nu_{D}$.
        \item $G_d\times \dd\times W_K$ muni de $\nu_{GDW}:=-\nu_{G}+\nu_{D}+\nu_{W}$.
\end{itemize}
Lorsqu'aucune ambigu\"it\'e n'en r\'esultera, on ommettra les indices et \'ecrira simplement $\nu$.

Les repr\'esentations consid\'er\'ees seront toujours \`a coefficients dans un corps alg\'ebriquement clos $C$ de caract\'eristique nulle. On notera $\Psi_H$ le groupe des $C$-caract\`eres de $H/H^0$, dits "non-ramifi\'es". Tout caract\`ere $\psi:\ZM \To{} C^\times$ d\'efinit un caract\`ere non ramifi\'e $\psi_H:=\psi\circ \nu_H$ (de m\^eme, on note $\psi_G$, $\psi_D$, etc ...).

\alin{Familles de repr{\'e}sentations irr{\'e}ductibles} \label{defun}
Le groupe $\Psi_H$  agit par torsion sur l'ensemble des classes
d'isomorphisme de repr{\'e}sentations irr{\'e}ductibles de $H$. Appelons
``orbites inertielles'' ses orbites. Soit $(\rho,V) \in \Irr{C}{H}$ telle
que $\rho_{|H^0}$ soit de longueur finie.


Si $\rho^0$  est un facteur
irr{\'e}ductible de $\rho_{|H^0}$, le normalisateur $H_\rho$ de la classe de
$\rho^0$ dans $H$ est d'indice fini $n_\rho$ et est d{\'e}termin{\'e} par cet
indice car on a alors $H_\rho=\nu_H^{-1}(n_\rho\ZM)$. Comme la notation le sugg{\`e}re, il ne d{\'e}pend pas de $\rho^0$
; il s'identifie au noyau commun des $\psi\in\Psi_H$ tels que
$\rho\simeq \rho\psi$.
Par cyclicit{\'e} de $H_\rho/H^0$, on 
peut {\'e}tendre $\rho^0$ en une repr{\'e}sentation $\wt{\rho^0}$ de
$H_\rho$ qui appara{\^\i}t dans $\rho_{|H_\rho}$.
La th{\'e}orie de Clifford montre  
alors que $\rho^0$, resp. $\wt{\rho^0}$, est de multiplicit{\'e} $1$
dans $\rho_{|H^0}$, resp $\rho_{|H_\rho}$ et que
 $\rho\simeq \cind{H_\rho}{H}{\wt{\rho^0}}$.

Consid{\'e}rons alors  la repr{\'e}sentation induite {\`a} supports compacts
$$\rho_{un}:=\cind{H^0}{H}{\rho^0}$$
dont on note $V_{un}$ l'espace.
On v{\'e}rifie facilement les propri{\'e}t{\'e}s suivantes :
\begin{enumerate}
\item Ses quotients irr{\'e}ductibles sont toutes les repr{\'e}sentations dans l'orbite inertielle de
$\rho$  et apparaissent avec multiplicit{\'e} $1$.
\item Elle est, {\`a} isomorphisme pr{\`e}s, ind{\'e}pendante du choix de
  $\rho^0$.
\end{enumerate}

\alin{Une caract\'erisation dans un cas particulier} \label{zetaiso}
Supposons de plus qu'il existe un sous-groupe central $Z\subset H$ tel que
$ZH^0$ soit d'indice fini dans $H$ et $Z^0:=H^0\cap Z$ soit
compact. 
\begin{nota} \label{covariants}
Pour un caract{\`e}re $\zeta:\; 
Z\To{} C^\times$ et une repr{\'e}sentation $(\rho,V)$ de $H$ nous noterons
$(\rho_\zeta,V_\zeta)$ la repr{\'e}sentation induite par $\rho$ sur les
$\zeta$-coinvariants de $V$ d{\'e}finis par
$V_\zeta:=V/<\rho(z)v-\zeta(z)v,v\in V>$. 
\end{nota}

Soit $\rho\in\Irr{C}{H}$ admettant un caract{\`e}re central $\omega_\rho$.
On v{\'e}rifie facilement qu'il existe un isomorphisme
$$ (\rho_{un})_\zeta \simeq \bigoplus_{\buildrel{\rho'\sim
  \rho}\over{(\omega_{\rho'})_{|Z}=\zeta}} \rho' $$
o{\`u} l'{\'e}quivalence $\sim$ est la torsion par un caract{\`e}re de $H/H^0$.
R{\'e}ciproquement, on a
\begin{lemme} \label{centun}
  Soit $\rho^0$ une $C$-repr{\'e}sentation admissible de $H^0$ telle qu'il existe
  $\rho\in\Irr{C}{H}$
  et un isomorphisme
$$ \beta_\zeta:\;\;{\cind{H^0}{H}{\rho^0}}_\zeta \simto \bigoplus_{\buildrel{\rho'\sim
  \rho}\over{(\omega_{\rho'})_{|Z}=\zeta}} \rho'$$
pour chaque caract{\`e}re $\zeta:\; A\To{} C^\times$.
Alors il existe un isomorphisme $\cind{H^0}{H}{\rho^0}\simto \rho_{un}$.
\end{lemme}
\begin{proof}
Soit $\zeta$ un caract{\`e}re tel que
$\zeta^0:=\zeta_{|Z^0}\neq(\omega_\rho)_{|Z^0}$. Dans ce cas, le terme de
droite de l'isomorphisme est nul.
Comme le foncteur $\cInd{H^0}{H}$ est fid{\`e}le, il s'ensuit que
$\rho^0_{\zeta^0}=0$. Comme $Z^0$ est compact, il s'ensuit que $Z^0$
agit sur $\rho^0$ via $\omega_\rho^0$.

Choisissons au contraire un $\zeta$ tel que
$\zeta_{|Z^0}=(\omega_\rho)_{|Z^0}$. Alors le terme de droite de
l'isomorphisme est non nul et on en d{\'e}duit par r{\'e}ciprocit{\'e} de
Frobenius un morphisme non nul $\rho^0 \To{} \rho_{|H^0}$. Comme
$\rho_{|H^0}$ est semi-simple, on en d{\'e}duit l'existence d'un
morphisme $\rho^0 \To{} \rho_{|H^0}$ d'image irr{\'e}ductible. En
induisant, on obtient donc un morphisme surjectif
$$ \alpha:\;\cind{H^0}{H}{\rho^0} \To{} \rho_{un}.$$
Pour tout $\zeta$,  $\alpha_\zeta$  factorise $\beta_\zeta$ et donc est un isomorphisme.
Soit $V$ le noyau du  morphisme $\alpha$. Montrons que $V$ est nul.
Soit $z\in Z$ tel que $\nu(z)$ engendre le groupe cyclique
$\nu(Z)$ ; alors   le sous-groupe $z^\ZM$ de $Z$ engendr{\'e} par
$z$ est discret et $Z=Z^0z^\ZM$. Toute $C$-repr{\'e}sentation de $H$ est
canoniquement munie d'une structure  
de $C[z^\ZM]$-module et tout morphisme de $CH$-repr{\'e}sentations commute {\`a}
ces structures (en d'autres termes $C[z^\ZM]$ s'identifie {\`a} un
sous-anneau du centre de la cat{\'e}gorie $\Mo{C}{H}$). Maintenant les
repr{\'e}sentations $\cind{H^0}{H}{\rho^0}$ et $\rho_{un}$ sont $C[z^\ZM]$-projectives
et $C[z^\ZM]$-admissibles. Le noyau $V$
est donc aussi $C[z^\ZM]$-projectif et $C[z^\ZM]$-admissible. Il s'ensuit
que pour tout caract{\`e}re $\zeta: Z\To{} C^*$ tel que
$\zeta_{|Z^0}=(\omega_\rho)_{|Z^0}$, la suite
$V_\zeta\To{} \cind{H^0}{H}{\rho^0}_\zeta \To{} (\rho_{un}) _\zeta$
reste exacte. Mais par hypoth{\`e}se, le second morphisme est un
isomorphisme. Donc $V_\zeta=0$. Par le lemme de Nakayama appliqu{\'e} aux
$C[z^\ZM]$-modules de type fini $V^J$ pour $J$ sous-groupe ouvert
compact de $H$, on en d{\'e}duit $V=0$.
\end{proof}

\begin{coro} \label{equivdefun} \label{defpiun}
Soit $\rho\in\Irr{\o\QM_l}{\dd}$ et $I\subseteq S_{d_\rho}$.
La repr\'esentation $\pi_{\rho,un}^I$ obtenue en appliquant la construction de \ref{defun} \`a $H \leftarrow G_d$ et $\rho \leftarrow \pi_\rho^I$ est isomorphe \`a la repr\'esentation suivante
$$\pi^I_{\rho,un} := \alpha_\rho^{-1}\left(\cind{(G'_{d_\rho})^0}{G'_{d_\rho}}{{\pi^{K_\rho}_{d_\rho,I}}_{|(G'_{d_\rho})^0}}\right)
$$
avec les notations de \ref{equivcat}.
\end{coro}

\begin{proof}
Lorsque $\rho=1$, l'assertion est tautologique. \`A partir de l\`a, pour $\rho$ quelconque, c'est une cons\'equence de la propri\'et\'e \ref{equivcat} iii) de l'\'equivalence de cat\'egories $\alpha_\rho$ et de la caract\'erisation donn\'ee par le lemme pr\'ec\'edent.
\end{proof}

\alin{Produits de groupes}
Donnons-nous maintenant deux
couples $(H_1,\nu_1)$ et $(H_2,\nu_2)$ comme en \ref{notationsnu}. 
On munit le
groupe produit $H:=H_1\times H_2$ du morphisme diff\'erence $\nu:=-\nu_1+\nu_2$ (la discussion qui suit s'applique aussi bien au morphisme somme, moyennant quelques changement \'evidents), ce qui nous donne un sous-groupe $H^0=\nu^{-1}(0)$ que l'on munit \`a son tour du morphisme $\nu_0:=\nu_1+0=0+\nu_2$. 

 
Soit $\rho_1$, resp. $\rho_2$, deux repr{\'e}sentations irr{\'e}ductibles de
$H_1$, resp. $H_2$, de restriction {\`a} $H_1^0$, resp. $H_2^0$,
artinienne. La restriction de $\rho:=\rho_1\otimes \rho_2$ {\`a}
$H^0$ est alors artinienne. On voit ais{\'e}ment que l'indice
$n_\rho=[H:H_\rho]$ est le p.g.c.d. de $n_{\rho_1}$ et $n_{\rho_2}$. 

\begin{lemme} \label{lemmeun} 
Supposons que $n_{\rho_1}=n_{\rho_2}$. Alors 
$$ \left(\rho_{un}\right)_{|H_1} \simeq \rho_{1,un}\otimes_C V_{\rho_2}. $$
\end{lemme}

\begin{proof}
Posons $n:=n_{\rho_1}=n_{\rho_2}$ et choisissons pour chaque $i=1,2$ un facteur irr\'eductible $\wt\rho_i^0$ de $(\rho_i)_{|\nu_i^{-1}(n\ZM)}$. Par d\'efinition, on a
$$\rho_{i,un} = \cind{H_i^0}{H_i}{(\wt\rho_i^0)_{|H_i^0}}.$$
Remarquons aussi que 
$$\wt\rho^0:=\cind{\nu_1^{-1}(n\ZM)\times \nu_2^{-1}(n\ZM)}{\nu^{-1}(n\ZM)}{\wt\rho_1^0\otimes\wt\rho_2^0} $$
est un facteur irr\'eductible de $\rho_{|\nu^{-1}(\ZM)}$, puisque pour chaque $i$, on a $\cind{\nu_{i}^{-1}(n\ZM)}{H_i}{\wt\rho_i^0}=\rho_i$.
Il s'ensuit que 
\begin{eqnarray*}
 \rho_{un} = \cind{H^0}{H}{\cind{\nu_1^{-1}(n\ZM)\times \nu_2^{-1}(n\ZM)}{\nu^{-1}(n\ZM)}{\wt\rho_1^0\otimes\wt\rho_2^0}_{|H^0}} 
 \simeq \cind{\nu_0^{-1}(n\ZM)}{H}{(\wt\rho_1^0\otimes\wt\rho_2^0)_{|\nu_0^{-1}(n\ZM)}}
\end{eqnarray*}
Par la formule de Mackey on obtient alors
$$ (\rho_{un})_{|H_1} \simeq \bigoplus_{h_2\in H_2/\nu_2^{-1}(n\ZM)} \cind{H_1^0}{H_1}{\wt\rho_1^0\otimes_C V_{\wt\rho_2^0}} 
\simeq \rho_{1,un}\otimes_C V_{\rho_2}. $$

\end{proof}

\section{Cohomologie \'equivariante des espaces de Berkovich} \label{berkovich}

Dans cette partie, le contexte est le suivant : on suppose qu'un
groupe  localement profini $G$ agit  sur un espace $K$-analytique $X$, au sens de
Berkovich. Celui-ci a d{\'e}fini une notion de
continuit{\'e} pour une telle action, et 
tout un formalisme cohomologique {\'e}tale-$G$-{\'e}quivariant 
"continu". C'est l'existence de ce formalisme qui nous pousse \`a utiliser les espaces de Berkovich et leur cohomologie \'etale plut\^ot que les espaces rigides, ou les espaces adiques de Huber. Il ne fait cependant aucun doute qu'un formalisme similaire existe dans le cadre des espaces de Huber.

Lorsque $\Lambda$ est un anneau de torsion premi\`ere \`a $p$,
le formalisme en question fournit, entre autres, un complexe canonique de la
cat{\'e}gorie d{\'e}riv{\'e}e born{\'e}e des $\Lambda G$-modules {\em lisses} dont la
cohomologie est la cohomologie {\'e}tale de $X$ {\`a} coefficients dans
$\Lambda$, munie de son action canonique de $G$.

Le but de cette partie est d'exposer ce formalisme en grande partie
non-publi{\'e} -- bien que d{\'e}ja utilis{\'e} dans \cite{HaCusp}, \cite{HaKyoto} et \cite{Hausb} -- puis d'en donner une variante $l$-adique  dans un langage inspir\'e de la cohomologie \'etale continue de Jannsen \cite{Jannsen} \cite{Huber} \cite{Fargues}, et enfin d'interpr\'eter certaines suites spectrales de type Hochschild-Serre construites par Fargues dans \cite{Fargues} comme des \'egalit\'es de complexes de cohomologie dans certaines cat\'egories d\'eriv\'ees.

\subsection{Coefficients de torsion}

Cette section est un bref expos{\'e} d'un manuscrit non publi{\'e} de Berkovich \cite{Bicunp}. Les
impr{\'e}cisions et erreurs {\'e}ventuelles ci-dessous sont de la seule
responsabilit{\'e} de l'auteur de ces lignes. Nous avons opt\'e pour une exposition plus \'el\'ementaire -- et beaucoup moins \'el\'egante -- que celle de \cite{Bicunp} ; nous inclurons quelques remarques \`a ce sujet.

\alin{Faisceaux {\'e}tales $G$-{\'e}quivariants} 
 Soit $X$ un espace $K$-analytique et $X_{et}$ son site
 {\'e}tale \cite{Bic2}. Tout endomorphisme analytique $\phi$ de $X$
 d{\'e}finit un endomorphisme $X_{et}\To{\phi} X_{et}$ du site $X_{et}$
 tel que pour tout morphisme
 {\'e}tale $U\To{f} X$, $\phi^{-1}(U,f)\simeq (U,f)\times_X (X,\phi)$. 
On a
 par cons{\'e}quent un couple d'endofoncteurs adjoints $(\phi^*,\phi_*)$ du topos associ{\'e}. 
Si  $\psi$ est un autre endomorphisme de $X$, on a des isomorphismes
 canoniques de foncteurs $(\phi\psi)^*\simto \psi^*\circ \phi^*$ et
 $(\phi\psi)_*\simto \phi_*\psi_*$. En
 particulier, si
$G$ est un groupe agissant sur $X$ par automorphismes
 analytiques, cette
 action induit aussi une action et une anti-action  sur le topos
 {\'e}tale de $X$.
Ces (anti-)actions ne sont pas {\em strictes}, mais les isomorphismes
de transitions {\'e}tant canoniques, nous ferons l'abus de les noter par un signe $=$.
 
L'action et l'anti-action sont reli{\'e}es par des isomorphismes
canoniques. Plus pr{\'e}cis{\'e}ment si $h\in G$, il y a des isomorphismes
canoniques de foncteurs $h^*\simeq 
h^!\simeq (h^{-1})_* \simeq (h^{-1})_!$. On a donc autant de fa{\c c}ons de d{\'e}finir un faisceau
$G$-{\'e}quivariant. Voici la d{\'e}finition que nous prendrons dans ce texte
(elle n'est pas standard mais est {\'e}quivalente aux autres) :
un faisceau $G$-{\'e}quivariant sur $X_{et}$ est un couple
$(\FC,\tau_\FC)$ o{\`u} $\FC$ est un faisceau sur $X_{et}$ et $\tau_\FC$
une famille d'isomorphismes $\tau_\FC(h):\; h_*\FC \simto \FC$, $h\in G$ v{\'e}rifiant
la condition de cocycle $\tau_\FC(h'h)=\tau_\FC(h')\circ
h'_*(\tau_\FC(h))$.

\alin{Faisceaux {\'e}tales $G$-{\'e}quivariants discrets/lisses} \label{fgeq}
Berkovich munit dans \cite[part 6]{Bic3} le groupe 
d'automorphismes analytiques  $\aut{}{X}$  de $X$ d'une certaine
topologie totalement discontinue.  
Par d{\'e}finition de cette topologie, tout
ouvert {\em distingu{\'e}} de $X$ ({\em i.e.} qui est diff{\'e}rence de deux domaines
analytiques compacts de $X$) est stabilis{\'e} par un sous-groupe
ouvert de $\aut{}{X}$. L'id\'ee ma\^itresse de Berkovich dans ce contexte est que {\em tout ouvert \'etale "distingu\'e" est aussi "stabilis\'e" par un sous-groupe ouvert de $\aut{}{X}$}.
Tentons d'expliquer ce que cela signifie. 

Par d\'efinition un ouvert \'etale $(U,f)$ est dit {\em distingu\'e} s'il peut se factoriser 
$U\To{i} \o{U} \To{f_{\o{U}}} X$ o{\`u} $\o{U}$ est quasi-{\'e}tale (au sens de \cite{Bic3}) et compact et $i$ fait de $U$ un ouvert distingu\'e de $\o{U}$. 
Alors d'apr{\`e}s \cite[Key Lemma
7.2.]{Bic3}, pour un tel ouvert \'etale 
 il existe un voisinage de $\id_X$ et, pour tout $\phi$ dans ce voisinage, un isomorphisme {\em canonique}
$i_\phi:\;\;(U,f) \simto \phi^{-1}(U,f)$. En particulier, il existe 
un sous-groupe ouvert
$\aut{U}{X}$ de $\aut{}{X}$ et
une action {\em canonique}  $\aut{U}{X}\To{\beta_U}
\aut{}{U}$ compatible avec $f$ : il suffit de poser
$\beta_U(\phi):\;\; (U,f) \To{i_\phi} \phi^{-1}(U,f) \To{can} (U,f)$.
C'est ce que nous entendons par la phrase "$\aut{U}{X}$ stabilise $(U,f)$".

Maintenant, si $G$ est un groupe topologique agissant sur $X$ par
automorphismes analytiques,
l'action est dite 
{\em continue} si le morphisme $G\To{}\aut{}{X}$ l'est. On notera
alors toujours $G_U$ l'image r{\'e}ciproque d'un groupe $\aut{U}{X}$
comme ci-dessus.
Si $(\FC,\tau)$ est un faisceau $G$-{\'e}quivariant,
et $(U,f)$ un ouvert {\'e}tale distingu{\'e}, l'action $\beta_U$
mentionn{\'e}e  ci-dessus
 munit $\FC(U)$ d'une action de $G_U$.
Berkovich dit alors qu'un faisceau {\'e}tale $G$-{\'e}quivariant en
$\Lambda$-modules  $\FC$
est {\em discret} si l'une des deux conditions {\'e}quivalentes suivantes
est v{\'e}rifi{\'e}e (voir aussi \cite[Part. 2]{HaKyoto}):
\begin{itemize}
\item
Pour tout morphisme {\'e}tale  $U\To{f} X$ 
et toute section $s\in
\FC(U)$, tout point $u\in U$ admet un voisinage distingu{\'e} $\UC$ tel
que  le stabilisateur de $s_{|\UC}$   dans $G_\UC$ soit ouvert (d{\'e}finition qui
ne d{\'e}pend pas du choix de $G_\UC$).
\item 
Pour tout morphisme quasi-{\'e}tale $\o{U}\To{f}X$ (voir \cite[part
3]{Bic3}) avec $\o{U}$ compact,
le $G_{\o{U}}$-module $f^*\FC(\o{U})$ est lisse.
\end{itemize}

La terminologie "faisceau $G$-discret" de Berkovich diff\`ere de celle de la th{\'e}orie des
repr{\'e}sentations qui emploierait plutot les termes {\em lisse} ou
{\em localement constant}
 (mais elle {\'e}vite les conflits terminologiques avec la th{\'e}orie
des faisceaux).
Pour \'eviter toute ambigu\"it\'e nous parlerons dor\'enavant de faisceaux $G$-\'equivariants {\em discrets/lisses}.

\bigskip
\noindent {\em Exemples :}
\begin{itemize}
\item Le groupe $G=PGL(V)$ agit contin{\^u}ment sur $\Omega$ et le faisceau
  constant $\Lambda_\Omega$ est $G$-discret.
\item Le groupe $\gal(K^a/K)$ agit contin{\^u}ment sur $X\hat\otimes
  \wh{K^a}$ pour tout $K$-espace analytique $X$.
\end{itemize}
\bigskip

Notons $\FC(X,\Lambda G_{disc})$, resp. $\FC(X,\Lambda G)$, la
cat{\'e}gorie des $\Lambda$-faisceaux {\'e}tales 
$G$-{\'e}quivariants, resp. $G$-{\'e}quivariants discrets/lisses
(attention aux notations! L'indice $disc$ signifie qu'on oublie la
topologie de $G$).
Ce sont deux cat{\'e}gories ab{\'e}liennes, on a un "plongement
canonique" (foncteur pleinement fid\`ele) $\FC(X,\Lambda G)\injo\FC(X,\Lambda G_{disc})$, et une suite
de morphismes dans $\FC(X,\Lambda G)$ est exacte \ssi\ elle l'est dans
$\FC(X,\Lambda G_{disc})$ \ssi\ la suite de morphismes de faisceaux
sous-jacente l'est.
Le plongement $\FC(X,\Lambda G)\injo\FC(X,\Lambda G_{disc})$ admet un adjoint \`a droite, dit  foncteur de 
``lissification''
$$\infty:\;\application{}{\FC(X,\Lambda G_{disc})}{\FC(X,\Lambda
  G)}{\FC}{\FC^\infty}$$
$\FC^\infty$ (ou $\infty(\FC)$ selon l'humeur) est le faisceau engendr\'e par le pr\'efaisceau sur les ouverts \'etales distingu\'es $U\mapsto \FC(U)^\infty$
 o{\`u} la notation
$\FC({U})^\infty$ repr{\'e}sente le sous-$G_{{U}}$-module lisse maximal de $\FC({U})$.
En particulier, pour tout ouvert \'etale $U\To{} X$, on a 
$$\Gamma_c(U,\FC^\infty)=\Gamma_c(U,\FC)^\infty.$$

\begin{boite}{Remarque :}
Une jolie construction de  \cite{Bicunp} d{\'e}finit  un site, not{\'e} $X(G)_{et}$ dont
la cat{\'e}gorie des $\Lambda$-faisceaux est isomorphe {\`a} $\FC(X,\Lambda
G)$. Il y a aussi un morphisme canonique de sites 
$X(G_{disc})_{et}\To{\iota} X(G)_{et}$ tel que $\iota^*$ s'identifie
au plongement $\FC(X,\Lambda G)\injo \FC(X,\Lambda
G_{disc})$ et par cons{\'e}quent le foncteur de lissification
$\infty$ s'identifie {\`a}  $\iota_*$. Ce langage est plus abouti que
celui, plus {\'e}l{\'e}mentaire, que nous suivons ici (voir remarque suivante...).
  \end{boite}

\alin{D{\'e}finition du $R\Gamma_c$}
D'apr{\`e}s Berkovich \cite{Bicunp} la cat{\'e}gorie $\FC(X,\Lambda G)$ des
faisceaux {\'e}tales $G$-{\'e}quivariants 
discrets/lisses sur $X$ contient suffisamment d'objets injectifs ; c'est aussi
une cons{\'e}quence de l'existence de l'induction, voir 
\ref{induction}.
 On notera
$D^+(X,\Lambda G)$ sa cat{\'e}gorie d{\'e}riv{\'e}e des complexes born{\'e}s inf{\'e}rieurement et
$D^b(X,\Lambda G)$ la sous-cat{\'e}gorie triangul{\'e}e des objets
cohomologiquement born{\'e}s de celle-ci. 

Si $Y\To{\phi} X$ est un morphisme $G$-{\'e}quivariant d'espaces
analytiques, une v\'erification \'el\'ementaire montre que le foncteur
$\phi_!$ image directe {\`a} supports propres\footnote{Ce serait faux pour
  $\phi_*$}
induit un foncteur $\phi_!^\infty : \FC(Y,\Lambda G)\To{}\FC(X,\Lambda G)$ exact {\`a}
gauche dont on note $R\phi^\infty_!:D^+(Y,\Lambda G)\To{} D^+(X,\Lambda G)$ le foncteur
d{\'e}riv{\'e}.

De m{\^e}me on obtient le foncteur d{\'e}riv{\'e} des sections {\`a} support
compact $$R\Gamma_c^\infty(X,.) : D^+(X,\Lambda G) \To{}
D^+(\Lambda G)$$ o{\`u} $D^+(\Lambda G)$ d{\'e}signe la cat{\'e}gorie d{\'e}riv{\'e}e des
$\Lambda G$-modules lisses (modules {\`a} gauche, de par nos conventions). 

La question naturelle qui se pose est alors de savoir si ces foncteurs
ont la bonne cohomologie, c'est-{\`a}-dire celle des foncteurs usuels
$R\phi_!$ et $R\Gamma_c$.
Pour cela notons $\iota_X :
\FC(X,\Lambda G) \To{} \FC(X,\Lambda)$ le foncteur d'oubli. On a
{\'e}videmment $\phi_!\circ \iota_Y =\iota_X \circ \phi^\infty_!$ et
$\Gamma_c(X,.)\circ \iota_X = \iota_{pt}\circ
\Gamma^\infty_c(X,.)$.
Berkovich montre alors
 \begin{prop} \label{torsion0}
On a $R\phi_!\circ \iota_Y =\iota_X \circ R\phi^\infty_!$ et
   $R\Gamma_c(X,.)\circ \iota_X = \iota_{pt}\circ R\Gamma^\infty_c(X,.)$
\end{prop}

Il suffit en fait de prouver que tout objet de $\FC(X,\Lambda G)$
admet une r{\'e}solution par des objets de $\FC(X,\Lambda G)$ dont les
faisceaux sous-jacents sont $\Gamma_c$-acycliques. Comme cette
propri{\'e}t{\'e} est tr{\`e}s importante (elle est aussi utilis{\'e}e dans
\cite{HaCusp} et \cite{HaKyoto}), et que nous aurons besoin d'une variante par la
suite, nous en donnons une preuve assez d{\'e}taill{\'e}e dans l'appendice
(le r{\'e}sultat original est le corollaire \ref{coroBerk} et la variante,
plus simple et suffisante ici, est la proposition \ref{alaBerk}), en
suivant les id{\'e}es de Berkovich, sous une forme plus {\'e}l{\'e}mentaire (nous
n'introduisons pas le site $X(G)$).

Ce r{\'e}sultat nous conduira {\`a} faire sans commentaires suppl{\'e}mentaires
l'abus de noter simplement 
$R\Gamma_c$ pour $R\Gamma_c^\infty$.
Appliquant le r{\'e}sultat de
finitude cohomologique de \cite[5.3.8]{Bic2} on obtient lorsque $\Lambda$
est {\em de 
  torsion} et $K$ est alg{\'e}briquement clos :

$$R\Gamma_c(X, .) : D^b(X,\Lambda G) \To{}
D^b(\Lambda G)$$

\subsection{Coefficients $l$-adiques} \label{coefladic}

 L'auteur ne connait pas de formalisme $l$-adique complet pour les
espaces analytiques, c'est-{\`a}-dire l'existence pour chaque espace
analytique $X$ d'une  cat{\'e}gorie
triangul{\'e}e $D^b_c(X,\o\QM_l)$ stable par les ``six
op{\'e}rations''. Berkovich a d{\'e}fini des groupes de cohomologie $l$-adique
{\`a} supports compacts pour un $\QM_l$-faisceau ``lisse'' 
(non publi{\'e}). Il  a aussi montr{\'e} que lorsqu'un 
tel faisceau est muni d'une action   
d'un groupe localement pro-$p$ et que l'espace analytique est
``quasi-alg{\'e}brique'', alors l'action obtenue en  cohomologie 
est
lisse, voir \cite[4.1.19]{Fargues}.

N{\'e}anmoins, ces groupes de cohomologie ne sont pas d{\'e}finis par des foncteurs d{\'e}riv{\'e}s. 
Or, pour
le pr{\'e}sent article, nous avons {\'e}videmment
besoin d'un formalisme en cat{\'e}gories d{\'e}riv{\'e}es.
Ainsi nous allons, en suivant des id{\'e}es originales de Jannsen
reprises par Huber puis Fargues, attacher {\`a} un espace
analytique raisonnable $X$ 
muni d'une action continue d'un groupe $G$ localement pro-$p$-fini,
un complexe
$R\Gamma_c^\infty(X,\o\QM_l) \in D^b(\o\QM_lG)$ canonique dans la cat{\'e}gorie d{\'e}riv{\'e}e des
$\o\QM_lG$-modules lisses et qui calcule la  cohomologie $l$-adique de
$X$ munie de l'action de $G$.

Notons que dans l'un des cas qui nous int{\'e}ressent pour cet article,
l'action de $GL_d(K)$ sur les rev{\^e}tements de $\Omega^{d-1}_K$,
Harris a d{\'e}ja exhib{\'e} dans \cite{HaCusp} un tel complexe de repr{\'e}sentations lisses
repr{\'e}sentant la cohomologie et l'a notamment utilis{\'e} pour d{\'e}finir une suite
spectrale d'uniformisation. Sa construction repose sur l'existence
d'un recouvrement par des ouverts distingu{\'e}s quasi-alg{\'e}briques
de nerf l'immeuble de Bruhat-Tits et sur la r{\'e}solution fonctorielle de
Berkovich \ref{coroBerk}. Il r{\'e}sultera des arguments qui suivent que
son complexe est un repr{\'e}sentant de notre $R\Gamma_c^\infty$.



\alin{La d{\'e}finition de Berkovich} \label{defBerk}
Fixons un anneau $\Lambda$ de valuation discr\`ete, complet 
 et de caract{\'e}ristique r{\'e}siduelle $\neq
p$, dont on note $\mG$ l'id\'eal maximal.
Soit $X$ un espace $K$-analytique. Nous conviendrons ici d'appeler un
syst{\`e}me projectif $(\FC_n)_n$ de $\Lambda$-faisceaux sur $X$ 
un ``$\Lambda$-syst{\`e}me local'' si
\begin{enumerate}
\item $\forall m\geq n,\; \FC_m\otimes \Lambda/\mG^n \simto \FC_n$.
\item $\forall n\in \NM,\;\; \FC_n$ est un faisceau {\'e}tale localement constant fini.
\end{enumerate}
Notons $\UM(X)$ l'ensemble des ouverts distingu{\'e}s de $X$. 
Berkovich d{\'e}finit  pour tout $q\in \NM$
$$ H^q_c(X,(\FC_n)_n):= \limi{U\in \UM(X)} \limp{n}
H^q_c(U,\FC_n). $$
Nous poserons aussi 
$$H^q_c(X,\Lambda):=H^q_c(X,(\Lambda/\mG^n)_n) \;\;\hbox{ et }\;\;
H^q_c(X,Q):= Q\otimes H^q_c(X,\Lambda)$$
 pour toute extension alg{\'e}brique $Q$ du corps des fractions de
$\Lambda$.

\alin{La d{\'e}finition de Jannsen-Huber-Fargues} \label{defladic}
Nous rappelons ici la d{\'e}finition de Huber \cite{Huber} dans le cadre
des espaces adiques, inspir{\'e}e de \cite{Jannsen} et reprise par Fargues
\cite[4.1]{Fargues} dans le contexte des espaces analytiques  de Berkovich.
 
On appelle $\Lambda_\bullet$-faisceau \'etale sur $X$ tout syst\`eme projectif $(\FC_n)_{n\in\NM}$ de $\Lambda$-faisceaux sur $X$ v\'erifiant $\mG^n\FC_n=0$ pour tout $n\in\NM$. Ces objets, munis d'une notion \'evidente de morphisme, forment une cat\'egorie ab\'elienne que nous noterons
$\FC(X,\Lambda_\bullet)$. Par des arguments
g{\'e}n{\'e}raux \cite[(1.1)]{Jannsen}, cette cat{\'e}gorie poss{\`e}de assez d'objets injectifs.
 Soit $\limproj$ le
foncteur
$$\application{\limproj:\;\;}{\FC(X,\Lambda_\bullet)}{\FC(X,\Lambda)}{(\FC_n)_{n\in
    \NM}}{\limp{n} \FC_n}, $$
on note exceptionnellement $\Gamma_!(X,.)$ le foncteur des sections {\`a}
supports compacts et on d{\'e}finit
$$ \Gamma_c(X,.):=\Gamma_!(X ,.)\circ \limproj :\FC(X,\Lambda_\bullet) \To{}
\Lambda-\hbox{mod}. $$
C'est un foncteur exact {\`a} gauche 
et on peut donc consid{\'e}rer son foncteur d{\'e}riv{\'e} $R\Gamma_c$. On sait
alors, \cite[4.1.9]{Fargues}, que lorsque $(\FC_n)_n$ est un
$\Lambda$-syst{\`e}me local, on a 
$$ R^q\Gamma_c(X,(\FC_n)_n) \simto H^q_c(X,(\FC_n)_n).$$
Rappelons aussi que d'apr{\`e}s \cite[4.1.4]{Fargues}, on a $R\Gamma_c
=R\Gamma_!\circ R\limproj$, ce qui permet de calculer la cohomologie
du syst{\`e}me $(\FC_n)_n$ comme l'hypercohomologie d'un complexe de $\Lambda$-faisceaux.

On d{\'e}finit enfin
pour toute extension alg{\'e}brique $Q$ du corps des fractions de
$\Lambda$ : 
 $$ R\Gamma_c(X,\Lambda):=R\Gamma_c(X,(\Lambda/\mG^n)_X) 
\;\;\hbox{ et }\;\;
 R\Gamma_c(X,Q):=Q\otimes_\Lambda R\Gamma_c(X,\Lambda)
$$

\alin{Actions d'un groupe discret} \label{defladicequ}
 Supposons maintenant donn{\'e}e une
action d'un groupe discret
$G$ sur $X$. On rappelle ici des constructions de \cite{Fargues}.
 Notons 
$\FC(X,\Lambda_\bullet G_{disc})$ la cat{\'e}gorie ab{\'e}lienne des
syst{\`e}mes projectifs $(\FC_n)_{n\in\NM}$ de $\Lambda$-faisceaux
$G$-{\'e}quivariants sur
$X$ tels que pour tout $n$ on a $\mG^n\FC_n=0$. De nouveau, on a des
foncteurs $\limproj^{disc}: \FC(X,\Lambda_\bullet G_{disc}) \To{} \FC(X,\Lambda
G_{disc})$ et $\Gamma_!^{disc}(X,.): \FC(X,\Lambda
G_{disc}) \To{} \Mo{\Lambda}{G_{disc}}$. 
On pose alors
$\Gamma_c^{disc}:=\Gamma_!^{disc} \circ \limproj^{disc}$ que l'on peut
d{\'e}river {\`a} droite et on montre comme pour les faisceaux non
$G$-{\'e}quivariants que $R\Gamma_c^{disc}=R\Gamma_!^{disc}\circ
R\limproj^{disc}$. D'apr{\`e}s \cite[lemme 4.2.6]{Fargues} le 
diagramme
$$\xymatrix{ D^+(X,\Lambda_\bullet G_{disc}) \ar[d] \ar[r]^{R\Gamma_c^{disc}}
  &  D^+(\Lambda G^{disc}) \ar[d] \\
D^+(X,\Lambda_\bullet)      \ar[r]^{R\Gamma_c} & D^+(\Lambda)}$$
dans lequel les fl{\`e}ches verticales sont les foncteurs d'oubli, est
commutatif.
En particulier le complexe $R\Gamma_c^{disc}(X,(\Lambda/\mG^n)_X)$
calcule bien la  cohomologie de $X$ {\`a} coefficients dans $\Lambda$ et
munie de son action de $G$.

Les choses se compliquent un peu lorsqu'on s'int{\'e}resse aux faisceaux
et modules lisses.

\alin{Actions lisses  d'un groupe topologique} \label{deflimpinf} on suppose dor{\'e}navant
que $G$ est localement pro-$p$ et agit contin{\^u}ment sur $X$.
 On note alors $\FC(X,\Lambda_\bullet G)$ la cat{\'e}gorie des
syst{\`e}mes projectifs $(\FC_n)_{n\in\NM}$ de faisceaux
$G$-{\'e}quivariants {\em discrets/lisses} tels que $\mG^n\FC_n=0$.
Pour un tel syst{\`e}me projectif, la limite projective $\limproj \FC_n$
est un $\Lambda$-faisceau $G$-\'equivariant qui n'est g\'en\'eralement pas discret/lisse. L'id{\'e}e na{\"\i}ve consistant
{\`a} lissifier la limite projective (ce qui revient {\`a} repr{\'e}senter la
limite projective dans la cat{\'e}gorie des faisceaux discrets/lisses) semble
{\^e}tre la bonne, comme on va l'expliquer maintenant.

Notons donc $\limproj^\infty$ le foncteur
$$ \application{\limproj^\infty:\;\;}{\FC(X,\Lambda_\bullet
  G)}{\FC(X,\Lambda G)}{(\FC_n)_n}{\infty(\limp{n}\FC_n)}$$
et $\Gamma_!^\infty : \FC(X,\Lambda G) \To{} \Mo{\Lambda}{G}$ le 
foncteur des sections {\`a} support compact. On pose 
$$\Gamma_c^\infty :=
\Gamma_!^\infty \circ \limproj^\infty.$$


	Par des arguments g\'en\'eraux (voir aussi \ref{preuveladic2}), la cat{\'e}gorie $\FC(X,\Lambda G)$
 poss{\`e}de  assez d'objets injectifs et l'on peut donc d{\'e}river les foncteurs exacts {\`a} gauche. On obtient ainsi le foncteur
 $$ R\Gamma_c^\infty : D^+(X,\Lambda G) \To{} D^+(\Lambda G).$$

Notons
 maintenant les plongements canoniques 
$ \iota_{disc,X} :\;\;\FC(X,\Lambda_\bullet G) \To{}
\FC(X,\Lambda_\bullet G_{disc}) $ et
 $\iota_{disc} :\;\; \Mo{\Lambda}{G}\To{} \Mo{\Lambda}{G_{disc}}$.
On a une tranformation naturelle 
$$ \Phi:\;\; \iota_{disc} \circ \Gamma_c^\infty \To{} \Gamma^{disc}_c \circ
\iota_{disc,X}.$$

\begin{lemme} \label{ladic0}
 On a $R(\Gamma_c^{disc}\circ
    \iota_{disc,X})=R\Gamma_c^{disc}\circ \iota_{disc,X}$.
\end{lemme}
\begin{proof}
Il s'agit de voir que $\iota_{disc,X}$ envoie suffisamment
d'objets de $\FC(X,\Lambda_\bullet G)$ sur des  objets
$\Gamma_c^{disc}$-acycliques de $\FC(X,\Lambda_\bullet G_{disc})$ (c'est-\`a-dire dont le faisceau sous-jacent est $\Gamma_c$-acyclique). La construction de tels objets repose sur la construction de Berkovich (pour des faisceaux de torsion) expliqu\'ee dans l'appendice et sur les techniques de cohomologie de syst\`emes projectifs de faisceaux \`a la Jannsen. Les d\'etails sont donn\'es dans l'appendice, paragraphe \ref{preuveladic2}
\end{proof}

\begin{prop} \label{ladic1}
Soit $R\Phi : \iota_{disc} \circ R\Gamma_c^\infty \To{} R\Gamma^{disc}_c \circ
\iota_{disc,X}$ la tranformation naturelle induite par $\Phi$ et le
lemme pr\'es\'edent. 
Pour tout $p\in \NM$, la transformation $R^p\Phi$ induit un
isomorphisme de foncteurs
$$R^p\Gamma_c^\infty \simto \infty\circ R^p\Gamma_c^{disc}\circ
\iota_{disc,X}.$$
\end{prop}

\begin{proof}
La preuve est un peu technique mais l'id{\'e}e est
simple : par d{\'e}finition on a $\Gamma_c^\infty=\Gamma^\infty_!\circ
\limproj^\infty = \infty\circ \Gamma_! \circ \limproj.$  
Le foncteur de lissification $\infty$ n'est pas exact sur la cat{\'e}gorie des $\Lambda
G_{disc}$-modules et d{\'e}river l'identit{\'e} ci-dessus devrait donc faire
intervenir le foncteur $R\infty$. Ce qu'affirme la proposition que nous
voulons montrer c'est qu'il n'en est rien. La raison est
essentiellement que le foncteur $\Gamma_!\circ \limproj$ envoie la
cat{\'e}gorie $\FC(X,\Lambda_\bullet G)$ dans une sous-cat{\'e}gorie
$\infty$-acyclique de $\Mo{\Lambda}{G_{disc}}$.

Pour formaliser tout {\c c}a, nous devons introduire auparavant une
certaine cat{\'e}gorie de $\Lambda G$-modules, interm{\'e}diaire entre la
cat{\'e}gorie lisse et la cat{\'e}gorie de tous les $\Lambda G$-modules.

\bigskip

\noindent{\sl Modules sur l'alg{\`e}bre de distributions :} Soit $G$ un groupe
localement pro-$p$ et 
$\Lambda$ un anneau tel que $p\in \Lambda^\times$. Notons $\DC_\Lambda(G)$, resp. $\HC_\Lambda(G)$,
l'alg{\`e}bre de convolution des mesures, resp. mesures  localement constantes, {\`a} support
compact et {\`a} valeurs dans $\Lambda$ sur $G$. 
L'anneau $\HC_\Lambda(G)$ n'a pas d'unit{\'e} mais assez d'idempotents  et 
on sait que la cat{\'e}gorie des $\Lambda G$-modules lisses est
canoniquement isomorphe {\`a} la cat{\'e}gorie des $\HC_\Lambda(G)$-modules
lisses (appel\'es aussi non-d{\'e}g{\'e}n{\'e}r{\'e}s ou unitaux). L'anneau $\DC_\Lambda(G)$ est
unitaire et contient les anneaux $\Lambda[G]$
et $\HC_\Lambda(G)$. On a $\DC_\Lambda(G)*\HC_\Lambda(G)=
\HC_\Lambda(G)$ de sorte que toute structure de
$\HC_\Lambda(G)$-module lisse s'{\'e}tend canoniquement en une structure
de $\DC_\Lambda(G)$-module : en effet, si $M$ est un $\HC_\Lambda(G)$-module lisse,
 $m\in M$ et $d\in \DC_\Lambda(G)$, on pose $d.m:=(d*e).m$, expression qui ne d\'epend pas du choix de l'idempotent $e$ de $\HC_\Lambda(G)$ fixant $m$.
On obtient le foncteur $\iota_\infty^D$ du  
 syst{\`e}me de foncteurs suivant :
$$\xymatrix{  & \DC_\Lambda(G)-\hbox{mod} \ar@<.5ex>[dl]^{\infty_D}
  \ar[dr]^{\iota_{D}^{disc}} & \\
\Mo{\Lambda}{G}  \ar@<.5ex>[ru]^{\iota_{\infty}^D}
\ar@<.5ex>[rr]^{\iota_{\infty}^{disc}}& & \Mo{\Lambda}{G_{disc}} 
\ar@<.5ex>[ll]^{\infty_{disc}} }$$
Les autres foncteurs $\iota$ sont des foncteurs d'oublis.
 Les foncteurs $\infty$ sont des foncteurs de lissification ; celui qui est
 not\'e $\infty_{disc}$ est le foncteur de lissification "habituel" que nous notons ailleurs simplement $\infty$ et qui {\`a} un $\Lambda G_{disc}$-module associe
 le sous-module de ses vecteurs lisses. On d{\'e}finit $\infty_D$
pour tout $\DC_\Lambda(G)$-module $M$ par la formule
$$ \infty_D(M): =\limi{H\subset G} e_H M,$$
la limite {\'e}tant prise sur un syst{\`e}me de voisinages de l'unit{\'e} form{\'e} de
pro-$p$-sous-groupes ouverts, et la notation $e_H$ d{\'e}signant
l'idempotent de $\HC_\Lambda(G)$ associ{\'e} {\`a} $H$. L'avantage d'avoir
introduit la cat{\'e}gorie $\DC_\Lambda(G)-\hbox{mod}$ 
est que le foncteur $\infty_D$ y est exact.
Cependant,
on prendra garde au
fait qu'en g{\'e}n{\'e}ral, $\infty_D(M) \neq \infty_{disc}(\iota_D^{disc}(M))$.
N{\'e}anmoins cet inconv{\'e}nient sera inoffensif pour l'utilisation que nous
ferons de ces objets. La motivation pour introduire l'anneau
$\DC_\Lambda(G)$ est la remarque suivante : soit $(M_n)_n$ un syst{\`e}me
projectif de $\Lambda G$-modules lisses. Le $\Lambda G$-module produit
$\limproj M_n$ n'est en g{\'e}n{\'e}ral pas lisse, mais est canoniquement muni
d'une structure de $\DC_\Lambda(G)$-module. De plus on a
$\infty_D(\limproj M_n)=\infty_{disc}(\iota_D^{disc}(\limproj M_n))$.

On obtient de cette mani{\`e}re
un foncteur de la cat{\'e}gorie des syst{\`e}mes projectifs de $\Lambda
G$-modules lisses vers celle des $\DC_\Lambda(G)$-modules.
En particulier, si $(\FC_n)_n$ est un syst{\`e}me projectif de
$\Lambda$-faisceaux $G$-{\'e}quivariants $G$-discrets/lisses, on a une
structure de $\DC_\Lambda(G)$-module sur $\limproj \,
\Gamma_c(X,\FC_n)$. Or, puisque les foncteurs de sections $\Gamma(U,.)$ commutent aux limites projectives, on a une inclusion
$$\Gamma_!(X,\limproj \FC_n) \subseteq \limproj \; \Gamma_c(X,\FC_n).$$
Puisque les mesures de $\DC_\Lambda(G)$ sont \`a supports compacts, le
sous-$\Lambda G$-module $\Gamma_!(X,\limproj \FC_n)$ est stable par
$\DC_\Lambda(G)$, de sorte qu'on obtient
 un foncteur 
$$\application{\Gamma_c^D:\;\;}{\FC(X,\Lambda_\bullet
  G)}{\DC_\Lambda(G)-\hbox{mod}}{(\FC_n)_n}{\Gamma_!(X,\limp{n} \FC_n)}$$
et une factorisation
$$\Phi :\;\; \iota_\infty^{disc}\circ \Gamma_c^\infty  \To{\Phi^{\infty,D}}
\iota_D^{disc} \circ \Gamma_c^D \To{\Phi^{D,disc}} \Gamma_c^{disc}\circ
\iota_{disc,X}. $$ 
Par les remarques pr{\'e}c{\'e}dentes, $\Phi^{\infty,D}$ et $\Phi^{D,disc}$ induisent des isomorphismes
$$\Gamma_c^\infty \simto \infty_D \circ \Gamma_c^D \;\;\hbox{  et }\;\;
 \iota_D^{disc}\circ \Gamma_c^D\simto \Gamma_c^{disc} \circ \iota_{disc,X} ,$$
 que l'on peut d{\'e}river en
$$
R\Gamma_c^\infty \simto \infty_D \circ R\Gamma_c^D \;\;\hbox{  et }\;\;
\iota_D^{disc}\circ R\Gamma_c^D\simto R\Gamma_c^{disc} \circ \iota_{disc,X}  ,$$
en utilisant l'exactitude de $\infty_D$ pour la premi{\`e}re et le lemme pr\'ec\'edent pour la deuxi{\`e}me.
Soit maintenant $p\in\NM$, en combinant les divers isomorphismes ci-dessus on obtient
$$ R^p\Gamma_c^\infty \simto \infty_D\circ R^p\Gamma_c^D \simeq \infty_{disc}\circ \iota_D^{disc} \circ R^p\Gamma_c^D \simto \infty \circ R^p\Gamma_c^{disc} \circ \iota_{disc,X}.$$

\end{proof}

La proposition \ref{ladic1} montre qu'en g{\'e}n{\'e}ral $R^p\Gamma_c^\infty
(X,(\FC_n)_n)$ est diff{\'e}rent de $R^p\Gamma_c (X,(\FC_n)_n)$ (c'est par
exemple le cas pour le syst{\`e}me $(j_!(\ZM/l^n))_n$ o{\`u} $j$ est
l'inclusion de $\Omega^{d-1}_K$ dans $\PM^{d-1}_K$).
 
Cependant on sait par Berkovich (voir
aussi \cite[4.1.19]{Fargues}) que sous l'hypoth{\`e}se technique de
quasi-alg{\'e}bricit{\'e} de $X$, {\em cf} \cite[4.1.11]{Fargues}, et parce
que nous avons 
suppos{\'e} $G$ localement pro-$p$ et $p\in \Lambda^\times$, l'action
de $G$ sur les groupes de 
cohomologie {\`a} supports compacts d'un $\Lambda$-syst{\`e}me local $G$-{\'e}quivariant
$(\FC_n)_n$ est lisse. 
On obtient alors
\begin{coro} \label{coroladic}
  Supposons $X$ quasi-alg{\'e}brique et soit $(\FC_n)_n$ un
  $\Lambda$-syst{\`e}me local $G$-{\'e}quivariant. Alors le morphisme
$$ R\Phi((\FC_n)_n):\;\;\iota_{disc}R\Gamma_c^\infty(X,(\FC_n)_n) \To{}
R\Gamma_c(X,(\FC_n)_n)$$
est un isomorphisme. En particulier, on a des isomorphismes canoniques
de $\Lambda G$-modules (lisses)
$$ R^q\Gamma_c^\infty(X,(\FC_n)_n) \simto H^q_c(X,(\FC_n)_n).$$
\end{coro}

 Nous ferons parfois l'abus de notation suivant :
 \begin{nota}
  Pour $X$
quasi-alg{\'e}brique et $Q$ une extension alg{\'e}brique de
$\hbox{Frac}(\Lambda)$, on pose
$$ R\Gamma_c(X,\Lambda):=R\Gamma_c^\infty(X,(\Lambda/\mG^n)_n) \;\;\hbox{ et }\;\;R\Gamma_c(X,Q):=Q\otimes_\Lambda R\Gamma_c(X,\Lambda).$$
\end{nota}
De cette mani{\`e}re on verra toujours $R\Gamma_c(X,\Lambda)$ comme un objet de
$D^b(\Lambda G)$ et pas seulement de $D^b(\Lambda G_{disc})$.

Signalons ici le r\'esultat suivant que nous n'utiliserons pas dans ce texte :
\begin{prop} \label{ladic2}
  On a $R\Gamma_c^\infty = R\Gamma_!^\infty \circ R\limproj^\infty$.
\end{prop}
La preuve en est donn{\'e}e dans l'appendice en \ref{preuveladic2}. Comme
pr{\'e}c{\'e}demment, l'int{\'e}r{\^e}t de cette proposition est qu'elle permet de
calculer la cohomologie munie de son action lisse comme
l'hypercohomologie d'un complexe de $\Lambda$-faisceaux
$G$-{\'e}quivariants discrets.

Pour terminer, on \'enonce et prouve le lemme suivant qui est utilis\'e dans la partie \ref{drinfeld} :
\begin{lemme} \label{lemmelimpinf}
Soit $Y\To{\psi} X$ un morphisme \'etale et $G$-\'equivariant.
On a un isomorphisme canonique de foncteurs $\FC(X,\Lambda_\bullet G) \To{} \FC(Y,\Lambda G)$.
$$ \psi^*\circ\limproj^\infty \simto \limproj^\infty \circ \psi^* .$$
\end{lemme}
\begin{proof}
Il suffit de v\'erifier que $\psi^*$ commute avec chacun des foncteurs $\limproj$ et $\infty$ (lissification). Pour le premier, c'est imm\'ediat. Pour le second, rappelons que $\psi^*$ poss\`ede un adjoint \`a gauche $\psi_!$. Ainsi avec les notations du paragraphe \ref{deflimpinf}, les foncteurs $\infty\circ \psi^*$ et $\psi^*\circ \infty$ sont respectivement adjoints \`a droite des foncteurs $\psi_!\circ \iota_{disc,Y}$ et $\iota_{disc,X}\circ \psi_!$. Mais ces deux derniers foncteurs sont bien-s\^ur canoniquement isomorphes.

\end{proof}

\begin{boite}{Remarque :}
En identifiant $\FC(X,\Lambda G)$ avec la cat{\'e}gorie
$\FC(X(G),\Lambda)$ des
$\Lambda$-faisceaux sur le site $X(G)_{et}$ ({\em cf} les remarques
pr{\'e}c{\'e}dentes), on constate  que le foncteur $\limproj^\infty$
s'identifie au foncteur limite projective usuel dans
$\FC(X(G),\Lambda)$.
\end{boite}

\subsection{Suites spectrales de Hochschild-Serre}

Dans cette section, on s'int\'eresse \`a la situation suivante : un groupe
discret $G$ agit librement sur un $K$-espace analytique $X$. On sait alors former le quotient $X/G$, {\em cf}  \cite[lemma 4]{Bic5},  et  l'application quotient $X\To{p} X/G$ est un rev{\^e}tement
analytique Galoisien. En particulier, $p$ est donc {\'e}tale. Dans ce
contexte, \cite[4.4]{Fargues} montre que pour un "faisceau $l$-adique"
ou sur un anneau Gorenstein 
fini, il existe une suite spectrale de
Hochschild-Serre, c'est-{\`a}-dire du type :
$$E_2^{pq}=\ext{p}{H^q_c(X,p^*\FC)}{1}{G} \Rightarrow H^{-q-p}_c(X/G,\FC)^*.$$
Nous avons besoin d'une interpr{\'e}tation en termes de cat{\'e}gories
d{\'e}riv{\'e}es de cette suite spectrale. Ceci nous conduira d'ailleurs {\`a}
donner une preuve diff{\'e}rente de celle de \cite[4.4]{Fargues}.
Le r\'esultat principal de cette section s'\'enonce ainsi (les notations seront expliqu\'ees au fil de la preuve) :
\begin{theo} \label{theoHS}
Soit $G$ un groupe discret agissant librement sur un $K$-espace analytique $X$ et $p:X\To{} X/G$ le quotient. On suppose $K$ alg\'ebriquement clos.
Si $\Lambda$ est un anneau de torsion premi\`ere \`a $p$ ou un anneau $l$-adique, alors il existe un isomorphisme de foncteurs
$D^b(X/G,\Lambda)\To{} D^-(\Lambda)$, resp. $D^b(X/G,\Lambda_\bullet)\To{} D^-(\Lambda)$ 
$$ \Phi:\;\;\Lambda \otimes_{\Lambda[G]}^L (R\Gamma_c^{eq}(X,p^*(.))) \To{}
R\Gamma_c(X/G,.).$$
\end{theo}

L'hypoth\`ese sur $K$ n'est pas incontournable : il
suffirait de supposer que la dimension cohomologique des
$\gal(K^{ca}/K)$-modules  de $l$-torsions est born{\'e}e pour obtenir les
r{\'e}sultats de cette section sur un anneau de $l$-torsion ou $l$-adique.
Quoiqu'il en soit, avec cette hypoth{\`e}se on sait que :
\begin{itemize}
\item $H^q_c(X,\FC)=0$ pour $q>2\dim(X)$ et tout faisceau ab{\'e}lien de
  torsion, \cite[5.3.8]{Bic2}.
\item $H^q_c(X,(\FC_n))=0$ pour $q>2\dim(X)+1$ pour  tout  $\Lambda_\bullet$-faisceau \'etale comme en \ref{defladic}, si $\Lambda$ est un anneau
  $l$-adique, $l\neq p$, \cite[4.1.9.(b)]{Fargues}.
\end{itemize}

\alin{Quelques sorites}
Nous d{\'e}signons par $\Lambda$ soit un anneau de torsion soit un anneau
$l$-adique. Nous consid{\`e}rerons alors exclusivement dans le premier cas les
cat{\'e}gories de $\Lambda$-faisceaux {\'e}tales et dans le second celles de
$\Lambda_\bullet$-faisceaux {\'e}tales d{\'e}finies en \ref{defladic} et \ref{defladicequ}.

Remarquons tout d'abord que comme l'action de $G$ sur $X/G$ est
triviale, la cat\'egorie $\FC(X/G,\Lambda G)$ est \'equivalente \`a la cat\'egorie
$\FC(X/G,\Lambda[G])$ des faisceaux
{\'e}tales en $\Lambda[G]$-modules. Ainsi le foncteur $\omega
:\Lambda-\hbox{mod} \To{} \Lambda G-\hbox{mod}$ qui munit un
$\Lambda$-module de l'action triviale de $G$ induit aussi un foncteur
$\omega:\FC(X/G,\Lambda) \To{} \FC(X/G,\Lambda G)$. Ces foncteurs
ont pour adjoints {\`a} gauche les foncteurs $M\mapsto
\Lambda\otimes_{\Lambda[G]} M$, resp. $\FC \mapsto
\Lambda\otimes_{\Lambda[G]} \FC$ o{\`u} le morphisme $\Lambda[G]\To{}
\Lambda$ est l'augmentation habituelle. La m{\^e}me discussion s'applique
verbatim aux cat{\'e}gories de  $\Lambda_\bullet$ faisceaux et modules.

{\em Notation} : pour {\'e}viter les confusions nous utiliserons les
notations $\Gamma_c$, resp. $\Gamma_c^{eq}$, pour d{\'e}signer le foncteur des
sections {\`a} support compact lorsqu'on le voit {\`a}
valeurs dans les $\Lambda$-modules, resp. dans les $\Lambda
G$-modules.
Par exemple, on a  $\Gamma^{eq}_c(X/G, \omega(.)) = \omega\circ
\Gamma_c(X/G, .)$ Les m{\^e}mes notations s'appliquent aux
$\Lambda_\bullet$-faisceaux et aux foncteurs $\Gamma_c=\Gamma_!\circ \limproj$.

\alin{Le morphisme $\Phi$}
Comme $p$ est un morphisme \'etale,
 le foncteur $p^*:\; \FC(X/G,\Lambda G)\To{} \FC(X,\Lambda G)$ a un
adjoint {\`a} gauche exact $p_!$, voir \cite[5.4.2(ii)]{Bic2}. 
Le morphisme
d'adjonction $p_!p^*\To{} \id$ induit donc un morphisme de foncteurs
$\FC(X/G,\Lambda) \To{} \Lambda G-\hbox{Mod}$ :
 $$ \Gamma_c^{eq}(X,p^*\omega(.)) \To{} \omega(\Gamma_c(X/G,.)).$$
 Lorsque $\Lambda$ est de torsion, par les r\'esultats de finitude cohomologique rappel\'es plus haut, 
on peut d{\'e}river ceci en un morphisme de foncteurs
$D^b(X/G,\Lambda)\To{} D^b(\Lambda G)$
$$ R\Gamma_c^{eq}(X,p^*\omega(.)) \To{} \omega(R\Gamma_c(X/G,.))$$
qui par adjonction
fournit
un morphisme canonique de
foncteurs $D^b(X/G,\Lambda) \To{} D^-(\Lambda)$
$$ \Phi:\;\;\Lambda \otimes_{\Lambda[G]}^L (R\Gamma_c^{eq}(X,p^*\omega(.))) \To{}
R\Gamma_c(X/G,.) $$
Lorsque $\Lambda$ est $l$-adique, on
construit $\Phi$ de la mani{\`e}re suivante : on v{\'e}rifie d'abord (voir
preuve de \cite[4.2.7]{Fargues}) que pour tout morphisme {\'e}tale
$U\To{f} V$ et tout $\Lambda_\bullet$-faisceau sur $V$, l'isomorphisme 
\'evident $f^*\limproj(\FC_n) \simto \limproj(f^*(\FC_n))$ se d\'erive en un isomorphisme
$f^*R\limproj(\FC_n) \simto R\limproj(f^*(\FC_n))$. On en d{\'e}duit alors un morphisme fonctoriel en $(\FC_n)_n$
\begin{eqnarray*} R\Gamma_c(V,(f^*\FC_n)_n) &\simeq &
  R\Gamma_!(V,R\limproj(f^*\FC_n)) \\
& \simeq & R\Gamma_!(V,f^*R\limproj(\FC_n)) \To{}
R\Gamma_!(U,R\limproj(\FC_n)) = R\Gamma_c(U,(\FC_n)_n),\end{eqnarray*}
la fl\`eche du milieu \'etant encore induite par l'adjonction $f_!f^*\To{}\id$.

Appliquant ceci  au morphisme $p$, et compte tenu de la finitude cohomologique rappel\'ee plus haut, on obtient formellement de la m{\^e}me
mani{\`e}re que  pour les
faisceaux  de torsion un morphisme de foncteurs
$D^b(X/G,\Lambda_\bullet)\To{} D^-(\Lambda)$ :
$$\Phi:\;\;\Lambda \otimes_{\Lambda[G]}^L (R\Gamma_c^{eq}(X,p^*\omega(.))) \To{}
R\Gamma_c(X/G,.) $$

\begin{prop} \label{HSder}
  Dans les deux cas consid{\'e}r{\'e}s, $\Phi$ est un isomorphisme de foncteurs.
\end{prop}

\begin{proof}
Nous dirons qu'un faisceau, resp. un $\Lambda_\bullet$-faisceau, $\FC$
sur ${X/G}_{et}$ est :
\begin{itemize}
\item $\Gamma_c$-acyclique s'il est $\Gamma_c(U,-)$-acyclique pour
  tout morphisme {\'e}tale $U\To{} X/G$.
\item $\Gamma_{c,et}^\vee$-acyclique, resp. $\Gamma_{c,loc}^\vee$-acyclique si pour tout morphisme surjectif $\xymatrix{U
  \ar@{->>}[r] & V}$ {\'e}tale, resp. localement iso, entre ouverts \'etales au-dessus de $X/G$, le complexe
de Cech ``dual''  
$$\xymatrix{ \cdots \ar[r]\ar@<.8ex>[r]\ar@<-.8ex>[r] &
  \Gamma_c(U\times_{V} U,\FC)
  \ar@<.5ex>[r] \ar@<-.5ex>[r] & \Gamma_c(U,\FC) \ar[r] &
  \Gamma_c(V,\FC) \ar[r] & 0 }$$
est acyclique.
\end{itemize}
Un morphisme $U\To{f} V$ est dit localement-iso si tout point de $U$
admet un voisinage ouvert $\UC$ tel que $f_{|\UC}:\UC\To{} f(\UC)$ est
un isomorphisme.
Un tel morphisme est en particulier {\'e}tale.

\begin{lemme} \label{resol}
  Tout faisceau de torsion, resp. $\Lambda_\bullet$-faisceau,  admet
  une r{\'e}solution born{\'e}e par des faisceaux $\Gamma_c$-acycliques et
  $\Gamma_{c,et}^\vee$-acycliques, resp. par des
  $\Lambda_\bullet$-faisceaux  $\Gamma_c$-acycliques et 
  $\Gamma_{c,loc}^\vee$-acycliques.  
\end{lemme}
Laissons ce lemme de c{\^o}t{\'e} et commen{\c c}ons la preuve de la proposition.
Fixons $\FC\in \FC(X/G,\Lambda)$, resp $\FC\in
\FC(X/G,\Lambda_\bullet)$,  et soit $\FC\injo\JC^\bullet$ une r{\'e}solution
born{\'e}e de $\FC$ comme dans le lemme ci-dessus.
 Par d{\'e}finition on a :
$$ R\Gamma_c(X/G,\FC)= \Gamma_c(X/G, \JC^\bullet) \;\;\hbox{ et }\;\;
R\Gamma^{eq}_c(X,p^*\omega \FC) = \Gamma^{eq}_c(X,p^*\omega \JC^\bullet)
$$
(en remarquant que $p^*\omega(\JC^\bullet)$ est une r{\'e}solution de
$p^*\omega(\FC)$ par des $G$-faisceaux $\Gamma_c^{eq}(X,.)$-acycliques).
La proposition est alors une cons{\'e}quence imm{\'e}diate du lemme suivant :
\begin{lemme}
Soit $\JC$ un faisceau de torsion, resp $\Lambda_\bullet$-faisceau,
$\Gamma_{c,loc}^\vee$-acyclique sur ${X/G}_{et}$. Alors 
\begin{enumerate}
\item $\Gamma_c^{eq}(X,p^*\omega \JC)$ est un $\Lambda[G]$-module
  acyclique pour le foncteur $\Lambda\otimes_{\Lambda[G]} -$.
\item Le morphisme canonique $\Lambda\otimes_{\Lambda[G]}
  \Gamma_c^{eq}(X,p^*\omega \JC) \To{} \Gamma_c(X/G,\JC)$ est un isomorphisme.
\end{enumerate}
\end{lemme}

Prouvons ce lemme.
Soit  $U\To{f} X/G$ un morphisme {\'e}tale. Consid{\'e}rons le
produit cart{\'e}sien $U\times_{X/G} X \To{\pi_X} X$ comme un ouvert {\'e}tale
$G$-{\'e}quivariant sur $X$ en le munissant de l'action $\id\times
g_X$. Le $\Lambda$-module $\Gamma_c(U\times_{X/G} X,p^*\omega\FC)$
est alors naturellement muni d'une action de $G$ et nous noterons
$\Gamma_c^{eq}(U\times_{X/G} X,p^*\omega\FC)$ pour pr{\'e}ciser que l'on
tient compte de cette action. La premi{\`e}re projection fournit un
morphisme {\'e}quivariant $\Gamma_c^{eq}(U\times_{X/G} X,p^*\omega\FC)
\To{} \omega(\Gamma_c(U,\FC))$. Dans le cas $U=X/G$, ce morphisme est
celui qui induit notre morphisme $\Phi$.

Supposons maintenant que le morphisme $f$ se factorise $U\To{f'}X\To{p}X/G$.
Alors le morphisme d'espaces analytiques au-dessus de $U$
$$\application{}{U\times G}{U\times_{X/G} X}{(u,g)}{(u,g.f'(u))}$$
est un isomorphisme $G$-{\'e}quivariant si l'on munit $U\times G$ de
l'action par translation {\`a} gauche sur le second terme, qui s'inscrit
dans le diagramme :
$$\xymatrix{ U \ar[r]^-{\id\times f'} \ar[rd]_{\id\times 1_G}  &
  U\times_{X/G}X  \ar[r]^-{\pi_U} & U \\  & U\times G \ar[u]^\sim
  \ar[ru]_{\pi_U} & }.$$
On en d{\'e}duit, dans le cas de torsion comme dans le cas des
$\Lambda_\bullet$-faisceaux,  un isomorphisme de $\Lambda[G]$-modules
fonctoriel en  $\FC$ :
$$ \Lambda[G]\otimes_\Lambda
\Gamma_c(U,\FC)\simto\Gamma_c^{eq}(U\times_{X/G} X,p^*\omega\FC) $$
qui s'inscrit dans un diagramme
$$ \xymatrix{ \Gamma_c(U,\FC) \ar[r]^-{(\id\times f')_!}
  \ar[rd]_{1_G\otimes \id}  &
  \Gamma_c^{eq}(U\times_{X/G}X,\FC)  \ar[r]^-{\pi_{U,!}} &
  \omega\Gamma_c(U,\FC) \\  & \Lambda[G]\otimes_\Lambda\Gamma_c(U,\FC) \ar[u]^\sim 
  \ar[ru]_{\hbox{aug}\otimes \id}& }$$
o{\`u} le triangle de droite est $G$-{\'e}quivariant mais pas celui de gauche
en g{\'e}n{\'e}ral.
Remarquons que l'isomorphisme ci-dessus d{\'e}pend de la
section $f'$ de $f$ contrairement {\`a} l'isomorphisme obtenu en composant
avec le foncteur $\Lambda\otimes_{\Lambda[G]} -$ 
$$ \Gamma_c(U,\FC) \simto \Lambda\otimes_{\Lambda[G]}\Gamma_c^{eq}(U\times_{X/G}
X,p^*\omega\FC) 
$$ 
qui n'est autre que l'inverse du morphisme adjoint au morphisme not{\'e} $\pi_{U,!}$
dans le diagramme.

Notons maintenant $X^{(0)}:=X/G$ et pour tout $i>0$ d{\'e}finissons par r{\'e}currence
$$ X^{(i)}:= X^{(i-1)} \times_{X/G} X \hbox{ et } p^{(i)}: X^{(i)}
\To{} X/G \;\hbox{la projection canonique}.$$
Appliquant la discussion pr{\'e}c{\'e}dente au morphisme {\'e}tale $f=p^{(i)}$,
 on obtient pour $i\geq 1$ des isomorphismes 
$$\Lambda[G]\otimes_\Lambda \Gamma_c(X^{(i)},\FC) \simto
\Gamma_c^{eq}(X^{(i+1)},p^*\omega\FC)  $$ (d{\'e}pendant du choix d'une
projection $X^{(i)}\To{} X$) et des isomorphismes
$$ \Lambda \otimes_{\Lambda[G]}
\Gamma_c^{eq}(X^{(i+1)},p^*(\omega\FC))  \simto \Gamma_c(X^{(i)},\FC)
$$ canoniques, donc en particulier compatibles aux diff{\'e}rentes
projections  $X^{(i)}\To{}X^{(i-1)}$. 

Par ailleurs, on peut {\'e}crire le complexe de Cech ``dual'' d'un faisceau $\FC$ sur $X/G_{et}$
associ{\'e} \`a un recouvrement 
$X\times_{X/G} U\To{\pi_U} U$ sous la forme
$$ \cdots \To{} \Gamma_c( X^{(i)}\times_{X/G} U,\FC) \To{}
\Gamma_c( X^{(i-1)}\times_{X/G} U,\FC)\To{} \cdots \To{}
\Gamma_c( X \times_{X/G}U,\FC) \To{\pi_{U,!}} \Gamma_c(U,\FC)$$
o{\`u} les  diff{\'e}rentielles sont donn{\'e}es par des sommes altern{\'e}es (dont la
forme pr{\'e}cise ne nous importe pas) de
morphismes associ{\'e}s aux projections $X^{(i)}\To{} X^{(i-1)}$.
Ce complexe est naturel par rapport aux $X/G$-morphismes $V\To{\phi}
U$. En particulier pour $V=X\To{\phi=p}U=X/G$, on obtient un morphisme
de complexes
$$\xymatrix{ \cdots \ar[r] & \omega\Gamma_c(X^{(i)},\JC) \ar[r] & \cdots  \ar[r] &
  \omega\Gamma_c(X,\JC) \ar[r]^{}  & \omega\Gamma_c(X/G,\JC) \\
 \cdots \ar[r] & \Gamma_c^{eq}( X^{(i+1)},p^*\omega\JC)
 \ar[r] \ar[u]^{(\id\times p)_!}  & \cdots  \ar[r] &
  \Gamma_c^{eq}(X^{(2)},p^*\omega\JC) \ar[r]^{} \ar[u]^{(\id\times p)_!} &
  \Gamma_c^{eq}(X,p^*\omega\JC) \ar[u]^{p_!} \\ 
}$$
o{\`u} on a rajout{\'e} les symboles $eq$ et $\omega$ pour tenir compte de
l'action de $G$.
Par hypoth{\`e}se, le faisceau $\JC$ est $\Gamma_{c,loc}^\vee$-acyclique,
de sorte que
les deux lignes du diagramme ci-dessus sont {\em exactes}.
D'apr{\`e}s la discussion pr{\'e}c{\'e}dente les termes
$\Gamma_c^{eq}(X^{(i)},p^*\omega\JC)$ pour $i\geq 2$ sont de la
forme $\Lambda[G]\otimes_\Lambda M$ pour un $\Lambda$-module $M$ et
sont donc acycliques pour le foncteur
$\Lambda\otimes_{\Lambda[G]}-$. On a donc des isomorphismes dans
$D^-(\Lambda)$ :
\begin{eqnarray*}
  & & \Lambda\otimes^L_{\Lambda[G]}
  \Gamma_c^{eq}(X,p^*\omega\JC) \\
  & \simeq & \left( \cdots \To{}
  \Lambda\otimes_{\Lambda[G]}\Gamma_c^{eq}(X^{(i+1)},p^*\omega\JC)\To{} \cdots  \To{} 
\Lambda\otimes_{\Lambda[G]}\Gamma_c^{eq}(X^{(2)},p^*\omega\JC) \right)\\
& \simeq & \left( \cdots \To{}
  \Gamma_c(X^{(i)},\JC) \To{} \cdots \To{} 
\Gamma_c(X^{},\JC)\right) \\
& \simeq & \Gamma_c(X/G,\JC).
 \end{eqnarray*}
D'o{\`u} le lemme, puis la proposition. 

\bigskip

Il reste maintenant {\`a} prouver le
lemme \ref{resol}.

\bigskip

\noindent{\em Cas $\Lambda$ de torsion} : dans ce cas, ceci doit {\^e}tre
bien connu. Tronquons une r{\'e}solution injective de $\FC$ en degr{\'e}
$2\dim(X)$ et notons $\JC^\bullet$ le complexe obtenu. Comme on sait
que les foncteurs $\Gamma_c(U,-)$ pour $U$ 
ouvert {\'e}tale de $X/G$ sont de dimension cohomologique $\leq 2\dim(X)$,
les faisceaux $\JC^i$ sont $\Gamma_c$-acycliques au sens donn{\'e}
plus haut. Remarquons maintenant que le complexe de Cech ``dual'' d{\'e}fini plus
haut associ{\'e} au recouvrement {\'e}tale $U\To{f} V$ n'est autre que l'image
par le foncteur $\Gamma_c(V,-)$ du complexe de faisceaux donn{\'e} par :
$$ \cdots \To{} f^{(i+1)}_!{f^{(i+1)}}^*\JC \To{d_{i+1}}
f^{(i)}_!{f^{(i)}}^*\JC \To{} \cdots \To{d_1} f_!f^*\JC  \To{d_0} \JC$$
o{\`u} $f^{(i)}$ d{\'e}signe la projection $U\times_V\cdots \times_V U\To{} V$
($i$ facteurs) et les diff{\'e}rentielles sont des sommes altern{\'e}es des
morphismes induits par les diff{\'e}rentes projections. Ce complexe est
{\em exact} comme on le v{\'e}rifie  
ais{\'e}ment sur les fibres. Montrons alors que le complexe de Cech
``dual'' est lui aussi
exact : soit $n\geq 0$, tronquons le complexe de faisceaux ci-dessus
en degr{\'e} $-n-2\dim(X)-1$ :
 $$ C_{-n-2\dim(X)-1} \To{} f^{(n+2\dim(X))}_!{f^{(n+2\dim(X))}}^*\JC
 \To{}\cdots\To{} \JC,$$
 on obtient une r{\'e}solution du premier terme $C_{-n-2\dim(X)-1}$ par
des faisceaux $\Gamma_c$-acycliques. En appliquant le foncteur
$\Gamma_c(V,-)$ on obtient un complexe dont la cohomologie est
$$\HC^q=H^{n+2\dim(X)+1+q}_c(V,C_{-n-2\dim(X)})$$ qui est nul pour $q\geq
-n$. Ceci montre l'exactitude du complexe de Cech ``dual'' en degr{\'e} $-n$.

\bigskip

\noindent{\em Cas $\Lambda$ anneau $l$-adique} : Soit $(\FC_n)_n\in
\FC(X/G,\Lambda_\bullet)$. Comme auparavant, tronquons une
r{\'e}solution de $(\FC_n)_n$ par des objets injectifs de
$\FC(X/G,\Lambda_\bullet)$ en degr{\'e} $2\dim(X)+1$. Le complexe  obtenu
$(\JC_n)_n^\bullet$ est form{\'e} de $\Lambda_\bullet$-faisceaux
$\Gamma_c$-acycliques pour la m{\^e}me raison que dans le cas de
torsion. On sait par ailleurs que le foncteur $\limproj :
\FC(X/G,\Lambda_\bullet) \To{} \FC(X/G,\Lambda)$ est aussi de dimension
cohomologique $\leq 2\dim(X) +1$. Une fa{\c c}on de le voir est la
suivante : en raisonnant sur les fibres on voit que
$$ \hbox{dim.coh.}(\limproj) \leq \sup_{U\To{} X/G}
\hbox{dim.coh.}(\Gamma(U,-)\circ \limproj)$$
et pour tout $U\To{} X/G$, on a $R(\Gamma(U,-)\circ
\limproj)=R(\limproj\circ \Gamma(U,-))=R\limproj\circ R\Gamma(U,-)$ o{\`u}
l'on sait que ce dernier $R\limproj : \Lambda_\bullet-\hbox{mod} \To{}
\Lambda-\hbox{mod}$ est de dimension cohomologique $\leq 1$ et que
$R\Gamma(U,-)$ est de dimension cohomologique $2\dim(X)$.

Ainsi les $\Lambda_\bullet$-faisceaux $(\JC_n)_n^i$ sont aussi
$\limproj$-acycliques, et par 
cons{\'e}quent les $\Lambda$-faisceaux $(\limproj \JC_n)^i$ sont
$\Gamma_!(U,-)$-acycliques pour tout morphisme {\'e}tale $U\To{} X$,
puisqu'on a la factorisation $R\Gamma_c(U,-)=R\Gamma_!(U,-)\circ
R\limproj$. On aimerait alors conclure par le m{\^e}me raisonnement que
dans le cas de torsion que les $\Lambda$-faisceaux $(\limproj\JC_n)^i$
sont $\Gamma_{!,et}^\vee$-acycliques (ce qui est {\'e}quivalent {\`a} dire que les
$\Lambda_\bullet$-faisceaux $(\JC_n)_n^i$ sont
$\Gamma_{c,et}^\vee$-acycliques).
Mais il nous manque la finitude cohomologique des foncteurs
$\Gamma_!(V,-)$ sur les cat{\'e}gories de $\Lambda$-faisceaux ($\Lambda$
n'est pas de torsion).

L'auteur ne sait pas si l'assertion du lemme est vraie pour les
rev{\^e}tements $U\To{f} V$ {\em {\'e}tales} quelconques. Mais dans le cas o{\`u}
$f$ est un isomorphisme local, on peut proc{\'e}der de la mani{\`e}re
suivante. Introduisons le site $V_{loc}$ dont les objets sont les
isomorphismes locaux $W\To{} V$, avec la d{\'e}finition habituelle d'un
recouvrement. 
On a des morphisme de sites {\'e}vidents $V_{et} \To{\rho^V} V_{loc} \To{\mu^V}
|V|$ o{\`u} $|V|$ est le site topologique de $V$. Les foncteurs $\mu^V_*$
et ${\mu^V}^*$ induisent des {\'e}quivalences r{\'e}ciproques de topos
$\wt{V_{loc}} \simeq \wt{|V|}$ et en particulier des {\'e}quivalences de
cat{\'e}gories $\FC(|V|,\Lambda)\simeq\FC(V_{loc},\Lambda)$. Il s'ensuit
que le foncteur $\Gamma_!^{loc}(V,-)$ des sections {\`a} support compact
sur $\FC(V_{loc},\Lambda)$ est de 
dimension cohomologique $\leq \dim(V)=\dim(X)$, et ceci {\em quel que
  soit} l'anneau $\Lambda$.

On a une factorisation $\Gamma_c(V,-)=
\Gamma_!^{loc}(V,-)\circ \rho^V_* \circ \limproj$.
Montrons maintenant que le foncteur exact {\`a} gauche $\rho^V_*\circ
\limproj : \FC(V_{et},\Lambda_\bullet)\To{} \FC(V_{loc},\Lambda)$ est
de dimension cohomologique finie : on remarque d'abord
que $\rho^V_*\circ \limproj =\limproj^{loc} \circ \rho^V_*$ o{\`u}
$\limproj^{loc}$ d{\'e}signe le foncteur limite projective
$\FC(V_{loc},\Lambda_\bullet) \To{} \FC(V_{loc},\Lambda)$. Ce foncteur
est de dimension cohomologique $\leq \dim(V)+1=\dim(X)+1$ pour la
m{\^e}me raison que plus haut.
En d{\'e}rivant on obtient
$$ R(\rho^V_*\circ \limproj)=R(\limproj^{loc} \circ \rho^V_*
)=R\limproj^{loc} \circ R\rho^V_*,$$
la seconde {\'e}galit{\'e} venant de ce que $\rho^V_*$ envoie injectifs
sur injectifs puisqu'il a un adjoint {\`a} gauche exact. Maintenant pour
tout faisceau ab{\'e}lien $\FC$, et tout $x\in V$ la fibre de
$R^q\rho^V_*(\FC)$ en $x$ est donn{\'e}e comme dans \cite[4.2.4]{Bic2} par 
$$R^q\rho^V_*(\FC)_x \simeq H^q(\gal(\HC(x)^{ca}/\HC(x)),\FC_x). $$
D'apr{\`e}s \cite[2.5.1]{Bic2} et notre hypoth{\`e}se $K$ alg{\'e}briquement clos,
on en d{\'e}duit que $\rho^V_*$ est de dimension cohomologique $\leq
\dim(X)$ sur la cat{\'e}gorie $\FC(V_{et},\Lambda_\bullet)$. On a donc
obtenu que le foncteur $\rho^V_*\circ \limproj$ 
est de dimension cohomologique $\leq 2\dim(X)+1$. Il s'ensuit que les
$\Lambda_\bullet$-faisceaux $(\JC_n)_n^i$ sont acycliques pour ce
foncteur et par cons{\'e}quent que les $\Lambda$-faisceaux
$\rho^V_*\circ\limproj(\JC_n)^i$ sur $V_{loc}$ sont
$\Gamma^{loc}_!(V,-)$-acycliques !

{\`A} partir de l{\`a}, on prouve par un raisonnement similaire au cas de
torsion que les faisceaux $\rho^V_*\circ\limproj(\JC_n)^i$ sur
$V_{loc}$ sont $\Gamma_{!,loc}^\vee$-acycliques, ce qui est encore
{\'e}quivalent {\`a} dire que les $\Lambda_\bullet$-faisceaux $(\JC_n)_n^i$
sont $\Gamma_{c,loc}^\vee$-acycliques au sens du lemme.

\end{proof}

\begin{coro}
  Soit $\FC \in \FC(X/G,\Lambda)$ avec $\Lambda$ de torsion ou
  $\FC\in\FC(X/G,\Lambda_\bullet)$ avec $\Lambda$ $l$-adique. Il existe une suite
  spectrale  
$$ E_2^{pq}=\tor{p}{\Lambda}{H^q_c(X,p^*\FC)}{G} \Rightarrow
H^{q-p}_c(X/G,\FC). $$
\end{coro}

Remarquons que lorsque $\Lambda$ est un corps ou un anneau Gorenstein,
ou $\FC$ est $l$-adique, on obtient en dualisant la suite spectrale
\cite[4.4.1]{Fargues}.

\section{Espaces modulaires de Drinfeld} \label{drinfeld}

Dans cette section, nous rappelons bri{\`e}vement la d{\'e}finition (d'une
version) des espaces de
modules de groupes formels que nous consid{\'e}rons et la d{\'e}finition usuelle
de leur cohomologie. R{\'e}sumons pour les
sp{\'e}cialistes : on ne fixe pas la hauteur de la quasi-isog{\'e}nie qui
rigidifie le probl{\`e}me de modules (suivant en cela Rapoport-Zink) et on
{\em ne fixe pas} de caract{\`e}re central, du moins pour l'{\'e}nonc{\'e} du
th{\'e}or{\`e}me principal.
Enfin on d{\'e}finit le complexe de cohomologie $R\Gamma_c$ qui joue
{\'e}videmment un r{\^o}le central dans ce texte.

\def\od{\OC_{D_d}}

On notera $\OC$ l'anneau des entiers de $K$, $\od$ celui de $D_d$ 
et $\wh\OC^{nr}$ celui de la compl{\'e}tion de la
sous-extension maximale non ramifi{\'e}e $K^{nr}$ de $K$ dans $K^{ca}$.
Le corps r{\'e}siduel de $K$ est $k$ et celui de $K^{nr}$ est une cl{\^o}ture
alg{\'e}brique $k^{ca}$ de $k$. Enfin $\varpi$ d{\'e}signe une uniformisante de
$\OC$ et $\wh\OC^{nr}$.

Les espaces qui nous int{\'e}ressent sont des espaces de modules de
$\OC$-modules formels. Si $B$ est une $\OC$-alg{\`e}bre, un $\OC$-module
formel sur $B$ est un groupe formel muni d'une action de $\OC$
relevant l'action naturelle sur l'alg{\`e}bre de Lie. Un $\od$-module
formel sur $B$ est un $\OC$-module formel muni d'une action de $\od$
qui {\'e}tend celle de $\OC$. Notons $\OC_d \subset \od$ l'anneau des
entiers d'une sous-extension non ramifi{\'e}e maximale de $K$ dans
l'alg{\`e}bre {\`a} division. Suivant Drinfeld, un $\od$-module formel
sur $B$ est dit sp{\'e}cial si son alg{\`e}bre de Lie est localement libre de rang
$1$ comme $\OC_d\otimes_\OC B$-module.

Une isog{\'e}nie de $\OC$-modules formels est une isog{\'e}nie des groupes
formels sous-jacents compatible aux actions de $\OC$. Sa
($\OC$-)hauteur est le quotient (entier) de la hauteur au sens des
groupes par le degr{\'e} de $k$ sur $\FM_p$. La hauteur d'un
$\OC$-module formel est la hauteur de l'isog{\'e}nie $X(\varpi)$.
Une quasi-isog{\'e}nie $X\To{} Y$ de $\OC$-modules formels est un
{\'e}l{\'e}ment de $\hom{X}{Y}{\OC-mf}\otimes_\OC K$ qui admet un inverse
dans $\hom{Y}{X}{\OC-mf}\otimes_\OC K$. On montre que c'est en fait une isog{\'e}nie {\`a}
multiplication par une puissance de $X(\varpi)$ pr{\`e}s, ce qui permet d'en
d{\'e}finir la hauteur.

\subsection{La tour de Drinfeld} \label{tourDr}

On sait qu'il existe sur $k^{ca}$ un $\od$-module formel $\XM$ sp{\'e}cial (donc de
dimension $d$) de hauteur $d^2$  et qu'il est unique {\`a} isog{\'e}nie
pr{\`e}s, \cite[Prop II.5.2]{BouCar}, \cite[3.60]{RZ}.  

\alin{Probl{\`e}me de modules}
\def\nilp{\hbox{Nilp}}
Soit $\nilp$ la cat{\'e}gorie des $\wh\OC^{nr}$-alg{\`e}bres o{\`u} l'image
de $\varpi$ est nilpotente.
 On consid{\`e}re le foncteur $\wt{G} : \nilp
\To{} \hbox {Ens}$ qui {\`a} $B$ associe l'ensemble des classes
d'isomorphisme de 
couples $(X,\rho)$ avec $X$
un $\od$-module formel sur $B$ et $\rho : \XM \otimes_{k^{ca}} (B/\varpi B)
\To{} X \otimes_B (B/\varpi B)$ une quasi-isog{\'e}nie.
On
a une d{\'e}composition {\'e}vidente $\wt{G} = \bigsqcup_{h\in \ZM} G^{(h)}$ o{\`u}
$G^{(h)}$ classifie les classes de couples $(X,\rho)$ avec $\rho$ de
hauteur $dh$. Chaque $G^{(h)}$ est non-canoniquement isomorphe {\`a}
$G^{(0)}$ et $G^{(0)}$ est le foncteur $G$ de Drinfeld
\cite{Drincov}. On sait que $G$ est repr{\'e}sentable ({\em cf}
\cite{Genestier}, \cite{RZ} pour une preuve {\em a priori}) par un
sch{\'e}ma formel localement de type fini sur $\wh\OC^{nr}$. On note
$\wh{\MC}_{Dr,0}^{(0)}$ ce sch{\'e}ma formel. De m{\^e}me on note $\wh{\MC}_{Dr,0}^{(h)}$ le
sch{\'e}ma formel repr{\'e}sentant $G^{(h)}$ et $\wh{\MC}_{Dr,0}$ celui qui repr{\'e}sente
$\wt{G}$. On a donc non-canoniquement $\wh{\MC}_{Dr,0}\simeq
\wh{\MC}_{Dr,0}^{(0)}\times \ZM$. 
Enfin on note sans $\;\wh{}\;$ les fibres g{\'e}n{\'e}riques au sens de
Raynaud-Berkovich de ces espaces : ce sont donc des $\knr$-espaces
analytiques au sens de \cite{Bic2}.



\alin{Structures de niveau}
Notons $({X}_u,\rho_u)$ l'objet universel au-dessus de $\wh\MC_{Dr,0}$. Le noyau
${X}_u[\varpi^n]$ de
la multiplication par $\varpi^n$ dans ${X}_u$ est un sch{\'e}ma formel en groupes
plat  fini de rang  $p^{nd^2}$ au-dessus de $\wh\MC_{Dr,0}$, et qui
est {\'e}tale en fibre g{\'e}n{\'e}rique. Plus 
pr{\'e}cis{\'e}ment, sa fibre g{\'e}n{\'e}rique est
localement pour la topologie {\'e}tale isomorphe {\`a} $\mdro \times
\od/\varpi^{n}\od$. Le
$(\od/\varpi^n\od)^\times$-torseur sur $\mdro$
$$ \underline{\hbox{Isom}_{\od}}\left((\varpi^{-n}\od/\od)_{\mdro},{X}_u[\varpi^n]\right)$$
est  donc repr{\'e}sent{\'e} par un rev{\^e}tement {\'e}tale de $\mdro$, galoisien  de
groupe $\od^\times /(1+\varpi^n\od)$, et que 
nous noterons $\mdrn$\footnote{Cette d{\'e}finition des rev{\^e}tements
  est emprunt{\'e}e {\`a} \cite{Genestier}}. Pour
$m\leq n$, l'inclusion $\varpi^{-m}\od/\od \subset \varpi^{-n}\od/\od$
induit un
morphisme $\mdrn\To{} \MC_{Dr,m}$. On obtient ainsi une ``tour'' de
rev{\^e}tements Galoisiens de $\mdro$ dont le pro-groupe de Galois est
$\od^\times$.

\alin{Actions des groupes}  \label{actionDr}
On sait  que le groupe des
quasi-isog{\'e}nies du $\od$-module formel 
$\XM$ s'identifie {\`a} $G_d=GL_d(K)$, \cite[3.60]{RZ} ou
\cite[II.5.2]{BouCar}. On r{\'e}cup{\`e}re donc une action ({\`a} gauche) de $G_d$ 
sur le foncteur $\wt{G}$ qui envoie le couple $(X,\rho)$ sur le couple
$(X,\rho \circ ({^tg})_\XM)$ o{\`u} $g_\XM$ d{\'e}signe la
quasi-isog{\'e}nie de $\XM$ associ{\'e}e {\`a} $g$ {\em et ${^tg}$ d\'esigne la transpos\'ee de $g$}\footnote{Nous nous \'ecartons ici des conventions usuelles. Il y a trois raisons \`a cela : d'une part, l'\'enonc\'e final (\ref{conj} i)) est plus joli car exempt de contragr\'ediente, d'autre part, cela rendra le morphisme des p\'eriodes de Drinfeld compatible avec l'action sur $\Omega$ de la partie \ref{dp}, et enfin, le th\'eor\`eme de Faltings fait intervenir un passage \`a la transpos\'ee.}. Cette action en induit
une sur $\wh{\MC}_{Dr,0}$ par Yoneda, et l'objet universel est muni
d'isomorphismes canoniques 
$({X}_u,\rho_u\circ g_\XM^{-1}) \simto g^*({X}_u,\rho_u) $. En particulier ${X}_u$ est
$G_d$-{\'e}quivariant au-dessus de $\wh\MC_{Dr,0}$.
Par cons{\'e}quent tous les $\mdrn$ sont
munis d'une action de $G_d$ (explicitement d{\'e}finie par
$$ X_u[\varpi^n] \simto g^*X_u[\varpi^n] = X_u[\varpi^n]\times_{\wh\MC_{Dr,0},
  g} \wh\MC_{Dr,0} \To{proj_1} X_u[\varpi^n] )$$
et les morphismes de transition sont
$G_d$-{\'e}quivariants. D'apr{\`e}s Berkovich, l'action de $G_d$ sur ces $\knr$-espaces
analytiques est continue au sens de \cite{Bic3} pour la topologie
naturelle de $G_d$.

\medskip

On d{\'e}finit maintenant l'action de $\dd$ sur $\MC_{Dr,0}$. Pour $d\in
\dd$ et $X$ un $\od$-module formel sur $B$, notons ${^dX}$ le
$\od$-module formel dont le $\OC$-module formel sous-jacent est
encore $X$ mais dont l'action de $\od$ est donn{\'e}e par
${^dX}(x):=X(d^{-1}xd)$, de sorte que $X(d^{-1})$ est une quasi-isog{\'e}nie
$X\To{} {^dX}$. Pour un couple $(X,\rho)$ on pose
$d(X,\rho):=({^dX}, \rho \circ \XM(d^{-1}))$. On obtient ainsi
une action {\`a} gauche de $\dd$  sur le foncteur
$\wt{G}$ qui en induit une 
sur $\wh\MC_{Dr,0}$ par Yoneda. L'objet universel est alors muni d'isomorphismes
canoniques $d^*(X_u,\rho_u) \simto (X_u^d,\rho_u\circ \XM(d^{-1}))$. 
On d{\'e}finit  la structure $\dd$-{\'e}quivariante  sur
$\mdrn$ au-dessus de $\mdro$ par la suite d'isomorphismes
\begin{eqnarray*}
    d^*  \underline{\hbox{Isom}}_{\od}\left(\varpi^{-n}\od/\od,
X_u[\varpi^n]\right) & \simeq &   \underline{\hbox{Isom}}_{\od}\left(\varpi^{-n}\od/\od,
d^*X_u[\varpi^n]\right) \\
& \simeq &  \underline{\hbox{Isom}}_{\od}\left(\varpi^{-n}\od/\od,
X_u^d[\varpi^n]\right) \\
& \simto &  \underline{\hbox{Isom}}_{\od}\left(\varpi^{-n}\od/\od,
X_u[\varpi^n]\right)
\end{eqnarray*}
o{\`u} les deux premiers isomorphismes sont canoniques et le troisi{\`e}me est
donn{\'e} par $\phi \mapsto \phi\circ\,\hbox{ad}(d^{-1})$, en notant
$\hbox{ad}(d^{-1})$ l'application $a\in \varpi^{-n}\od/\od \mapsto
d^{-1}ad \in \varpi^{-n}\od/\od$.   
Comme pr{\'e}c{\'e}demment pour l'action de $G_d$, ces structures
$\dd$-{\'e}quivariantes induisent (et 
{\'e}quivalent {\`a}) des actions de $\dd$ sur les $\mdrn$ compatibles avec les
morphismes de transition.
Notons que si
$d\in \od^\times$, l'action de $d$ est triviale sur $\mdro$ et co{\"\i}ncide
sur $\mdrn$ avec l'action naturelle sur le torseur qui d{\'e}finit
$\mdrn$. En particulier l'action de $\dd$ se factorise sur chaque $\mdrn$
par un quotient discret.

\medskip

Passons {\`a} l'action de Galois : fixons un g{\'e}n{\'e}rateur topologique $\phi$
de $\gal(K^{nr}/K)=\gal(k^{ca}/k)$, on peut munir les $\knr$-espaces
analytiques $\mdrn$ d{\'e}finis ci-dessus de structures
$\phi^\ZM$-{\'e}quivariantes au-dessus de $\knr$\footnote{Remarquons
  que la donn{\'e}e d'une structure $\phi^\ZM$-{\'e}quivariante sur un $\knr$-espace
  analytique  {\'e}quivaut {\`a} celle d'une action de $\phi^\ZM$ sur
  l'``espace analytique au-dessus de $K$'' sous-jacent ({\em cf}
    \cite{Bic2} pour la terminologie), compatible {\`a} l'action de
    $\phi^\ZM$ sur $\MC(\knr)$} 
(appel{\'e}es ``donn{\'e}es de descente {\`a} la Weil'' dans
\cite[3.45]{RZ}). Commen{\c c}ons par $\wh\MC_{Dr,0}$ : si  $(X,\rho)$ est un
$\od$-module formel rigidifi{\'e}
au-dessus de $B$, alors on en d{\'e}finit un autre $(\phi^*X,{^\phi\rho})$
au-dessus de $\phi^*B :=B\otimes_{\OC^{nr},\phi} \OC^{nr}$ en posant
$\phi^*X:=X\otimes_{B} \phi^*B$ et 
$$^\phi\rho:\; \XM \otimes_{k^{ca}} \phi^*(B/\varpi) \To{Frob\otimes \id}
\phi^*\XM \otimes_{k^{ca}} \phi^*(B/\varpi) \To{\phi^*\rho} \phi^*X
\otimes_{\phi^* B} \phi^* (B/\varpi).$$
On obtient ainsi un morphisme de foncteurs $\wt{G} \To{} \phi_*\wt{G}$
dont on voit imm{\'e}diatement qu'il est inversible, ce qui nous donne une
structure $\phi^\ZM$-{\'e}quivariante sur $\wt{G}$, et donc sur $\wh\MC_{Dr,0}$, au-dessus 
de $\hbox{Spf}(\wh\OC^{nr})$. Nommons momentan{\'e}ment $\tau_\phi:
\phi^*\wh\MC_{Dr,0}\simto \wh\MC_{Dr,0}$ le $\hbox{Spf}(\wh\OC^{nr})$-isomorphisme 
associ{\'e}  {\`a} cette $\phi^\ZM$-structure {\'e}quivariante (adjoint du pr{\'e}c{\'e}dent), on obtient pour
l'objet universel l'existence d'un isomorphisme canonique
$\sigma_\phi:\;\;(\phi^*(X_u),{^\phi\rho_u})\simto
\tau_\phi^*(X_u,\rho_u) $  au-dessus de $\phi^*\wh\MC_{Dr,0}$, d'o{\`u}, par
composition, une structure $\phi^\ZM$-{\'e}quivariante
$$ \phi^*(X_u) \To{\sigma_\phi} \tau_\phi^*X_u =
X_u\times_{\wh\MC_{Dr,0}} \phi^*\wh\MC_{Dr,0} \To{proj_1} X_u $$
au-dessus de $\hbox{Spf}(\wh\OC^{nr})$.
En passant aux fibres g{\'e}n{\'e}riques, on obtient sur les $\mdrn$ les structures
$\phi^\ZM$-{\'e}quivariantes au-dessus de $\MC(\knr)$ (notation de
\cite{Bic2}) cherch{\'e}es. Par construction, elles sont compatibles aux
morphismes de transition. Remarquons que l'action de $\phi$ sur
$\mdrn$ envoie $\mdrn^{(h)}$ sur $\mdrn^{(h+1)}$.

{\'E}tendons maintenant les scalaires {\`a} $\ka$, en suivant la d\'efinition de \cite[1.4]{Bic2}.
 On obtient le $\ka$-espace analytique 
 $$\mdrn^{d/K}:=\mdrn
\wh\otimes_{\knr} \ka$$ 
de l'introduction. La structure
$I_K$-{\'e}quivariante sur $\mdrn^{d/K}$ au-dessus de $\ka$ induite par ce changement
de base et la structure $\phi^\ZM$-{\'e}quivariante pr{\'e}c{\'e}dente se
``recollent''  en une structure $W_K$-{\'e}quivariante au-dessus de
$\MC(\ka)$. De mani{\`e}re {\'e}quivalente, on obtient une action de
$W_K$ {\em \`a droite} sur l'``espace analytique   
au-dessus de $K$'' $\mdrn^{d/K}$, compatible {\`a} l'action naturelle de
$W_K$ sur $\ka$.  

\medskip 

Les actions de $G_d, \dd$ et $W_K$ qu'on vient de d{\'e}finir sur les espaces
analytiques $\mdrn^{d/K}$ commutent. On v{\'e}rifie facilement que le
stabilisateur de $\MC^{d/K,(h)}_{Dr,n}$ 
est $(G_d\times \dd\times W_K)^0$, avec la notation introduite en \ref{defun}.
Remarquons pour terminer que l'action de $G_d\times \dd\times W_K$ est triviale sur
le sous-groupe $\Delta$ de $G_d\times \dd$ 
d\'efini par 
\ini
\begin{equation} \label{defDelta}
\Delta=\hbox{Im}\left(\application{}{K^\times}{G_d\times\dd}{z}{(z,z)}\right)
\end{equation}


\alin{Le morphisme de p{\'e}riodes} \label{perDr}
Soit $\Omega^{d-1}_K$ 
l'espace sym{\'e}trique
de Drinfeld de la partie \ref{dp}.
Nous noterons $\Omega^{d-1,nr}_K := \Omega^{d-1}_K\wh\otimes_K \knr$ son
changement de base {\`a} $\knr$. Drinfeld a construit dans
\cite{Drincov} un isomorphisme $\xi^{(0)}:\;\;\wh\MC_{Dr,0}^{(0)}\To{}
\wh{\Omega}^{d-1,nr}_K$ o{\`u} $\wh{\Omega}^{d-1,nr}_K$ d{\'e}signe un certain
mod{\`e}le formel de ${\Omega}^{d-1,nr}_K$ d{\'e}fini auparavant par
Deligne.
Nous n'avons besoin ici que de la fibre g{\'e}n{\'e}rique de ce
morphisme $\xi$, mais nous avons par contre besoin  de l'{\'e}tendre {\`a} tout
$\MC_{Dr,0}$. Pour cela nous utilisons l'isomorphisme de $\knr$-espaces
analytiques $\tau_{\phi^h}:\;\;(\phi^h)^*\mdro \To{} \mdro$ qui induit
un isomorphime $(\phi^h)^*\MC_{Dr,0}^{(0)} \simto \MC_{Dr,0}^{(h)}$
et nous d{\'e}finissons $\xi^{(h)}$ par le diagramme commutatif
$$ \xymatrix{ (\phi^h)^* \mdro^{(0)} \ar[d]_{(\phi^h)^*\xi^{(0)}} \ar[r]^\sim & \mdro^{(h)}
  \ar@{-->}[d]^{\xi^{(h)}} \\ (\phi^h)^*\Omega^{d-1,nr}_K \ar[r]^\sim
  & \Omega^{d-1,nr}_K }$$
o{\`u} la fl{\`e}che du bas est la donn{\'e} de descente naturelle sur
$\Omega^{d-1,nr}_K$. Enfin on pose
$\xi := \bigsqcup_{h\in\ZM} \xi^{(h)} : \MC_{Dr,0} \To{} \Omega^{d-1,nr}_K$.

Une construction directe (et plus ``simple'' que celle de Drinfeld) de
ce $\xi$ est propos{\'e}e dans le 
chapitre 5 de Rapoport-Zink \cite{RZ}, mais n'est r{\'e}dig{\'e}e qu'en
in{\'e}gales caract{\'e}ristiques (techniques cristallines).

\begin{fact} Le morphisme
$$ \xi:\;\; \mdro \To{} \Omega^{d-1,nr}_K $$  est
$G_d\times \dd\times W_K$ {\'e}quivariant pour l'action naturelle de $G_d$ sur
$\Omega^{d-1,nr}_K$ (qui se factorise donc par $PGL_d(K)$), l'action
triviale de $\dd$ et l'action naturelle de $W_K$ (celle donn{\'e}e par
l'extension des scalaires).
\end{fact}

\begin{proof}
Tout est contenu dans les r\'ef\'erences usuelles \cite{Drincov}, \cite{BouCar}, \cite{Genestier} et \cite{RZ}. Commentons seulement la compatibilit\'e \`a l'action de $G_d$ : dans les r\'ef\'erences pr\'ec\'edentes  l'action de $G_d$ sur les $\MC_{Dr,n}$ est l'action "naturelle" dont la notre se d\'eduit par $g\mapsto {^tg}^{-1}$. Le morphisme de p\'eriode y est \'equivariant pour l'action naturelle sur le mod\`ele "dual" de $\Omega^{d-1}_K$, \`a savoir l'ensemble des hyperplans de $K^d$ ne contenant pas de droites rationnelles. C'est pourquoi, avec notre normalisation diff\'erente, le morphisme de p\'eriodes reste $G_d$-\'equivariant \`a condition d'utiliser le mod\`ele de la partie \ref{dp} pour $\Omega_K^{d-1}$.

\end{proof}

\subsection{La tour de Lubin-Tate} \label{tourLT}
Soit $\XM$ un $\OC_K$-groupe formel sur $\o\FM_p$ de dimension $1$ et
hauteur $d$, au sens de \cite{Drinell} ou \cite{HG}. On sait qu'un tel
objet existe et est unique {\`a} isomorphisme pr{\`e}s, \cite[prop
1.6-7]{Drinell}.  
\alin{Modules des d{\'e}formations de $\XM$} On peut comme dans le livre
de Rapoport-Zink exprimer le probl{\`e}me de modules sur la cat{\'e}gorie
$\nilp$ et de mani{\`e}re semblable au cas de la tour de Drinfeld. Nous
donnons malgr{\'e} tout la d{\'e}finition historique (et la seule {\'e}crite en
{\'e}gales caract{\'e}ristiques) en termes de d{\'e}formations.
 Soit $\CC$ la
cat{\'e}gorie des $\wh\OC^{nr}$-alg{\`e}bres locales de corps r{\'e}siduel
$\o\FM_p$ et compl{\`e}tes pour leur 
topologie adique. Une {\em d{\'e}formation par quasi-isog{\'e}nie} de $\XM$ sur une
telle alg{\`e}bre $R$ est une paire $(X,\rho)$ form{\'e}e d'un $\OC_K$-module
formel $X$ sur $R$ et d'une quasi-isog{\'e}nie $\XM\To{\rho}X\otimes_R
R/\MG_R$.  Notons $\hbox{Def}$ le foncteur de   
$\CC$ dans les ensembles qui 
{\`a} $R$ associe l'ensemble des classes d'isomorphisme de d{\'e}formations
par quasi-isog{\'e}nie $(X,\rho)$ de $\XM$ sur $R$. Ce foncteur est une
r{\'e}union disjointe de sous-foncteurs $\hbox{Def}^{(h)}$ classifiant les
couples $(X,\rho)$ avec $\rho$ de hauteur $h$. Tous les $\hbox{Def}^{(h)}$
sont non-canoniquement isomorphes {\`a} $\hbox{Def}^{(0)}$. De plus, comme une
quasi-isog{\'e}nie de hauteur nulle entre deux $\OC_K$-modules formels de
dimension $1$ sur $\o\FM_p$ est un isomorphisme (c'est une
cons{\'e}quence de \cite{Drinell},props 1.6-2 et 1.7), on voit que
$\hbox{Def}^{(0)}$ est le foncteur de d{\'e}formations {\'e}tudi{\'e}
par Drinfeld.  
D'apr{\`e}s \cite[4.2]{Drinell} on sait que $\hbox{Def}^{(0)}$ est
repr{\'e}sentable 
par l'alg{\`e}bre $R^{(0)}:=\wh\OC^{nr}[[T_1,\cdots,T_{d-1}]]$  des s{\'e}ries
formelles {\`a} $d-1$ variables sur $\wh\OC^{nr}$. Nous noterons alors
$\MC_{LT,0}^{(0)}$ le $\knr$-espace analytique de Berkovich
associ{\'e} au sch{\'e}ma 
formel $\wh\MC_{LT,0}^{(0)}:=\hbox{Spf}(R^{(0)})$ (suivant la proc{\'e}dure de
Raynaud-Berthelot d{\'e}crite dans \cite{Bic4}) : c'est la boule
unit{\'e} ouverte de dimension $d-1$. 
Il s'ensuit que les foncteurs $\hbox{Def}^{(h)}$ pour $h\in \NM$ sont
aussi repr{\'e}sentables et nous noterons $\MC_{LT,0}^{(h)}$ les espaces
analytiques associ{\'e}s. Nous posons enfin
$\wh\MC_{LT,0}:=\bigsqcup_{h\in\NM} \wh\MC_{LT,0}^{(h)}$ et notons sans
$\;\wh{}\; $ l'espace analytique associ{\'e}.

\alin{Structures de niveau} On peut suivre  la m{\^e}me proc{\'e}dure que dans
le cas $Dr$. En notant encore $(X_u,\rho_u)$ l'objet universel sur
$\wh\MC_{LT,0}$, on d{\'e}finit $\mltn$ comme le rev{\^e}tement {\'e}tale
galoisien de
groupe $GL_d(\OC/\varpi^n\OC)$ de $\mlto$ repr{\'e}sentant le torseur
$$ \underline{\hbox{Isom}}((\varpi^{-n}\OC/\OC)^d,X_u[\varpi^n]).$$
Il se trouve que dans ce cas $LT$ on peut faire mieux en interpr{\'e}tant 
$\MC_{LT,n}$ comme la fibre g{\'e}n{\'e}rique d'un probl{\`e}me de modules
classifiant les ``structures de niveau de Drinfeld'' : cela  a
l'avantage de simplifier la description de l'action de
$G_d$\footnote{On peut n{\'e}anmoins se passer d'une telle
  interpr{\'e}tation modulaire comme il est esquiss{\'e}
  dans \cite[5.34]{RZ} et expliqu{\'e} dans \cite[2.3.8.3]{Fargues}}.

Soit  $\Lambda \subset K^d$ un $\OC$-r{\'e}seau.
Si $R\in \CC$ et $(X,\rho)$ est une d{\'e}formation par quasi-isog{\'e}nie
de $\XM$ au-dessus
de $R$, alors une $\Lambda$-structure de niveau $n$ de Drinfeld sur $X$ est un
morphisme de $\OC$-modules $\psi$ de  $\varpi^{-n}\Lambda/\Lambda$ vers l'id{\'e}al
maximal $\MG_R$ muni de la structure de $\OC$-module d{\'e}finie par $X$
et tel que
$$ \prod_{x\in \varpi^{-1}\Lambda/\Lambda} (T-\psi(x)) \;\;\hbox{ divise }\;\;
  X_\varpi(T) $$
o{\`u} $X_\varpi(T)$ d{\'e}signe la s{\'e}rie formelle donnant l'action de
$\varpi$ sur $X$.
Les foncteurs $\hbox{Def}^{(h)}_{\Lambda,n}$ classifiant les triplets $(X,\rho,\psi)$ {\`a}
isomorphisme pr{\`e}s sont repr{\'e}sentables \cite[4.3]{Drinell} par des
alg{\`e}bres $R^{(h)}_{\Lambda,n}$ finies et 
plates sur les $R^{(h)}$ et nous noterons $\wh\MC_{LT,\Lambda,n}:=\bigsqcup_{h\in
  H} \hbox{Spf} R^{(h)}_{\Lambda,n}$. Lorsque $\Lambda$ est le
r{\'e}seau ``canonique'' $\OC^d\subset K^d$, la fibre g{\'e}n{\'e}rique de
$\wh\MC_{LT,\Lambda,n}$ s'identifie canoniquement {\`a} $\MC_{LT,n}$. 
{\'E}videmment pour $n$ fix{\'e}, tous les $\wh\MC_{LT,\Lambda,n}$ sont
(non-canoniquement) isomorphes.



\alin{Action des groupes} \label{actionLT}
On sait que le groupe des quasi-isog{\'e}nies du
$\OC$-module formel $\XM$ s'identifie {\`a} $\dd$. Ceci nous permet de
d{\'e}finir une action \`a gauche de $\dd$ sur $\mlto$, puis sur les
$\mltn$
d'une mani{\`e}re exactement analogue {\`a}  celle par laquelle on a d{\'e}fini
l'action de $G_d$ sur $\mdro$ et les $\mdrn$. Explicitement, on envoie
un triplet $(X,\rho,\psi)$ sur le triplet $(X, \rho\circ d_{\XM}^{-1},\psi)$.

\medskip

L'action de $G_d$ est plus d{\'e}licate. Elle est soigneusement d{\'e}finie
dans \cite[II.2]{HaTay} 
(o{\`u} les d{\'e}formations sont cependant par isomorphismes) et
\cite{Strauch} (qui traite les d{\'e}formations par
quasi-isog{\'e}nies). On rappelle ici bri{\`e}vement cette
d{\'e}finition, en l'exposant un peu diff\'eremment des r\'ef\'erences
usuelles. Pour deux r\'eseaux $\Lambda, \Lambda'$ dans $K^d$, nous
noterons $d(\Lambda,\Lambda')$ la distance combinatoire entre les
points de l'immeuble de $PGL_d(K)$ associ\'es. C'est le minimum de la
somme $r+k$ pour tous les couples d'entiers
 $r,k\in\ZM$  tels
tels que 
\ini\begin{equation}\label{condreseaux}
\Lambda\subseteq
\varpi^{-r} \Lambda'\subseteq \varpi^{-r-k}\Lambda.
\end{equation} 
\begin{lemme} \label{blob} Il existe une famille de morphismes
  $\dd$-\'equivariants de foncteurs
$$
  \alpha_{\Lambda'|\Lambda}^{n'|n}:\;\;\hbox{Def}_{\Lambda,n} 
\To{} \hbox{Def}_{\Lambda',n'}
$$
indic\'ee par les quadruplets $(\Lambda,\Lambda',n,n')$ o\`u
  $\Lambda, \Lambda'$ sont des r\'eseaux de $K^d$ et $n,n'\in\NM$ sont
  tels que $n-n'\geq d(\Lambda,\Lambda')$, et telle que 
\begin{enumerate}
\item si $\Lambda=\Lambda'$, alors $
  \alpha_{\Lambda'|\Lambda}^{n'|n}$ est le morphisme de restriction de
  la structure de niveau.
\item pour tout autre $\Lambda'', n''$ avec $n'-n''\geq
    d(\Lambda',\Lambda'')$, on a $\alpha_{\Lambda''|\Lambda}^{n''|n}=
    \alpha_{\Lambda''|\Lambda'}^{n''|n'} \circ
    \alpha_{\Lambda'|\Lambda}^{n'|n}$.
\item les fibres g\'en\'eriques des morphismes de sch\'emas formels
  $\wh\MC_{LT,\Lambda,n} \To{}\wh\MC_{LT,\Lambda',n'}$ associ\'es sont
  \'etales, finies et surjectives (au sens de Berkovich, par exemple).
\end{enumerate}
\end{lemme}

\begin{proof}
Remarquons que dans l'\'egalit\'e ii), le quadruplet
$(\Lambda,\Lambda'',n,n'')$ v\'erifie bien $n-n''\geq
d(\Lambda,\Lambda'')$ par l'in\'egalit\'e triangulaire.
 
Fixons maintenant $(\Lambda,\Lambda',n,n')$ et choisissons $(k,r)$
v\'erifiant \ref{condreseaux} et tels que $n-n' \geq k+r$.
On d\'efinit un morphisme $(\alpha_{\Lambda'|\Lambda}^{n'|n})_{k,r}$
comme dans l'\'enonc\'e du lemme
en associant {\`a} $(X,\rho,\psi)$ au-dessus de $R\in\CC$ le triplet
$(X'_{k,r},\rho'_{k,r},\psi'_{k,r})$ au-dessus de $R$ d{\'e}fini par
\begin{itemize}
\item $X'_{k,r}:=X/\psi(\varpi^{-r}\Lambda'/\Lambda)$, dont l'existence (et
  la d{\'e}finition pr{\'e}cise) est assur{\'e}e par \cite[Lemma
  4.4]{Drinell}, ou \cite[Lemma II.2.4]{HaTay}.
\item $\rho'_{k,r} :\XM \To{\varpi^{-r}} \XM \To{\rho} X\otimes_R R/\MG_R
  \To{can} X'\otimes_R R/\MG_R $.
\item $\psi'_{k,r} :\;\; \varpi^{-n'} \Lambda'/\Lambda' \To{\times
    \varpi^{-r}} \varpi^{-r-n'} \Lambda'/\varpi^{-r}\Lambda' \injo
  \xymatrix{ \varpi^{-n}\Lambda/\varpi^{-r}\Lambda' \ar@{-->}[r] &
    (\MG_R,X') \\ \varpi^{-n}\Lambda/\Lambda \ar@{->>}[u]
    \ar[r]^{\psi} & (\MG_R,X) \ar[u]^{can}}$
\end{itemize}
Il est clair que ce morphisme est $\dd$-\'equivariant.
Il nous faut voir qu'il ne d\'epend pas du couple $(k,r)$ tel que
$n-n'\geq k+r$. Soit $(k',r')$ un autre tel couple, et supposons que
$r' \geq r$.
Dans ce cas la multiplication par
  $\varpi^{r'-r}$ dans le $\OC_K$-module formel
  $X/\psi(\varpi^{-r}\Lambda'/\Lambda)$ induit un isomorphisme
  $X/\psi(\varpi^{-r'}\Lambda'/\Lambda) \simto
  X/\psi(\varpi^{-r}\Lambda'/\Lambda)$. On v\'erifie alors ais\'ement
  que cet isomorphisme induit un
  isomorphisme de triplets $(X_{k',r'}',\rho'_{k',r'},\psi'_{k',r'})
  \simto (X_{k,r}',\rho'_{k,r},\psi'_{k,r})$.

La propri\'et\'e i) est imm\'ediate. Pour la propri\'et\'e ii),
choisissons $(k,r)$ comme ci-dessus et $(k',r')$ tel que $n'-n''\geq
k'+r'$ et $\Lambda'\subseteq \varpi^{-r'} \Lambda''\subseteq \varpi^{-r'-k'}\Lambda'$. Alors
soit $(k'',r''):=(k'+k,r'+r)$, on a bien $n-n''\geq k''+r''$ et
$\Lambda\subseteq \varpi^{-r''} \Lambda''\subseteq
\varpi^{-r''-k''}\Lambda$ et on v\'erifie sur la d\'efinition que 
$(\alpha_{\Lambda''|\Lambda}^{n''|n})_{k'',r''}=
    (\alpha_{\Lambda''|\Lambda'}^{n''|n'})_{k',r'} \circ
    (\alpha_{\Lambda'|\Lambda}^{n'|n})_{k,r} $.

Il reste \`a montrer l'\'etale finitude et la surjectivité des morphismes induits sur les
fibres g\'en\'eriques. Pour cela, on utilise le fait que lorsque on a
deux morphismes analytiques  tels que $g\circ f$ est \'etale, alors
($g$  \'etale) $\Rightarrow$ ($f$ étale), et ($f$ étale surjectif) $\Rightarrow$  ($g$ étale). Rappelons maintenant que les
morphismes de  restriction de la structure de niveau sont \'etales surjectifs (et
finis) en fibre g\'en\'erique. Ainsi, \'etant donn\'e un quadruplet
$(\Lambda,\Lambda',n,n')$ tel que $n-n'\geq d(\Lambda,\Lambda')=:\delta$,
l'\'egalit\'e $\alpha_{\Lambda'|\Lambda'}^{0|n'}\alpha_{\Lambda'|\Lambda}^{n'|n}=
\alpha_{\Lambda'|\Lambda}^{0|\delta}\alpha_{\Lambda|\Lambda}^{\delta|n}$
montre que $\alpha_{\Lambda'|\Lambda}^{n'|n}$ est \'etale  en
fibre g\'en\'erique \ssi\ $\alpha_{\Lambda'|\Lambda}^{0|\delta}$
l'est, et donc \ssi\ $\alpha_{\Lambda'|\Lambda}^{2\delta|3\delta}$
l'est !
Ce dernier est le premier morphisme de la suite
$$ \hbox{Def}_{\Lambda,3\delta} \To{}\hbox{Def}_{\Lambda',2\delta}
\To{}\hbox{Def}_{\Lambda,\delta} \To{} 
\hbox{Def}_{\Lambda',0}. $$
En fibre g\'en\'erique, le caract\`ere \'etale de la compos\'ee des deux
derni\`eres fl\`eches
$\alpha_{\Lambda'|\Lambda}^{0|\delta} \circ
\alpha_{\Lambda|\Lambda'}^{\delta|2\delta} =
\alpha_{\Lambda'|\Lambda'}^{0|2\delta}$ implique que celle du milieu,
$\alpha_{\Lambda|\Lambda'}^{\delta|2\delta}$, est non-ramifi\'ee. Mais
alors, toujours en fibre g\'en\'erique, le caract\`ere \'etale et fini de 
 la compos\'ee des deux premi\`eres fl\`eches implique par \cite[3.2.9]{Bic2} que
$\alpha_{\Lambda'|\Lambda}^{2\delta|3\delta}$ est \'etale et finie.
La surjectivité se voit de la m\^eme mani\`ere.

\end{proof}

{\em Remarque :} Les $\alpha_{\Lambda'|\Lambda}^{n'|n}$ {\em ne commutent
  pas}, en g\'en\'eral, avec les morphismes d'oubli $\hbox{Def}_{\Lambda,n}\To{}
\hbox{Def}$.

\medskip

Pour d{\'e}finir l'action de $G_d$ sur le syst{\`e}me projectif des
$\wh\MC_{LT,n}$, remarquons que tout {\'e}l{\'e}ment $g\in G_d$ induit un
isomorphisme {\'e}vident $\hbox{Def}_{\Lambda,n} \simto
\hbox{Def}_{g\Lambda,n}$ en envoyant le triplet $(X,\rho,\psi)$ sur le
triplet $(X,\rho,\psi\circ g^{-1})$. Nous posons alors pour tous $n,n' \in \NM$
``ad{\'e}quates'', c'est-\`a-dire tels que $n-n'\geq d(g\OC^d,\OC^d)$
\ini\begin{equation}
  \label{acG}
 g^{n'|n} :\;\; \hbox{Def}_{\OC^d,n} \simto \hbox{Def}_{g\OC^d,n}
 \To{\alpha_{\OC^d|g\OC^d}^{n'|n}} \hbox{Def}_{\OC^d,n'}.
\end{equation}
On note de la m\^eme mani\`ere les morphismes de sch\'emas formels associ\'es.
Par le lemme pr\'ec\'edent, ils sont $\dd$-{\'e}quivariants,  se
composent comme on l'imagine $ (g'g)^{n''|n} =g'^{n''|n'}\circ
g^{n'|n}$ pour des $n,n',n''$ ad{\'e}quates, et lorsque $g=1$
on retrouve le morphisme de restriction de la structure de niveau.
En particulier, lorsque $(m,m')$ est un autre couple ``ad\'equate''
pour $g$, avec $m\leq n$ et $m'\leq n'$, les morphismes $
g^{n'|n}$ et $g^{m'|m}$ 
commutent aux morphismes de restriction du niveau
$\wh\MC_{LT,n} \To{}  \wh\MC_{LT,m}$ et $\wh\MC_{LT,n'}
\To{}  \wh\MC_{LT,m'}$, par l'\'egalit\'e $
1^{m'|n'}g^{n'|n}= g^{m'|m} 
1^{m|n}$.
Il ne commutent aux morphismes d'oubli vers $\wh\MC_{LT,0}$ que lorsque
$g\OC^d=\OC^d$ et $g$ agit trivialement sur $\wh\MC_{LT,0}$, {\em
  i.e.} lorsque $g\in GL_d(\OC)$.

Heuristiquement, on obtient une action \`a gauche de $G_d$ sur la
``limite projective''  des $\wh\MC_{LT,n}$.
Remarquons que l'action du centre $K^\times$ de $G_d$ est d\'efinie sur chaque
$\MC_{LT,n}$  (par les morphismes $z^{n|n}$)
et que cette action est inverse de celle du
 centre $K^\times$ de
$\dd$. L'action de $\dd\times G_d$ sur le syst\`eme projectif est donc
triviale sur le sous-groupe $\Delta$ d\'efini en \ref{defDelta}.

Enfin, d'apr\`es le lemme pr\'ec\'edent, 
 les fibres g{\'e}n{\'e}riques des morphismes $g^{n'|n}$ 
  sont des morphismes {\'e}tales
  finis et surjectifs ({\em i.e.} des revêtements étales) d'espaces analytiques. 
  Leurs degrés ne dépendent pas de $g$ et sont égaux à $[1+\varpi^{n'}M_d(\OC):1+\varpi^{n}M_d(\OC)]=q^{d^2(n-n')}$.

\medskip

Pour terminer, d\'efinissons l'espace de l'introduction
par
$$\mltn^{d/K}:=\mltn\wh\otimes_\knr \ka.$$
 L'action de $W_K$ sur les espaces analytiques
 est d{\'e}finie de mani{\`e}re analogue
au cas $Dr$ par recollement d'une donn\'ee de descente pour un
Frobenius et de l'action de 
l'inertie donn\'ee par l'extension des scalaires
(c'est un cas particulier de l'action d{\'e}finie par
Rapoport-Zink \cite[3.48]{RZ} dans un contexte plus g{\'e}n{\'e}ral).
C'est une action {\em \`a droite}.

\alin{Morphisme de p{\'e}riodes} \label{perLT}
Dans le cas $LT$, le morphisme de p{\'e}riodes n'est d{\'e}fini qu'en
fibre g{\'e}n{\'e}rique. Il a {\'e}t{\'e} d'abord d{\'e}fini dans
\cite[par. 23]{HG} {\`a} partir de calculs explicites. La d{\'e}finition
la plus visiblement intrins{\`e}que 
est celle de  \cite[ch. 5]{RZ} mais elle fait appel aux techniques
cristallines des groupes $p$-divisibles et n'est {\'e}crite qu'en
in{\'e}gales caract{\'e}ristiques\footnote{Cependant, A. Genestier sait
  adapter au cas d'{\'e}gales caract{\'e}ristiques, par exemple en
  utilisant le module de coordonn{\'e}es}.

Dans chacune des d{\'e}finitions, la construction repose sur l'existence
d'un $\knr$-espace vectoriel $\MM$ de dimension $d$  muni d'une
action de $\dd$ et  attach{\'e} au $\OC$-module formel $\XM$ (son
module de Dieudonn{\'e} dans le cas d'in{\'e}gales caract{\'e}ristiques ou
son module de coordonn{\'e}es en {\'e}gales caract{\'e}ristiques) et d'un
isomorphisme $\dd$-{\'e}quivariant de $\OC_{\mlto}$-modules
$$ \MM \otimes_\knr \OC_{\mlto} \simto M_{X_u} $$
o{\`u} $M_{X_u}$ est un fibr{\'e} vectoriel au-dessus de $\mlto$ attach{\'e}
{\`a}  l'objet
universel $X_u$ (obtenu comme  fibre g{\'e}n{\'e}rique de l'alg{\`e}bre de Lie
de l'extension additive universelle de l'objet
universel $X_u$ au-dessus de 
$\wh\MC_{LT,0}$, resp. de son module de coordonn{\'e}es en {\'e}gales
caract{\'e}ristiques). Notons alors $\hbox{Lie}_{X_u}$ le 
$\OC_{\mlto}$-module inversible obtenu comme  fibre 
g{\'e}n{\'e}rique de l'alg{\`e}bre de Lie de $X_u$
au-dessus de $\wh\MC_{LT,0}$ ; la compos{\'e}e
$$ 
 \MM \otimes_\knr \OC_{\mlto} \simto M_{X_u} \To{} \hbox{Lie}_{X_u}$$
d{\'e}finit un morphisme de $\knr$-espaces analytiques
$$ \xi:\;\;\mlto \To{} \PM(\MM) \simeq \PM^{d-1,nr}$$
qui est le morphisme de p{\'e}riodes voulu.
Par ailleurs, les isomorphismes $\MM(\XM) \To{Frob}
\MM(\phi^*\XM)=\phi^*\MM(\XM)$ munissent $\MM$ et donc $\PM(\MM)$
d'une structure $\phi^\ZM$-{\'e}quivariante au-dessus de $\knr$ (ou
donn{\'e}e de descente {\`a} la Weil). 

\begin{fact} \label{equperLT}
  Le morphisme $\xi$ est {\'e}tale surjectif et $\dd\times W_K$-{\'e}quivariant. Il est aussi
  $G_d$-{\'e}quivariant au sens suivant : pour tout $g\in G_d$ et tous $n,n'$
  ad{\'e}quates, le diagramme
$$\xymatrix{ \mltn \ar[r]^{g^{n'|n}} \ar[d]_{\xi_n} &
  \MC_{LT,n'} \ar[ld]^{\xi_{n'}} \\  \PM(\MM)
  & \\ }$$
est commutatif (bien-s{\^u}r on a not{\'e} $\xi_n$ la compos{\'e}e de $\xi$
avec le morphisme d'oubli de la structure de niveau).
\end{fact}

\begin{proof}
Le fait que $\xi$ est étale au sens rigide-analytique est prouvé dans
\cite[5.17]{RZ}. Il se trouve qu'il est aussi étale au sens de
Berkovich car son bord relatif est vide. La surjectivité de $\xi$ est prouvée, au moins pour les points classiques (rigides analytiques)  dans
\cite[23.5]{HG}, mais la preuve de {\em loc.cit} montre aussi la surjectivité pour les points de Berkovich. Les propriétés d'équivariance se trouvent dans
\cite[5.37]{RZ} ou \cite[23.28]{HG}.   
\end{proof}

La donn{\'e}e de descente sur $\MM$ n'est pas effective mais l'est sur
$\PM(\MM)$ : notons en effet $S_K^{1/d}$ le $K$-espace analytique
obtenu par analytification de la $K$-vari{\'e}t{\'e} de Severi-Brauer
associ{\'e}e {\`a} l'alg{\`e}bre centrale simple sur $K$ d'invariant
$1/d$. Alors 

\begin{fact}
  Il existe un isomorphisme $\dd\times W_K$-{\'e}quivariant de
  $\knr$-espaces analytiques $$\PM(\MM) \simto S^{1/d}_K\wh\otimes_K
  \knr.$$
\end{fact}
\begin{proof}
Dans le cas où $K$ est de caractéristique nulle, le $\OC_K$-groupe formel $\XM$ sur $\o\FM_p$ est $p$-divisible et $\MM$ est l'isocristal associé. Sa dimension sur $\wh{K^{nr}}$ est la hauteur $d$ de $\XM$, et son unique pente est $1/d$. En particulier $\MM$ est irréductible et donc "défini sur $\FM_p$" au sens suivant : il existe une base de $\MM$ dans laquelle le Frobenius est donné par la matrice 
$\Phi=\left(\begin{array}{cccc} 0 & 1 & 0 & 0 \\ 0 & \ddots & \ddots  & 0 \\ 0 & 0 & \ddots & 1 \\ \varpi & 0 & 0 & 0 
 \end{array}\right)$. Ainsi, la forme rationnelle de $\PM^{d-1}$ correspondant à la donnée de descente sur $\PM(\MM)$ est celle associée au $1$-cocycle $\gal(\wh{K^{nr}}/K)=\phi^\ZM \To{} \aut{}{\PM^{d-1}}=PGL_{d-1}(\wh{K^{nr}})$ qui envoie $\phi$ sur la matrice $\Phi$. Cette forme rationnelle est la variété de Severi-Brauer d'invariant $1/d$.

Dans le cas où $K$ est de caractéristique positive, le même raisonnement s'applique au module de coordonnées $\MM$.
\end{proof}


\subsection{D{\'e}finition  du
  $R\Gamma_c$} \label{defRG} 

Dans cette section, nous proposons une d{\'e}finition des complexes de
cohomologie des tours $\mlt^{d/K}$ et $\mdr^{d/K}$. On les obtient comme
{\'e}valuation en un faisceau constant du foncteur d{\'e}riv{\'e} d'un
certain foncteur ``sections {\`a} 
support compact sur la tour''.  Comme on veut que ce foncteur soit {\`a} valeurs dans la
cat{\'e}gorie des $G_d\times \dd\times W_K$-modules, le plus facile (surtout
du c{\^o}t{\'e} $LT$)  est de
prendre pour source de ce foncteur la cat{\'e}gorie des faisceaux
{\'e}quivariants sur (la descente de) l'espace des p{\'e}riodes.

Nous ferons usage du formalisme d{\'e}crit dans la
partie \ref{berkovich}. 

\alin{Coefficients de torsion} \label{cohoDr}
Dans cette section,  
pour faire une construction ou exprimer une propri{\'e}t{\'e} suppos{\'e}e
s'appliquer de la m{\^e}me mani{\`e}re {\`a} la tour de 
Drinfeld $\mdr^{d/K}$  ou {\`a} celle de Lubin-Tate $\mlt^{d/K}$, nous ommettrons
simplement la notation $Dr$ ou $LT$.
Par exemple,
 pour tout $n\in \NM$, $\xi_n$ d{\'e}signera soit la compos{\'e}e
$$ \xi_{Dr,n} :\;\; \mdrn^{d/K} \To{} \mdro^{d/K} \To{} \Omega^{d-1,nr}_K \To{}
\PC_{Dr}:=\Omega^{d-1}_K $$
soit la compos{\'e}e
$$ \xi_{LT,n} :\;\; \mltn^{d/K} \To{} \mlto^{d/K} \To{} \PM^{d-1,nr}_K \To{}
\PC_{LT}:=S^{1/d}_K. $$
Il s'agit de morphismes d'espaces ``analytiques au-dessus de $K$'' au
sens de \cite{Bic2}, les deux premiers espaces \'etant
$\ka$-analytiques, le troisi\`eme $\knr$-analytique, et le dernier $K$-analytique. 
Dans chacun des cas, l'espace des p{\'e}riodes $\PC$ est muni d'une action
 continue du groupe $J$ en posant  $J_{Dr}:=G_d$ et $J_{LT}:=\dd$ et d'une
 action triviale de $W_K$. Le morphisme $\xi_{n}$ est
 formellement {\'e}tale et $J \times W_K$-{\'e}quivariant. Dans le cas $Dr$,
 il est aussi $\dd$-{\'e}quivariant pour l'action triviale de $\dd$ sur
 $\Omega^{d-1}_K$ et l'action d{\'e}finie en \ref{actionDr} sur les $\mdrn^{d/K}$. Dans le
 cas $LT$, c'est plut{\^o}t le syst{\`e}me inductif des $\xi_{LT,n}$ qui est
 $G_d$-{\'e}quivariant au sens de \ref{equperLT}.

Dans chacun des cas, si
 $\FC$ est un faisceau {\'e}tale de torsion sur $\PC$
qui est $J$-{\'e}quivariant
discret/lisse au sens de \ref{fgeq}, alors
pour tout $n$, le groupe ab\'elien
$\Gamma_c(\MC_n^{d/K},\xi_n^*(\FC))$ est muni d'une action \`a gauche
{\em lisse} de $J$ et d'une action  \`a droite de $W_K$.
Si $n' \leq n$, le morphisme $\MC_n^{d/K} \To{\pi_{n,n'}}
\MC_{n'}^{d/K}$ est fini, de sorte que $\pi_{n,n',!}=\pi_{n,n',*}$ et
on a un morphisme naturel
$$ \pi_{n,n'}^*:\;\;\Gamma_c(\MC_{n'}^{d/K},\xi_{n'}^*(\FC)) \To{}
\Gamma_c(\MC_{n'}^{d/K},\pi_{n,n',!}\pi_{n,n'}^*\xi_{n'}^*(\FC)) =
\Gamma_c(\MC_n^{d/K},\xi_n^*(\FC))$$
qui est $ J\times W_K$-\'equivariant.
Alors 
le groupe ab{\'e}lien
$$ \Gamma_c(\MC^{d/K},\FC):= \limi{n\in \NM}
\Gamma_c(\MC_n^{d/K},\xi_n^*(\FC)) $$
est muni d'une action lisse de $J \times W_K$.
\begin{fact} \label{actiontours}
  Le groupe ab\'elien $\Gamma_c(\MC^{d/K},\FC)$ est muni d'une action
  de $G_d\times \dd\times W_K$ lisse sur les deux premiers facteurs et
  triviale  sur le
sous-groupe $\Delta\times 1$ (rappelons que $\Delta$ a \'et\'e
d\'efini en \ref{defDelta}).
\end{fact}
\begin{proof}
Le cas $Dr$ est facile et laiss\'e au lecteur ; en fait, on a d\'eja
une action de $G_d\times \dd\times W_K$ sur chaque
$\Gamma_c(\MC_{Dr,n}^{d/K},\xi_n^*(\FC))$.
On s'occupe donc du cas $LT$, o\`u on a d\'eja l'action de $\dd\times
W_K$ par ce qui pr\'ec\`ede.

Fixons $g\in G_d$ et $n \geq n'$ tels que
le morphisme $g^{n'|n} : \mltn^{d/K} \To{} \MC_{LT,n'}^{d/K}$
soit d\'efini. Alors ce morphisme est fini de sorte que 
$(g^{n'|n})_!=(g^{n'|n})_*$ et on
obtient le morphisme
\begin{eqnarray*}
g^{n'|n,*}:\;\;  \Gamma_c(\MC_{LT,n'}^{d/K},\xi_{LT,n'}^*(\FC)) &  \To{} &
\Gamma_c(\MC_{LT,n'}^{d/K},({g^{n'|n}})_! (g^{n'|n})^*
\xi_{LT,n'}^*(\FC)) \\&=  & \Gamma_c(\MC_{LT,n}^{d/K},\xi_{LT,n}^*(\FC))
\end{eqnarray*}

Les propri\'et\'es des morphismes $g^{n'|n}$ mentionn\'ees lors de leur
construction en \ref{actionLT} assurent la commutativit\'e des
diagrammes du type 
$$\xymatrix{\Gamma_c(\MC_{LT,n'}^{d/K},\xi_{LT,n'}^*(\FC)) \ar[d]_{\pi_{m',n'}^*}
  \ar[r]^{g^{n'|n,*}} &
\Gamma_c(\MC_{LT,n}^{d/K},\xi_{LT,n}^*(\FC)) \ar[d]^{\pi_{m,n}^*} \\
\Gamma_c(\MC_{LT,m'}^{d/K},\xi_{LT,m'}^*(\FC)) \ar[r]^{g^{m'|m,*}} &
\Gamma_c(\MC_{LT,m}^{d/K},\xi_{LT,m}^*(\FC))
}$$
o\`u $m\geq n$ et $m'\geq n'$ sont tels que $g^{m'|m}$ est
d\'efini. Ainsi on obtient \`a la limite inductive un morphisme
$g^*:\; \Gamma_c(\MC^{d/K},\FC) \To{} \Gamma_c(\MC^{d/K},\FC)$. Par
d\'efinition il commute \`a l'action de $\dd\times W_K$, et si $g'$
est un autre \'el\'ement de $G_d$, alors $(gg')^*={g'}^*g^*$ de sorte
qu'on obtient une action {\em \`a droite} de $G_d$ sur
$\Gamma_c(\MC^{d/K},\FC)$. 
{\em Pour nous ramener \`a une action \`a gauche,
 nous ferons agir $g$ sur $\Gamma_c(\MC^{d/K},\FC)$ par son inverse
${g}^{-1}$.}

Soit maintenant $z\in K^\times$ et $z_G$, resp. $z_D$, son image dans
le centre de $G_d$, resp. dans celui de $\dd$. 
Alors les actions g\'eom\'etriques de
$z_{G}$ et $z_D$ sur la tour stabilisent chaque
$\MC^{d/K}_{LT,n}$ et sont inverses l'une de l'autre. 
Il en est donc de m\^eme sur les sections $\Gamma_c(\MC^{d/K}_{LT,n},\FC)$,
d'o\`u la
trivialit\'e de l'action de $\Delta\times 1$.

\end{proof}

Si $\Lambda$ est un anneau de torsion, on d{\'e}finit ainsi un foncteur
exact {\`a} gauche 
$$\Gamma_c(\MC^{d/K},-):\;\; \FC(\PC,\Lambda J) \To{}
\Lambda(GD\times W_K)-\hbox{mod},$$
en notant $\FC(\PC,\Lambda J) $ la cat\'egorie des $\Lambda$-faisceaux
\'etales $J$-\'equivariants discrets/lisses sur $\PC$, et
$\Lambda(GD\times W_K)-\hbox{mod}$ 
la cat{\'e}gorie ab{\'e}lienne des
$\Lambda$-modules munis d'une action lisse de $GD=(G_d\times \dd)/\Delta$ et d'une action
commutante de $W_K$.

Remarquons maintenant que la construction du syst\`eme inductif 
$(\Gamma_c(\MC^{d/K}_n,\xi_n^*(\FC)))_{n\in\NM}$ et de l'action de
$G_d\times \dd\times W_K$ sur icelui reposaient sur la finitude (en
fait la propret\'e) des morphismes $\pi_{n,n'}$ et $\alpha_g^{n'|n}$,
et que ces m\^emes propri\'et\'es permettent de d\'efinir de la
m\^eme mani\`ere pour chaque $q\in\NM^*$ un syst\`eme inductif
$(H^q_c(\MC^{d/K}_n,\xi_n^*(\FC)))_{n\in\NM}$ muni d'une action de
$G_d\times \dd\times W_K$.

\begin{fact}
 On a pour tout $q\in\NM$ un isomorphisme
canonique de foncteurs :
$$ R^q\Gamma_c(\MC^{d/K},-) \simto \limi{n\in\NM} H^q_c(\MC_n^{d/K},\xi_n^*-).$$
\end{fact}
\begin{proof}
En effet, par commutation de la cohomologie aux limites inductives
filtrantes, on a $$R^q\Gamma_c(\MC^{d/K},-) =\limi{}
R^q\left(\Gamma_c(\MC_n^{d/K},-)\circ \xi_n^*\right).$$ 
Il nous faut donc montrer que $R^q(\Gamma_c(\MC_n^{d/K},-)\circ
\xi_n^*)=R^q\Gamma_c(\MC_n^{d/K},-)\circ \xi_n^*$ pour tout $n$. Or on a une factorisation $\xi_n^* : (\PC)_{et}  \To{\beta^*} \left(\PC{\wh\otimes}_K \wh{K^{ca}}\right)_{et} \To{{\xi_n^{ca}}^*} (\MC_n^{d/K})_{et} $ o\`u $\xi_n^{ca}$ est un morphisme \'etale de $\wh{K^{ca}}$-espaces analytiques et $\beta$ est le changement de base (ce n'est pas un morphisme d'espaces analytiques mais un morphisme des sites \'etales concern\'es). Ainsi ${\xi_n^{ca}}^*$ envoie injectifs sur injectifs puisqu'il est adjoint \`a droite du foncteur exact ${\xi_n^{ca}}_!$, et $\beta^*$ envoie suffisamment d'injectifs sur des injectifs, par exemple tous ceux de la forme $\beta_*(\IC)$, avec $\IC$ injectif.
\end{proof}

\alin{Coefficients $l$-adiques} \label{defrgadic}
Fixons maintenant un anneau de valuation discr\`ete complet $\Lambda$ de
caract{\'e}ristique r{\'e}siduelle $\neq p$, et
$\FC_\bullet=(\FC_m)_{m\in \NM}$
un $\Lambda_\bullet$-faisceau $J$-{\'e}quivariant discret/lisse sur $\PC$,
au sens
de la section \ref{coefladic}.
Pour $n\in\NM$, nous notons
$$ \xi_n^{*,\infty}(\FC_\bullet) :=
\mathop{\limp{m\in\NM}}\nolimits^\infty(\xi_n^{*}\FC_m)$$ 
le  $\Lambda$-faisceau $J$-{\'e}quivariant discret/lisse sur
$\MC^{d/K}_n$ associ{\'e} comme en
\ref{deflimpinf} au $\Lambda_\bullet$-faisceau $J$-{\'e}quivariant discret/lisse
$(\xi_n^*\FC_m)_{m\in\NM}$ sur $\MC^{d/K}_n$. 
D'apr\`es \ref{lemmelimpinf}, les morphismes $\pi_{n,n'}$ pour $n\geq n'$ \'etant \'etales, on a 
$$ \xi_n^{*,\infty}(\FC_\bullet) \simeq \pi_{n,n'}^*\xi_{n'}^{*,\infty}(\FC_\bullet). $$
La finitude de ces m\^emes $\pi_{n,n'}$ permet donc comme dans le cas
de torsion de d\'efinir des morphismes de $\Lambda$-modules
$$ \pi_{n,n'}^*:\;\; \Gamma_c(\MC^{d/K}_{n'},\xi_{n'}^{*,\infty}(\FC_\bullet)) \To{}
\Gamma_c(\MC^{d/K}_n,\xi_{n}^{*,\infty}(\FC_\bullet)) $$
et de poser 
$$ \Gamma_c(\MC^{d/K},(\FC_m)_{m\in\NM}):= \limi{n\in \NM}
\Gamma_c\left(\MC_n^{d/K},\xi_{n}^{*,\infty}(\FC_\bullet)\right).$$
Ce groupe est muni par d\'efinition d'une action {\em lisse} de $J$
\`a gauche et d'une action commutante de $W_K$ \`a droite.
On construit l'action du troisi\`eme groupe ($\dd$ pour le cas $Dr$ ou
$G_d$ pour le cas $LT$) exactement comme dans le cas de torsion. En
particulier dans le cas $LT$, on utilise le fait que les morphismes
$\alpha_g^{n'|n}$ sont \'etales pour pourvoir \'ecrire
$$ {\alpha_g^{n'|n}}^*\xi_{n'}^{*,\infty}(\FC_\bullet) \simeq
\xi_n^{*,\infty}(\FC_\bullet) $$
gr\^ace \`a \ref{equperLT} et \ref{lemmelimpinf}


On d\'efinit ainsi un foncteur 
$$ \Gamma_c(\MC^{d/K},-):\;\; {\FC(\PC,\Lambda_\bullet J)} \To{}
\Lambda(GD\times W_K)-\hbox{mod} $$
qui est exact {\`a} gauche. 

Supposons maintenant que le $\Lambda_\bullet$-faisceau $\FC_\bullet$ soit un
$\Lambda$-syst\`eme local au sens de \ref{defBerk}, et notons 
$H^q_c(\MC_n^{d/K},\xi_n^*((\FC_m)_m))$ la cohomologie \`a supports
compacts de son image inverse sur $\MC_n^{d/K}$, cohomologie
dont la d\'efinition par Berkovich est rappel\'ee en \ref{defBerk}.
Toujours les m\^emes propri\'et\'es de propret\'e des $\pi_{n,n'}$ et
$\alpha_g^{n'|n}$ permettent de d\'efinir un syst\`eme inductif
$(H^q_c(\MC_n^{d/K},\xi_n^*((\FC_m)_m)))_{n\in\NM}$ muni d'une
action du groupe $G_d\times\dd\times W_K$.

\begin{fact}
  Pour tout $q\in \NM$ et tout $\Lambda$-syst{\`e}me local $J$-{\'e}quivariant
$(\FC_m)_{m\in\NM}$,
 on a un isomorphisme canonique
  $(G_d\times \dd\times W_K)$-{\'e}quivariant :
$$ R^q\Gamma_c(\MC^{d/K},\FC_\bullet) \simto \limi{n\in \NM}
H^q_c\left(\MC_n^{d/K},\xi_n^*((\FC_m)_m)\right) $$
\end{fact}
\begin{proof}
  Comme dans le cas de torsion, la commutation de la cohomologie aux limites
  inductives filtrantes et le fait que $\xi_n^*$ envoie suffisamment d'injectifs sur des injectifs montrent que 
\begin{eqnarray*}
  R^q\Gamma_c(\MC^{d/K},\FC_\bullet) & \simeq & \limi{n\in\NM}
\left(R^q\left(\Gamma_c(\MC_n^{d/K},-) \circ \limproj^\infty \circ
    \xi_n^* \right)(\FC_\bullet) 
 \right) \\
& \simeq &  \limi{n\in\NM} \left(
R^q\left(\Gamma_c(\MC_n^{d/K},-) \circ \limproj^\infty  \right)
(\xi_n^*(\FC_\bullet)) \right)
\end{eqnarray*}
Mais d'apr{\`e}s le corollaire \ref{coroladic} appliqu{\'e} aux espaces analytiques
quasi-alg{\'e}briques (car lisses) $\MC_n^{d/K}$, on a pour tout $\Lambda$-syst{\`e}me
local
$$ R^q\left(\Gamma_c(\MC_n^{d/K},-) \circ \limproj^\infty
\right)(\xi_n^*(\FC_\bullet)) \simeq
H^q_c\left(\MC_n^{d/K},\xi_n^*(\FC_m)_m\right).$$
 \end{proof}

Prenant $\Lambda:=\ZM_l$ et $(\FC_m)_{m\in\NM}= (\ZM/l^m\ZM)_{m\in\NM}$
dans les d{\'e}finitions pr{\'e}c{\'e}dentes,
nous posons
$$  R\Gamma_c(\MC^{d/K},\ZM_l):=R\Gamma_c(\MC^{d/K},(\ZM/l^m\ZM)_{m\in\NM}).$$
C'est un objet de la cat{\'e}gorie
d{\'e}riv{\'e}e born{\'e}e $D^b_{\ZM_l}(GD)$ de la cat\'egorie
ab\'elienne $\Mo{\ZM_l}{GD}$ des $\ZM_l(GD)$-modules
lisses, qui est muni d'une action de $W_K$ et dont la cohomologie est $G_d\times
\dd\times W_K$-isomorphe {\`a} la limite inductive $H^q_c(\MC^{d/K},\ZM_l)$ des
$H^q_c(\MC_n^{d/K},\ZM_l)$, $n\in\NM$. 
Nous posons enfin 
$$ R\Gamma_c(\MC^{d/K},\o\QM_l):= \o\QM_l\otimes_{\ZM_l}
R\Gamma_c(\MC^{d/K},\ZM_l) \;\; \in D^b_{\o\QM_l}(GD). $$
Par commutation du produit tensoriel aux limites inductives, la
cohomologie de ce
complexe est donn{\'e}e par 
$$ R^q\Gamma_c(\MC^{d/K},\o\QM_l) \simeq H^q_c(\MC^{d/K},\o\QM_l) :=
\limi{n\in\NM} H^q_c(\MC_n^{d/K},\o\QM_l).$$

\subsection{Variantes et décompositions}

Dans cette section, on discute deux variantes de la construction de la section précédente dont on tire quelques propriétés des complexes $R\Gamma_c(\MC^{d/K},\Lambda)$ où $\Lambda$  désigne un anneau commutatif o{\`u} $p$ est inversible, et qui sera soit de torsion, soit local complet, soit {\'e}gal {\`a} $\o\QM_l$. 
On s'intéresse ensuite à l'action du centre de la catégorie des $\Lambda(GD)$-modules lisses sur ces complexes et on montre une propriété de finitude dans le cas $LT$.

\alin{Premi\`ere variante} \label{variante}
Dans \cite{CarAA}, \cite{HaTay} et \cite{Boyer1}, les auteurs consid{\`e}rent
plut{\^o}t la cohomologie des espaces
$\MC_{?,n}^{d/K,(0)}:=\MC_{?,n}^{(0)}\wh\otimes_\knr \ka$ pour $?=LT$
ou $Dr$,  rencontr\'es
dans les paragraphes \ref{tourLT} et \ref{tourDr}.
Ces espaces ne sont plus munis
d'une action de $G_d\times \dd\times W_K$ mais plut{\^o}t du sous-groupe 
 $(G_d\times \dd\times W_K)^0$ (la notation $^0$ a \'et\'e introduite en \ref{notationsnu}). L'action de ce sous-groupe
 consid{\'e}r{\'e}e 
dans les articles sus-cit{\'e}s est simplement la restriction \`a
$\MC^{d/K,(0)}$ de l'action
qu'on a d{\'e}finie en \ref{actiontours} sur les espaces $\MC^{d/K}$. 
{
On laisse au lecteur le soin de vérifier que le lien entre ces espaces de cohomologie et ceux de
ce texte est le suivant : le morphisme naturel  $
H^q_c(\MC^{d/K,(0)},\Lambda) \To{}  H^q_c(\MC^{d/K},\Lambda)$ est
$(G_d\times\dd\times W_K)^0$-\'equivariant et induit un
isomorphisme  $G_d\times \dd\times W_K$-{\'e}quivariant
$$ \cind{(G_d\times \dd\times
  W_K)^0}{G_d\times \dd\times W_K}{ H^q_c(\MC^{d/K,(0)},\Lambda)} \simto
H^q_c(\MC^{d/K},\Lambda)$$ 
o{\`u} {\sl ind} d{\'e}signe l'induction {\`a} supports compacts.
Plus généralement, soit $\FC$ un faisceau étale abélien de torsion
$J$-équivariant discret/lisse  sur  $\PC$ (avec les notations de
\ref{cohoDr}),
et d\'esignons par $\xi_n^{(0)}$ la restriction de $\xi_n$ \`a
$\MC_n^{d/K,(0)}$.
Alors on peut munir comme dans les paragraphes pr\'ec\'edents le
groupe ab\'elien $\limi{n\in \NM} 
\Gamma_c(\MC_n^{d/K,(0)},\xi_n^{(0),*}(\FC))$ d'une action
du groupe $(G_d\times\dd\times W_K)^0$ et on a un
isomorphisme $G_d\times \dd \times W_K$-équivariant 
$$\cind{(G_d\times\dd \times W_K)^0}{G_d\times\dd\times W_K}{\limi{n\in \NM}
\Gamma_c(\MC_n^{d/K,(0)},\xi_n^{(0),*}(\FC))} \simto  \limi{n\in \NM}
\Gamma_c(\MC_n^{d/K},\xi_n^*(\FC)).$$

Soit $\psi :\ZM \To{} \Lambda^\times$ un caract\`ere et
$\psi_{GDW}:= \psi\circ\nu_{GDW}$
le caract\`ere de $G_d\times \dd\times W_K$ associ\'e comme au paragraphe \ref{notationsnu}.
L'isomorphisme pr\'ec\'edent implique l'existence d'un isomorphisme 
$G_d\times \dd\times W_K$-\'equivariant et fonctoriel en le
$\Lambda$-faisceau $\FC$ 
$$a_\psi(\FC):\;\;  \limi{n\in \NM}
\Gamma_c(\MC_n^{d/K},\xi_n^*(\FC)) \simto \psi_{GDW}\otimes \limi{n\in \NM}
\Gamma_c(\MC_n^{d/K},\xi_n^*(\FC)).$$
La m\^eme discussion s'applique aux $\Lambda_\bullet$-faisceaux
lorsque $\Lambda$ est local complet, de sorte que dans chaque cas, on obtient
un isomorphisme
\ini\begin{equation}\label{indui}
Ra_\psi(\Lambda):\;\;  R\Gamma_c(\MC^{d/K},\Lambda) \simto
\psi_{GDW}\otimes R\Gamma_c(\MC^{d/K},\Lambda)
\end{equation}
Plus pr\'ecis\'ement, on a obtenu

\begin{lemme} \label{torsion} Il existe un isomorphisme dans $D^b_\Lambda(GD)$
$$ Ra_\psi(\Lambda):\;\; R\Gamma_c(\MC^{d/K},\Lambda) \simto
\psi_{GD}\otimes R\Gamma_c(\MC^{d/K},\Lambda)$$
qui est $\psi_{W}$-\'equivariant, au sens o\`u le diagramme suivant
est commutatif
$$\xymatrix{
  \endo{D^b_\Lambda(GD)}{
    R\Gamma_c(\MC^{d/K},\Lambda)}^\times  \ar[r]^{Ra_\psi(\Lambda)}
  & \endo{D^b_\Lambda(GD)}{\psi_{GD} \otimes R\Gamma_c(\MC^{d/K},\Lambda)}^\times \\
W_K \ar[ru]_{\psi_W\otimes\gamma} \ar[u]^{\gamma} }
$$
\end{lemme}
On a not\'e dans cet énoncé $\gamma :\;W_K \To{} \endo{D^b_\Lambda(GD)}{
    R\Gamma_c(\MC^{d/K},\Lambda)}^\times$ le morphisme de groupes
  donnant l'action de $W_K$ sur $R\Gamma_c(\MC^{d/K},\Lambda)$.

\alin{Deuxi\`eme variante} \label{variante2}
Rappelons que $\varpi$ d\'esigne une
uniformisante de $\OC$. Les \'el\'ements $(\varpi,1,1)$ et $(1,\varpi^{-1},1)$ de
$G_d\times\dd\times W_K$ agissent de la m\^eme mani\`ere sur les
$\knr$-espaces analytiques $\MC_n$ et l'action
du groupe libre \`a un g\'en\'erateur $\varpi^\ZM$ qui s'en d\'eduit 
est clairement libre. Le quotient $\MC_n/\varpi^\ZM$ est encore un
$\knr$-espace analytique : il est (non-canoniquement) isomorphe \`a
$\bigsqcup_{h=0}^{d-1} \MC_n^{(h)}$.
Comme $\varpi^\ZM$ est central dans $G_d\times\dd\times W_K$, ces
quotients sont encore munis d'une action du produit triple.
Les morphismes de transition $\pi_{n,n'}$ passent au quotient ainsi
que les morphismes de p\'eriodes, de sorte qu'on
peut appliquer les constructions des paragraphes \ref{cohoDr} et \ref{defrgadic} aux deux tours
$(\MC_{?,n}^{d/K}/\varpi^\ZM)_{n\in\NM}$.

\begin{nota}
Pour all\'eger un peu les notations, nous surlignerons les quotients par $\varpi^\ZM$. En particulier nous noterons $\o{GD}:=GD/\varpi^\ZM$  et nous utiliserons la notation
$ \o\MC^{d/K}:=\MC^{d/K}/\varpi^\ZM $
que nous d\'ecorerons d'indices $Dr$, $LT$ ou $n$ selon les besoins.
\end{nota} 
 
On obtient en particulier un complexe (deux complexes, en fait)
$$R\Gamma_c(\o\MC^{d/K},\Lambda) \in D^b_\Lambda(\o{GD}).$$
Pour le comparer au complexe $R\Gamma_c(\MC^{d/K},\Lambda)$,
introduisons le foncteur exact \`a droite des ``$\varpi^\ZM$-coinvariants'' :
$$\application{}{\Mo{\Lambda}{GD}}{\Mo{\Lambda}{\o{GD}}}{V}
{\Lambda\otimes_{\Lambda[\varpi^\ZM]}V}$$ 
o\`u le produit tensoriel est pris par rapport au morphisme
d'augmentation $\Lambda[\varpi^\ZM]\To{}\Lambda$. Ce foncteur est
adjoint \`a gauche du foncteur d'inclusion $\omega :\;\;
\Mo{\Lambda}{\o{GD}} \To{} \Mo{\Lambda}{GD}$. 

\begin{lemme} \label{quotient}
Il existe un isomorphisme ``canonique'' dans $D^b_\Lambda(\o{GD})$ 
$$ \Lambda {\mathop\otimes\limits^L}_{\Lambda[\varpi^\ZM]} R\Gamma_c(\MC^{d/K},\Lambda)
\simto R\Gamma_c(\o\MC^{d/K},\Lambda).$$ 
\end{lemme}
\begin{proof}
Soit $\FC$ un $\Lambda$-faisceau, resp. $\FC_\bullet$ un
$\Lambda_\bullet$-faisceau,  $J$-\'equivariant discret/lisse sur
$\PC$. Soit $p_n$ le morphisme de projection $\MC_n^{d/K} \To{}
\MC_n^{d/K}/\varpi^\ZM$ et $\o\xi_n :\MC_n^{d/K}/\varpi^\ZM \To{} \PC$
la descente de $\xi_n$. On a $\xi_n=\o\xi_n\circ p_n$. 
Comme $p_n$ est \'etale, le morphisme ${p_n}_!$ est adjoint \`a gauche
de $p_n^*$. 
On a donc une application $J\times
W_K$-\'equivariante
$${p_n}_!:\;\; \Gamma_c(\MC_n^{d/K},\xi_n^*(\FC)) \To{}
\Gamma_c(\MC_n^{d/K}/\varpi^\ZM,\o\xi_n^*(\FC)). $$
Utilisant \ref{lemmelimpinf} on a de m\^eme une application 
$$ {p_n}_!:\;\; \Gamma_c(\MC_n^{d/K},\xi_n^{*,\infty}(\FC_\bullet)) \To{}
\Gamma_c(\MC_n^{d/K}/\varpi^\ZM,\o\xi_n^{*,\infty}(\FC_\bullet)). $$
Par ailleurs, le morphisme $p_n$ est ``restreint'' au sens de
\cite[5.3.6]{Bic2}, de sorte que par {\em loc.cit}, le carr\'e
cart\'esien o\`u $n'\geq n$
$$\xymatrix{ \MC_n^{d/K} \ar[r]^{p_n} & \MC_n^{d/K}/\varpi^\ZM \\
\MC_{n'}^{d/K}  \ar[r]_{p_{n'}} \ar[u]^{\pi_{n',n}} &
\MC_{n'}^{d/K}/\varpi^\ZM \ar[u]_{\o\pi_{n',n}} }
$$
assure que le morphisme naturel $\pi_{n',n}^*\circ {p_n}_! \To{}
{p_{n'}}_!\circ \o\pi_{n',n}^*$ est un isomorphisme.
Ceci signifie que les applications $({p_n}_!)_{n\in\NM}$ ci-dessus
sont compatibles aux transitions dans les syst\`emes inductifs.
On v\'erifie de m\^eme qu'elles sont compatibles \`a l'action du
``troisi\`eme groupe'', de sorte 
qu'on obtient  des transformations naturelles 
$$ \Gamma_c(\MC^{d/K},-) \To{} \Gamma_c(\o\MC^{d/K},-)$$
entre foncteurs $\FC(\PC,\Lambda J)\To{} \Lambda(GD)\times
W_K\hbox{-Mod}$, resp. $\FC(\PC,\Lambda_\bullet J)\To{} \Lambda(GD)\times
W_K\hbox{-Mod}$.
En d\'erivant et en \'evaluant convenablement, on obtient un morphisme
$W_K$-\'equivariant dans $D^b_\Lambda(GD)$
$$ R\Gamma_c(\MC^{d/K},\Lambda)
\To{} \omega R\Gamma_c(\o\MC^{d/K},\Lambda)$$
o\`u $\omega$ est le foncteur $D^b_\Lambda(\o{GD})\To{}
D^b_\Lambda(GD)$ d'infation.
  Le morphisme de l'\'enonc\'e est alors obtenu par adjonction \`a partir
  du morphisme ci-dessus.
Pour montrer que c'est un isomorphisme, il suffit de voir qu'il induit
des isomorphismes au niveau des groupes de cohomologie (o\`u l'on peut
oublier l'action des groupes). Or cela est \'evident puisqu'on a des
isomorphismes de $\ka$-espaces analytiques \'equivariants sous $\varpi^\ZM$ (mais pas
\'equivariants sous $G_d\times \dd\times W_K$) 
$$ \MC^{d/K} \simto \left(\bigsqcup_{h=0}^{d-1}
  \MC^{d/K,(h)}\right)\times \varpi^\ZM  \simto
\left(\o\MC^{d/K}\right)\times \varpi^\ZM.$$
\end{proof}
\alin{D{\'e}composition selon le centre de la cat{\'e}gorie des $\Lambda
  \dd$-modules lisses} 
 Notons $\ZG_\Lambda( GD)$ le centre de la
cat{\'e}gorie ab{\'e}lienne $\Mo{\Lambda}{GD}$ des $\Lambda GD$-modules lisses. 
 Comme le centre d'une cat{\'e}gorie
ab{\'e}lienne agit encore sur la cat{\'e}gorie d{\'e}riv{\'e}e born{\'e}e
associ{\'e}e (et s'identifie m{\^e}me au centre de celle-ci), on obtient une
action canonique de $\ZG_\Lambda( GD)$ sur
$R\Gamma_c(\MC^{d/K},\Lambda)$. En particulier tout idempotent de
$\ZG_\Lambda( GD)$ fournit un facteur direct de $R\Gamma_c(\MC^{d/K},\Lambda)$.
De plus, en utilisant les m{\^e}mes notations pour les groupes $\dd$ et $G_d$, on a des morphismes
canoniques $\ZG_\Lambda(\dd) \To{} \ZG_\Lambda( GD)$ et $\ZG_\Lambda( G_d)
\To{} \ZG_\Lambda( GD)$ induits par
les inclusions $\dd \injo GD$ et $G_d\injo GD$.

La suite de sous-ensembles  $1+\varpi^n\od$, $n\in\NM$ de $\dd$ est un syst{\`e}me
fondamental de voisinages de l'unit{\'e} form{\'e} de pro$-p$-sous-groupes
ouverts compacts distingu{\'e}s de $\dd$. Comme $p$ est inversible dans
$\Lambda$, ces pro-$p$-sous-groupes ouverts d{\'e}finissent des idempotents de l'alg{\`e}bre
de Hecke $\HC_\Lambda(\dd)$. Comme ils sont distingu{\'e}s, ces
idempotents sont centraux : on note
$\varepsilon_{n,D}$ leurs images dans $\ZG_\Lambda( \dd)$ et
$\varepsilon_{n,D}R\Gamma_c(\MC^{d/K},\Lambda) \in D^b_\Lambda( GD)$ le facteur direct 
de $R\Gamma_c(\MC^{d/K},\Lambda)$ associ{\'e}.

{\em Interpr{\'e}tation g{\'e}om{\'e}trique :}
Du côté $Dr$, on prouve facilement (comme dans le lemme \ref{niveauLT} ci-dessous, en plus simple) que l'inclusion de foncteurs
$\Gamma_c(\mdrn^{d/K},-) \circ \xi_{Dr,n}^* \injo \Gamma_c(\mdr^{d/K},-)$
induit un isomorphisme $W_K$-{\'e}quivariant dans $D^b_\Lambda (GD)$
$$ R\Gamma_c(\MC_{Dr,n}^{d/K},\Lambda) \simto \varepsilon_{n,D}R\Gamma_c(\mdr^{d/K},\Lambda).$$
Par contre, il n'y a pas d'interpr{\'e}tation g{\'e}om{\'e}trique {\'e}vidente
dans le cas $LT$.

\alin{D{\'e}composition selon le centre de la cat{\'e}gorie des $\Lambda
  G_d$-modules lisses} 
On suppose toujours que 
$\Lambda$ est une anneau o{\`u} $p$ est inversible. Il existe une
d{\'e}composition ``par le niveau''  similaire {\`a} celle que nous avons
utilis{\'e}e pour $\dd$ ci-dessus, mais beaucoup moins triviale. Notons
$\Mo{\Lambda}{G_d}_n$ la sous-cat{\'e}gorie pleine de $\Mo{\Lambda}{G_d}$ formée des objets engendr{\'e}s
par leur vecteurs invariants sous le pro-$p$-sous-groupe ouvert
$H_n:=1+\varpi^nM_d(\OC)$ de $G_d=GL_d(K)$, et notons $\HC_\Lambda(G_d,H_n)$
l'alg{\`e}bre de Hecke de la paire $(G_d,H_n)$ {\`a} coefficients dans $\Lambda$.
\begin{fact} \label{hn}
  La sous-cat{\'e}gorie $\Mo{\Lambda}{G_d}_n$ est ``facteur direct'' de la
  cat{\'e}gorie $\Mo{\Lambda}{G_d}$. Les foncteurs
$$
{V}\mapsto {V^{H_n}}
\;\;\hbox{ et } \;\;
{M}\mapsto{\cind{H_n}{G_d}{1}\otimes_{\HC_\Lambda(G_d,H_n)} M}
$$
sont des {\'e}quivalences ``inverses''  l'une de l'autre entre
$\Mo{\Lambda}{G_d}_n$ et la cat{\'e}gorie des $\HC_\Lambda(G_d,H_n)$-modules.
\end{fact}
\begin{proof} elle est donn\'ee
  dans l'appendice \ref{decomposition}.
\end{proof}

On associe ainsi {\`a} tout entier $n\in\NM$ un idempotent central
$\varepsilon_{n,G}$ de $\ZG_\Lambda( G_d)$ qui "projette" la catégorie $\Mo{\Lambda}{G_d}$ sur sa sous-catégorie $\Mo{\Lambda}{G_d}_n$. 
Faisant agir cet idempotent sur la catégorie $\Mo{\Lambda}{GD}$, on en obtient un "facteur direct" $\varepsilon_{n,G}\Mo{\Lambda}{GD}$. On peut interpréter ce facteur direct en termes de modules lisses sur l'algèbre de Hecke $\HC_\Lambda(GD,H_n)$ des mesures localement constantes à supports compacts et à valeurs dans $\Lambda$ qui sont invariantes à droite et à gauche sous l'image de $H_n\To{} GD$
(algèbre qui n'a pas d'unité, mais "suffisamment d'idempotents"). En effet, en notant toujours $\cInd{H_n}{GD}$ l'induction lisse à supports compacts, on déduit du fait précédent que les foncteurs 
$$
{V}\mapsto {V^{H_n}}
\;\;\hbox{ et } \;\;
{M}\mapsto{\cind{H_n}{GD}{1}\otimes_{\HC_\Lambda(GD,H_n)} M}
$$
sont inverses l'un de l'autre.
Cette fois-ci l'interpr{\'e}tation
g{\'e}om{\'e}trique se fait du c{\^o}t{\'e} $LT$ :
\begin{lemme}\label{niveauLT}
 L'inclusion de foncteurs
$\Gamma_c(\mltn^{d/K},-) \circ \xi_{LT,n}^{*}\injo \Gamma_c(\mlt^{d/K},-)$ induit
un isomorphisme $W_K$-{\'e}quivariant dans
$D^b( \HC_\Lambda(GD,H_n))$
$$ R\Gamma_c(\mltn^{d/K},\Lambda) \simto
R\Gamma_c(\mlt^{d/K},\Lambda)^{H_n}$$
qui induit {\`a} son tour un isomorphisme dans $D^b_\Lambda( GD)$
$$ \cind{H_n}{GD}{1}\otimes_{\HC_\Lambda(GD,H_n)}R\Gamma_c(\mltn^{d/K},\Lambda) \simto 
\varepsilon_{n,G} R\Gamma_c(\mlt^{d/K},\Lambda).$$
\end{lemme}
\begin{proof}
Nous traitons seulement le cas $\Lambda$ complet de valuation discr\`ete et laissons le cas de torsion au lecteur. 
Il nous suffit de montrer que pour tout $\Lambda_\bullet$-faisceau étale $\dd$-équivariant
discret/lisse $\FC_\bullet$ sur $\PC_{LT}=S^{1/d}_K$, le morphisme canonique  
$\Gamma_c(\mltn^{d/K},\xi_{LT,n}^{*,\infty}(\FC_\bullet)) \To{} \Gamma_c(\mlt^{d/K},\FC_\bullet)$
induit un isomorphisme $\dd\times W_K$-équivariant
\ini
\begin{equation}
\label{inv}
 \Gamma_c(\mltn^{d/K},\xi_{LT,n}^{*,\infty}(\FC_\bullet)) \simto \Gamma_c(\mlt^{d/K},\FC_\bullet)^{H_n}. 
 \end{equation}
En effet, cela montrera d'abord  l'existence de l'action de $\HC_\Lambda(GD,H_n)$ sur le terme de gauche et permettra donc de définir $R\Gamma_c(\mltn^{d/K},\Lambda)$ comme un objet $W_K$-équivariant de $D^b( \HC_\Lambda(GD,H_n))$. Puis cela montrera aussi les autres assertions du lemme, puisque le foncteur des points fixes sous $H_n$ est bien-sûr exact.

Pour prouver \ref{inv}, rappelons que  $\Gamma_c(\mlt^{d/K},\FC_\bullet)=\limi{m\in \NM}
\Gamma_c(\MC_{LT,m}^{d/K},\xi_{LT,m}^{*,\infty}(\FC_\bullet)) $ 
et que  l'action du groupe $GL_d(\OC_K)$ 
est simplement donnée par la structure de pro-revêtement étale galoisien de groupe $GL_d(\OC_K)$ du système $(\MC_{LT,m}^{d/K})_{m\in\NM}$ 
au-dessus de $\MC_{LT,0}^{d/K}$. Le groupe $GL_d(\OC_K)$ agit donc sur chaque 
$\Gamma_c(\MC_{LT,m}^{d/K},\xi_{LT,m}^{*,\infty}(\FC_\bullet))$ et il s'agit de montrer que pour $m\geq n$, l'application 
$\Gamma(\mltn^{d/K},\xi_{LT,n}^{*,\infty}(\FC_\bullet)) \To{} \Gamma(\MC_{LT,m}^{d/K},\xi_{LT,m}^{*,\infty}(\FC_\bullet))$ induit un isomorphisme 
$$ \Gamma_c(\mltn^{d/K},\xi_{LT,n}^{*,\infty}(\FC_\bullet)) \simto \Gamma_c(\MC_{LT,m}^{d/K},\xi_{LT,m}^{*,\infty}(\FC_\bullet))^{H_n}.$$
Rappelons que le morphisme 
$\pi_{m,n}:\;\MC_{LT,m}^{d/K}\To{} \MC_{LT,n}^{d/K}$ est étale galoisien de groupe $H_n/H_m$ et que par le lemme \ref{lemmelimpinf}, on a un isomorphisme canonique
$\pi_{m,n}^*(\xi_{LT,n}^{*,\infty}(\FC_\bullet))\simto \xi_{LT,m}^{*,\infty}(\FC_\bullet)$.
Ainsi l'isomorphisme cherch\'e est tautologique si on enl\`eve l'indice $c$ (supports compacts). On en d\'eduit l'isomorphisme avec indice $c$, par propret\'e et surjectivit\'e de $\pi_{m,n}$.

 \end{proof}

\alin{Propri{\'e}t{\'e}s de finitude} Dans ce paragraphe,
on rappelle une propriété de finitude de la cohomologie dans le cas $LT$ et on en tire une conséquence sur le $R\Gamma_c$ qui sera importante plus tard pour l'étude de l'action de l'inertie. La propriété analogue dans le cas $Dr$ demeure un
myst{\`e}re  si on en cherche une explication directe ne
faisant pas appel {\`a} la comparaison de Faltings.


\begin{fact} \label{admicohoent}
  Les $\Lambda$-modules $H^q_c(\o\MC_{LT}^{d/K},\Lambda)$, $q\in \NM$, 
  sont {\em $\o{GD}$-admissibles} (et m\^eme $G_d$-admissibles).
\end{fact}
\begin{proof}
 Rappelons que 
  $H_n=1+\varpi^nM_d(\OC)$. Par une variante évidente du lemme \ref{inv}, on a un isomorphisme $H^q_c(\mlt^{d/K}/\varpi^\ZM,\Lambda)^{H_n} \simeq H^q_c(\mltn^{d/K}/\varpi^\ZM,\Lambda)$.  Il nous faut donc voir que pour tout
  $n\in\NM$, le $\Lambda$-module $H^q_c(\mltn^{d/K}/\varpi^\ZM,\Lambda)$ est de type
  fini, ou de manière équivalente que $H^q_c(\mltn^{d/K,(0)},\Lambda)$ est de type fini. 
Pour cela, rappelons que $\mltn^{(0)}$ est la fibre g{\'e}n{\'e}rique d'un
  sch{\'e}ma formel $\wh\MC_{LT,n}^{(0)}$. On sait par diverses versions du th\'eor\`eme de Serre-Tate, {\em cf} \cite[II.2.7]{HaTay} et \cite[7.4.4]{Boyer1}, que ce dernier est
  isomorphe au compl{\'e}t{\'e} formel d'un sch{\'e}ma alg{\'e}brique $S_n$ propre sur $\OC^{nr}$
  (une certaine vari{\'e}t{\'e} de Shimura ou une vari{\'e}t{\'e} de
  ``Drinfeld-Stuhler'') en un point ferm{\'e} $x_n$ de sa fibre sp{\'e}ciale. Notons
  $S_n^{an}$  le $\ka$-espace analytique associ{\'e} (par GAGA ou par fibre
  g{\'e}n{\'e}rique de la compl{\'e}tion formelle $\wh{S_n}$ le long de la fibre sp{\'e}ciale :
  c'est la m{\^e}me chose car $S_n$ est propre). On a un plongement
  canonique $\mltn^{(0)} \injo S_n^{an}$ qui identifie $\mltn^{(0)}$ avec l'image
  r{\'e}ciproque de $x_n$ par le morphisme de sp{\'e}cialisation $S_n^{an}
  \To{} \o{S_n}$. Cette image r{\'e}ciproque est ouverte dans $S_n^{an}$
  et son compl{\'e}mentaire s'identifie {\`a} la fibre g{\'e}n{\'e}rique du sous-sch{\'e}ma
  formel ouvert de $\wh{S_n}$ d{\'e}fini par $\o{S_n}\setminus \{x\}$. Ce
  compl{\'e}mentaire est donc un domaine analytique de $S^{an}_n$.  Comme
  $S^{an}_n$ est compact, il s'ensuit que $\mltn^{(0)}$ est un ouvert {\em
    distingu{\'e}} de $S^{an}_n$ (diff\'erence de deux espaces
  compacts). 
Dans le cas de torsion, la cohomologie d'un tel ouvert distingu\'e est
de type fini par la suite exacte associ\'ee \`a la diff\'erence des
deux espaces compacts et \cite[cor. 5.6]{Bic3}. Le cas $\Lambda$ complet (et $\o\QM_l$) est expliqu\'e dans
 \cite[4.1.15 et  4.1.17]{Fargues}.

{\em Remarque :} bien que cela n'apparaisse pas clairement (c'est
cach{\'e} dans les r{\'e}f{\'e}rences {\`a} \cite{Bic3} et \cite{Fargues}), tout cela repose sur la
comparaison des cycles {\'e}vanescents de Berkovich pour les sch{\'e}mas
formels  et les cycles {\'e}vanescents alg{\'e}briques.

\end{proof}

\begin{nota} \label{notaidemp} On notera
  $\varepsilon_{n,GD}:=\varepsilon_{n,G}\varepsilon_{n,D}=
  \varepsilon_{n,D}\varepsilon_{n,G} \in \ZG_\Lambda(GD)$ et $\Mo{\Lambda}{GD}_n$ la
  sous-cat\'egorie facteur direct  de $\Mo{\Lambda}{GD}$ associ\'ee
  \`a cet idempotent.  
\end{nota}

\begin{coro} \label{finit} 
Le $\Lambda$-module
$\endo{D^b(\o{GD})}{\varepsilon_{n,GD}R\Gamma_c(\o\MC_{LT}^{d/K},\Lambda)}$
est de type fini. 
\end{coro}

\begin{proof}
D'apr\`es le lemme \ref{niveauLT} et le fait rappel\'e ci-dessus, 
il suffit de prouver le r\'esultat suivant, de nature "th\'eorie des repr\'esentations" :

\medskip

\noindent{\em Soit $\Lambda$ un anneau noeth\'erien o{\`u} $p$ est inversible.
  Soit $A^\bullet, B^\bullet \in D^b_\Lambda(\o{GD})$ deux complexes {\`a}
  cohomologie {\em admissible et de niveau fini}. 
  Alors  $\hom{A^\bullet}{B^\bullet}{D^b_\Lambda(\o{GD})}$ est un
  $\Lambda$-module de type fini.}

\medskip

En utilisant les suites spectrales habituelles, on constate que
$\hom{A^\bullet}{B^\bullet}{D^b_\Lambda(\o{GD})}$ est un sous-quotient
de
$$ \bigoplus_{q,k\in\ZM}
\ext{q-k}{\HC^q(A^\bullet)}{\HC^k(B^\bullet)}{\Lambda\o{GD}} $$
et il suffit donc de montrer que si $V$ et $W$ sont deux
$\Lambda\o{GD}$-modules lisses admissibles et de niveau $\leq n\in
\NM$, alors pour tout $i\in\NM$, 
$\ext{i}{V}{W}{\Lambda\o{GD}}$ est de type fini sur $\Lambda$.
Gr{\^a}ce {\`a} la suite exacte 
$$1\To{} G_d/\varpi^\ZM \To{} \o{GD} \To{} \dd/K^\times \To{} 1$$
et la compacit{\'e} de $\dd/K^\times$, il suffit de prouver que
les $\Lambda$-modules $\ext{i}{V}{W}{\Lambda \o{G_d}}$ 
sont de type fini.
Comme  le $\Lambda$-module $\hom{P}{W}{\o{G_d}}$ est de type fini d\`es que $P$ est projectif de type fini dans $\Mo{\Lambda}{\o{G_d}}$,
il suffit  de prouver que $V$ admet une r\'esolution par des objets projectifs de type fini.
Nous prouvons ceci dans l'appendice, {\em cf} \ref{corofinit}, suivant des arguments de Schneider-Stuhler \cite{SS2} prolong\'es par Vign{\'e}ras dans
\cite{Vigsheaves}. 
\end{proof}

\alin{Lien avec les cycles \'evanescents, dans le cas $LT$}
Dans \cite{CarAA}, \cite{Boyer1}, \cite{Boyer2} et \cite{HaTay}, les
auteurs non-seulement consid\`erent plut\^ot la ``sous-tour''
$(\MC_{LT,n}^{d/K,(0)})_{n\in\NM}$ comme dans la premi\`ere variante
\ref{variante}, mais consid\`erent aussi les cycles \'evanescents
plut\^ot que la cohomologie \`a supports compacts. Ils consid\`erent
donc les $\Lambda$-modules
$$  R^q\Psi_\eta(\wh\MC^{(0)}_{LT},\Lambda):= \limi{n\in\NM}
R^q\Psi_\eta(\wh\MC^{(0)}_{LT,n},\Lambda) $$     
qu'ils munissent d'une action de $(G_d\times\dd\times W_K)^0$.
Les cycles évanescents sont ceux que Berkovich associe aux schémas
formels localement formellement de type fini  dans  \cite{Bic4}. Ce
sont en général des faisceaux étales sur la fibre spéciale géométrique
d'un tel schéma formel. Dans le cas présent, la fibre spéciale étant
un point, ce sont seulement des $\Lambda$-modules. 

Nous devons pour la suite les comparer aux $\Lambda$-modules
$H^q_c(\MC^{d/K, (0)}_{LT},\Lambda)$ de ce texte.

\begin{lemme} \label{comparcyclesev}
Supposons que $\Lambda=\ZM/n\ZM$ avec $(n,p)=1$ ou $\Lambda=\o\QM_l$.
  Il existe des isomorphismes $(G_d\times \dd\times
  W_K)^0$-\'equivariants $$H^q_c(\MC^{d/K,
    (0)}_{LT},\Lambda)\simto
  R^{2d-2-q}\Psi_\eta(\wh\MC^{(0)}_{LT,n},\Lambda)^\vee\otimes|-|^{1-d}$$ 
o\`u le $^\vee$ d\'esigne le $\Lambda$-module lisse contragr\'edient.
\end{lemme}

\begin{proof}
La preuve repose essentiellement sur les faits suivants : les cycles
\'evanescents co\"incident avec la cohomologie sans supports, et 
celle-ci est duale de la cohomologie \`a supports compacts. 
La seule chose qui demandera une v\'erification est la compatibilit\'e
avec l'action de $G_d$. Il y a aussi une
complication technique dans le cas $l$-adique, $\Lambda=\o\QM_l$, o\`u
la cohomologie sans supports n'a pas \'et\'e d\'efinie en g\'en\'eral,
mais dans le cas pr\'esent, une d\'efinition ad hoc et un formalisme
de dualit\'e de Poincar\'e sont donn\'es dans \cite[5.9]{Fargues}.

{\em Premi\`ere \'etape :} Par d\'efinition des cycles \'evanescents,
lorsque $\Lambda$ est un anneau de torsion,
on a d'apr\`es \cite[Cor 2.3.ii)]{Bic4} des isomorphismes canoniques
$$  H^q(\MC_{LT,n}^{d/K, (0)},\Lambda) \simto
R^{q}\Psi_\eta(\wh\MC^{(0)}_{LT},\Lambda)$$
compatibles aux applications de transitions $\pi_{n,m}^*$, pour 
$n\leq m\in\NM$.
Dans le cas $\Lambda=\o\QM_l$, de tels isomorphismes existent encore
en d\'efinissant la cohomologie (sans supports) convenablement, {\em
  cf} \cite[5.9.1,5.9.4]{Fargues}.

{\em Deuxi\`eme \'etape :} Lorsque $\Lambda=\ZM/n\ZM$ avec $(n,p)=1$, le
formalisme de la dualit\'e de Poincar\'e tel qu'il est d\'evelopp\'e
dans \cite[7.3-7.4]{Bic2} fournit des accouplements parfaits pour
$0\leq q\leq 2d-2$
$$ \la\, ,\, \ra_n :\;\; H^q_c(\MC_{LT,n}^{d/K,(0)},\Lambda) \times
H^{2d-2-q}(\MC_{LT,n}^{d/K,(0)},\Lambda) \To{} \Lambda(1-d) $$ qui v\'erifient
la compatibilit\'e suivante aux morphismes de transition pour $n\leq
m$ :
$$ \forall \alpha\in H^q_c(\MC_{LT,m}^{d/K,(0)},\Lambda), \forall \beta
\in H^q(\MC_{LT,n}^{d/K,(0)},\Lambda),\;\; \la \alpha, \pi_{m,n}^* \beta
\ra_m = \la \pi_{m,n,!}\alpha,\beta\ra_n. $$
Modifions alors ces accouplements en posant pour tout $n\in\NM$ 
$$ (\alpha,\beta)_n:= [H_1:H_n]^{-1} \la \alpha,\beta\ra_n $$
qui est bien d\'efini puisque $H_1$ et $H_n$ sont des pro-$p$-groupes
et $p$ est inversible dans $\Lambda$.
En appliquant la compatibilit\'e ci-dessus et en tenant compte de ce
que l'endomorphisme $\pi_{m,n,!}\circ \pi_{m,n}^*$ de
$H^q_c(\MC_{LT,n}^{d/K,(0)},\Lambda)$ est la multiplication par
$[H_n:H_m]$, on v\'erifie que pour tous $n\leq m$, on a
$$ \forall \alpha\in H^q_c(\MC_{LT,n}^{d/K,(0)},\Lambda), \forall \beta
\in H^q(\MC_{LT,n}^{d/K,(0)},\Lambda),\;\; (\pi_{m,n}^*\alpha, \pi_{m,n}^* \beta
)_m = ( \alpha,\beta)_n. $$
En d'autres termes, les accouplements $(\,,\,)_n$ sont compatibles aux morphismes
de transition et d\'efinissent un accouplement entre les limites
inductives :
$$ (\,,\,) :\;\;  H^q_c(\MC_{LT}^{d/K,(0)},\Lambda) \times
H^{2d-2-q}(\MC_{LT}^{d/K,(0)},\Lambda) \To{} \Lambda(1-d). $$ 

Dans le cas $l$-adique, la dualité de Poincaré nécessaire est donnée par \cite[5.9.2]{Fargues}. 

{\em Troisi\`eme \'etape :} Il faut v\'erifier la compatibilit\'e \`a
l'action de $(G_d\times \dd \times W_K)^0$. Pour cela, il est plus
commode de consid\'erer le probl\`eme analogue pour l'espace
$\o\MC_{LT}^{d/K}$, muni de
 l'action de $G_d\times \dd \times W_K$ (la discussion précédente sur la dualité s'y applique verbatim). En effet dans ce cas, le groupe $\dd\times W_K$ agit
 sur chaque $\o\MC_{LT,n}^{d/K}$ et la dualit\'e de Poincar\'e est
 bien-s\^ur compatible aux automorphismes rationnels et l'action de
 Galois. Il suffit donc de consid\'erer l'action de $G$.
Pour cela fixons $g\in G$ et $m,n\in \NM$ tels que $m\geq n+d(\OC^d,g\OC^d)$, de sorte que le morphisme $g^{n|m} : \o\MC_{LT,m}^{d/K}\To{}\o\MC_{LT,n}^{d/K}$ soit bien défini. Rappelons qu'il s'agit d'un revêtement étale de degré $[H_n:gH_mg^{-1}]=[H_n:H_m]$.
On calcule alors pour tous $\alpha\in H^q(\o\MC_{LT,n}^{d/K},\Lambda), \beta\in H^q_c(\o\MC_{LT,n}^{d/K},\Lambda)$
\begin{eqnarray*}
	\left(g^{n|m,*}\alpha,g^{n|m,*}\beta\right)_m = [H_n:H_m]^{-1}\left(\alpha,g^{n|m}_!g^{n|m,*}\beta\right)_n 
	= \left(\alpha,\beta\right)_n 
\end{eqnarray*}

Ceci montre que $H^q(\o\MC_{LT}^{d/K},\Lambda)$ est la contragrédiente de $H^q_c(\o\MC_{LT}^{d/K},\Lambda)$. On en déduit le lemme.

\end{proof}

\subsection{Comparaison de Faltings}
Dans \cite{FaltDrin}, Faltings a esquiss\'e les arguments menant au th\'eor\`eme ci-dessous, au moins dans le cas o\`u $K$ est $p$-adique. Ces arguments ont \'et\'e compl\'et\'es et \'etendus au cas d'\'egales caract\'eristiques par L. Fargues, A. Genestier et V. Lafforgue, lors d'un groupe de travail \`a l'IHES s'\'etalant sur un semestre. Ils devraient \^etre prochainement \'ecrits puis, esp\'erons-le,  publi\'es. 
\begin{theo} (Faltings, Fargues, Genestier, Lafforgue) \label{faltings}
Il existe des isomorphismes $G_d\times \dd\times W_K$-\'equivariants (compte tenu de la normalisation des actions que nous avons choisie)
$$ H^q_c(\MC^{d/K}_{LT},\Lambda) \simto H^q_c(\MC^{d/K}_{Dr},\Lambda) $$
pour tout $q$ et tout anneau $\Lambda$ de torsion ou complet de valuation discr\`ete tel que $p\in \Lambda^\times$.
\end{theo}

Comme on l'a dit dans l'introduction, on esp\`ere que la r\'edaction pr\'ecise des arguments montrera aussi :

\begin{conj}
Il existe un isomorphisme $W_K$-\'equivariants dans $D^b_\Lambda({GD})$ 
$$ R\Gamma_c\left( \MC^{d/K}_{LT},\Lambda\right) \simto R\Gamma_c\left( \MC^{d/K}_{Dr},\Lambda\right) $$
pour tout anneau $\Lambda$ comme au-dessus.
\end{conj}

\subsection{Premi\`eres r\'eductions vers la conjecture \ref{conj}}

Dans cette section, on suppose $\Lambda=\o\QM_l$. 

\alin{\'Equivalence des deux points de la conjecture \ref{conj}} \label{reduc1}
Soit $\pi \in \Irr{\o\QM_l}{G_d}$ de caract\`ere central $\omega$ et
$\rho \in \Mo{\omega}{\dd}$, c'est-\`a-dire une repr\'esentation de
$\dd$ de caract\`ere central $\omega$.
On a une factorisation
$$ \hom{-}{\pi\otimes\rho^\vee}{GD}:\;\;\Mo{\o\QM_l}{GD}
\To{\hom{-}{\pi}{G_d}} \Mo{\omega}{\dd} 
\To{\hom{\rho}{-}{\dd}} \o\QM_l-\hbox{e.v.} $$
o\`u les deux derni\`eres cat\'egories sont semi-simples.
En particulier, on a un isomorphisme
$$  \hom{\rho}{\HC^*(R\hom{R\Gamma_c}{\pi}{G_d})}{\dd} \simto 
\HC^*(R\hom{R\Gamma_c}{\pi\otimes\rho^\vee}{GD}) $$
qui bien-s\^ur est $W_K$-\'equivariant (la d\'efinition de $\HC^*$ a
\'et\'e donn\'ee dans l'introduction).
Il est alors clair que dans la conjecture \ref{conj}, le point i)
implique le point ii). Compte tenu du fait que la famille de
foncteurs $\hom{\rho}{-}{\dd}$, $\rho \in \Irr{\omega}{\dd}$
est fid\`ele sur la cat\'egorie $\Mo{\omega}{\dd}$, la
r\'eciproque est encore vraie.

\alin{R\'eduction aux repr\'esentations de caract\`ere central trivial sur
  $\varpi$} \label{reduc2}
Le lemme suivant montre que pour \'etudier la conjecture \ref{conj},
on peut se ramener \`a l'\'etude du complexe
$R\Gamma_c(\o\MC^{d/K},\o\QM_l)$. Plus pr\'ecis\'ement il
montre qu'il nous suffira de montrer les points i) ou ii) de
\ref{conj} en y rempla{\c c}ant $\MC^{d/K}$ par
$\o\MC^{d/K}=\MC^{d/K}/\varpi^\ZM$ et en y supposant les
repr\'esentations $\pi$ et $\rho$  triviales sur l'uniformisante
$\varpi$.
Rappelons qu'\'etant donn\'e un caract\`ere $\psi$ de $\ZM$, on lui
associe des caract\`eres non ramifi\'es $\psi_{G}$, $\psi_{D}$, $\psi_{W}$
etc... des groupes $G_d$, $\dd$, $W_K$, etc..., cf \ref{notationsnu}.

\begin{lemme}
Soit $\pi\in\Irr{\o\QM_l}{G_d}$  et
$\rho\in \Irr{\o\QM_l}{\dd}$ de caract\`eres centraux ``inverses l'un de l'autre''
$\omega_\pi=\omega_\rho^{-1} :\; 
K^\times\To{} \o\QM_l^\times $.
Notons  $\psi : z\in\ZM \mapsto \omega_\pi(\varpi)^{z}\in\o\QM_l^*$.
Alors la repr\'esentation $\psi_{GD}(\pi\otimes\rho) =
(\psi_{G_d}^{-1}\pi)\otimes (\psi_{\dd} \rho)$ de $GD$ est triviale sur
$\varpi^\ZM$ et on a un isomorphisme \'equivariant
\begin{eqnarray*}
 & & \HC^*\left(R\hom{R\Gamma_c(\MC^{d/K},\o\QM_l)}{\pi\otimes
   \rho}{D^b_{\o\QM_l}(GD)}\right)  
 \\ &   & \simto_{W_K}  \psi_{W_K}.\HC^*\left(R\hom{R\Gamma_c(\o\MC^{d/K},\o\QM_l)}
{\psi_{GD}(\pi\otimes \rho)}{D^b_{\o\QM_l}(\o{GD})}\right) 
\end{eqnarray*}
\end{lemme}
\begin{proof} Pour alléger les notations,
on ommetra la notation $\HC^*$, les coefficients $\o\QM_l$, et le signe $D^b$ dans la preuve. 
Le lemme \ref{torsion} fournit le premier $W_K$-isomorphisme   
\begin{eqnarray*}
R\hom{R\Gamma_c(\MC^{d/K})}{\pi\otimes
   \rho}{GD} & \simto & \psi_{W_K}R\hom{\psi_{GD}^{-1}\otimes R\Gamma_c(\MC^{d/K})}{\pi\otimes
   \rho}{GD} \\ & \simto & \psi_{W_K}R\hom{R\Gamma_c(\MC^{d/K})}{\psi_{GD}(\pi\otimes
   \rho)}{GD} \\
       & \simto & \psi_{W_K}R\hom{R\Gamma_c(\o\MC^{d/K})}{\psi_{GD}(\pi\otimes
   \rho)}{\o{GD}}  
\end{eqnarray*}
et le troisième isomorphisme est donné par le lemme \ref{quotient} et la propriété d'adjonction du foncteur $\o\QM_l \otimes^L_{\o\QM_l[\varpi^\ZM]} -$.
  
\end{proof}

\alin{D\'ecomposition par les idempotents centraux primitifs} \label{pardecomp}
Nous avons d\'eja rappel{\'e} en \ref{decbernstein} la
d{\'e}composition par Bernstein de la cat{\'e}gorie ab{\'e}lienne $\Mo{\o\QM_l}{G_d}$ en
blocs. Une d{\'e}composition similaire mais beaucoup plus simple existe
pour $\Mo{\o\QM_l}{\dd}$. Notons $\Mo{\rho}{\dd}$ la sous-cat{\'e}gorie pleine de
$\Mo{\o\QM_l}{\dd}$ dont tous les sous-quotients irr{\'e}ductibles sont dans
l'orbite inertielle de la repr{\'e}sentation irr{\'e}ductible $\rho$ ({\em i.e.} isomorphes \`a} $\rho\psi$ pour
un caract{\`e}re non ramifi{\'e} $\psi$ de $\dd$). Alors on a la d{\'e}composition 
$$\Mo{\o\QM_l}{\dd} \simeq \prod_{\rho\in \Irr{\o\QM_l}{\dd}/\sim} \Mo{\rho}{\dd}$$
o{\`u} $\sim$ est l'{\'e}quivalence inertielle. Cette d{\'e}composition de la
cat{\'e}gorie  correspond {\`a} une {d{\'e}compo}sition de
son centre $\ZG_{\o\QM_l}(D)$: les idempotents primitifs de celui-ci sont donc en
bijection avec les classes d'inertie de repr{\'e}sentations
irr{\'e}ductibles. Notant $[\rho]$ l'idempotent associ\'e \`a la classe de $\rho$, et 
appliquant ceci au complexe de cohomologie de $\MC^{d/K}$, on obtient une d\'ecomposition
\ini\begin{equation} \label{rhoiso1}
R\Gamma_c(\MC^{d/K},\o\QM_l) \simeq \bigoplus_{\rho\in\Irr{\o\QM_l}{\dd}/\sim}
R\Gamma_c(\MC^{d/K},\o\QM_l)[\rho] \;\;\hbox{ dans } D^b_{\o\QM_l}({GD}).
\end{equation}

Maintenant, soit $\rho$ une $\o\QM_l$-repr\'esentation irr\'eductible du groupe
$\o\dd:=\dd/\varpi^\ZM$. Comme ce groupe est compact, $\rho$ d\'efinit un
idempotent primitif de $\ZG_{\o\QM_l}(\o\dd)$ que nous noterons
 $[\rho]$. Remarquons que si $\rho$ se factorise par $\dd/(1+\varpi^n\od)$, alors
$\varepsilon_{n,D}[\rho]=[\rho]$.
En d\'ecomposant la cat\'egorie $\Mo{\o\QM_l}{\o\dd}$ suivant
les composantes $\rho$-isotypiques, on 
 en d{\'e}duit la d{\'e}composition canonique
\ini\begin{equation} \label{rhoiso}
R\Gamma_c(\o\MC^{d/K},\o\QM_l) \simeq \bigoplus_{\rho\in\Irr{\o\QM_l}{\o\dd}}
R\Gamma_c(\o\MC^{d/K},\o\QM_l)[\rho] \;\;\hbox{ dans } D^b_{\o\QM_l}(\o{GD}).
\end{equation}
On prendra garde aux sens diff\'erents que prend le symbole $[\rho]$ dans les deux d\'ecompositions ci-dessus. Le lien entre ces deux d\'ecompositions est le suivant :
$$ \o\QM_l \otimes^L_{\o\QM_l[\varpi^\ZM]} \left(R\Gamma_c(\MC^{d/K},\o\QM_l)[\rho]\right) \simto \bigoplus_{\rho'\sim\rho,\omega_{\rho'}(\varpi)=1} R\Gamma_c(\o\MC^{d/K},\o\QM_l)[\rho'].
$$

Pour ce qui est de l'action du groupe $G_d$, on peut utiliser la
description de $\ZG_{\o\QM_l}(G_d)$ donn\'ee par Bernstein. Comme on l'a rappel{\'e} en
\ref{decbernstein}, les idempotents centraux primitifs sont dans ce
cas en bijection avec les
classes inertielles de paires Levi-cuspidale  $(M,\tau)$. Nous
noterons $[M,\tau]$ l'idempotent associ{\'e}. Nous dirons que $(M,\tau)$
est de niveau $\leq n$ si le bloc $\MC_{M,\tau}$ est ``inclus'' dans
$\Mo{\o\QM_l}{G_d}_n$, ou de mani{\`e}re {\'e}quivalente si $\varepsilon_{n,G}
[M,\tau]=[M,\tau]$ dans $\ZG_{\o\QM_l}(G_d)$.
La d{\'e}composition \ref{rhoiso} se raffine en une d{\'e}composition
$$R\Gamma_c(\o\MC^{d/K},\o\QM_l) \simeq
\bigoplus_{[\rho],[M,\tau]} 
R\Gamma_c(\o\MC^{d/K},\o\QM_l)[\rho][M,\tau] \;\;\hbox{ dans } D^b_{\o\QM_l}(\o{GD}).$$
En fait nous verrons que la description de Harris/Boyer implique que
le facteur associ{\'e} {\`a} $[\rho][M,\tau]$ est nul sauf si $[\rho]$ et
$[M,\tau]$ sont en ``correspondance de Jacquet-Langlands'' (dans un sens
assez {\'e}vident, en tout cas laiss{\'e} {\`a} la sagacit{\'e} du
lecteur).
Une preuve purement locale de ce fait est peut-{\^e}tre
envisageable en suivant \cite{Faltrace} ou \cite{Strauch}.

Fixons un idempotent primitif $[\rho]$ de $\ZG_{\o\QM_l}(\o\dd)$,
resp. $[M,\tau]$ de $\ZG_{\o\QM_l}( G_d)$ et supposons $\rho$ et $\tau$ de
niveau $\leq n$. 

\begin{prop} \label{continuite} L'image du
  morphisme 
$$W_K \To{\gamma} \endo{D^b_{\o\QM_l}(\o{GD})}{R\Gamma_c(\o\MC^{d/K},\o\QM_l)[\rho][M,\tau]}$$
donnant l'action de l'inertie est contenue dans un sous
$\ZM_l$-module de type fini.
\end{prop}

\begin{proof} 
Par d\'efinition des complexes de cohomologie, l'action de $W_K$ se
factorise
$$ W_K \To{} \endo{D^b_{\ZM_l}(\o{GD})}{R\Gamma_c(\o\MC^{d/K},\ZM_l)}
 \To{\otimes_{\ZM_l} \o\QM_l}  \endo{D^b_{\o\QM_l}
   (\o{GD})}{R\Gamma_c(\o\MC^{d/K},\o\QM_l)}.
$$
Mais comme $\rho$ et $\tau$ sont suppos{\'e}es de niveau $\leq n$,
l'action de $W_K$ sur le facteur associ{\'e} {\`a} $[\rho][M,\tau]$
se factorise plus pr\'ecis\'ement 
$$ W_K \To{}\endo{D^b_{\ZM_l}(\o{GD})}{\varepsilon_{n,GD}R\Gamma_c(\o\MC^{d/K},\ZM_l)}
 \To{\otimes_{\ZM_l} \o\QM_l}  \endo{D^b_{\o\QM_l}
   (\o{GD})}{R\Gamma_c(\o\MC^{d/K},\o\QM_l)[\rho][M,\tau]}.$$
La  proposition est alors une cons\'equence du corollaire \ref{finit} du c\^ot\'e $LT$, ainsi que du th\'eor\`eme \ref{faltings} du c\^ot\'e $Dr$.

\end{proof}

\section{R{\'e}alisation cohomologique des correspondances} \label{real}

\subsection{Description de la cohomologie}

\alin{Repr{\'e}sentations elliptiques ``cohomologiques''}
Soit $\rho \in \Irr{\o\QM_l}{\dd}$, on a d{\'e}fini et classifi{\'e} en \ref{classrepell} les
repr{\'e}sentations de $G_d$ elliptiques $\pi_\rho^I$, o{\`u} $I\subseteq S_{d_\rho}$,
associ{\'e}es {\`a} $\rho$. Parmi celles-ci seules certaines apparaissent dans
la cohomologie des espaces modulaires, et m\'eritent donc d'\^etre appel\'ees  
``cohomologiques''. Pour les d{\'e}crire on utilise la bijection
$S_{d_\rho}\simto \I{d_\rho-1}$ d{\'e}crite en \ref{rappelclassif} et déja utilisée en \ref{corlanell}. On pose alors
pour $0\leq i\leq d_\rho -1$ :
$$ I(i):=\I{i}\subseteq \I{d_\rho-1} \;\; \hbox{ et } \;\;
\pi^i_\rho:=\pi^{I(i)}_\rho .$$
Soulignons que pour $i=0$, la d{\'e}finition donne $I(i)=\emptyset$ et par
cons{\'e}quent $\pi^0_\rho=\pi^\emptyset_\rho=JL_d(\rho)$ est la s{\'e}rie
discr{\`e}te de $G_d$ associ{\'e}e {\`a} $\rho$ par la correspondance de Jacquet-Langlands.

\alin{Description de Harris-Boyer et variantes }
Nous donnons ici la description, conjecturée par Harris du côté $Dr$ (communication priv\'ee) et annoncée par Boyer du côté $LT$ (au moins sur un corps de fonctions) des divers  espaces de cohomologie que nous
avons d{\'e}finis dans la partie \ref{drinfeld}
pour les espaces analytiques $ \MC_?^{d/K}$ et $\o\MC_?^{d/K}$.

Commen{\c c}ons par quelques remarques générales : tout d'abord, les espaces 
$H^i_c(\MC^{d/K}_{?},\o\QM_l)$ ne sont non-nuls que pour
$d-1\leq i\leq 2d-2$. En effet les espaces analytiques 
$\MC^{d/K}_{?,n}$ sont quasi-affines au sens de \cite[5]{Bic4} et on peut appliquer \cite[6.2]{Bic4}.
 Ensuite 
on calcule facilement en plus haut degré
$$H_c^{2d-2}(\MC^{d/K},\o\QM_l)=(|-|^{1-d})_{un} \;\hbox{ et }\;
H_c^{2d-2}(\o\MC^{d/K},\o\QM_l)=\bigoplus_{\psi^d=1} \psi_{GDW}\otimes
|-|^{1-d}$$ 
o{\`u} on utilise la notation $un$  introduite au paragraphe \ref{defun}
et on a \'etendu trivialement \`a $G_d\times \dd$ 
le caract{\`e}re non ramifi\'e $|-|$ de $W_K$ qui envoie un Frobenius g\'eom\'etrique sur $q^{-1}$, puis la notation
$\psi_H$ introduite en \ref{notationsnu} pour tout caractère
$\ZM\To{}\o\QM_l^\times$ (ici d'ordre $d$). 


Jusqu'{\`a} pr{\'e}sent, les seules informations que l'on avait en toute
g{\'e}n{\'e}ralit{\'e} concernaient la 
contribution de la partie cuspidale de $G$ (qui n'appara{\^\i}t que pour
$i=d-1$ dans le cas $Dr$, voir Harris \cite{HaCusp} et Hausberger
\cite{Hausb}, et dans  le cas $LT$, voir la remarque de Carayol
apr{\`e}s \cite[5.6]{CarBou}) ainsi que la partie
associ{\'e}e {\`a} $\mdro$, calcul{\'e}e par Schneider et Stuhler dans \cite{SS1}, {\em cf} partie \ref{dp}.


Voici la description promise. Par prudence nous l'\'enon{\c c}ons sous la forme d'une conjecture :

\begin{conj} \label{conjHarris}
  Soit $i\geq 0$ et $n\geq 0$, alors on a un isomorphisme de $G_d\times \dd\times W_K$-modules 
$$ H^{d-1+i}_c(\MC^{d/K},\o\QM_l) \simeq \bigoplus_{\rho\in \Irr{\o\QM_l}{\dd}/\sim}
\left(\pi^i_{\rho^\vee}\otimes\rho\otimes \sigma_{d/d_\rho}(\tau_\rho^0)|-|^{{{d_\rho-d}\over 2}-i}\right)_{un}$$
o{\`u} la somme est {\'e}tendue aux repr{\'e}sentations irr{\'e}ductibles de
$\dd$ prises {\`a} torsion pr{\`e}s par un caract{\`e}re 
non-ramifi\'e 
et la notation $un$ a {\'e}t{\'e} d{\'e}finie en \ref{defun}, tandis que les notations $d_\rho$ et $\tau_\rho^0$ sont celles de \ref{notationsrho}.
En particulier, on a aussi un isomorphisme 
$$ H^{d-1+i}_c(\o\MC^{d/K},\o\QM_l) \simeq \bigoplus_{\rho\in \Irr{\o\QM_l}{\dd}, \omega_\rho(\varpi)=1}
\pi^i_{\rho^\vee}\otimes\rho\otimes \sigma_{d/d_\rho}(\tau_\rho^0)|-|^{{{d_\rho-d}\over 2}-i}.$$

\end{conj}

L'{\'e}nonc{\'e} suivant est une cons{\'e}quence des r{\'e}sultats annonc{\'e}s par Pascal 
Boyer dans \cite{Boyer2} :
\begin{theo}(Boyer) La conjecture \ref{conjHarris} est vraie dans le
  cas Lubin-Tate $\MC^{d/K}=\MC_{LT}^{d/K}$ lorsque $K$ est d'{\'e}gales caract{\'e}ristiques.
\end{theo}
\begin{proof}
  L'{\'e}nonc{\'e} pr{\'e}cis que l'on trouve dans \cite{Boyer2} concerne les
  cycles \'evanescents  de $\wh\MC_{LT}$ o{\`u} l'on a fix{\'e} un caract{\`e}re
  central. Nous devons expliquer ici comment cet {\'e}nonc{\'e} entraine
  formellement celui de \ref{conjHarris}.
    Pour commencer, comparons les notations de {\em loc.cit} avec les notres :
  on a les correspondances suivantes 
\def\llr{\longleftrightarrow}
  $$ \tau_\rho^0 \llr \pi_0, \,\,\,\,\rho \llr JL(St_s(\pi_0)),\,\,\, 
    \sigma_{d/d_\rho}(\pi_0)\llr L_g(\pi_0),\;\;\;\; 
    d_\rho \longleftrightarrow s $$
    et en ce qui concerne les repr\'esentations elliptiques,
$$    \pi^i_\rho =\pi_\rho^{\{1,\cdots,i\}} \llr 
    [\overrightarrow{i},\overleftarrow{s-i-1}]_{\pi_0} \;\;\hbox{ et }\;\;
    \pi_\rho^{\{i+1,\cdots, d_\rho-1\}} \llr [\overleftarrow{i},\overrightarrow{s-i-1}]_{\pi_0}
$$ 

Notons temporairement $H^i_c:=H^i_c(\mlt^{d/K},\o\QM_l)$, puis
$h^i_c:=H^i_c(\mlt^{d/K,(0)},\o\QM_l)$ et enfin $\psi^i:=R^i\Psi_\eta(\wh\MC_{LT}^{(0)},\o\QM_l)$. Rappelons qu'on a 
$H^i_c \simeq
\cind{(G_d\times \dd\times W_K)^0}{G_d\times \dd\times W_K}{h^i_c}$, voir \ref{variante}, et que $h^i_c\simeq (\psi^{2d-2-i})^\vee \otimes |-|^{-d+1}$ 
en notant avec un $ ^\vee$ la 
 contragr\'ediente de $\psi^{2d-2-i}$, par \ref{comparcyclesev}.

Soit $\rho\in \Irr{\o\QM_l}{\dd}$, un des r\'esultats principaux de \cite{Boyer2} est :
$$ \hom{\rho}{\psi^{d-d_\rho+i}}{\OC_D^\times} \mathop{\simeq}\limits_{G_d\times W_K} \pi_{\rho^\vee}^{\{i+1,\cdots, d_\rho-1 \}} \otimes \sigma_{d/d_\rho}(\tau_\rho^0)|-|^{-(d-d_\rho+2i)/2}. $$ 
Rappelons comment le terme de gauche fournit une repr\'esentation de $G_d\times W_K$ : on peut prolonger $\rho$ trivialement au produit triple puis la restreindre au sous-groupe $(G_d\times \dd\times W_K)^0$, on utilise alors le fait que $\OC^\times_D$ est un sous-groupe normal de ce dernier et que  la projection canonique induit un isomorphisme $(G_d\times \dd\times W_K)^0/\OC_D^\times\simto G_d\times W_K$. Soulignons aussi que Boyer normalise l'action de $G_d$ de la m\^eme mani\`ere que dans ce texte (la torsion par $g\mapsto {^tg}^{-1}$ du c\^ot\'e $Dr$ correspond \`a l'action naturelle du c\^ot\'e $LT$, via l'isomorphisme de Faltings).

Par dualit\'e et en posant $j=d_\rho-1-i$, on obtient
$$\hom{h^{d-1+j}_c}{\rho^\vee}{\OC_D^\times}\mathop{\simeq}\limits_{G_d\times W_K} 
\pi_{\rho^\vee}^{\{d_\rho-j,\cdots, d_\rho-1 \}} \otimes \sigma_{d/d_\rho}(\tau_{\rho}^0)|-|^{\frac{d-d_\rho}{2}+j}.$$ 

Par ailleurs, si $\zeta$ est un caract\`ere de $K^\times$, la partie $\zeta_{|\OC_K^\times}$-covariante $(h^i_c)_{\zeta_{\OC_K^\times}}$, {\em cf} \ref{covariants}, se prolonge en une repr\'esentation $(h^i_c)_\zeta$ du groupe
$$(\nu_{GDW})^{-1}(d\ZM)= (1\times K^\times \times 1)(G_d\times\dd\times W_K)^0$$  de $G_d\times\dd\times W_K$ en posant $(h^i_{c})_\zeta(k):=\zeta(k)$ pour tout $k\in K^\times$.
On a alors 
$$ (H^i_c)_\zeta \simeq \cind{(\nu_{GDW})^{-1}(d\ZM)}{G_d\times\dd\times W_K}{(h^i_c)_\zeta},$$
et par suite,
lorsque $\zeta=\omega_\rho$, 
$$  \hom{h^{i}_c}{\rho}{\OC_D^\times} \mathop{\simeq}\limits_{G_d\times W_K} \hom{(h^i_{c})_\zeta}{\rho}{\OC^\times_D\varpi^\ZM} \mathop{\simeq}\limits_{G_d\times W_K}
\hom{(H^i_c)_\zeta}{\rho}{\dd},
$$
le second isomorphisme est une version \'equivariante de la r\'eciprocit\'e de Frobenius.
La d\'ecomposition de la repr\'esentation $(H_c^i)_\zeta$  selon les idempotents centraux primitifs de $\Mo{\zeta}{\dd}$ peut alors s'\'ecrire sous la forme
\begin{eqnarray*} (H^{d-1+j}_c)_\zeta & \simeq & \bigoplus_{\rho\in \Irr\zeta\dd}  \rho\otimes \hom{(H^{d-1+j}_c)_\zeta}{\rho}{\dd}^\vee \\
& \simeq & \bigoplus_{\rho\in \Irr\zeta\dd}  \rho\otimes\pi^j_{\rho^\vee}\otimes \sigma_{d/d_\rho}(\tau_{\rho}^0)|-|^{\frac{d_\rho-d}{2}-j}. 
\end{eqnarray*}
On a utilis\'e le calcul de contragr\'ediente 
$(\pi_{\rho}^{\{d_\rho-j,\cdots, d_\rho-1 \}})^\vee\simeq \pi_{\rho^\vee}^j$ de \ref{dualiteell}.


Rappelons maintenant qu'{\`a} (la classe d'inertie de) $\rho$
est attach{\'e} un idempotent primitif $[\rho]$ de
$\ZG_{\o\QM_l}(\dd)$. 
Faisant agir ce dernier sur la d\'ecomposition pr\'ec\'edente, on obtient
un isomorphisme
$$ H^{d-1+j}_c[\rho]_\zeta \simto \bigoplus_{\buildrel{\rho'\sim
    \rho}\over{{\omega_{\rho'}}_{|K^\times} = \zeta}} \pi^j_{{\rho'}^\vee}\otimes \rho'
\otimes \sigma_{d/d_\rho}(\tau_{\rho'}^0) |-|^{\frac{d_\rho-d}{2}-j}.$$
Mais d'apr{\`e}s la compatibilit{\'e} des correspondances de Langlands et
Jacquet-Langlands {\`a} la torsion par les caract{\`e}res non ramifi{\'e}s,
le terme de droite est  encore isomomorphe {\`a}
$$ \bigoplus_{\buildrel{\tau \sim
    \pi^i_{\rho^\vee}\otimes \rho\otimes
    \sigma_{d/d_\rho}(\tau_\rho^0)}\over{{\omega_{\tau}}_{|1\times K^\times\times 1} = \zeta}} \tau$$
D'apr{\`e}s le lemme \ref{centun}, cela implique le premier isomorphisme de
\ref{conjHarris}, puis le second, compte tenu de \ref{quotient}.

\end{proof}

\begin{coro}(Faltings-Fargues-Genestier-Lafforgue) La conjecture \ref{conjHarris} est vraie dans le cas Drinfeld $\MC^{d/K}=\MC_{Dr}^{d/K}$ lorsque $K$ est d'{\'e}gales caract{\'e}ristiques.
\end{coro}

C'est en effet une cons\'equence du th\'eor\`eme de Boyer et du th\'eor\`eme \ref{faltings}.

\subsection{Scindage du $R\Gamma_c$}
\`A partir de maintenant, nous supposerons {\em toujours} que la cohomologie est d\'ecrite selon la conjecture \ref{conjHarris}. Nous ne le rappellerons pas syst\'ematiquement.

\begin{lemme} \label{corodesc}
Soit $\pi \in \Irr{\o\QM_l}{G_d}$. 
\begin{enumerate}
        \item Si $\pi$ n'est pas elliptique, alors $R\hom{R\Gamma_c(\MC^{d/K},\o\QM_l)}{\pi}{D^b(G_d)}=0$.
        \item Si $\pi=\pi_\rho^I$ pour $\rho\in\Irr{\o\QM_l}{\o{\dd}}$ et $I\subseteq S_{d_\rho}$, alors 
        $$ R\hom{R\Gamma_c(\MC^{d/K},\o\QM_l)}{\pi_\rho^I}{D^b(G_d)} \simeq R\hom{R\Gamma_c(\o\MC^{d/K},\o\QM_l)[\rho^\vee]}{\pi_\rho^I}{D^b(G_d)}$$
        avec la notation de \ref{rhoiso}.
\end{enumerate}

\end{lemme}

\begin{proof}
D'apr\`es la suite spectrale
$$ \ext{p}{H^q_c(\MC^{d/K},\o\QM_l)}{\pi}{G_d} \Rightarrow R^{p-q}\hom{R\Gamma_c(\MC^{d/K},\o\QM_l)}{\pi}{D^b(G_d)}, $$
la description de \ref{conjHarris}, le lemme \ref{lemmeun} et la remarque \ref{equivdefun}, 
le complexe $R\hom{R\Gamma_c}{\pi}{G_d}$ est non nul seulement s'il existe $\rho\in\Irr{\o\QM_l}{\dd}$ telle que 
$\ext{*}{\pi^*_{\rho,un}}{\pi}{G_d}\neq 0$. Ceci implique alors que $\pi$ est dans le m\^eme bloc de Bernstein que $\pi_\rho^*$. En appliquant l'\'equivalence de cat\'egorie $\alpha_\rho$ de \ref{equivcat}, on voit que $\ext{*}{\pi^{K_\rho}_{d_\rho,*}}{\alpha_\rho(\pi)}{G^0_{d_\rho}}\neq 0$, ce qui ne peut se produire que si $\pi$ a le m\^eme support cuspidal que $\pi^*_\rho$, c'est-\`a-dire si $\pi$ est elliptique.

Supposons maintenant $\pi=\pi_\rho^I$ avec $\omega_\rho(\varpi)=1$. D'apr\`es \ref{quotient} et \ref{rhoiso} on a
\begin{eqnarray*}
R\hom{R\Gamma_c(\MC^{d/K},\o\QM_l)}{\pi_\rho^I}{D^b(G_d)}  &\simeq & R\hom{R\Gamma_c(\o\MC^{d/K},\o\QM_l)}{\pi_\rho^I}{D^b(G_d)} \\
& \simeq & \bigoplus_{\rho'\in\Irr{\o\QM_l}{\o{\dd}}} R\hom{R\Gamma_c(\o\MC^{d/K},\o\QM_l)[\rho']}{\pi_\rho^I}{D^b(G_d)}.
\end{eqnarray*}
Mais d'apr\`es la description \ref{conjHarris}, le m\^eme raisonnement que ci-dessus montre que le complexe  $$R\hom{R\Gamma_c(\o\MC^{d/K},\o\QM_l)[\rho']}{\pi_\rho^I}{D^b(G_d)}$$ est non-nul seulement si $\rho'\simeq \rho^\vee$.

\end{proof}

\begin{prop} (Prop. \ref{prop1}) \label{prop1bis}
Supposons valide la description de \ref{conjHarris}. 
 Alors il existe des isomorphismes (pas uniques) dans $D^b(GD)$, resp. dans $D^b(\o{GD})$,
 $$R\Gamma_c(\MC^{d/K},\o\QM_l)
    \simto \bigoplus_{i\in \NM} H^i_c(\MC^{d/K},\o\QM_l)[-i],$$
 $$\hbox{resp. } R\Gamma_c(\o\MC^{d/K},\o\QM_l)
    \simto \bigoplus_{i\in \NM} H^i_c(\o\MC^{d/K},\o\QM_l)[-i].$$
\end{prop}

\begin{proof} 
Nous traitons le cas de $\MC^{d/K}$ et laissons l'adaptation, \'evidente, au cas $\o\MC^{d/K}$ au lecteur.
La d\'ecomposition canonique 
\ref{rhoiso1}
nous ram\`ene \`a prouver que pour toute (classe d'inertie de) repr\'esentation $\rho\in\Irr{\o\QM_l}{\dd}$, le complexe $R\Gamma_c(\MC^{d/K},\o\QM_l)[\rho]$ est scindable.
Pour all\'eger les notations, nous noterons (dans cette preuve seulement) $R\Gamma_c[\rho]$ ce complexe et $H^i_c[\rho]:=H^i_c(\MC^{d/K},\o\QM_l)[\rho]$.
La conjecture \ref{conjHarris} et le lemme \ref{lemmeun} nous donnent alors la cohomologie de
$R\Gamma_c[\rho]$ :
$$ \HC^{d-1+i}(R\Gamma_c[\rho])= H^{d-1+i}_c[\rho] \simeq_{GD} 
\left\{\begin{array}{ll} (\pi^i_{\rho^\vee}\otimes
\rho)_{un} \otimes_{\o\QM_l} V_{\sigma(\tau_\rho^0)} & \hbox{si } 0\leq
i\leq d_\rho-1 \\
0& \hbox{sinon} \end{array}\right. 
$$

Pour montrer que $R\Gamma_c[\rho]$ est scind{\'e} nous voulons appliquer
le crit{\`e}re donn{\'e} par le lemme \ref{scindage}. Selon ce crit{\`e}re il nous
suffit de montrer que $\ext{k}{H^{d-1+i}_c[\rho]}{H^{d-1+j}_c[\rho]}{GD}$ est
nul d{\`e}s que $k=i-j+1$. 

 La suite exacte $1\To{} G_d \To{}  GD \To{} \dd/K^\times \To{} 1$ et la compacit\'e de $\dd/K^\times$ 
montrent que pour $V,W$ deux $\o\QM_l(GD)$-modules lisses, on a
$$ \ext{k}{V}{W}{GD} \simeq \left(\ext{k}{V}{W}{G_d}\right)^{\dd/K^\times}.
$$  Or pour $0\leq i\leq d_\rho-1$,
on a par le lemme \ref{lemmeun} et la remarque \ref{equivdefun} $$H^{d-1+i}_c[\rho] \simeq_{G_d} \pi^i_{\rho,un} \otimes_{\o\QM_l}
(V_\rho\otimes V_{\sigma(\tau_\rho^0)}).$$
On est donc amen\'e \`a calculer $\ext{k}{\pi^i_{\rho,un}}{\pi^j_{\rho,un}}{G_d}$.
Pour cela on utilise l'\'equivalence de cat\'egories $\alpha_\rho$ et les notations de \ref{equivcat}, puis le lemme de Shapiro qui nous donnent
\begin{eqnarray*}
\ext{k}{\pi^i_{\rho,un}}{\pi^j_{\rho,un}}{G_d} & \simeq & \ext{k}{\pi_{d_\rho,i}^{K_\rho}}{\cind{(G'_{d_\rho})^0}{G'_{d_\rho}}{\pi_{d_\rho,j}^{K_\rho}}}{(G'_{d_\rho})^0} \\
& \simeq  & \ext{k}{\pi_{d_\rho,i}^{K_\rho}}{\pi_{d_\rho,j}^{K_\rho}}{(G'_{d_\rho})^0}\otimes_{\o\QM_l} \o\QM_l[G'_{d_\rho}/(G'_{d_\rho})^0].
\end{eqnarray*}
Mais le groupe $G'_{d_\rho}/SL_{d_\rho}(K_\rho)$ est isomorphe \`a $\OC_{K_\rho}^\times$ donc est compact, de sorte que
$$ \ext{k}{\pi_{d_\rho,i}^{K_\rho}}{\pi_{d_\rho,j}^{K_\rho}}{(G'_{d_\rho})^0} \simeq
\ext{k}{\pi_{d_\rho,i}^{K_\rho}}{\pi_{d_\rho,j}^{K_\rho}}{SL_{d_\rho}(K_\rho)}
.$$
On peut alors appliquer \ref{theoext} i) qui nous dit que le $\hbox{Ext}^k$ est nul sauf si $k=|i-j|$. En
particulier, il est nul si $k=i-j+1$.
\end{proof}

Dor\'enavant, et motiv\'es par le paragraphe \ref{reduc2} et le corollaire \ref{corodesc}, nous nous restreindrons \`a l'\'etude du complexe
$R\Gamma_c(\o\MC^{d/K},\o\QM_l)$ de $D^b_{\o\QM_l}(\o{GD})$, et plus pr\'ecis\'ement de sa composante $R\Gamma_c(\o\MC^{d/K},\o\QM_l)[\rho]$ ({\em cf} la d\'ecomposition \ref{rhoiso}) o\`u $\rho\in\Irr{\o\QM_l}{\o{\dd}}$. Nous all\`egerons les notations en notant cette composante $R\Gamma_c[\rho]$ et sa cohomologie  $H^i_c[\rho]$. Il sera aussi plus agr\'eable de poser 
\ini\begin{equation}\label{defsigmaprime}
 \sigma'(\rho):=\sigma_{d/d_\rho}(\tau_\rho^0)|-|^{\frac{d_\rho-d}{2}} 
 \end{equation}
et nous noterons $V_{\sigma'(\rho)}$ l'espace de cette repr\'esentation de $W_K$.
Nous travaillons sous l'hypoth\`ese de validit\'e de \ref{conjHarris} et nous supposerons fix\'es des isomorphismes $G_d\times\dd\times W_K$-\'equivariants
\ini\begin{equation}\label{isocoho}  H^{d-1+i}_c[\rho] \simto \pi^i_{\rho^\vee}\otimes\rho\otimes\sigma'(\rho)|-|^{-i},\;\hbox{ pour } i=0,\cdots, d_\rho.
\end{equation}

Introduisons enfin l'objet $\CC_\rho$ de $D^b_{\o\QM_l}(\o{GD})$
$$ \CC_\rho := \bigoplus_{i=0}^{d_\rho} \left(\pi^i_{\rho^\vee}\otimes \rho\right)[-(d-1+i)]. $$
\begin{rema}\label{triv}
Lorsque $\DC$ est une cat\'egorie $R$-lin\'eaire, $R$ \'etant un anneau commutatif unitaire, et $L$ est un $R$-module libre de type fini, on a un (une classe d'isomorphisme d') endofoncteur $C\mapsto (C \otimes_R L)$ de $\DC$. On v\'erifie imm\'ediatement qu'on a un isomorphisme canonique de $R$-alg\`ebres
$$ \endo{\DC}{C\otimes L} \simto \endo{\DC}{C} \otimes_R \endo{R}{L} . $$
\end{rema}

Par la suite, nous appellerons {\em scindage} de $R\Gamma_c[\rho]$ tout isomorphisme
\ini\begin{equation}\label{defscin}
 \alpha:\;\; R\Gamma_c[\rho] \simto \CC_\rho\otimes_{\o\QM_l} V_{\sigma'(\rho)} \;\;\hbox{ dans }\;\; D^b_{\o\QM_l}(\o{GD}).
 \end{equation}
 qui induit les isomorphismes $\ref{isocoho}$ sur les groupes de cohomologie ; c'est une l\'eg\`ere variante de la notion utilis\'ee dans la section \ref{sectionscin}.
Par la remarque ci-dessus, un tel scindage induit donc un isomorphisme
$$\alpha_*:\;\; \endo{D^b(\o{GD})}{R\Gamma_c[\rho]} \simto \endo{D^b(\o{GD})}{\CC_\rho} \otimes_{\o\QM_l} \endo{\o\QM_l}{V_{\sigma'(\rho)}}. $$
Par ailleurs, on a la description de l'alg\`ebre des endomorphismes de $\CC_\rho$ 
$$ \endo{D^b(\o{GD})}{\CC_\rho} = \bigoplus_{i\geq
  j} \ext{i-j}{\pi^i_{\rho^\vee}\otimes \rho}{\pi^j_{\rho^\vee}\otimes \rho}{\o{GD}}, $$ 
le produit sur l'espace de droite {\'e}tant donn{\'e} par le $\cup$-produit. Or, pour tous $i\geq j$, on a des isomorphismes compatibles aux $\cup$-produits
\begin{eqnarray*}
\ext{i-j}{\pi^i_{\rho^\vee}\otimes \rho}{\pi^j_{\rho^\vee}\otimes \rho}{\o{GD}}
& \simto  & \ext{i-j}{\pi^i_{\rho^\vee}}{\pi^j_{\rho^\vee}}{\o{G_d}}
\end{eqnarray*}
induits par la suite exacte $\o{\dd}\injo \o{GD} \To{} \o{G_d}$.
Nous avons d\'eja expliqu\'e en \ref{beta} comment la table de multiplication du $\cup$-produit de la proposition \ref{theoextell} permet de construire un (des) isomorphisme(s) de $\o\QM_l$-alg\`ebres  
$$ \beta :\;\; \TC_{d_\rho} \simto \bigoplus \ext{i-j}{\pi^i_{\rho^\vee}}{\pi^j_{\rho^\vee}}{\o{G_d}}, $$
(o\`u $\TC_{d_\rho}$ d\'esigne toujours l'alg\`ebre des matrices triangulaires sup\'erieures de taille $d_\rho$) qui par composition nous donne(nt) un (des) isomorphisme(s)
$ \beta :\;\; \TC_{d_\rho} \simto \endo{D^b(\o{GD})}{\CC_\rho}.$ 
 Comme dans \ref{beta}, soulignons :
\begin{rema}\label{choixgen}
La construction de $\beta$ d\'epend du choix de g\'en\'erateurs des $\o\QM_l$-droites $$\ext{1}{\pi^{i}_{\rho^\vee}}{\pi^{i-1}_{\rho^\vee}}{\o{G_d}}$$
pour $i=1,\cdots, d_\rho-1$. Changer ce choix de g\'en\'erateurs revient simplement \`a composer $\beta$ avec la conjugaison par une matrice diagonale.
\end{rema}
L'isomorphisme $\beta$ \'etant fix\'e, nous noterons par le m\^eme symbole l'isomorphisme d'alg\`ebres 
$$ \beta :\;\; \endo{\o\QM_l}{V_{\sigma'(\rho)}}\otimes_{\o\QM_l} \TC_{d_\rho} \simto \endo{D^b(\o{GD})}{\CC_\rho\otimes_{\o\QM_l} V_{\sigma'(\rho)}}.$$
que l'on en d\'eduit gr\^ace \`a la remarque \ref{triv}.
R\'ecapitulons tout cela :

\begin{prop} \label{propnoneq} 
\`A tout scindage $\alpha$ comme en \ref{defscin} et tout choix de g\'en\'erateurs 
\ref{choixgen} est associ\'e un
 isomorphisme de $\o\QM_l$-alg{\`e}bres $\beta^{-1}\circ \alpha_*$ s'inscrivant dans les diagrammes commutatifs 
$$\xymatrix{ \endo{D^b(\o{GD})}{R\Gamma_c[\rho]} \ar[r]^\sim_{\beta^{-1}\circ \alpha_*} \ar[d]^{can}
& \endo{\o\QM_l}{V_{\sigma'(\rho)}} \otimes_{\o\QM_l} \TC_{d_\rho}
\ar[d]^{\id\otimes E_{ii}} \\
 \endo{\o{GD}}{H^{d-1+i}_c[\rho]} \ar[r]^{(\ref{isocoho})}  &
 \endo{\o\QM_l}{V_{\sigma'(\rho)}} },$$
o\`u  $i\in\{0,\cdots, d_\rho\}$ et $E_{ii}$ d\'esigne la coordonn\'ee diagonale $(i,i)$ de $\TC_{d_\rho}$
 \end{prop}

Cela montre en particulier la proposition \ref{prop2} annonc\'ee dans l'introduction, et la pr\'ecise sensiblement.

\subsection{Action de $W_K$ sur $R\Gamma_c[\rho]$}
Le but de ce paragraphe est de comprendre le morphisme canonique
$$ \gamma_\rho : \;\; W_K \To{}\aut{D^b(\o{GD})}{R\Gamma_c[\rho]} $$
qui d{\'e}crit l'action de $W_K$ sur $R\Gamma_c[\rho]$.
\'Evidemment, nous travaillons toujours sous l'hypoth\`ese que {\em la cohomologie est d\'ecrite par  \ref{conjHarris}}. 
Le r{\'e}sultat final que nous visons est une version $W_K$-{\'e}quivariante
de la proposition \ref{propnoneq}. La premi\`ere \'etape rec\`ele encore une inconnue :

\begin{prop} \label{propeq}
Pour tout rel\`evement de Frobenius g\'eom\'etrique $\phi$, il existe un unique scindage $\alpha_\phi$ comme en \ref{defscin}, un unique sous-ensemble $I_\phi\subseteq S_{d_\rho}=\{1,\cdots,d_\rho-1\}$, et un isomorphisme $\beta_\phi$ associ\'e \`a un certain
choix de g\'en\'erateurs \ref{choixgen},
  tels que le diagramme suivant commute :
$$\xymatrix{\endo{D^b(\o{GD})}{R\Gamma_c[\rho]} \ar[r]^\sim_{ \beta_\phi^{-1}\circ \alpha_{\phi*}} 
& \endo{\o\QM_l}{V_{\sigma'(\rho)}} \otimes_{\o\QM_l} \TC_{d_\rho} \\
W_K \ar[u]^{\gamma_\rho}  \ar[ur]_{\sigma'(\rho)\otimes \tau^\phi_{d_\rho,I_\phi}} & }.$$
\end{prop}

\noindent {\bf Remarque :} Dans cette proposition et dans la suite, on omettra g\'en\'eralement  de pr\'eciser la d\'ependance des constructions en le choix de $\mu$, g\'en\'erateur de $\ZM_l(1)$. Le lecteur retrouvera cette d\'ependence sans difficult\'es.

Le reste de ce paragraphe est consacr{\'e} {\`a} la preuve de cette
proposition. Dans le paragraphe suivant nous d{\'e}terminerons l'ensemble
$I_\phi$ qui est la derni{\`e}re inconnue, du c\^ot\'e $Dr$ et sous une hypoth{\`e}se de validit{\'e} de la
conjecture monodromie-poids. Nous montrerons que $I_\phi=\emptyset$ (et donc est ind\'ependant de $\phi$), et ceci impliquera ({\em cf} remarque \ref{remaunibeta}) l'{\em unicit\'e} du choix de g\'en\'erateurs n\'ecessaire \`a la d\'efinition de $\beta_\phi$ comme en \ref{choixgen}.
Comme dans la partie \ref{indepbeta}, on peut alors v\'erifier que $\beta_\phi$ est en fait ind\'ependant de $\phi$.

\begin{lemme} (monodromie quasi-unipotente) \label{monqunip}
  Il existe un unique endomorphisme nilpotent $$N_\rho\in
  \endo{D^b_{\o\QM_l}(\o{GD})}{R\Gamma_c[\rho]}$$ tel que, si 
  $I_\rho\subseteq I_K$ d\'esigne le noyau de $\sigma'(\rho)_{|I_K}$, alors
 $$\forall i\in I_\rho,\;\; \gamma_\rho(i) = \hbox{exp}(N_\rho t_\mu(i)) .$$
 De plus on a 
 \ini
\begin{equation}
  \label{eqNrho}
\forall w\in W_K,\;\;  \gamma_\rho(w)N_\rho\gamma_\rho(w)^{-1} = q^{-\nu(w)} N_\rho=|w|N_\rho.
\end{equation}
\end{lemme}
\begin{proof}
Soit $\NC(\rho)$ le noyau du morphisme canonique
$$  \endo{D^b(\o{GD})}{R\Gamma_c[\rho]} \To{can} \prod_{i}
\endo{\o{GD}}{H^{d-1+i}_c[\rho]}. $$
Par le lemme \ref{actionphi}, $\NC(\rho)$ est form{\'e}
d'{\'e}l{\'e}ments nilpotents d'ordre $\leq d_\rho$. 
De plus,
 la proposition \ref{continuite} affirme que l'intersection $\NC(\rho)^0$ de $\NC(\rho)$ avec l'image de l'application canonique
$$ \endo{D^b_{\ZM_l}(\o{GD})}{R\Gamma_c(\o\MC,\ZM_l)} \To{} \endo{D^b_{\o\QM_l}(\o{GD})}{R\Gamma_c[\rho]} $$
est un $\ZM_l$-module de type fini et par cons\'equent est $l$-adiquement complet.

Le groupe $I_\rho$ de l'\'enonc\'e est un sous-groupe d'indice fini de
$I_K$ qui est aussi dans le noyau de l'action de $W_K$ sur chaque $H^i_c[\rho]$, par la description de la cohomologie  \ref{conjHarris}. Par cons\'equent
$\gamma_\rho$ envoie $I_\rho$ dans le sous-groupe $\id+\NC(\rho)$ de
$\aut{}{R\Gamma_c[\rho]}$, et m\^eme dans le sous-groupe $\UC(\rho)^0:=1+\NC(\rho)^0$ puisque l'action de $W_K$ sur $R\Gamma_c[\rho]$ est induite par celle sur $R\Gamma_c(\o\MC,\ZM_l)$. 
Mais puisque $\NC(\rho)^0$ est $l$-adiquement complet, l'application canonique
$$ \UC(\rho)^0 \To{} \limp{n} \UC(\rho)^0/(\id+l^n\NC(\rho)^0) $$
est un isomorphisme de groupe qui munit $\UC(\rho)^0$ d'une structure de pro-$l$-groupe.
Comme un morphisme de groupes entre groupes profinis est automatiquement continu, le morphisme $I_\rho \To{} \UC(\rho)^0$ se factorise par le plus grand pro-$l$-quotient de $I_\rho$.
Mais celui-ci n'est autre que l'image de $I_\rho$ par
le morphisme  $t_\mu : I_K \To{} \ZM_l$. 
Cette image est d'indice fini, donc de la forme $l^m\ZM_l$, et l'on peut par cons\'equent \'ecrire  
$$ \forall i\in I_\rho,\;\; \gamma_\rho(i)= u^{t_\mu(i)/l^m} $$
o\`u $u \in \UC(\rho)^0$ est l'image par $\gamma_\rho$ de n'importe quel \'el\'ement de $I_\rho$ s'envoyant sur $l^m$ par $t_\mu$.

Posant alors $N_\rho:= l^{-m} \hbox{log}(u) \in \NC(\rho)$, (le logarithme \'etant bien d\'efini puisque $u$ est unipotent),
on obtient
$$  \gamma_\rho (i) = \hbox{exp}(t_\mu(i).N_\rho) ,\;\;\forall i\in I_\rho. $$
Puisque pour tous $(w,i)\in W_K\times I_\rho$, on a $t_\mu(wiw^{-1})=q^{-\nu(w)}t_\mu(i)$, l'endomorphisme $N_\rho$ doit satisfaire
l'{\'e}quation \ref{eqNrho}.
Enfin, son unicit\'e est \'evidente.




\end{proof}

Soit $\phi$ un rel\`evement de Frobenius g\'eom\'etrique. Consid\'erons 
l'application
\ini
\begin{equation}\label{gammaphi} 
\application{\gamma^\phi_\rho(w):\;\;}{W_K}{\endo{D^b(\o{GD})}
  {R\Gamma_c[\rho]}^\times}{w}{\gamma_\rho(w) \hbox{exp}(-N_\rho t_\mu(i_\phi(w))) \hbox{ o{\`u} } w=\phi^{\nu(w)}i_\phi(w)}.
\end{equation}
L'{\'e}quation \ref{eqNrho} assure que l'application $\gamma^\phi_\rho$  est  un morphisme de groupes $W_K \To{} 
\aut{}{R\Gamma_c[\rho]}$. Par construction, celui-ci se factorise par
$W_K\To{} W_K/I_\rho$ qui est discret.

\begin{lemme} \label{lemmeeq}
Il existe un unique scindage $\alpha_\phi$ comme en \ref{defscin}  tel que le diagramme suivant commute :
$$\xymatrix{\endo{D^b(\o{GD})}{R\Gamma_c[\rho]} \ar[r]^\sim_{\beta^{-1}\circ \alpha_{\phi*}} 
& \endo{\o\QM_l}{V_{\sigma'(\rho)}} \otimes_{\o\QM_l} \TC_{d_\rho} \\
W_K \ar[u]^{\gamma^\phi_\rho}  \ar[ur]_{\sigma'(\rho)\otimes \tau_{S_{d_\rho}}} & }.$$
Rappelons que $\tau_{S_{d_\rho}}$ est la repr\'esentation semi-simple de $W_K$ correspondant, \`a torsion pr\`es, \`a la repr\'esentation triviale de $G_{d_\rho}$. En cons\'equence, la commutation du diagramme est ind\'ependante du choix de g\'en\'erateurs de \ref{choixgen} pour d\'efinir $\beta$.
\end{lemme}

\begin{proof}
Comme dans le lemme \ref{monqunip}, notons $I_\rho\subset I_K$ le noyau de $\sigma'(\rho)_{|I_K}$, qui est aussi celui de $(\gamma^\phi_\rho)_{|I_K}$, par construction. 
Puisque $I_K/I_\rho$ est fini, on peut trouver un entier
$n\in\NM$ tel que (l'image de) $\phi^n$ soit {\em central} dans le groupe
$W_K/I_\rho$.
Par le lemme de Schur, l'automorphisme $\sigma'(\rho)(\phi^n)$ de $V_{\sigma'(\rho)}$ est un scalaire que l'on notera $\omega\in\o\QM_l^\times$.
Par la description \ref{isocoho} de la cohomologie, on voit que l'endomorphisme $\HC^{d-1+j}(\gamma^\phi_\rho(\phi^n))$ de $H^{d-1+j}_c$ ({\em i.e} l'action de $\phi^n$ sur $H^{d-1+j}_c$) est annul\'e par le polyn\^ome $X-q^{nj}\omega$. Par le lemme \ref{actionphi} ii), il existe donc un unique scindage $\alpha_\phi$ tel que (pour tout choix de $\beta$)  
$$ \beta^{-1}\circ \alpha_{\phi*}(\gamma^\phi_\rho(\phi^n)) = \omega \otimes \hbox{Diag}(1,q^n,\cdots,q^{nd_\rho}) \in \endo{\o\QM_l}{V_{\sigma'(\rho)}}\otimes \TC_{d_\rho}.$$
Or, le commutant de $\omega \otimes \hbox{Diag}(1,q^n,\cdots,q^{nd_\rho})$ dans  $\endo{\o\QM_l}{V_{\sigma'(\rho)}}\otimes \TC_{d_\rho}$ est $\endo{\o\QM_l}{V_{\sigma'(\rho)}}\otimes \o\QM_l^{d_\rho}$ o\`u $\o\QM_l^{d_\rho}$ d\'esigne ici la sous-alg\`ebre des matrices diagonales de $\TC_{d_\rho}$. Ainsi, puisque $\phi^n$ est central dans $W_K/\ker(\gamma^\phi_\rho)$, on a 
$$ \hbox{Im} (\beta^{-1}\circ \alpha_{\phi*}\circ \gamma^\phi_\rho) \subset \endo{\o\QM_l}{V_{\sigma'(\rho)}}\otimes \o\QM_l^{d_\rho}. $$
Mais par les diagrammes commutatifs de la proposition \ref{propnoneq} et la $W_K$-\'equivariance des isomorphismes \ref{isocoho}, on en d\'eduit que 
\begin{eqnarray*}
 \forall w\in W_K,\;\; (\beta^{-1}\circ \alpha_{\phi*})(\gamma^\phi_\rho(w)) & = & \sum_{i=0}^{d_\rho} \sigma'(\rho)(w)|w|^{-i} \otimes E_{ii} \\
 & = & \sum_{i=0}^{d_\rho} \sigma'(\rho)(w) \otimes |w|^{-i}E_{ii}
\end{eqnarray*}
o\`u $E_{ii}$ d\'esigne toujours la matrice \'el\'ementaire de coordonn\'ees $(i,i)$ de $\TC_{d_\rho}$. Le scindage $\alpha_\phi$ fait donc bien commuter le diagramme du lemme. Pour l'unicit\'e, remarquons que toute autre solution $\alpha'$ v\'erifie la conclusion du point ii) de \ref{actionphi} pour $\phi^n$, et co\"incide donc avec $\alpha$ par l'assertion d'unicit\'e de \ref{actionphi} ii).

\end{proof}

\alin{Preuve de la proposition \ref{propeq}}
On garde les notations pr\'ec\'edentes ; en particulier, $\alpha_\phi$ est le scindage de $R\Gamma_c[\rho]$ du lemme \ref{lemmeeq}, et
$N_\rho$ 
l'endomorphisme nilpotent de $R\Gamma_c[\rho]$ donn{\'e} par le lemme
\ref{monqunip}. On voudrait montrer que pour un bon $\beta=\beta_\phi$ (correspondant \`a un choix de g\'en\'erateurs des $\ext{i}{\pi^i_{\rho^\vee}}{\pi^{i-1}_{\rho^\vee}}{\o{G_d}}$), il existe un certain sous-ensemble $I\subseteq S_{d_\rho}$ tel que
$$(\beta_\phi^{-1}\circ\alpha_{\phi*})(N_\rho) = (\id_{V_{\sigma'(\rho)}}\otimes N_{I}) \in
\endo{\o\QM_l}{V_{\sigma'(\rho)}} \otimes_{\o\QM_l} 
\TC_{d_\rho}$$
o\`u $N_I$ a \'et\'e d\'efini en \ref{deftauI}.
Choisissons un $\beta$ arbitraire pour commencer.
Remarquons que par construction, $N_\rho$ induit l'endomorphisme nul en cohomologie et par cons\'equent on a
$$(\beta^{-1}\circ\alpha_{\phi*})(N_\rho) \in \sum_{i>j} \endo{\o\QM_l}{V_{\sigma'(\rho)}}
\otimes E_{ij}, $$
en notant $(E_{ij})_{i\geq j}$ la base "canonique" de $\TC_{d_\rho}$. 
D'autre part, on d\'eduit de \ref{eqNrho} la relation 
\ini\begin{equation}\label{rel}
 \gamma^\phi_\rho(w) N_\rho \gamma^\phi_\rho(w)^{-1} = |w| N_\rho 
\end{equation}
pour tout $w\in W_K$. Appliqu{\'e}e {\`a} l'{\'e}l{\'e}ment $\phi^n$ de la preuve du
lemme \ref{monqunip}, pour lequel on a $$(\beta^{-1}\circ\alpha_{\phi*})(\phi^n)= \sum_{i=0}^{d_\rho-1} \omega q^{ni}(\id \otimes E_{ii}),$$
 cette relation implique que 
$$(\beta^{-1}\circ\alpha_{\phi*})(N_\rho) \in \sum_{i=1}^{d_\rho-1} \endo{\o\QM_l}{V_{\sigma'(\rho)}}
\otimes E_{i-1,i}.$$ 
  {\'E}crivons donc $(\beta^{-1}\circ\alpha_{\phi*})(N_\rho)= \sum_i M_i\otimes E_{i-1,i}$ avec $M_i\in
\endo{\o\QM_l}{V_{\sigma'(\rho)}}$. 
Par d{\'e}finition de $\alpha_\phi$, on a pour tout $w\in W_K$ 
$$ (\beta^{-1}\circ\alpha_{\phi*})(\gamma^\phi_\rho(w)) = \sum_{i=0}^{d_\rho-1} 
\sigma'(\rho)(w) \otimes |w|^{-i}E_{ii}. $$
La relation \ref{rel} implique donc que chaque $M_i$ commute avec
$\sigma'(\rho)(w)$ pour tout $w$. Par le lemme de Schur, on a donc
$M_i\in\o\QM_l$. Quitte \`a changer $\beta$, on peut alors supposer que $M_i$ est soit nul, soit \'egal \`a $1$. Il n'y a plus qu'\`a poser $I:=\{i\in\{1,\cdots,d_\rho-1\}, M_{d-i}=0\}$.

\begin{rema} \label{remaunibeta}
Supposons que $I=\emptyset$, {\em i.e.} que tous les $M_i$ sont non-nuls. Alors le changement de $\beta$ est {\em unique}.
\end{rema}

\subsection{Uniformisation $p$-adique et conjecture de puret{\'e}} \label{unifpurete}

Le but de cette section est de prouver que l'ensemble $I_\phi\subseteq S_{d_\rho}$ de la proposition \ref{propeq} est l'ensemble vide, ce qui ach\`evera la description explicite de l'action de $W_K$ sur le complexe $R\Gamma_c[\rho]$. Par construction de cet ensemble $I_\phi$, il est \'equivalent de montrer que l'endomorphisme $N_\rho$ de $R\Gamma_c[\rho]$ d\'efini dans le lemme \ref{monqunip} (le logarithme de la partie unipotente de la monodromie) est d'ordre exactement $d_\rho$.

Jusqu'ici notre \'etude concernait indiff\'eremment les c\^ot\'es $LT$ et $Dr$ (sous la description \ref{conjHarris} plus le th\'eor\`eme de Faltings \ref{faltings} dans le cas $Dr$), 
mais la seule strat\'egie que nous avons trouv\'ee pour estimer l'ordre de nilpotence de $N_\rho$  concerne pour l'instant le c\^ot\'e $Dr$ et utilise la th\'eorie de l'uniformisation $p$-adique de Drinfeld. Elle suppose aussi la validit\'e de la conjecture dite "monodromie-poids" de Deligne. 
Rappelons de quoi il s'agit.

\alin{Puret\'e de la filtration de monodromie} \label{monpds}
 On sait depuis Grothendieck que sur toute repr\'esentation   $l$-adique continue de dimension finie $(\sigma,V)$ de $W_K$, l'inertie $I_K$ agit de mani\`ere quasi-unipotente, 
c'est-\`a-dire qu'il existe 
un unique endomorphisme nilpotent $N_\sigma$ de $V$ tel qu'il existe un sous-groupe $I_\sigma \subset I_K$ d'indice fini dont l'action est d\'ecrite par
$$ \forall i\in I_\sigma,\;\; \sigma(i)=\hbox{exp}(N_\sigma t_\mu(i)).$$
L'op\'erateur $N_\sigma$ est donc le logarithme de la partie
unipotente de la monodromie de $\sigma$, mais par abus de langage,
nous l'appelerons simplement "op\'erateur de monodromie de
$\sigma$"\footnote{Bien-s\^ur $N_\sigma$ d\'epend du g\'en\'erateur
  $\mu$ de $\ZM_l(1)$, mais $N_\sigma\otimes \mu^* : \;
  V\otimes_{\ZM_l} \ZM_l(1) \To{} V$ n'en d\'epend pas.}. 
Il v\'erifie n\'ecessairement l'\'equation habituelle
$wN_\sigma w^{-1}=|w|N_\sigma$ de sorte que pour tout rel\`evement de
Frobenius g\'eom\'etrique $\phi$, 
l'application
\ini\begin{equation}
\label{sigmalisse}
w\mapsto  \sigma^\phi(w):=\sigma(w)\hbox{exp}(-N_\sigma t_\mu(i_\phi(w)))\;\;\hbox{ o\`u } \; w=\phi^{\nu(w)}i_\phi(w) 
 \end{equation}
d\'efinit une repr\'esentation {\em lisse} de $W_K$ sur $V$. Rappelons aussi que par
\cite[8.4.2]{DelAntwerp}, la classe d'isomorphisme de $\sigma^\phi$ ne d\'epend pas du choix de $\phi$.

Ceci s'applique en particulier aux espaces de cohomologie $l$-adique $H^i(X_{K^{ca}},\o\QM_l)$ d'une vari\'et\'e $X$ propre sur $K$.
\`A l'op\'erateur de monodromie $N$ est associ\'ee une filtration croissante stable sous $W_K$ de l'espace $V$ dite "filtration de monodromie"\footnote{qui elle est bien ind\'ependante du choix de $\mu$} $\cdots \subseteq M_iV\subseteq M_{i+1}V \subseteq\cdots $ de longueur finie et caract\'eris\'ee par les propri\'et\'es que $N(M_iV)\subseteq M_{i-2}V$ pour tout $i$, et $N$ induit des isomorphismes $N^i:\;\; \hbox{Gr}_M^{i}V \otimes |-|^i\simto \hbox{Gr}_{M}^{-i}V$ pour tout $i\geq 0$.
Par ailleurs, Deligne d\'efinit une autre filtration croissante stable sous $W_K$ de $V$, dite "filtration par les poids", $\cdots \subseteq W_iV \subseteq W_{i+1}V\subseteq \cdots$ caract\'eris\'ee par la propri\'et\'e que les valeurs propres de {\em tout} rel\`evement  de Frobenius g\'eom\'etrique sont des entiers alg\'ebriques dont tous les conjugu\'es complexes sont de norme complexe $q^{i/2}$.

\begin{conj} \label{conmp}(Monodromie-Poids) Si $X$ est propre et lisse sur $K$, 
alors pour tout $i\in \NM$, on a
$$ M_i\left(H^j(X^{ca},\o\QM_l)\right)= W_{i+j}\left(H^j(X^{ca},\o\QM_l)\right) .$$
\end{conj}

Lorsque $K$ est d'\'egales caract\'eristiques, l'\'enonc\'e est essentiellement contenu et d\'emontr\'e dans les travaux de Deligne sur les conjectures de Weil. Le cas d'in\'egales caract\'eristiques est tr\`es peu avanc\'e, m\^eme dans les cas de r\'eduction semi-stable.

Remarquons que pour tout $i\in\ZM$ on a $N(W_iV)\subseteq W_{i-2}V$, de sorte que par la caract\'erisation de la filtration de monodromie, la conjecture ci-dessus est \'equivalente \`a l'assertion : {\em Pour tout $i\geq 0$, $N$ induit un isomorphisme  
\ini\begin{equation}\label{equiMP}
N^i:\;\; \hbox{Gr}_W^{i+j}\left(H^j(X^{ca},\o\QM_l)\right)\otimes |-|^i \simto \hbox{Gr}_{W}^{-i+j}\left(H^j(X^{ca},\o\QM_l)\right).
\end{equation}
}

\alin{Uniformisation $p$-adique} Nous n'aurons besoin dans ce texte que d'une version simple de la propri\'et\'e d'uniformisation $p$-adique, notamment nous n'utiliserons pas le langage des vari\'et\'es de Shimura. 

Commen{\c c}ons par rappeler que les espaces $\o\MC_{Dr,n}^{d/K}$ sont munis d'une donn\'ee de descente sous le groupe de Weil $W_K$, qui est d\'eduite de celle que nous avons d\'efinie sur les espaces $\MC_{Dr,n}^{d/K}$ en \ref{actionDr} par passage au quotient. Il se trouve que cette donn\'ee de descente n'est pas effective sur $\MC_{Dr,n}^{d/K}$, mais l'est sur $\o\MC_{Dr,n}^{d/K}$. Ceci est d\'emontr\'e dans un contexte plus g\'en\'eral dans \cite[Thm 3.49]{RZ}.
On peut expliciter la descente de $\o\MC_{Dr,0}^{d/K}$ ainsi : notons $K_d$ l'extension non-ramifi\'ee de degr\'e $d$ de $K$. Alors on a
$$ \o\MC_{Dr,0}^{d/K} \simeq \wh{K^{ca}} \wh\otimes \left(\res{K_d}{K}{\Omega^{d-1}_K\wh\otimes_K K_d} \right) $$
o\`u l'action de $\o{GD}$ sur la restriction des scalaires de $K_d$ \`a $K$ de $\Omega^{d-1}_K\wh\otimes_K K_d$ est le produit de l'action naturelle de $\o{GD}/\o{\dd}\simeq PG_d$ sur $\Omega^{d-1}_K$ et de l'action de $\o{GD}$ sur $K_d$ donn\'ee par $(g,d)\mapsto \phi^{\nu(g)+\nu(d)}$.
Nous noterons ici $\o\MC_{Dr,n}$ les $K$-espaces analytiques  "descendus" des $\wh{K^{ca}}$-espaces $\o\MC_{Dr,n}^{d/K}$. Il sont toujours munis d'actions du groupe  $\o{GD}$, les morphismes de restriction de la structure de niveau se descendent aussi et on obtient une tour de   rev\^etements \'etales de groupe $\OC_D^\times$ du $K$-espaces analytique  $\o\MC_{Dr,0}=\res{K_d}{K}{\Omega^{d-1}_K\wh\otimes_K K_d}$.

Soit $\Gamma$ un sous-groupe {\em discret, cocompact et sans torsion} de $\o{G_d}$. 
On sait que l'action d'un tel sous-groupe  sur l'espace analytique $\o\MC_{Dr,0}$ est {\em libre}, et il en est donc de m\^eme de l'action sur les rev\^etements $\o\MC_{Dr,n}$.
Par \cite[lemma 4]{Bic5}, l'espace annel\'e quotient est alors muni d'une structure de $K$-espace analytique quotient que nous noterons $\o\MC_{Dr,n}/\Gamma$. 
Le th\'eor\`eme suivant affirme que ces espaces sont alg\'ebrisables.

\begin{theo} (Mustafin, Cherednik, Drinfeld, Rapoport-Zink, Varshawski...) Il existe une vari\'et\'e alg\'ebrique $S_{\Gamma,0}$, propre et lisse sur $K$
et une tour $(S_{\Gamma,n})_{n\in\NM}$ de rev\^etements \'etales de $S_{\Gamma,0}$ de groupe $\OC_D^\times/(1+\varpi^n\OC_D)$, munies d'une tour d'isomorphismes compatibles
$$ \theta_n:\;\;S_{\Gamma,n}^{an} \simto \o\MC_{Dr,n}/\Gamma.$$
Toutes ces vari\'et\'es sont uniques \`a isomorphisme canonique pr\`es et sont munies d'une action de $\o{\dd}$ prolongeant l'action galoisienne de $\OC_D^\times$.
\end{theo}

\begin{proof}
Dans le cas o\`u $K$ est $p$-adique, la r\'ef\'erence la plus commode est le th\'eor\`eme 6.36 de \cite{RZ}. L'uniformisation y est obtenue par voie modulaire et contient beaucoup plus d'informations que ce dont nous avons besoin ici. 
Dans le cas d'\'egales caract\'eristiques, le r\'esultat dans la g\'en\'eralit\'e demand\'ee est d\^u \`a Drinfeld, {\em cf} les commentaire de la partie III.5 de \cite{BouCar}.
\end{proof}

Appliquons maintenant le th\'eor\`eme de comparaison du type GAGA de Berkovich \cite[7.1]{Bic2} :
\begin{coro}
Il y a des isomorphismes $\dd\times W_K$-\'equivariants
$$ H^i(S_{\Gamma,n}\otimes_K K^{ca},\o\QM_l) \simto H^i_c(\o\MC_{Dr,n}^{d/K}/\Gamma,\o\QM_l) $$
o\`u le terme de gauche d\'esigne la cohomologie \'etale $l$-adique au sens des vari\'et\'es alg\'ebriques. En particulier, si $\rho\in \Irr{\o\QM_l}{\dd/\varpi^\ZM(1+\varpi^n\OC_D)}$, on en d\'eduit des isomorphismes $W_K$-\'equivariants
$$ H^i(S_{\Gamma,n}^{ca},\o\QM_l)[\rho]\simto H^i_c(\o\MC^{d/K}_{Dr,n}/\Gamma,\o\QM_l)[\rho].$$
\end{coro}

\alin{Application}
Fixons $\rho$ de niveau $\leq n$, c'est-\`a-dire $\rho\in\Irr{\o\QM_l}{\dd/\varpi^\ZM(1+\varpi^n\OC_D)}$.
On a alors un premier isomorphisme dans $D^b(\o{GD})$, \'evident par construction de $R\Gamma_c(\MC^{d/K}_{Dr},\o\QM_l)$, et qui est $W_K$ \'equivariant  :
$$ R\Gamma_c(\o\MC_{Dr,n}^{d/K},\o\QM_l)[\rho] \simto R\Gamma_c(\o\MC_{Dr}^{d/K},\o\QM_l)[\rho]=R\Gamma_c[\rho]. $$
Par la version "foncteurs d\'eriv\'es" de la suite spectrale de Hochshild-Serre de \ref{theoHS}, on a un isomorphisme $W_K$-\'equivariant dans $D^b(\o{\dd})$ :
\ini\begin{equation}\label{isoHS}
 \o\QM_l \otimes^L_{\o\QM_l[\Gamma]} R\Gamma_c(\o\MC_{Dr,n}^{d/K},\o\QM_l)[\rho]
\simto R\Gamma_c(\o\MC_{Dr,n}^{d/K}/\Gamma,\o\QM_l)[\rho]. 
\end{equation}
Dans l'expression de gauche, le complexe est vu \`a travers le foncteur d'oubli $D^b(\o{GD})\To{} D^b(\Gamma\times \o\dd)$ associ\'e \`a l'inclusion de groupes $\Gamma\times \o\dd \To{} \o{GD}$.
Soit alors $N_{\rho,\Gamma,j}$ l'endomorphisme du $\o\QM_l$-espace vectoriel $H^{j}(S_{\Gamma,n}^{ca},\o\QM_l)[\rho]$ fonctoriellement induit par $N_\rho$ via les deux isomorphismes pr\'ec\'edents, le passage \`a la cohomologie en degr\'e $j$ et le th\'eor\`eme de comparaison GAGA. Comme tous les isomorphismes utilis\'es sont $W_K$-\'equivariants, la d\'efinition de $N_\rho$ montre que le sous-groupe $I_\rho\subset I_K$ du lemme \ref{monqunip} agit sur $H^{j}(S_{\Gamma,n}^{ca},\o\QM_l)[\rho]$ par $i\mapsto \hbox{exp}(N_{\rho,\Gamma,j}t_\mu(i))$, de sorte que $N_{\rho,\Gamma,j}$ est l'op\'erateur de monodromie de la repr\'esentation de $W_K$ sur $H^{j}(S_{\Gamma,n}^{ca},\o\QM_l)[\rho]$.

Fixons dor\'enavant un rel\`evement de Frobenius g\'eom\'etrique $\phi$ et 
notons $H^{j,\phi}(S_{\Gamma,n}^{ca},\o\QM_l)[\rho]$ la repr\'esentation lisse de $W_K$ associ\'ee, comme en \ref{sigmalisse}.

Combinons alors l'isomophisme \ref{isoHS} avec le scindage $\alpha_\phi$ de \ref{lemmeeq} : la suite spectrale de Hochschild-Serre d\'eg\'en\`ere en des isomorphismes $\o\dd\times W_K$-\'equivariants
\begin{eqnarray*} 
H^{j,\phi}(S_{\Gamma,n}^{ca},\o\QM_l)[\rho] & \simto  & \bigoplus_{i=0}^{d_\rho-1}
\tor{j-d+1-i}{\o\QM_l}{\pi_{\rho^\vee}^i}{\Gamma}\otimes \rho\otimes\sigma'(\rho)|-|^{-i} \\
& \simeq & \bigoplus_{i=0}^{d_\rho-1} \ext{d-1-j+i}{\pi^i_{\rho^\vee}}{\o\QM_l}{\Gamma}^* \otimes \rho\otimes\sigma'(\rho)|-|^{-i}\\
& \simeq & \bigoplus_{i=0}^{d_\rho-1}
\ext{d-1-j+i}{\pi_{\rho^\vee}^i}{\CC^\infty(\o{G_d}/\Gamma,\o\QM_l)}{\o{G_d}}^* \otimes \rho\otimes\sigma'(\rho)|-|^{-i}
\end{eqnarray*}
La repr\'esentation $\CC^\infty(\o{G_d}/\Gamma,\o\QM_l)$ de $\o{G_d}$ est admissible et semi-simple, avec constituents "unitarisables" : en particulier, les seules repr\'esentations elliptiques qui peuvent y apparaitre sont les s\'eries discr\`etes et les repr\'esentations de Speh locales.
Soit $m_{\rho,\Gamma}^\emptyset$, resp. $m'_{\rho,\Gamma}$, la multiplicit\'e de la s\'erie discr\`ete $\pi^\emptyset_{\rho^\vee}$, resp. de la repr\'esentation de Speh $\pi^{S_{d_\rho}}_{\rho^\vee}$  dans cette repr\'esentation $\CC^\infty(\o{G_d}/\Gamma,\o\QM_l)$. Compte tenu du calcul de $\hbox{Ext}$ entre repr\'esentations elliptiques (le calcul partiel de Schneider-Stuhler suffit ici), on trouve la description suivante :

\noindent{\em Pour $j=d-1$, on a
$$ {H}^{d-1,\phi}(S_{\Gamma,n}^{ca},\o\QM_l)[\rho] \mathop{\simeq}\limits_{\dd\times W_K}
 \left\{ \begin{array}{ll}
 \left(\mathop{\bigoplus}\limits_{i=0}^{d_\rho-1} \rho\otimes \sigma'(\rho)|.|^{-i}\right)^{m_{\rho,\Gamma}^\emptyset} & \hbox{si }
d_\rho \hbox{ est pair} \\
\left(\mathop{\bigoplus}\limits_{i=0}^{d_\rho-1} \rho\otimes \sigma'(\rho)|.|^{-i}\right)^{m_{\rho,\Gamma}^\emptyset} \oplus \left(\rho\otimes\sigma'(\rho)|.|^{-k}\right)^{m'_{\rho,\Gamma}} & \hbox{si }
d_\rho=1+2k
\end{array}\right.
$$
et pour $j\neq d-1$, on a 
$$ {H}^{j,\phi}(S_{\Gamma,n}^{ca},\o\QM_l)[\rho] \mathop{\simeq}\limits_{\dd\times W_K}
 \left\{ \begin{array}{ll}
 0 & \hbox{si }
j+d_\rho-d \hbox{ est impair} \\
\left(\rho\otimes\sigma'(\rho)|.|^{-k}\right)^{m'_{\rho,\Gamma}} & \hbox{si }
j+d_\rho-d=2k\geq 0
\end{array}\right.
$$
}
Observons en particulier que l'action de $W_K$ sur ${H}^{j,\phi}(S_{\Gamma,n}^{ca},\o\QM_l)[\rho]$ est semi-simple. On en d\'eduit la remarque suivante :

\begin{rema} (sous la description \ref{conjHarris})
L'action d'un rel\`evement de Frobenius $\phi$ sur les espaces de cohomologie $H^i(S_{\Gamma,n}^{ca},\o\QM_l)$ est semi-simple.
\end{rema}

Revenons \`a notre probl\`eme initial ; 
par le th\'eor\`eme de "multiplicit\'es limites" de \cite[1.3]{Rog}, on sait que, quitte \`a remplacer $\Gamma$ par un sous-groupe d'indice fini, on peut supposer $m_{\rho,\Gamma}^\emptyset >0$. 
On peut maintenant \'enoncer
\begin{prop} \label{propequMP}
Avec les notations ci-dessus, supposons $\Gamma$ "assez petit" pour que $m_{\rho,\Gamma}^\emptyset\neq 0$. Alors les propri\'et\'es suivantes sont \'equivalentes (toujours sous la description \ref{conjHarris}) :
\begin{enumerate}
\item L'endomorphisme nilpotent $N_\rho$ de $R\Gamma_c[\rho]$ d\'efini dans le lemme \ref{monqunip} est d'ordre $d_\rho$ ({\em i.e.} v\'erifie $N_\rho^{d_\rho-1}\neq 0$).
\item L'op\'erateur de monodromie $N_{\rho,\Gamma,d-1}$ de $H^{d-1}(S_{\Gamma,n}^{ca},\o\QM_l)[\rho]$ est d'ordre $d_\rho$.
\item La conjecture monodromie-poids est v\'erifi\'ee pour la partie $\rho$-covariante de la cohomologie $H^{j}(S_{\Gamma,n}^{ca},\o\QM_l)[\rho]$ pour tout $j$. 
\item Pour tout rel\`evement de Frobenius g\'eom\'etrique $\phi$, le sous-ensemble $I_\phi\subseteq S_{d_\rho}$ de la proposition \ref{propeq} est l'ensemble vide $I_\phi=\emptyset$.
\end{enumerate}
\end{prop}
\begin{proof}
On a d\'eja remarqu\'e que $i)$ et $iv)$ sont \'equivalents. 
Une autre implication facile est $ii)\Rightarrow i)$, puisque $N_\rho$ est d'ordre au plus $d_\rho$ et induit $N_{\rho,\Gamma,d-1}$. 

Rappelons maintenant que la repr\'esentation $\sigma'(\rho)=\sigma_{d/d_\rho}(\tau^0_\rho)|-|^{\frac{d_\rho-d}{2}}$ est pure de poids $d-d_\rho$. 
Ainsi la description de la cohomologie ci-dessus montre que $\hbox{Gr}^{d-d_\rho}_W(H^{d-1}(S_{\Gamma,n}^{ca},\o\QM_l)[\rho])\neq 0$, et donc la conjecture monodromie-poids dans sa version \ref{equiMP} implique que $N_{\rho,\Gamma,d-1}^{d_\rho-1}\neq 0$. On a donc $iii)\Rightarrow iv)$. 

Il nous reste \`a prouver que $iv)\Rightarrow iii)$. En fait le seul espace de cohomologie qui peut poser probl\`eme pour la conjecture monodromie-poids est celui de degr\'e m\'edian $j=d-1$, les autres \'etant purs. Il nous faut alors expliciter l'action de $N_{\rho,\Gamma,d-1}$ sur 
$H^{d-1}(S_{\Gamma,n}^{ca},\o\QM_l)[\rho]$. Mais si on suppose la propri\'et\'e $iv)$, alors dans l'isomorphisme
$$H^{d-1}(S_{\Gamma,n}^{ca},\o\QM_l)[\rho] \simeq 
\bigoplus_{i=0}^{d_\rho-1}
\ext{i}{\pi_{\rho^\vee}^i}{(\pi_{\rho^\vee}^\emptyset)^{m_{\rho,\Gamma}^\emptyset}\oplus (\pi_{\rho^\vee}^{S_{d_\rho}})^{m'_{\rho,\Gamma}}}{\o{G_d}}^* \otimes \rho\otimes\sigma'(\rho)|-|^{-i},
$$
induit par le scindage $\alpha_\phi$, l'op\'erateur $N_\rho$ agit par $\cup$-produit et la description \ref{theoextell} ii) de ce $\cup$-produit montre que $N_\rho$ induit des isomorphismes
$$ \ext{i}{\pi_{\rho^\vee}^i}{(\pi_{\rho^\vee}^\emptyset)^{m_{\rho,\Gamma}^\emptyset}}{\o{G_d}}^* \otimes \rho\otimes\sigma'(\rho)|-|^{-i} \simto
\ext{i-1}{\pi_{\rho^\vee}^{i-1}}{(\pi_{\rho^\vee}^\emptyset)^{m_{\rho,\Gamma}^\emptyset}}{\o{G_d}}^* \otimes \rho\otimes\sigma'(\rho)|-|^{-i+1}
$$ pour tout $i\in \{1,\cdots, d_\rho-1\}$.
Or, toujours par la description de la cohomologie donn\'ee plus haut, on a pour 
$i\neq \frac{d_\rho-1}{2}$  
$$ \ext{i}{\pi_{\rho^\vee}^i}{(\pi_{\rho^\vee}^\emptyset)^{m_{\rho,\Gamma}^\emptyset}}{\o{G_d}}^* \otimes \rho\otimes\sigma'(\rho)|-|^{-i} \simto \hbox{Gr}^{(d-1)+(1-d_\rho+2i)}_W(H^{d-1}(S_{\Gamma,n}^{ca},\o\QM_l)[\rho] ) $$
et pour $i=\frac{d_\rho-1}{2}$ (lorsque $d_\rho$ est impair !) on a
$$ \ext{i}{\pi_{\rho^\vee}^i}{(\pi_{\rho^\vee}^\emptyset)^{m_{\rho,\Gamma}^\emptyset} \oplus (\pi_{\rho^\vee}^{S_{d_\rho}})^{m'_{\rho,\Gamma}}}{\o{G_d}}^* \otimes \rho\otimes\sigma'(\rho)|-|^{-i} \simto \hbox{Gr}^{d-1}_W(H^{d-1}(S_{\Gamma,n}^{ca},\o\QM_l)[\rho] ). $$
Dans tous les cas, la description de l'action de $N_\rho$ par $\cup$-produit montre que pour tout $k\in\ZM$, $N^k$ induit un isomorphisme
$$ N^k :\;\; \hbox{Gr}_W^{d-1+k}\left(H^{d-1}(S_{\Gamma,n}^{ca},\o\QM_l)[\rho]\right) \simto \hbox{Gr}_W^{d-1-k}\left(H^{d-1}(S_{\Gamma,n}^{ca},\o\QM_l)[\rho]\right) $$
et par \ref{equiMP}, la conjecture monodromie-poids pour la partie $\rho$-covariante de la cohomologie en d\'ecoule.

\end{proof}

\begin{rema} \label{preuveprop3}
On a ainsi montr\'e la proposition \ref{prop3}, mais bien-s\^ur, les r\'esultats de cette section et de la pr\'ec\'edente sont bien plus pr\'ecis.
\end{rema}

\subsection{Preuve du th{\'e}or{\`e}me \ref{thprin}}
Fixons $\rho\in \Irr{\o\QM_l}{\o\dd}$ et
 $I\subseteq S_{d_\rho}$, et choisissons un rel\`evement de Frobenius g\'eom\'etrique $\phi$. Nous noterons dans cette section
$$ \HC_\rho^I:= \HC^*\left( R\hom{R\Gamma_c[\rho^\vee]}{\pi_\rho^I \otimes \rho^\vee}{D^b(\o{GD})}   \right). $$
C'est un $\o\QM_l$-espace vectoriel gradu\'e de dimension finie muni
de l'action par automorphismes de degr\'e $0$ de $W_K$ induite par
l'action $\gamma_{\rho^\vee} :\; W_K\To{}
\aut{D^b(\o{GD})}{R\Gamma_c[{\rho^\vee}]}$. Nous noterons aussi
$\HC_\rho^{I,\phi}$ la repr\'esentation gradu\'ee lisse de $W_K$
associ\'ee \`a $\HC_\rho^I$ comme en \ref{sigmalisse} ; elle est aussi
induite par l'action $\gamma^\phi_{\rho^\vee}:\; W_K\To{}
\aut{D^b(\o{GD})}{R\Gamma_c[{\rho^\vee}]}$ d\'efinie en \ref{gammaphi}. Enfin,
nous noterons $N^I_\rho$ l'endomorphisme de degr\'e $0$ (d'espace vectoriel
gradu\'e) de $\HC_\rho^I$ induit fonctoriellement par l'endomorphisme
$N_{\rho^\vee}$ de $R\Gamma_c[{\rho^\vee}]$ du lemme \ref{monqunip}. L'endomorphisme $N_\rho^I$ est
l'op\'erateur de monodromie de la $W_K$-repr\'esentation gradu\'ee $\HC^I_\rho$
et la connaissance de cette derni\`ere \'equivaut \`a celle du couple
$(\HC_\rho^{I,\phi},N_\rho^I)$.

En utilisant le scindage $\alpha_\phi$ de \ref{lemmeeq}, on obtient
imm\'ediatement la premi\`ere ligne de
\begin{eqnarray*}
 \HC_\rho^{I,\phi} & \simeq_{W_K} & \bigoplus_{i=0}^{d_\rho-1} \HC^*\left(R\hom{\pi_\rho^i\otimes {\rho^\vee}\otimes \sigma'({\rho^\vee})|-|^{-i}}{\pi_\rho^I \otimes {\rho^\vee}}{D^b(\o{GD})} \right)[d-1+i] \\
 & \simeq_{W_K} & \bigoplus_{i=0}^{d_\rho-1}
 \ext{\delta(i,I)}{\pi_\rho^i\otimes {\rho^\vee}}{\pi_\rho^I\otimes
   {\rho^\vee}}{\o{GD}} \otimes_{\o\QM_l}
 \left(\sigma'(\rho^\vee)|-|^{-i}\right)^\vee [d-1+i-\delta(i,I)] \\
 & \simeq_{W_K} & \sigma'(\rho^\vee)^\vee \otimes\left(\bigoplus_{i=0}^{d_\rho-1}
 \ext{\delta(i,I)}{\pi_\rho^i}{\pi_\rho^I}{\o{G_d}} \otimes_{\o\QM_l}
 |-|^{i} [d-1+i-\delta(i,I)]\right)
 \end{eqnarray*}
 et la deuxi\`eme ligne vient des isomorphismes canoniques 
$$ \ext{*}{\pi^I_\rho\otimes {\rho^\vee}}{\pi^J_\rho\otimes {\rho^\vee}}{\o{GD}}
\simto \ext{*}{\pi^I_\rho}{\pi^J_\rho}{\o{G_d}} $$ 
et du calcul d'extensions de \ref{theoextell}.

On veut maintenant expliciter l'endomorphisme $N_\rho^I$. On proc\`ede
exactement comme dans le paragraphe \ref{preuvetheodp}.
Notons $\beta_{i-1,i} \in \ext{1}{\pi_\rho^{i}
  }{\pi_\rho^{i-1}}{\o{G_d}}$, o\`u $i\in
\{1,\cdots,d_\rho-1\}$ les g\'en\'erateurs qui d\'efinissent l'isomorphisme $\beta_\phi$ de  la
proposition \ref{propeq} (uniques par la remarque \ref{remaunibeta}) et choisissons 
pour chaque $j\in \{0,\cdots,d_\rho-1\}$ un
g\'en\'erateur de la droite
$\ext{\delta(i,I)}{\pi^i_\rho}{\pi^I_\rho}{\o{G_d}}$.
Alors l'action (\`a droite) de $\alpha(N_{\rho^\vee})$ sur $\HC_\rho^I$ est
donn\'ee par 
$$ N.e_i = e_i \cup \beta_{i,i+1} \in \o\QM_l e_{i+1}. $$
La formule pour le $\cup$-produit de \ref{theoextell} ii) montre que
$e_i\cup \beta_{i,i+1}$ est non-nul \ssi\
$\delta(i+1,I)=\delta(i,I)+1$ et des consid\'erations
\'el\'ementaires montrent que cela se produit \ssi\ $i+1\in I^c$
(compl\'ementaire de $I$ dans $S_{d_\rho}=\{1,\cdots,d_\rho-1\}$).

Oublions maintenant la graduation de $\HC^I_\rho$. La
discussion de \ref{preuvetheodp} montre alors que 
$$ \HC^I_\rho \simeq_{W_K} \sigma'(\rho^\vee)^\vee \otimes \tau_{d_\rho,I}
|-|^{d_\rho-1}, $$
et compte tenu de la d\'efinition \ref{defsigmaprime} de $\sigma'(\rho)$, de l'\'egalit\'e $\tau_{\rho^\vee}^0=(\tau_\rho^0)^\vee$ et de la compatibilit\'e de la correspondance de Langlands aux contragr\'edientes, on obtient
\begin{eqnarray*}
 \HC^I_\rho & \simeq_{W_K}  & \sigma_{d/d_\rho}(\tau_\rho^0)|-|^{\frac{d-d_\rho}{2}}\otimes \tau_{d_\rho,I}|-|^{d_\rho-1} \\
& \simeq_{W_K} & \sigma_{d/d_\rho}(\tau_\rho^0)|-|^{\frac{d-1}{2}}\otimes \tau_{d_\rho,I}|-|^{\frac{d_\rho-1}{2}} \simeq \sigma_d(\pi_\rho^I)|-|^{\frac{d-1}{2}} 
\end{eqnarray*}
par la description \ref{corell} de la correspondance de Langlands pour les repr\'esentations elliptiques.

En  vertu de \ref{corodesc} et \ref{reduc1}, nous avons donc achev\'e la d\'emonstration du th\'eor\`eme \ref{thprin}.

\section{Retour sur le demi-plan} \label{retour}

Nous reprenons les notations de la partie \ref{dp}.

\subsection{Un analogue du th\'eor\`eme principal sur un anneau de coefficients fortement
  banal}

Fixons un nombre premier $l\neq p$ fortement banal pour $PGL_d(K)$, au sens de \ref{bonbanal}.
Dans cette section, on \'etudie la cohomologie de $\Omega_K^{d-1}$ \`a coefficients dans un anneau du type $\Lambda=\OC_{E}/\lambda^n\OC_E$ o\`u $E$ est une extension finie de $\QM_l$ et $\lambda$ une uniformisante de $E$. 

On reprend les notations $\pi^\Lambda_I$, pour $I\subseteq S$ de \ref{defrep} et la num\'erotation de $S$ de \ref{rappelclassif}. On sait par \ref{cohodp} ou \cite{O2} que la cohomologie de $\Omega_K^{d-1}$ est encore donn\'ee par
$$ H^{d-1+i}_c(\Omega_K^{d-1,ca},\Lambda) \simto \pi^\Lambda_{\{1,\cdots,i\}}\otimes
|.|^{-i}.$$
Comme $l$ est banal, on a $\Lambda[X]=(X-q^i)\Lambda[X]+(X-q^j)\Lambda[X]$ pour tous $0\leq i < j\leq d-1$ et on en d\'eduit gr\^ace au lemme \ref{actionphi} l'existence pour tout rel\`evement de Frobenius g\'eom\'etrique $\phi$ d'un unique scindage 
$$\alpha_\phi :\; R\Gamma_c(\Omega_K^{d-1,ca},\Lambda)\simto \bigoplus_{i=0}^{d-1} \pi^\Lambda_{\{1,\cdots,i\}}[d-1+i]$$
 comme en \ref{scfrob}. Comme $l$ est fortement banal, le th\'eor\`eme \ref{theoext} permet d'appliquer la m\^eme discussion que celle pr\'ec\'edant la proposition \ref{propdp2} ; \`a tout choix de g\'en\'erateurs des $\Lambda$-modules libres de rang $1$ $\ext{1}{\pi^\Lambda_{\{1,\cdots,i\}}}{\pi^\Lambda_{\{1,\cdots,i\}}}{\Lambda G}$ est associ\'e un isomorphisme 
 $$\beta^{-1}\circ \alpha_{\phi *}:\;\;\endo{D^b_\Lambda(PG_d)}{R\Gamma_c(\Omega_K^{d-1,ca},\Lambda)} \simto \TC_d^\Lambda $$
o\`u $\TC_d^\Lambda$ d\'esigne la $\Lambda$-alg\`ebre des matrices triangulaires sup\'erieures.

Notons comme d'habitude $\gamma: W_K \To{} \endo{D^b_\Lambda(PG_d)}{R\Gamma_c(\Omega_K^{d-1,ca},\Lambda)}^\times $ le morphisme de groupes donnant l'action de $W_K$ sur ${R\Gamma_c(\Omega_K^{d-1},\Lambda)}$.
Comme le groupe des matrices unipotentes sup\'erieures est un $l$-groupe, on constate que $\gamma_{|I_K}$ se factorise par le $l$-quotient de $I_K$. Comme $l$ est banal pour $PGL_d(K)$, on a en particulier $l>d$ donc le logarithme, resp. l'exponentielle, de tout \'el\'ement unipotent, resp. nilpotent, d'ordre $\leq d$ est d\'efini. Cela permet de d\'efinir l'\'el\'ement nilpotent $N_\mu \in \endo{D^b_\Lambda(PG_d)}{R\Gamma_c(\Omega_K^{d-1},\Lambda)}$ comme en \ref{definitionN}.
Le morphisme $\gamma$ est alors enti\`erement d\'etermin\'e par l'\'equation 
$$ (\beta\circ \alpha_{\phi *})(\gamma(w)) = \hbox{Diag}(1,q^{\nu(w)},\cdots, q^{(d-1)\nu(w)}) \times \hbox{exp}\left(t_\mu(i_\phi(w)).(\beta\circ\alpha_{\phi *})(N_\mu)\right). $$

Comme en \ref{diag}, si $(E_{i,j})_{i\leq j}$ d\'esigne la base canonique de $\TC_d^\Lambda$, alors il existe $a_1,\cdots, a_{d-1} \in \Lambda$ tels que 
$$ \beta\circ\alpha_{\phi *}(N_\mu)= \sum_{i=1}^{d-1} a_i E_{i-1,i}. $$
En conjuguant par une matrice diagonale, ce qui revient \`a changer le choix de g\'en\'erateurs qui d\'efinit  $\beta$, on peut mettre les $a_i$ sous la forme $a_i=\lambda^{n_i}$. Le r\'esultat suivant demandera plus d'efforts que  son analogue dans le cas $\o\QM_l$ :
\begin{prop} \label{proplambda}
On a
$ \lambda^{n-1} . N_\mu^{d-1} \neq 0 $ dans $\endo{D^b_\Lambda(G)}{R\Gamma_c(\Omega_K^{d-1},\Lambda)}$.
(Rappelons que $\Lambda=\OC_E/\lambda^n\OC_E$).
\end{prop}

D\'efinissons maintenant des repr\'esentations
$\tau_I^{\phi,\mu}(\Lambda)$ de $W_K$ \`a coefficients dans $\Lambda$
par les m\^emes formules que dans le paragraphe \ref{deftauI} en
rempla{\c c}ant simplement $\o\QM_l$ par $\Lambda$.
Le corollaire suivant est le strict analogue de la proposition \ref{propdp2} :

\begin{coro} Pour tout rel\`evement de Frobenius g\'eom\'etrique $\phi$,
il existe un unique scindage $\alpha_\phi$ de $R\Gamma_c(\Omega_K^{d-1,ca},\o\QM_l)$ et un unique choix de g\'en\'erateurs \ref{generateurs} d\'efinissant un isomorphisme $\beta_\phi$ comme en \ref{beta} tels
que le diagramme suivant soit commutatif
$$\xymatrix{   \endo{D^b_\Lambda(PG_d)}{R\Gamma_c(\Omega_K^{d-1,ca},\Lambda)}
  \ar[r]^-{\beta_\phi^{-1}\alpha_{\phi*}} & \TC_d \\ W_K \ar[u]^{\gamma}
  \ar[ur]_{\tau_\emptyset^{\phi,\Lambda}} & }
$$
De plus, le choix de g\'en\'erateurs (et donc $\beta_\phi$) est en fait ind\'ependant de $\phi$.
\end{coro}

\`A partir de l\`a on d\'emontre \`a partir du calcul d'extensions de \ref{theoext}, et exactement comme dans le paragraphe \ref{preuvetheodp}, le r\'esultat suivant :
\begin{theo} 
  Pour tout $I\subseteq S$, il existe  un isomorphisme $W_K$-{\'e}quivariant
$$ \HC^*(R\hom{R\Gamma_c(\Omega_K^{d-1,ca},\Lambda)}{\pi^\Lambda_I}{D^b_\Lambda(PG_d)})
\simto \tau_I^{\Lambda}.$$
\end{theo}
Lorsque $\Lambda = \o\FM_l$ (donc avec $l$ fortement banal), on obtient une r\'ealisation de la correspondance de Langlands-Vign\'eras, puisque, comme dans le cas classique, on a 
$$ \tau_I^{\o\FM_l}\simeq \sigma_d(\pi^{\o\FM_l}_I)\otimes |.|^{\frac{d-1}{2}} .$$

\alin{Preuve de la proposition \ref{proplambda}}
Dans le cas $l$-adique, la preuve  de la
non-nullit\'e de $N^{d-1}$ en \ref{Nnonnul}
reposait sur la suite spectrale de Rapoport-Zink associ\'ee \`a un
quotient du sch\'ema formel $\wh\Omega_K^{d-1}$ de Deligne par un
sous-groupe $\Gamma$ de $PG_d$ discret, sans torsion et
cocompact. Elle utilisait deux propri\'et\'es de cette suite spectrale
vraiment li\'ees aux coefficients $l$-adiques : la
d\'eg\'en\'er\'escence en $E_2$ et l'existence d'un accouplement
``d\'efini positif'' sur $E_1^{1-d,2d-2}$. Dans le cas que l'on
envisage maintenant, la d\'eg\'en\'er\'escence en $E_2$ sera encore
v\'erifi\'ee, toujours par un argument de poids, {\em cf} \ref{degener} ci-dessous. Par
contre, pallier l'absence de positivit\'e semble plus difficile : cela
demande des informations pr\'ecises sur $\Gamma$, par exemple son
covolume ou la combinatoire de son action sur l'immeuble, etc... 

La strat\'egie que nous adoptons ici consiste \`a uniformiser
``partiellement'' : nous quotienterons simplement par un sous-groupe discret
cocompact et sans torsion {\em d'un tore maximal de $PG_d$}. Le
sch\'ema formel obtenu $\wh\Omega/\Gamma$ n'est plus alg\'ebrisable,
mais nous allons v\'erifier que le formalisme de Rapoport-Zink s'y applique sans probl\`emes. Les d\'efinitions et r\'esultats qui suivent ne sont cetainement pas optimaux, mais suffiront pour nos affaires.

\begin{DEf} \label{defsest}
Un sch\'ema formel $\XG$ quasi-s\'epar\'e au-dessus de $\hbox{Spf}(\OC_K)$ sera dit {\em fortement semi-stable} si :
\begin{enumerate}
	\item localement pour la topologie de Zariski, $\XG$ est isomorphe \`a $\hbox{Spf}(\OC_K\la T_1,\cdots, T_n\ra /(T_1\times\cdots\times T_m -\varpi))$
	\item Les composantes irr\'eductibles de la fibre sp\'eciale $Y=\XG\times_{\OC_K} k$ sont propres et lisses de dimension constante sur $k$. 
\end{enumerate}
\end{DEf}

Soit $I$ l'ensemble des composantes irr\'eductibles de $Y$. On met une structure simpliciale sur $I$ en d\'eclarant qu'un sous-ensemble $J\subset I$ est un simplexe si $Y_J:=\bigcap_{j\in J} Y_j \neq \emptyset$. Dans ce cas, si $d-1$ d\'esigne la dimension commune des composantes irr\'eductibles de $Y$, alors $Y_J$ est une sous-vari\'et\'e lisse de pure dimension $d-|J|$. On note $I^{(m)}$ l'ensemble des simplexes de dimension $m-1$ (donc de cardinal $m$) de $I$.

Notons maintenant $X$ la fibre g\'en\'erique au sens de Berkovich de $\XG$ : c'est un espace $K$-analytique lisse.

\begin{prop} \label{RZschfor}
Soit $\XG$ un sch\'ema formel fortement semi-stable sur $\OC_K$. Alors les modules $H^q_c(X^{ca},\Lambda)$ sont mod\'er\'ement ramifi\'es et 
il existe une suite spectrale
$$ E_1^{-r,q+r}= \bigoplus_{k \geq \hbox{sup}(0,-r)} \bigoplus_{J\in I^{(r+2k+1)}} H^{q+1-|J|}(Y_J^{ca},\Lambda(-r-k)) \Rightarrow H^q_c(X^{ca},\Lambda)$$
dont la diff\'erentielle $d_1$ de degr\'e $(-1,1)$ en $(r,q)$ s'\'ecrit
\begin{eqnarray*}
	d_1^{-r,q+r} & = & \sum_{k\geq \hbox{sup}(o,-r)} (-1)^k \sum_{J\in I^{(r+2k+1)}} \left(\sum_{K\supset J} \varepsilon_{JK} \hbox{Res}_{JK}\right)   \\
	 & + & \sum_{k\geq \hbox{sup}(o,-r+1)} (-1)^{r+k} \sum_{J\in I^{(r+2k+1)}} \left( \sum_{K\subset J} \varepsilon_{JK} \hbox{Gys}_{JK}\right)
\end{eqnarray*}
o\`u
$$ \hbox{Res}_{JK} :\;\; H^{q+1-|J|}(Y^{ca}_J,\Lambda) \To{} H^{q+1-|J|}(Y^{ca}_K,\Lambda)=H^{(q+1)+1-|K|}(Y^{ca}_K,\Lambda) $$
est le morphisme de restriction associ\'e \`a $Y_K\injo Y_J$ lorsque $J\subset K$, 
$$ \hbox{Gys}_{JK} :\;\; H^{q+1-|J|}(Y^{ca}_J,\Lambda) \To{} H^{q+1-|J|+2}(Y^{ca}_K,\Lambda)=H^{(q+1)+1-|K|}(Y^{ca}_K,\Lambda) $$
est le morphisme de Gysin associ\'e \`a $Y_J\injo Y_K$ lorsque $K\subset J$, et les $\varepsilon_{J,K}$ sont des signes donn\'es par le choix d'une orientation de l'ensemble simplicial $I^{(\bullet)}$.

De plus, si $\mu$ est un g\'en\'erateur topologique de $\ZM_l(1)$,
l'action de $\mu-1$   sur l'aboutissement est induite par l'endomorphisme $\nu$ de degr\'e $(-2,0)$ en $(r,q)$  d\'efini par
$$ \nu^{-r,q+r} =\sum_{k \geq \hbox{sup}(0,-r+1)} \sum_{J\in I^{(r+2k+1)}} \id_{H^{q+1-|J|}(Y_J^{ca},\Lambda(-r-k)))} \otimes \mu$$
\end{prop}

La preuve de cette proposition consiste \`a v\'erifier que le formalisme de Rapoport-Zink dans \cite{RZss} s'applique encore aux cycles \'evanescents d\'efinis par Berkovich. Nous la reportons \`a l'appendice.

On suppose toujours que la dimension commune des composantes de $Y$ est $d-1$.
Dans les cas particuliers $(r,q)=(d-1,d-1)$ et $(r,q)=(1-d,d-1)$, 
la description de la diff\'erentielle et de l'op\'erateur $\nu$ s'ins\`erent dans 
le diagramme suivant :
$$\xymatrix{E_1^{1-d,2d-2} = \Lambda[I^{(d)}](1-d) \ar[d]_{d_1^{1-d,2d-2}}  \ar[r]^{\nu^{d-1}=\mu^{d-1}} &   \Lambda[I^{(d)}]= E_1^{d-1,0} \\ E_1^{2-d,2d-2} = \Lambda[I^{(d-1)}](1-d) & \Lambda[I^{(d-1)}]=E_1^{d-2,0} \ar[u]_{d_1^{d-2,0} }
}$$	
o\`u $\Lambda[?]$ d\'esigne le $\Lambda$-module des fonctions \`a support fini sur l'ensemble $?$ et 
$$ d_1^{1-d,2d-2}(f)(K) = \sum_{J\supset K} \varepsilon_{JK}.f(J) \;\;\hbox{ et }\;\;
 d_1^{d-1,0}(g)(J) = \sum_{K\subset J} \varepsilon_{JK}.g(K). $$
Munissons alors $\Lambda[?]$ de la forme bilin\'eaire "canonique" pour laquelle la base $?$ est autoduale. On constate que $d_1^{1-d,2d-2}$ et $d_1^{d-1,0}$ sont adjointes l'une de l'autre (en n\'egligeant la torsion \`a la Tate pour simplifier l'exposition). 
En particulier, on a $\im(d_1^{d-1,0}) \subset \ker(d_1^{1-d,2d-2})^\perp$.
On en d\'eduit le lemme suivant
\begin{lemme}
Si $\lambda^{n-1} \ker(d_1^{1-d,2d-2}) \subsetneq \ker(d_1^{1-d,2d-2})^\perp$, alors le morphisme induit par $\nu^{d-1}$
$$ \lambda^{n-1}\nu^{d-1} :\;\; E_2^{1-d,2d-2} \To{} E_2^{d-1,0} $$
est non-nul.
\end{lemme}

Venons-en maintenant \`a la situation qui nous int\'eresse, \`a savoir $\XG=\wh\Omega_K^{d-1}/\Gamma$ o\`u $\Gamma$ est un sous-groupe discret sans torsion de $PGL_d(K)$.  Lorsque $\Gamma=\{1\}$, on sait que $\wh\Omega_K^{d-1}$ est fortement semi-stable au sens de \ref{defsest}.
L'ensemble $I$ des composantes irr\'eductibles est en bijection avec l'ensemble $BT$ des sommets de l'immeuble de Bruhat-Tits de $PGL_d(K)$, et la structure simpliciale $BT^{(\bullet)}$ qui en d\'ecoule est celle de Bruhat-Tits. 
Par exemple, $BT^{(d)}$ est l'ensemble des chambres de l'immeuble et $BT^{(d-1)}$ celui des murs (facettes de codimension $1$).
En g\'en\'eral, pour que $\XG=\wh\Omega_K^{d-1}/\Gamma$ soit fortement semi-stable, il faut et il suffit  que l'\'etoile d'un sommet soit contenue dans un domaine fondamental de $\Gamma$ dans $BT^{(\bullet)}$ (sinon il y a des points doubles). 
On a alors $I^{(m)}=\Gamma\ba BT^{(m)}$. Nous fixons une orientation $SL_d(K)$-\'equivariante de $BT^{(\bullet)}$ de la mani\`ere suivante : on choisit d'abord un ordre sur les sommets d'une chambre $\Delta$ et on le transporte \`a toute autre chambre par l'action transitive de $SL_d(K)$, sans ambigu\"it\'e puisque le stabilisateur et le fixateur de $\Delta$ co\"incident dans $SL_d(K)$ ; on obtient donc une orientation sur chaque simplexe maximal et on v\'erifie que ces orientations se recollent bien le long des murs.
L'orientation de $BT^{(\bullet)}$ obtenue est \'equivariante par tout $\Gamma$ discret sans-torsion et induit donc une orientation de $\Gamma\ba BT^{(\bullet)}$.

Notons $A$ l'ensemble des sommets de 
l'appartement de $BT$ associ\'e au tore maximal $T$ de $PGL_d(K)$ que l'on a choisi il y a d\'eja longtemps, et $A^{(\bullet)}$ l'ensemble simplicial associ\'e. 
Fixons aussi une chambre $\Delta$ dans $A^{(d)}$. Nous noterons $A^{(d)}_{+}$, resp. $A^{(d)}_{-}$, l'ensemble des chambres de $A^{(d)}$ qui sont \`a distance {\em paire}, resp. {\em impaire} de $\Delta$.

\begin{lemme} Soit $\Gamma$ un sous-groupe discret sans torsion et cocompact de $T$ tel que le quotient $\wh\Omega_K^{d-1}/\Gamma$ soit fortement semi-stable. 
Alors les fonctions caract\'eristiques $1_{A^{(d)}_{\pm}}$  de $A^{(d)}_{\pm}$ dans $BT^{(d)}$ d\'efinissent deux \'el\'ements de $E_1^{1-d,2d-2}$ tels que,
 avec notre choix d'orientation 
\begin{enumerate}
	\item $1_{A^{(d)}_+} -  1_{A^{(d)}_-} \in \ker d_1^{1-d,2d-2}$
  \item $\lambda^{n-1} (1_{A^{(d)}_+} -1_{A^{(d)}_-}) \notin (\ker d_1^{1-d,2d-2})^\perp$ d\`es que $|\Gamma \ba A^{(d)}|$ est inversible dans $\Lambda$. 
\end{enumerate}
\end{lemme}
\begin{proof}
Le fait que ces fonctions caract\'eristiques induisent des \'el\'ements de $E_1^{1-d,2d-2} = \Lambda[\Gamma\ba BT^{(d)}]$ d\'ecoule de 
\begin{itemize}
	\item ${A^{(d)}_+}$ et ${A^{(d)}_-}$ sont stables par $\Gamma$ puisque un domaine fondamental de $\Gamma$ contient l'\'etoile d'un sommet.
	\item $\Gamma$ a un nombre fini d'orbites sans $A^{(d)}$ puisqu'il est cocompact dans $T$.
\end{itemize}

Montrons maintenant le point i). Soit $M$ un mur de $ BT^{(d-1)}$. L'ensemble des chambres de $\Gamma\ba A^{(d-1)}$ contenant $\Gamma.M$ est non vide seulement si $M\in A^{(d-1)}$, auquel cas il contient deux \'el\'ements, $\Gamma\Delta_+$ et $\Gamma\Delta_-$ o\`u
$\Delta_+$ et $\Delta_-$ sont les chambres mitoyennes de $A^{(d)}$ contenant $M$, 
la premi\`ere paire et la seconde impaire. On a alors
$$ d_1^{1-d,2d-2}(1_{A^{(d)}_+}-1_{A^{(d)}_-})(M) = \varepsilon_{\Gamma\Delta_+,\Gamma.M}-\varepsilon_{\Gamma\Delta_-,\Gamma.M}.$$
Mais par d\'efinition de notre orientation, puisque $\Delta_+$ se d\'eduit de $\Delta_-$ par un \'el\'ement de $SL_d(K)$ laissant fixe le mur mitoyen, on a $\varepsilon_{\Gamma\Delta_+,\Gamma.M}=\varepsilon_{\Delta_+,M}=\varepsilon_{\Delta_-,M}=\varepsilon_{\Gamma\Delta_-,\Gamma.M}$. 

Passons maintenant au point iii). En notant $\la .,. \ra$, la forme bilin\'eaire canonique de $E^{1-d,2d-2}_1$, on a
$$ \left\la \lambda^{n-1} (1_{A^{(d)}_+}-1_{A^{(d)}_-}), (1_{A^{(d)}_+}-1_{A^{(d)}_-}\right\ra = \lambda^{n-1} |\Gamma\ba A^{(d)}|, $$ 
d'o\`u l'assertion.
\end{proof}

Remarquons qu'il est facile de choisir $\Gamma$ tel que $|\Gamma\ba A^{(d)}|$ soit inversible dans $\Lambda$. C'est par exemple le cas pour le r\'eseau engendr\'e par les $\alpha^\vee(\varpi_K)$ o\`u $\alpha^\vee$ parcourt les racines simples de $T$ : dans ce cas on a en effet $|\Gamma\ba A^{(d)}|=|W_G|=d!$ qui est inversible dans l'anneau banal $\Lambda$.

Pour appliquer les deux lemmes pr\'ec\'edents \`a notre probl\`eme, il nous faut montrer la d\'eg\'eneresence de la suite spectrale.

\begin{lemme} \label{degener}(On rappelle que $\Lambda$ est en particulier banal pour $PGL_d(K)$).
Pour tout $\Gamma$ discret sans torsion dans $PGL_d(K)$, la suite spectrale de Rapoport-Zink de la proposition \ref{RZschfor} associ\'ee \`a $\XG=\wh\Omega_K^{d-1}/\Gamma$ d\'eg\'en\`ere en $E_2$. 
\end{lemme}
\begin{proof}
Il s'agit moralement, comme souvent dans ces situations, d'un argument de poids. Pour un $\XG$ g\'en\'eral, il n'y a pas de th\'eorie des poids "\`a valeurs dans $\Lambda$", mais dans notre cas les vari\'et\'es $Y_J$ ont une cohomologie particuli\`erement simple : nous allons expliquer (rappeler) pourquoi pour tout $J$, on a
\begin{enumerate}
	\item $H^{*}(Y^{ca}_J,\Lambda)$ est non nul seulement si $*\in 2\NM$ et $0\leq *\leq 2d-2|J|$.
	\item Le Frobenius g\'eom\'etrique $\sigma$ agit sur $H^{2i}(Y^{ca}_J,\Lambda)$ par multiplication par $q^i$.
\end{enumerate}

Rappelons tout d'abord la structure de $Y_J$ : soit $B^n$ la $k$-vari\'et\'e obtenue en \'eclatant $\PM^n$ le long du  produit des (id\'eaux d\'efinissant ses) sous-vari\'et\'es lin\'eaires $k$-rationnelles. Alors pour tout $J$, les composantes connexes de $Y_J$ sont de la forme $ B^{d_1}\times\cdots \times B^{d_{|J|}}$ avec $d_1+\cdots +d_{|J|}=d-|J|$, {\em cf} \cite[par. 6]{Ito}. Il suffit donc de prouver les deux assertions ci-dessus pour les vari\'et\'es $B^n$. Mais on peut obtenir ces vari\'et\'es selon la suite d'\'eclatements plus agr\'eables suivante ({\em cf } \cite[par. 4]{Ito}) :
$$ B^n=Y_{n-1}\To{\psi_{n-1}} Y_{n-2}\To{} \cdots \To{} Y_1\To{\psi_1} Y_0=\PM^n $$
o\`u $Y_{m+1}$ est l'\'eclat\'e de $Y_m$ le long de la r\'eunion $Z_m$ des {\em transform\'es stricts}  des
 sous-vari\'et\'es lin\'eaires de dimension $m$ de $\PM^n$ dans $Y_m$. 
L'avantage est qu'\`a chaque \'etape, la sous-vari\'et\'e $Z_m$ le long de laquelle on \'eclate est lisse de pure dimension $m$.

Alors d'apr\`es [SGA5], VII, 8.5 : on a pour tout $m$ des suites exactes :
$$ 0\To{} H^{*+2m-2n}(Z^{ca}_m,\Lambda) \To{} H^{*-2}(\psi_m^*(Z_m)^{ca},\Lambda) \oplus H^*(Y_m^{ca},\Lambda) \To{} H^*(Y_{m+1}^{ca},\Lambda) \To{} 0 .$$
Or, par construction, les composantes connexes de $Z_m$ sont toutes isomorphes \`a $B^m$. De plus, $\psi_m:\;\psi_m^*(Z_m) \To{} Z_m$ est un $\PM^{n-m-1}$-fibr\'e. On montre alors les deux assertions annonc\'ees par double r\'ecurrence sur $n$ et $m$.

\medskip

Appliquons ceci \`a la suite spectrale de Rapoport-Zink : on obtient que $E_1^{-r,m+r}$ est non-nul seulement si $m+r$ est pair auquel cas l'action du Frobenius g\'eom\'etrique $\sigma$ sur $E_1^{-r,m+r}$ est la multiplication par $q^{(m+r)/2}$. Rappelons que $0\leq m+r \leq 2d-2$ et, sous l'hypoth\`ese $\Lambda$ banal, on a $\Lambda[X]=(X-q^i)\Lambda[X]+(X-q^j)\Lambda[X]$ pour tous $0\leq i, j\leq d-1$. Comme les diff\'erentielles $d_n^{-r,m+r}$ sont des applications $\sigma$-\'equivariantes $E_n^{-r,m+r} \To{}E_n^{-r+n,m+r-n+1}$, elles sont nulles d\`es que $n\geq 2$.

\end{proof}

D'apr\`es les trois lemmes pr\'ec\'edents, si $\Gamma$ est le sous-groupe discret sans torsion et cocompact de $T$ engendr\'e par les $\alpha^\vee(\varpi_K)$, alors le morphisme induit par $\nu^{d-1}$
$$ \lambda^{n-1}\nu^{d-1} :\;\; E_\infty^{1-d,2d-2} \To{} E_\infty^{d-1,0} $$
dans la suite spectrale de Rapoport-Zink du quotient $\wh\Omega_K^{d-1}/\Gamma$ est non-nul.
D'apr\`es la fin de la proposition \ref{RZschfor}, cela implique que l'action de $(1-\mu)^{d-1}$  sur $H^{d-1}(\Omega^{d-1,ca}/\Gamma,\Lambda)$ n'est pas annul\'ee par $\lambda^{n-1}$. 
D'apr\`es le th\'eor\`eme \ref{theoHS}, cette action n'est donc pas annul\'ee non plus sur $R\Gamma_c(\Omega^{d-1,ca}_K,\Lambda)$. On en d\'eduit la proposition \ref{proplambda}.

\subsection{La suite spectrale de Rapoport-Zink pour les sch\'emas formels fortement semi-stables}

Le but de ce paragraphe est de prouver la proposition \ref{proplambda}.

\section{Appendice 1 : un calcul de la cohomologie du demi-plan}

Nous calculons la cohomologie \`a support compact des espaces sym\'etriques $\Omega_K^{d-1}$ de Drinfeld par une m\'ethode diff\'erente de celle de Schneider et Stuhler, qui pr\'esente l'avantage de fonctionner aussi pour les coefficients $l$-adiques. Au-lieu d'utiliser l'arrangement "\`a l'infini" des sous-vari\'et\'es lin\'eaires rationnelles de $\PM^{d-1}_K$, dont le nerf est l'immeuble de Tits, nous utilisons un recouvrement de $\Omega_K^{d-1}$ par des ouverts "distingu\'es" dont le nerf s'identifie  \`a l'immeuble de Bruhat-Tis. Cette approche est similaire \`a celle que Drinfeld a utilis\'ee pour calculer la cohomologie de $\Omega^1_K$. Nous tomberons sur des probl\`emes combinatoires qui ont d\'eja \'et\'e r\'esolus par Schneider et Stuhler : en cela, notre m\'ethode n'est pas tout-\`a-fait ind\'ependante de \cite{SS1}.

On peut remarquer qu'il y aurait une troisi\`eme approche naturelle pour calculer ces espaces de cohomologie : via les cycles \'evanescents du mod\`ele formel $\wh\Omega_K^{d-1}$ de Deligne, et par une des deux suites spectrales disponibles. \`A premi\`ere vue, les probl\`emes combinatoires qui se posent dans ces approches ont l'air encore diff\'erents.

\subsection{Le recouvrement ouvert de Drinfeld}

Nous consid\'erons ici encore $\Omega_K^{d-1}$ comme un espace de Berkovich. Rappelons  sa d\'efinition, au moins en termes ensemblistes : 
nommons l'espace vectoriel $V:=K^d$ 
et  $V^*$ son dual.
Notons respectivement $\AM(V)$ et $\PM(V)$ les $K$-espaces
analytiques associ{\'e}s par Berkovich {\`a} $V$, et
 $S(V^*)=\bigoplus_{n\in \NM} S_n(V^*)$ l'alg{\`e}bre sym{\'e}trique
associ{\'e}e {\`a} $V$. Alors on a  
les descriptions ensemblistes  
$$ \AM(V) =\{K-\hbox{seminormes multiplicatives}\;\; S(V^*)\To{}\RM_+ \}, $$
$$ \PM(V)= (\AM(V)\setminus\{0\})/\sim $$
$ \hbox{ o{\`u} } \; x\sim y
 \hbox{ \ssi\ } \exists \lambda>0, \forall n\in \NM, \forall f\in S_n(V^*),\;
|f(x)|=\lambda^n|f(y)| $ (suivant Berkovich, nous notons indiff{\'e}remment $x(f)=|f(x)|$ pour ``la semi-norme $x$ {\'e}valu{\'e}e en la fonction $f$''.). On renvoie au livre de Berkovich pour la d\'efinition des  topologies et des faisceaux d'anneaux topologiques correspondants. 
Posons maintenant
$$ \wt\Omega(V) : = \{x\in \AM(V),\;\; x_{|V^*} \hbox{ est une $K$-norme sur
  $V^*$} \}, $$
alors $\Omega^{d-1}_K$ s'identifie au sous-espace analytique $\Omega(V)$ de $\PM(V)$ dont l'ensemble sous-jacent est l'image de $\wt\Omega(V)$.  
Puisqu'on a une base de $V$, on a une identification $PGL(V)=PGL_{d}(K)$ et on peut utiliser les notations de \ref{defrep} et \ref{ss}.
Nous voulons prouver le r\'esultat suivant :
\begin{theo}  \label{cohodp} Soit $\Lambda$ un anneau de torsion premi\`ere \`a $p$ ou un anneau de valuation discr\`ete complet et de caract\'eristiques diff\'erentes de $p$.
  Pour $i=0,\cdots, d-1$, il existe des isomorphismes $PGL_d(K)\times W_K$-{\'e}quivariants
$$ H^{d-1+i}_c(\Omega_K^{d-1,ca},\Lambda) \simto \pi^\Lambda_{\{1,\cdots,i\}}(-i)$$
\end{theo}
Par la suite nous abr\`egerons $\Omega:=\Omega(V)=\Omega^{d-1}_K$ et $G:=PGL(V)$.

\ali On notera $BT$ l'immeuble de Bruhat-Tits associ{\'e} {\`a} $PGL(V^*)$
consid{\'e}r{\'e} comme ensemble simplicial et $|BT|$ sa r{\'e}alisation
g{\'e}om{\'e}trique, consid{\'e}r{\'e} comme immeuble euclidien. Rappelons la
description ensembliste
 $$ |BT| = \{ \hbox{classes d'homoth{\'e}tie de $K$-normes sur $V^*$} \}
 $$
qui permet de d{\'e}finir {\`a} moindres frais l'application de Drinfel'd 
$$\application{\tau:\;}{\Omega}{|BT|}{x}{x_{|V}:\; v\mapsto |v(x)|} .$$
D'apr\`es \cite{BicCRAS}, cette application est surjective et continue.

\ali
Pour $q=0,\cdots,d$, on note $BT_q$ l'ensemble des $q$-simplexes de $BT$.
Si $F\subset BT$ est un simplexe, on notera  $|F|$ la facette de
$|BT|$ associ{\'e}e : on a pour toute paire de simplexes $F'\subset F
\equ |F'|\subset \o{|F|}$ o{\`u} $\o{|F|}$ d{\'e}signe l'adh{\'e}rence de
$|F|$.
Pour un simplexe $F$, on notera aussi
 $|F|^*\subset |BT|$ son lien
: c'est la r{\'e}union des facettes $|F'|$ avec $F'\supset F$. Si
$F=\{s_1,\cdots,s_n\}$ avec $s_i\in BT$ des sommets, on a bien-s{\^u}r
$|F|^*=\cap_i |s_i|^*$. Il s'ensuit que la famille $(\tau^{-1}(|s|^*))_{s\in BT}$ de
sous-ensembles de $\Omega$ est un recouvrement ouvert analytique de
$\Omega$ dont le nerf est justement $BT$.
De plus, chaque $\tau^{-1}(|s|^*)$ est distingu\'e au sens de Berkovich, {\em i.e.} peut s'\'ecrire comme diff\'erence de deux domaines analytiques compacts.

 Le groupe $G=PGL(V)$ agit par automorphismes analytiques sur
$\Omega$ et simpliciaux sur $|BT|$. L'application $\tau$ est
$PGL(V)$-{\'e}quivariante. Notre but est \'evidemment de calculer la cohomologie au moyen du recouvrement par les $\tau^{-1}(|s|^*)$. Comme on veut un calcul \'equivariant et comme il n'existe pas d'orientation $PGL(V)$-invariante sur $BT$, on modifie le complexe de Cech usuel selon la proc\'edure introduite par Schneider-Stuhler dans \cite{SS1} et \cite{SS2}.
En fait, pour nous raccrocher \`a certains r\'esutats combinatoires de \cite{SS1}, nous utiliserons leur notion de syst\`emes de coefficients que nous rappelons ci-dessous.

\subsection{Syst{\`e}mes de coefficients sur l'immeuble}

\alin{Syst\`emes de coefficients} (voir Schneider-Stuhler \cite{SS1},\cite{SS0},\cite{SS2}).  Si $\CC$ est une cat{\'e}gorie, un syst{\`e}me de coefficients {\`a} valeurs dans
$\CC$ sur $BT$ est la donn{\'e}e  
\begin{enumerate}
\item d'objets $X_F\in\CC$ pour chaque simplexe $F\subset BT$.
\item de morphismes $r^F_{F'} : X_F \To{} X_{F'}$ pour chaque paire de
  facettes telle que $F'\subset F$ satisfaisant la relation de
  transitivit{\'e} $r^{F'}_{F''}\circ r^{F}_{F'}=r^{F}_{F''}$.
\end{enumerate}
En d'autres termes, notons $\BC\TC$ la cat{\'e}gorie associ{\'e}e au
complexe simplicial $BT$ dont les objets sont les simplexes et les
morphismes sont induits par les inclusions de simplexes : un syst{\`e}me
de coefficients {\`a} valeurs dans $\CC$ est simplement un foncteur
contravariant $\BC\TC \To{} \CC$. Ils forment de mani{\`e}re claire une
cat{\'e}gorie, qui est ab{\'e}lienne lorsque $\CC$ l'est.

Supposons maintenant que $\CC$ est munie d'une action de $G=PGL(V)$ ;
pour tout $g\in G$, on note $g_\CC$ l'endofoncteur de $\CC$ associ{\'e}.
Une structure  $G$-{\'e}quivariante sur  un objet $X$ de $\CC$ est 
 une famille d'isomorphismes $\tau(g) : g_\CC(M)
\simto M$ satisfaisant la condition de cocycle usuelle.

\smallskip

\noindent{\em Exemples :}
\begin{itemize}
\item Si $\CC$ est la cat{\'e}gorie des $\Lambda$-modules munie de l'action
  triviale de $G$, un objet $G$-{\'e}quivariant est simplement un
  $\Lambda G$-module.
\item Si $\CC = \FC(\Omega^{ca},\Lambda)$, on fait agir $G$ via l'image directe
  $g\mapsto g_*$ et un objet $G$-{\'e}quivariant est simplement un
  faisceau $G$-{\'e}quivariant au sens usuel.
\end{itemize}

\smallskip

Nous appellerons syst{\`e}me de coefficients $G$-{\'e}quivariant {\`a}
valeurs dans $\CC$
tout syst{\`e}me de coefficients muni des donn{\'e}es suppl{\'e}mentaires
suivantes :
\begin{enumerate}
\item pour tout $g\in G$ et tout simplexe $F\subset BT$, un
  isomorphisme $g^X_F : g_\CC(X_F)\simto X_{g(F)}$ tel que
\item si $h,g\in G$ et $F$ est un simplexe, $h^X_{g(F)}\circ h_\CC(g^X_F)=
  (h\circ g)^X_{F}$ et
\item si $F'\subset F$, alors $r^{g(F)}_{g(F')}\circ g^X_F =
  g^X_{F'} \circ g_\CC(r^F_{F'})$.
\end{enumerate}
En d'autres termes, c'est un objet $G$-{\'e}quivariant de la
cat{\'e}gorie des syst{\`e}mes de coefficients {\`a} valeurs dans $\CC$
munie de l'action de $G$ suivante : $X \mapsto g_{\BC\TC}\circ
X\circ {g^{-1}}_{\CC}$.

\alin{Complexes de cha\^ines} 
Schneider et Stuhler ont d\'efini dans \cite[II.1]{SS2} une notion de facette orient\'ee de l'immeuble (poly)-simplicial associ\'e \`a un groupe $p$-adique quelconque.
Dans le cas $PGL(V)$ qui nous int\'eresse, la d\'efinition est plus \'el\'ementaire, {\em cf} \cite[par. 3]{SS0} :
\'etant donn\'es deux ordres totaux sur un simplexe de $BT$, on les d\'eclare \'equivalents s'ils se d\'eduisent l'un de l'autre par une permutation {\em paire} des sommets. 
Alors un {\em simplexe orient\'e} est par d\'efinition un couple $(F,c)$ form\'e d'un simplexe et d'une classe d'\'equivalence d'ordres totaux sur ce simplexe.
 On note $BT_{(q)}$ l'ensemble des facettes orient{\'e}es
  de dimension $q$. Si $F'\subset F$, on note $\partial^F_{F'}(c)$ l'orientation de $F'$ induite par l'orientation $c$ de $F$.

Supposons que $\CC$ est ab{\'e}lienne avec limites directes exactes.
Si $X$ est un syst{\`e}me de coefficients {\`a} valeurs dans $\CC$, on lui
associe le complexe born{\'e} d'objets de $\CC$, dit complexe de cha{\^\i}nes
de $X$
$$ C^{or}_c(BT,X) := \;\; C^{or}_c(BT_{(d)},X) \To{d_{d-1}} \cdots \To{d_0}
C^{or}_c(BT_{(0)},X) $$
o{\`u} 
$$C^{or}_c(BT_{(q)},X) := \ker (\id + \alpha)$$
avec la d{\'e}finition suivante de $\alpha$ :
 posons pour toute facette orient{\'e}e $X_{(F,c)}:=X_F$, alors 
  $$\alpha : \bigoplus_{(F,c)\in
    BT_{(q)}} X_{(F,c)} \To{} \bigoplus_{(F,c)\in
    BT_{(q)}} X_{(F,c)}$$
est le morphisme induit par $X_{(F,c)} \To{\id} X_{(F,-c)}$. 
La diff{\'e}rentielle est induite par les morphismes
$$ \bigoplus_{F'\subset F} r^F_{F'} :\;\; X_{(F,c)} \To{}
\bigoplus_{F'\subset F} X_{(F',\partial^F_{F'}(c))}.$$
Nous noterons $\HC_*(BT,X)$ les objets d'homologie du complexe
$C^{or}_c(BT,X)$.

Rappelons que $G$ agit sur les facettes orient{\'e}es de $BT_{(q)},
q\geq 0$, de sorte que si $X$ est un syst{\`e}me de coefficients
$G$-{\'e}quivariant, alors pour tout $q\geq 0$, l'objet
$M:=C^{or}_c(BT_{(q)},X)$ est par construction un objet $G$-{\'e}quivariant de $\CC$.
De plus les diff{\'e}rentielles sont compatibles {\`a} ces structures
$G$-{\'e}quivariantes, de sorte que les objets d' homologie
$\HC_*(BT,X)$ sont
eux-aussi $G$-{\'e}quivariants.


\ali
Pour tout simplexe $F$
dans $BT$, on note 
$$U_F^{ca}:=\tau^{-1}(|F|^*)\hat\otimes\wh{K^{ca}} \;\; \hbox{ et } j_F
: U_F^{ca} \injo \Omega^{ca}$$ l'immersion ouverte associ{\'e}e. 
Pour $p\in\NM$, on v\'erifie que les applications
\begin{itemize}
\item  $(F\in BT^{(\bullet)}) \mapsto H^p_c(U_F^{ca},\Lambda)$ et
\item $(F'\subset F) \mapsto (H^p_c(U_{F}^{ca},\Lambda)\To{}H^p_c(U_{F'}^{ca},\Lambda))$, (le morphisme induit par l'immersion ouverte $U_{F}^{ca}\subset  U_{F'}^{ca}$)
\end{itemize}
d\'efinissent un syst\`eme de coefficients $G$-\'equivariant  \`a valeurs dans les $\Lambda W_K$-modules que nous noterons simplement $F \mapsto H^p_c(U_F^{ca},\Lambda)$.


\begin{prop}
 \label{sspec}
Il existe une suite spectrale $G\times W_K$-\'equivariante
$$E_1^{pq} = C^{or}_c(BT_{(q)},F \mapsto H^p_c(U_F^{ca},\Lambda)) \Rightarrow H^{p+q}_c(\Omega^{ca},\Lambda) 
$$
dont la diff\'erentielle $d_1^{pq}$ est celle du complexe de cha\^ines du syst\`eme de coefficients $F\mapsto H^p_c(U_F^{ca},\Lambda)$.
\end{prop}

\begin{proof}
Seul le cas o\`u $\Lambda$ est de valuation discr\`ete complet demande peut-\^etre une preuve d\'etaill\'ee. Nous traitons tout-de-m\^eme le cas de torsion en m\^eme temps.

Si $\FC$ est un faisceau {\'e}tale ab\'elien sur $\Omega^{ca}$, on lui associe un syst{\`e}me de coefficients $\o\FC$ sur $BT$ {\`a} valeurs dans
$\FC(\Omega^{ca},\ZM)$ en posant :
\begin{itemize}
\item pour toute facette $F$, $\o\FC_F:={j_F}_!j_F^*(\FC)$.
\item pour tout couple de facettes $F'\subset \o{F}$, le morphisme
  $\o\FC_F \To{r^F_{F'}} \o\FC_{F'}$ est induit par l'immersion
  ouverte $U_F^{ca}\injo U_{F'}^{ca}$.
\end{itemize}
De plus, toute structure de faisceau $G$-{\'e}quivariant sur $\FC$
induit une structure $G$-{\'e}quivariante sur le syst{\`e}me de
coefficients $\o\FC$.

Le foncteur $\FC \mapsto C^{or}_c(BT,\o\FC)$ est visiblement exact, 
et on peut  augmenter le complexe $C^{or}_c(BT,\o\FC)$ par le morphisme 
${\epsilon_\FC:\;}{C^{or}_c(BT,\o\FC)} \To{} {\FC}$ somme des
morphismes canoniques ${j_s}_!j_s^* \To{} \id$ pour $s$ sommet de $BT$.

\begin{lemme}
 Soit $\FC$ un faisceau \'etale ab{\'e}lien (resp. ab{\'e}lien $G$-{\'e}quivariant) sur
 $\Omega^{ca}$. Alors 
 \begin{enumerate}
 \item le complexe augment{\'e}
 $C^{or}_c(BT,\o\FC)\To{\epsilon} \FC$ est une r{\'e}solution de $\FC$
 (resp. une r{\'e}solution $G$-{\'e}quivariante de $\FC$).
\item si $\FC_{|U}$ est $\Gamma_!(U,-)$-acyclique (sections \`a supports compacts) pour tout ouvert $U$ de $\Omega^{ca}$, alors $\CC^{or}_c(BT_{(q)},\o\FC)$ est $\Gamma_!(\Omega^{ca},-)$-acyclique.
\end{enumerate}
\end{lemme}
\begin{proof}
On laisse le lecteur  v\'erifier le i) sur les fibres. Pour le ii), on v\'erifie qu'une somme directe d'objets $\Gamma_!$-acycliques est $\Gamma_!$-acyclique. On est ramen\'e \`a prouver que pour $F$ facette de $BT_q$, ${j_F}_!j_F^*(\FC)$ est $\Gamma_!(\Omega^{ca},-)$-acyclique. Mais ceci vient de l'isomorphisme canonique $R\Gamma_!(U_F^{ca},-)\simeq R\Gamma_!(\Omega^{ca},-)\circ {j_F}_!$. 
\end{proof}



Si on applique le foncteur $\Gamma_!(\Omega^{ca},-)$ au 
syst{\`e}me  de
coefficients $\o\FC$, on obtient le syst\`eme de coefficients $\Gamma_!\o\FC$ en groupes ab\'eliens $F\mapsto \Gamma_!(U_F^{ca},\FC)$.
Si de plus $\FC$ est muni d'une structure $G$-{\'e}quivariante, alors le
syst{\`e}me de coefficients $\Gamma_!\o\FC$ l'est  aussi.

\def\u{\underline}

 Choisissons maintenant un complexe 
$$\IC^*:= \; \IC^0 \To{} \IC^1 \To{}  \cdots \To{} \IC^n \To{} \cdots $$
d'objets injectifs de $\FC(\Omega^{ca}, \Lambda (W_K\times G)_{disc})$ (notations de \ref{defladic}) 
qui soit quasi-isomorphe 
\begin{itemize}
	\item \`a $\underline\Lambda$ lorsque $\Lambda$ est de torsion, et
	\item \`a $R\limproj(\underline{\Lambda/\lambda^n\Lambda})$ lorsque $\Lambda$ est de valuation discr\`ete complet d'uniformisante $\lambda$.
\end{itemize}
On sait alors, {\em cf} \ref{defladic}, que l'objet $R\Gamma_c(\Omega^{ca},\Lambda)$ de  $D^+(\Lambda G_{disc})$ est repr\'esent\'e dans chacun des cas par le complexe $\Gamma_!(\Omega^{ca},\IC^*)$.

Consid\'erons maintenant le bicomplexe $C^{or}_c(BT,\o\IC^*)$. Celui-ci s'augmente dans le complexe simple $\IC^*$ et cette augmentation est un quasi-isomorphisme par le point i) du lemme pr\'ec\'edent. De plus, c'est un bicomplexe de faisceaux $\Gamma_!(\Omega^{ca},-)$-acycliques par le point ii)  donc $R\Gamma_c(\Omega^{ca},\Lambda)$, dans chacun des cas, est aussi repr\'esent\'e par le complexe simple associ\'e au bicomplexe 
$C^{or}_c(BT,\Gamma_!\IC^*)$ :
$$ {\mathbf s}\left( C^{or}_c(BT,\Gamma_!\IC^*) \right) \simto R\Gamma_c(\Omega^{ca},\Lambda)$$
En filtrant le bicomplexe $C^{or}_c(BT,\Gamma_!\IC^*)$  par l'indice $*$ et parce que le foncteur "complexe de cha\^ines" est exact en la variable "syst\`eme de coefficients", on obtient une suite spectrale
$$ E_1^{pq} = C^{or}_c(BT_{(q)},\HC^p(\Gamma_!\IC^*)) \Rightarrow H^{p+q}_c(\Omega^{ca},\Lambda) $$
dont la diff\'erentielle $d_1^{pq}$ est celle du complexe de cha\^ines du syst\`eme de coefficients $F\mapsto \HC^p(\Gamma_!(U_F^{ca},\IC^*))$ (cohomologie du complexe $\Gamma_!(U_F^{ca},\IC^*)$).
Il ne nous reste plus qu'\`a v\'erifier qu'on a des isomorphismes 
$$ \Gamma_!(U_F^{ca},\IC^*) \simto R\Gamma_c(U_F^{ca},\Lambda) $$
fonctoriels en les immersions $U_F\injo U_{F'}$. Ceci est \'evident dans le cas de torsion mais demande une explication dans le cas $\lambda$-adique. Dans ce cas,
 on sait que $R\Gamma_c(U_F^{ca},-)=R\Gamma_!(U_F^{ca},-) \circ R\limproj $, donc il nous suffit de v\'erifier que pour tout ouvert $U$, on a un isomorphisme de foncteurs $D^+(\Omega^{ca},\Lambda_\bullet G_{disc}) \To{} D^+(U^{ca},\Lambda G_{disc})$ 
$$ j_U^*\circ R\limproj \simto R\limproj \circ j_U^*. $$
Mais ceci est prouv\'e dans \cite[4.2.7]{Fargues} (c'est d'ailleurs tr\`es simple : l'isomorphisme sans le $R$ (pour des faisceaux) est \'el\'ementaire et $j_U^*$ envoie injectifs sur injectifs donc on peut le d\'eriver).


\end{proof}


\ali Rappelons maintenant l'exemple essentiel de syst{\`e}me de coefficients de
$\Lambda$-modules introduit par Schneider et Stuhler. Pour $F$
simplexe, on note  $\wh{G}_F$ son stabilisateur, $G_F$ son fixateur et
$G_F^+$ le pro-$p$-radical de $G_F$ (voir par exemple \cite[part 6]{SS1} o\`u
le groupe $G_F^+$ est not\'e $U_\sigma$). Partant alors d'un
$\Lambda G$-module lisse $V$ on d{\'e}finit le syst{\`e}me de coefficients
$\gamma_0(V)$ par :
\begin{itemize}
\item Pour tout simplexe $F$, ${\gamma_0(V)}_F:=V^{G_F^+}$.
\item Pour toute paire $F'\subset F$, $r^F_{F'}$ est l'inclusion
  $V^{G_F^+}\subset V^{G_{F'}^+}$ dont l'existence est assur{\'e}e par
  l'inclusion $G_{F'}^+\subset G_F^+$.
\end{itemize}
Bien s{\^u}r, l'action de $G$ sur $V$ induit une structure
$G$-{\'e}quivariante sur $\gamma_0(V)$.

Le point clef de la preuve du th\'eor\`eme \ref{cohodp} est la proposition suivante qui sera prouv\'ee dans la prochaine section :
\begin{prop} \label{cle}
Pour $0\leq i\leq d-1$, 
il existe un isomorphisme de syst\`emes de coefficients $G$-\'equivariants en $\Lambda W_K$-modules 
$$ \left(F\mapsto  H^{d-1+i}_c(U_F^{ca},\Lambda)\right) \simto \gamma_0(\pi^\Lambda_{\{1,\cdots,i\}}(-i)). $$
De plus, pour $p \notin \{d-1,\cdots,2d-2\}$, on a $H^p_c(U_F^{ca},\Lambda)=0$ pour toute $F$.
\end{prop}

\`A partir de l\`a, la preuve de \ref{cohodp} est facile : on sait par \cite[Thm 6.8]{SS1} que les complexes $C^{or}_c(BT,\gamma_0(\pi^\Lambda_I))$ sont des r\'esolutions des repr\'esentations $\pi^\Lambda_I$. La suite spectrale \ref{sspec} d\'eg\'en\`ere donc en des isomorphismes
$$ E_2^{*0} \simto H^*_c(\Omega^{ca},\Lambda)$$
avec $E_2^{*0}=  \pi_{\{1,\cdots,i\}}^\Lambda(-i)$ si $*=d-1+i$ pour $0\leq i\leq d-1$ et $E_2^{*0}=0$ sinon.

\subsection{Arrangements}

Nous commen{\c c}ons par quelques pr{\'e}liminaires simpliciaux et homologiques.

\alin{Ensembles partiellement ordonn{\'e}s} 
\label{poset}
Soit $(\XG,\leq)$ un ensemble
muni d'un ordre partiel. On peut lui associer une
cat{\'e}gorie not{\'e}e encore $\XG$ dont les objets sont les {\'e}l{\'e}ments de
$\XG$ et les ensembles de morphismes sont singletons ou vides selon
que $X \leq X'$ ou non, avec une loi de composition {\'e}vidente. On
appelle syst{\`e}me de coefficients sur $X$ {\`a} valeurs dans une cat{\'e}gorie
ab{\'e}lienne $\CC$ tout foncteur $\XG \To{} \CC$.

Supposons que $\XG$ admet un plus petit {\'e}l{\'e}ment $Z$ et notons
$\XG^{(0)}:=\XG\setminus\{Z\}$.  
On peut associer \`a $\XG$ 
l'ensemble simplicial dont les sommets sont les \'el\'ements de
$\XG^{(0)}$ et les $n$-simplexes sont les sous-ensembles
totalement ordonn{\'e}s {\`a} $n+1$ {\'e}l{\'e}ments de $\XG^{(0)}$. Pour un tel simplexe $S$ on
pose $X_S:=\hbox{max}(S)$ l'{\'e}l{\'e}ment maximal de $S$. On notera aussi
$\XG^{(n)}$ l'ensemble de ces $n$-simplexes. 
Un syst{\`e}me de
coefficients $\VC :\; X\mapsto \VC_X$ sur $\XG$ induit un syst{\`e}me de coefficients sur
l'ensemble simplicial associ{\'e} {\`a} $\XG$ en posant :
\begin{itemize}
\item $\VC_S:=\VC_{X_S}$ et
\item pour tout $S\subset S'$, la fl{\`e}che $\iota_S^{S'}: \VC_S\To{}
  \VC_{S'}$ est donn{\'e}e par la 
  fl{\`e}che associ{\'e}e par $\VC$ {\`a} $X_S\leq X_{S'}$.
\end{itemize}

Soient $S\subset S'$ deux simplexes tels que $|S'|=|S|+1$. On d{\'e}finit
un signe
$\epsilon_{S,S'}:=(-1)^{\hbox{rang}_S(X')}$ o{\`u} $S'=S \cup
  \{X'\}$ et $\hbox{rang}_S(X')=|\{X\in S, X\leq X'\}|$. 
{\`A} un syst{\`e}me de coefficients $\VC$ sur $\XG$ on associe alors le
complexe de cocha{\^\i}nes suivant :
$$ (\CC^*(\XG,\VC),\partial^*) :\;\; \VC_Z  \To{\partial^0}
\bigoplus_{S\in \XG^{(0)}}
\VC_S \To{\partial^1} \cdots \To{\partial^{k}} 
\bigoplus_{S\in\XG^{(k)}} \VC_S \To{\partial^{k+1}} \cdots $$
o{\`u} la diff{\'e}rentielle est d{\'e}finie par $\partial^0=\sum_{X\in\XG^{(0)}}
\iota_Z^{X}$ et
$$\partial^{k+1}=\bigoplus_{S\in\XG^{(k)}}\left(\sum_{S'\in\XG^{(k+1)},S'\supset S}
  \epsilon_{S,S'} \iota_S^{S'} \right) \;\;\hbox{ pour } k\geq 0.$$

\alin{Arrangements} \label{arr}
Soit $Z$ un $K$-espace analytique compact. On appellera arrangement
dans $Z$ tout ensemble fini $\XG$ de sous-espaces analytiques compacts
stable par intersection et contenant $Z$. En particulier $\XG$ est un ensemble
partiellement ordonn{\'e} pour la contenance (c'est {\`a} dire $X\leq X'$
\ssi\ $X\supseteq X'$) muni d'un plus petit {\'e}l{\'e}ment, $Z$. 
On note $d_\XG$ la
dimension ({\em i.e. } son cardinal moins $1$) d'un simplexe maximal.

On s'int{\'e}resse {\`a} la cohomologie de l'ouvert compl{\'e}mentaire $U(\XG):=Z\setminus
\bigcup_{X\in \XG^{(0)}} X$. Notons $j_U:U(\XG)\injo Z$ et $i_X: X\injo Z , X\in \XG$ les inclusions canoniques. 
Comme dans la section pr\'ec\'edente, on dispose de  syst\`emes de coefficients 
\begin{itemize}
\item $(X\in \XG)\mapsto H^p_c(X^{ca},\Lambda)$ et
\item $(X'\subset X) \mapsto \iota_X^{X'} : 
 (H^p_c(X^{ca},\Lambda)\To{} H^p_c({X'}^{ca},\Lambda)$, la fl{\`e}che \'etant induite par
  l'inclusion $X'\injo X$.
\end{itemize}
On supposera pour simplifier que les espaces compacts $Z$ et $X\in\XG$ sont "quasi-alg\'ebriques" ({\em i.e.} satisfont les hypoth\`eses de \cite[Cor 5.6]{Bic3}). On sait alors que
dans le cas o\`u $\Lambda$ est $\lambda$-adique, on a simplement $H^p_c(X^{ca},\Lambda)=\limproj\, H^p(X^{ca},\Lambda/\lambda^n)$ et les fl\`eches de transition ci-dessus se d\'eduisent des fl\`eches naturelles en cohomologie (sans supports) de torsion.

\begin{prop}
 \label{sspec2}
Il existe une suite spectrale $W_K$-\'equivariante
$$E_1^{pq}(\XG) = \CC^{q}(\XG,X\mapsto H^p_c(X^{ca},\Lambda)) \Rightarrow H^{p+q}_c(U(\XG)^{ca},\Lambda) 
$$
dont la diff\'erentielle $d_1^{pq}$ est celle du complexe de cocha\^ines du syst\`eme de coefficients $X\mapsto H^p_c(X^{ca},\Lambda)$.
\end{prop}

\begin{proof}
Encore une fois, ceci est compl\`etement standard dans le cas de torsion. Par prudence, nous expliquons les arguments dans le cas $\lambda$-adique.

Comme dans la section pr\'ec\'edente, tout faisceau ab{\'e}lien $\FC$ sur $Z$ d{\'e}finit un syst{\`e}me de
coefficients $\o\FC$ sur $\XG$ {\`a} valeurs dans $\FC(Z,\ZM)$ en posant 
\begin{itemize}
\item $\o\FC_X:=i_{X,*}\circ i_X^* (\FC)$.
\item la fl{\`e}che $\iota_X^{X'}$ pour $X\leq X'$ est donn{\'e}e par
  l'inclusion $X'\injo X$.
\end{itemize}
Le complexe $\CC^*(\XG,\o\FC)$ s'augmente alors du terme
${j_U}_!j_U^*(\FC)$. On obtient ainsi un complexe
$$ \CC^*(\XG,\o\FC)^+:\;\;{j_U}_!j_U^*(\FC) \injo \FC \To{\partial^0}
 \bigoplus_{S\in\XG^{(0)}}
i_{X_S,*}i_{X_S}^* (\FC) \To{\partial^1} \cdots \To{\partial^{d_X}} 
\bigoplus_{S\in\XG^{(d_X)}} i_{X_S,*}i_{X_S}^* (\FC) 
$$


\begin{lemme}
  Le complexe $\CC^*(\XG,\o\FC)^+$ ci-dessus est exact. 
\end{lemme}
\begin{proof}
On le montre 
sur les fibres en remarquant que si $x\in Z$, le sous ensemble simplicial $\XG_x:=\{X\in\XG, x\in X\}$ est contractile car le sous-ensemble ordonn\'e sous-jacent poss\`ede un \'el\'ement maximal, {\em cf} \cite{Quillen}.
\end{proof}


Choisissons maintenant pour tout $n\in\NM^\times$ une r\'esolution $\JC_n^*$ du faisceau constant $\Lambda/\lambda^n$ par des faisceaux injectifs dans $\FC(Z,\Lambda/\lambda^n)$.
Posons alors $\IC_n^*:=\bigoplus_{0\leq i\leq n} \JC_n^*$ et $\pi_n : \IC_n^*\To{} \IC_{n-1}^*$ la projection sur les $(n-1)$-premi\`eres composantes. On construit ainsi une r\'esolution injective du $\Lambda_\bullet$-faisceau $(\Lambda/\lambda^n)_n$ sur $Z$ telle que pour tout $X\in \XG$ et tout $p\in\NM$,
\begin{enumerate}
	\item  $i_{X,*}i_X^* (\IC^p_n)_n$ est un objet $\limproj \circ \Gamma(Z^{ca},-)$-acyclique.
	\item $\HC^p(\limproj\circ \Gamma(Z^{ca},i_{X,*}i_X^*\IC_n^*)) = H^p(X^{ca},\Lambda)$.
\end{enumerate}
En effet, on a des suites exactes
$$ R^1\limproj \, H^{q-1}(X^{ca},i_X^*(\IC^p_n)) \injo R^q(\limproj\circ \Gamma(Z^{ca},-))(i_{X,*}i_X^*(\IC^p_n)) \twoheadrightarrow \limproj \, H^q(X^{ca},i_X^*(\IC^p_n)), $$ qui entrainent le point i) puisque pour tous $p,n$ et $q>0$, on a  $H^q(X^{ca},i_X^*(\IC^p_n))=0$ et pour tout $p\in\NM$, le syst\`eme projectif $(\Gamma(X^{ca},i_X^*(\IC^p_n)))_n$ est \`a transitions surjectives, donc $\limproj$-acyclique.
Ce point i) implique la premi\`ere \'egalit\'e ci-dessous
\begin{eqnarray*} 
 \limproj\circ \Gamma(Z^{ca},i_{X,*}i_X^*\IC_n^*) & = & R(\limproj\circ \Gamma(Z^{ca},-))(i_{X,*}(\Lambda/\lambda^n)) \\
 & \simeq & R(\Gamma(Z^{ca},-)\circ \limproj)(i_{X,*}(\Lambda_\lambda^n)_n) \\
 &  \simeq & R(\Gamma(Z^{ca},-)\circ i_{X,*} \left(R\limproj((\Lambda/\lambda^n)_n) \right) \\
 & \simeq & R(\Gamma(X^{ca},-)\circ R\limproj ((\Lambda/\lambda^n)_n) 
 \end{eqnarray*}
d' o\`u le point ii). 

Ainsi, en appliquant le foncteur $\limproj\circ \Gamma(Z^{ca},-)$ au bicomplexe $\CC^{\#}(\XG,(\o\IC_n^*)_n)$, on obtient un quasi-isomorphisme de complexes
$$ R\Gamma_c(Z^{ca},(j_{U!}(\Lambda/\lambda^n))_n) \simto {\mathbf s}\left(\CC^{\#}\left(\XG, X\mapsto\limproj\circ \Gamma(X^{ca}, \IC^*_n)\right)\right). $$
En filtrant le bicomplexe par l'indice $*$, et par exactitude du foncteur "complexe de cocha\^ines", on obtient donc une suite spectrale
$$E_1^{pq}= \CC^q(\XG, X\mapsto H^p_c(X^{ca},\Lambda)) \Rightarrow H_c^{p+q}(Z^{ca},j_{U!}(\Lambda/\lambda^n)_n)$$
dont la diff\'erentielle est celle du complexe de cocha\^ines du syst\`eme de coefficients 
$X\mapsto H_c^p(X^{ca},\Lambda)$.

Il nous reste \`a montrer que l'application canonique
$$ R\Gamma_c(U(\XG)^{ca},(\Lambda/\lambda^n)_n) \To{} 
R\Gamma_c(Z^{ca},(j_{U!}(\Lambda/\lambda^n))_n)$$
est un quasi-isomorphisme. Ceci n'est pas vrai pour un ouvert quelconque (par exemple pour $U=\Omega$), mais l'est dans notre cas, car $U(\XG)$ est distingu\'e. En effet 
par \cite[4.1.9]{Fargues} et \cite[4.1.12]{Fargues},
 on a pour tout $q$ 
$$ R^q\Gamma_c(Z^{ca},j_{U,!}(\Lambda/\lambda^n))) \simto \limproj \, H^q_c(Z^{ca},j_{U!}(\Lambda/\lambda^n)) $$
et le terme de droite est aussi isomorphe \`a $H^q_c(U(\XG)^{ca},\Lambda)$, toujours car $U(\XG)$ est distingu\'e et par les m\^emes r\'ef\'erences.


\end{proof}

\alin{Morphismes d'arrangements} \label{morarr}
Supposons maintenant donn{\'e} un second arrangement $\XG'$. Par d\'efinition un {\em morphisme d'arrangements} $\XG\To{}\XG'$ est un couple $(f,\phi)$ o\`u 
\begin{itemize}
\item $f\in\aut{}{Z}$ est un automorphisme analytique de $Z$ et
\item  $\phi$ une
application strictement croissante $\XG \To{\phi}\XG'$ telle que pour tout $X\in
\XG$, on ait $\phi(X)\supseteq f(X)$.
\end{itemize}
Nous dirons aussi parfois que "$\phi$ est un morphisme d'arrangements au-dessus de $f$".
Notons que sur les ouverts compl\'ementaires, on a  $U(\XG') \subseteq f(U(\XG))$.

L'application $\phi$ induit aussi un morphisme d'ensembles
simpliciaux qui permet de construire des morphismes de complexes de cocha\^ines pour des syst\`emes de coefficients ad\'equates. Par exemple,pour tout faisceau ab\'elien $\FC$ on d{\'e}finit un morphismes de complexes
augment{\'e}s 
%
$\CC^*(f,\phi) : f^*\CC^*(\XG',\o\FC)^+ \To{}
\CC^*(\XG,\o{f^*\FC})^+$  par adjonction du morphisme suivant

$$\xymatrix{f_!{j_U}_!j_U^*(f^*\FC) \ar[r] & f_*f^*\FC \ar[r]  &
  \cdots \ar[r] & 
\bigoplus_{S\in\XG^{(k)}} f_*i_{X_S,*}i_{X_S}^*(f^*\FC) \ar[r] &
\cdots \\
{j_{U'}}_!j_{U'}^*(\FC) \ar[r] \ar[u] & \FC \ar[r] \ar@{=}[u] &
  \cdots \ar[r] & 
\bigoplus_{S\in{\XG'}^{(k)}} i_{X_S,*}i_{X_S}^*(\FC)
\ar[r]\ar[u]^{\CC^k(\phi)}  &
\cdots
}$$o{\`u} la fl{\`e}che de gauche est induite par l'inclusion $U'\injo f(U)$ et
 $$\CC^k(\phi) = \bigoplus_{S'\in {\XG'}^{(k)}}\sum_{S\in \phi^{-1}(S')}
i_{S'S} $$
avec $ i_{S'S}:  i_{X_S,*} i_{X_S}^*(\FC) \To{} f_*i_{X_{S'},*}i_{X_{S'}}^*(f^*\FC) $ le morphisme induit par l'inclusion $f(X_S)\injo X_{S'}$.


De m\^eme, pour tout $p$ on a un morphisme 
$$   \CC^*(f,\phi) :\;\;\CC^*(\XG',X'\mapsto H^p_c(X',\Lambda)) \To{}
\CC^*(\XG,X\mapsto H^p_c(X,\Lambda)) 
$$  
Nous laisserons le lecteur se convaincre de 
\begin{prop} Tout morphisme d'arrangements $\XG\To{(f,\phi)}\XG'$
induit un morphisme de suites spectrales $(E_r^{pq})(\XG')_{p,q,r} \To{} (E_r^{pq})(\XG)_{p,q,r}$ donn\'e par le morphisme $\CC^*(f,\phi)$ ci-dessus sur les termes $E_1^{p*}(\XG)$ et compatible  sur l'aboutissement au morphisme $H^{p+q}_c(U(\XG')^{ca},\Lambda) \To{f^{-1}_!} H^{p+q}_c(U(\XG)^{ca},\Lambda)$ induit par l'inclusion $U(\XG')\injo f(U(\XG))$.
\end{prop}

Les arrangements de $Z$ munis des morphismes que l'on vient de d\'efinir forment une cat\'egorie $\hbox{Arr}(Z)$. Le groupe des automorphismes analytiques $\aut{}{Z}$ agit strictement sur cette cat\'egorie.

\alin{Fonctions associ{\'e}es {\`a} un sommet}
Apr{\`e}s ces pr{\'e}liminaires simpliciaux, revenons au demi-plan. Notre but
est d'exhiber pour chaque facette $F$ de l'immeuble un arrangement
convenable $\XG_F$ dans $\PM(V)$ tel que $U_F=\PM(V)\setminus
\cup_{X\in \XG^{(0)}_F} X$.

Soit $s\in BT$ un sommet de l'immeuble. Choisissons un $\OC_K$-r{\'e}seau
g{\'e}n{\'e}rateur $\VC^*\subset V^*$ tel que la norme 
$$\application{|.|_{\VC^*} :\;}{ V^*}{\RM_+}{f}{\hbox{max}\{|x|_K, x\in
  K^* \hbox{ t.q. } x^{-1}f\in \VC^*\}}$$
repr{\'e}sente le sommet $s$. Posons alors $x_{\VC^*}:=\sup_{f\in
  \VC^*}|f(x)|$ pour tout $x\in \AM(V)$. Alors on v{\'e}rifie facilement
que la fonction 
$$\application{[.,.]_s:\;}{(V^*\setminus\{0\})\times
  (\AM(V)\setminus\{0\})}{\RM_+}{(f,x)}{\frac{|f(x)|}{x_{\VC^*}|f|_{\VC^*}}} $$
ne d{\'e}pend que de $s$ (et pas du choix du r{\'e}seau  $\VC^*$), et
se descend en une fonction  
$$[.,.]_s:\;P(V^*)\times \PM(V)\To{}\RM_+.$$

\begin{lemme} \label{l2}
  On a $\tau^{-1}(s^*) = \{x\in \PM(V),\;\forall f\in P(V^*), \;[f,x]_s
  > q^{-1}\}$.
\end{lemme}
\begin{proof}
(Indications.) Si $x\in\PM(V)$, on a $\tau(x)\in s^*$ \ssi\ il existe $k\in \NM$ et des sommets voisins de $s$ et deux-\`a-deux voisins et distincts $s_1,\cdots s_k$ tels que $\tau(x)$ soit dans l'enveloppe convexe stricte (la facette) des points $s,s_1,\cdots,s_k$.
En \'ecrivant les normes des $s_i$ en coordonn\'ees, on constatera sans difficult\'es que l'in\'egalit\'e de l'\'enonc\'e est v\'erifi\'ee.

R\'eciproquement, \'etant donn\'e $x$ satisfaisant l'in\'egalit\'e de l'\'enonc\'e, on r\'ecup\`ere les $s_i$ par l'\'egalit\'e
$$ \{s,s_1,\cdots, s_k\} =\{\hbox{sommets}\; t\;\hbox{ tels que} \;\; \forall f\in P(V^*),\; [f,x]_t>q^{-1}\}.$$

\end{proof}

\alin{Immeuble de Tits} Notons $\XG^{(0)}$ l'ensemble de tous les
  sous-$K$-espaces
vectoriels  de $V^*$ propres ({\em i.e.} non nuls et diff{\'e}rents de
$V^*$) et $\XG:=\XG^{(0)}\cup\{0\}$.
 Ces ensembles sont partiellement ordonn{\'e}s par l'inclusion et le
complexe simplicial associ{\'e} {\`a} $\XG^{(0)}$ est bien connu : il s'agit de l'immeuble
de Tits de $PGL(V^*)$. Il a {\'e}t{\'e} muni dans \cite{SS2} d'une topologie
``pro-sph{\'e}rique'' pour laquelle l'action naturelle de $PGL(V^*)$ est
continue. 
Pour $W\in \XG$ on d{\'e}finit
$$W^\perp :=\{x\in \PM(V),\;\forall f\in W,\;|f(x)|=0\}.$$
C'est un sous-espace $K$-analytique projectif de $\PM(V)$ et c'est
$\PM(V)$ tout entier si $W=0$. Pour un
sommet $s$, on d{\'e}finit aussi
$$\HC_s(W) := \{x\in \PM(V),\;\forall f\in W,\;[f,x]_s\leq q^{-1}\}$$
Si $W=0$ c'est encore $\PM(V)$ et si $W\neq 0$, c'est un domaine analytique compact (au sens de Berkovich) auquel on peut penser comme \`a un certain ``voisinage
tubulaire'' de $W^\perp$. 
Nous voulons aussi attacher \`a une facette $F$ des voisinages
tubulaires $\HC_F(W)$ de $W^\perp$. Pour cela, rappelons qu'une telle
facette correspond \`a une suite infinie 
$$ \left(\cdots \subset \VC_i^*\subset
\VC_{i+1}^* \subset  \cdots  \right)_{i\in\ZM}$$
 de $\OC$-r\'eseaux dans $V^*$ telle que
pour tout $i\in \ZM$, $\VC_i^*=\varpi_K\VC_{i+q+1}^*$
 o\`u $q$ est la dimension de $F$. Par ailleurs, un \'el\'ement $x\in
\PM(V)$ est repr\'esent\'e par une semi-norme multiplicative $\wt{x}$
sur $S(V^*)$. Notons $\kappa(x)$ le corps des fractions de l'anneau
quotient de $S(V^*)$ par l'id\'eal $\{f\in S(V^*),\;|f(x)|=0\}$. Ce
corps est muni d'une norme r\'esiduelle et on a
une factorisation :
$$\wt{x} :\; S(V^*)\To{} \kappa(x) \To{| |_{\kappa(x)}} \RM_+$$
dont nous noterons encore $\wt{x}$ la fl\`eche de gauche.
Nous posons alors 
\ini\begin{equation}\label{defHF}
 \HC_F(W):=\{x\in \PM(V),\;\; \forall i\in \ZM,\;
\wt{x}(W\cap\VC^*_{i+1})\subseteq  \wt{x}(\VC^*_i) \}, 
\end{equation}
le choix de $\wt{x}$ au-dessus de $x$ \'etant clairement indiff\'erent.
Lorsque $F=s=( \cdots \varpi^{-i}_K\VC^*\subset
\varpi^{-i-1}_K\VC^*\subset\cdots)$ est un sommet, la d\'efinition co\"incide
avec la 
pr\'ec\'edente en vertu des \'equivalences
\begin{eqnarray*}
  \forall i\in \ZM,\wt{x}(W\cap\VC^*_{i+1})\subseteq  \wt{x}(\VC^*_i)
& \equ & \wt{x}(W\cap\VC^*)\subseteq  \varpi_K\wt{x}(\VC^*) \\
& \equ & \forall f\in W \hbox{ t.q. } |f|_{\VC^*}= 1, |f(x)|\leq
q^{-1} x_{\VC^*}
\end{eqnarray*}

Par ailleurs, il est clair que si $F'\subset F$, alors pour tout $W$,
on a $\HC_{F'}(W) \subset \HC_F(W)$.

\begin{prop} \label{geom}
Notons $\XG_F$ l'ensemble de tous les ferm{\'e}s $\HC_F(W)$ obtenus
pour $W$ parcourant $\XG.$
  \begin{enumerate}
  \item Pour toute facette $F$, $\XG_F$ est un arrangement dans
    $\PM(V)$ dont l'ouvert compl{\'e}mentaire est 
     $U_F=\tau^{-1}(F^*)$.
\item
 L'application $W\mapsto
  \HC_F(W)$ est strictement croissante et induit pour tout $k$ une bijection
 :
 $$\XG^{(k)}/G_F^+ \simto \XG^{(k)}_F.$$
En particulier lorsque $F'\subset F$, on obtient une application croissante $\XG_F\To{} \XG_{F'}$ qui est un morphisme d'arrangement au-dessus de l'identit\'e ({\em cf }\ref{morarr}).

\item Pour tout $g\in PGL(V)$ et toute facette $F$, l'application $W\mapsto g^*(W)$ induit un isomorphisme d'arrangements $\XG_F\simto \XG_{gF}$ au dessus de $g$.

\item Pour toute facette $F$ et tout $W\in \XG_+$, l'application
  canonique $H^q(\HC_F(W)^{ca},\Lambda) \To{} H^q({W^\perp}^{ca},\Lambda)$ est un
  isomorphisme pour tout $q\geq0$.

\end{enumerate}
\end{prop}






\begin{proof}
Soit $F$ une facette repr\'esent\'ee par la suite de r\'eseaux $ \left(\cdots \subset \VC_i^*\subset
\VC_{i+1}^* \subset  \cdots  \right)_{i\in\ZM}$. 
De la  d\'efinition \ref{defHF}, on d\'eduit que
pour $W,W'\in \XG^{(0)}$, on a
$$ \HC_F(W)=\HC_F(W') \hbox{ \ssi\ } \forall i\in \ZM,\; W\cap \VC^*_i +
\VC^*_{i-1} = W'\cap \VC^*_i +  \VC^*_{i-1}. $$
Ceci a deux cons\'equences :
\begin{enumerate}
\item Si $F'\subset F$, on a $\HC_{F'}(W)=\HC_{F'}(W')
  \Rightarrow \HC_F(W)=\HC_F(W')$, ou, en d'autres termes,
  l'application $W\mapsto \HC_F(W)$ se factorise par l'application
  $W\mapsto \HC_{F'}(W)$.
\item  D'apr\`es l'\'egalit\'e 
 $ G_F^+ = \{g\in G,\;\; \forall i\in \ZM,\; g\VC^*_i\subseteq \VC^*_i
  \hbox{ et } g_{|\VC^*_i/\VC^*_{i-1}} = \id \} $, l'application
  $W\mapsto \HC_F(W)$ se factorise par  $\XG \To{}
  \XG/G_F^+ \To{} \XG_F$. En particulier, $\XG_F$ est fini.
\end{enumerate}
Pour $W,W'\in \XG^{(0)}$, on veut maintenant trouver $W''$ tel que 
$$ \HC_F(W)\cap \HC_F(W') = \HC_F(W'').$$
Pour cela, le sous-espace $W''\subset V^*$ doit v\'erifier 
$$\forall i\in\ZM,\;\;  W''\cap\VC^*_i +\VC^*_{i-1} = W\cap\VC^*_i + W'\cap\VC^*_i
+\VC^*_{i-1}. $$
Il suffit bien-s\^ur que l'\'egalit\'e ci-dessus soit v\'erifi\'ee
pour $i\in \{0,\cdots,q+1\}$ o\`u $q$ est la dimension de $F$. Plus
g\'en\'eralement, on a :
\begin{lemme}
  Dans le contexte ci-dessus, soit $(\o{W_i})_{i=0,\cdots, q+1}$ une
  famille de sous-$k$-espaces vectoriels respectifs des espaces
  $\o{V_i}:=\VC^*_i/\VC^*_{i-1}$. Alors il existe un sous-espace
  $W\subset V^*$ tel que pour tout $i$, l'image de $W\cap \VC^*_i$
  dans $\o{V_i}$ soit $\o{W_i}$.
\end{lemme}
La preuve de ce lemme  se
ram\`ene  \`a un simple exercice d'alg\`ebre lin\'eaire sur
le $k$-espace vectoriel $\varpi^{-1}\VC^*_0/\VC^*_0$ que nous laissons
au lecteur... Moyennant quoi, l'ensemble $\XG_F$ est bien un
arrangement au sens de \ref{arr}.

Dans le cas o\`u $F=s$ est un sommet, la premi\`ere d\'efinition que
nous avons donn\'ee de $\HC_s(W)$ et le lemme \ref{l2} montrent que
l'ouvert compl\'ementaire de $\XG_s$ est bien $U_s=\tau^{-1}(s^*)$.
Pour une facette g\'en\'erale, on a $\tau^{-1}(F^*)=\bigcap_{s\in F}
\tau^{-1}(s^*)$ et  les inclusions $\HC_{s}(W)\subset
\HC_F(W)$ pour chaque sommet $s\in F$ montrent que l'ouvert
compl\'ementaire de $\XG_F$ est inclus dans $\tau^{-1}(F^*)$. Pour
montrer l'autre inclusion, c'est un peu plus d\'elicat ; fixons  $W\in\XG^{(0)}$ et convenons de noter  $\o{W_i}$ l'image de $W\cap \VC^*_i$ dans
$\o{V_i}$. En appliquant le lemme pr\'ec\'edent, on trouve des
sous-espaces $W_i\subset V^*$ tels que :
$$ \forall j\in \ZM,\;\; \cas{\hbox{l'image de } W_i\cap
  \VC^*_{j}\hbox{ dans } \VC^*_j/\VC^*_{j-1}}{\o{W_j}}{j=i \hbox{ mod }
  q}{0}{j\neq i\hbox{ mod } q}$$ 
On v\'erifie alors sur les d\'efinitions que $\HC_F(W)=\bigcap_{i=0,\cdots, q+1}
\HC_F(W_i)$ et que d'autre part, si $s_i$ d\'esigne le sommet de $F$
repr\'esent\'e par le r\'eseau $\VC^*_i$, alors
$$ \HC_F(W_i) = \HC_{s_{i-1}}(W_i) .$$
Il s'ensuit que $\HC_F(W)$ est inclus dans le compl\'ementaire de
$\tau^{-1}(F^*)$, et on obtient ainsi la deuxi\`eme inclusion
cherch\'ee. Enfin, la discussion pr\'ec\'edente montre aussi que les $\HC_F(W)$ sont des domaines analytiques compacts de $\PM(V)$ (puisqu'on le sait pour les $\HC_{s_{i-1}}(W_i)$) et la preuve du i) est maintenant achev\'ee.

Passons \`a ii) : la croissance de $W\mapsto \HC_F(W)$ est
claire (rappelons que $\XG$ est ordonn\'e par inclusion et $\XG_F$ par contenance !). Pour voir la stricte croissance, il suffit de remarquer que
$\dim_K(W)= \sum_{i=0}^{q+1} \dim_{k}(\o{W_i})$ et de rappeler que
$\HC_F(W)=\HC_F(W')$ \ssi\ $\o{W_i}=\o{W'_i}$ pour tout $i\in \ZM$.
On a d\'eja remarqu\'e que la surjection $\XG^{(k)} \To{} \XG_F^{(k)}$
se factorise par $\XG^{(k)}\To{} \XG^{(k)}/G_F^+$. Le fait que
l'application $\XG^{(k)}/G_F^+
\To{} \XG_F^{(k)}$ est une bijection a \'et\'e  montr\'e par
Schneider-Stuhler \cite[Corollary 6.5.]{SS1}.

Le point iii) est une cons\'equence \'evidente des d\'efinitions.

Pour la preuve de iv), on s'inspire de l'id\'ee de \cite[Prop
1.6]{SS1} : traduite dans nos notations, elle consiste \`a exhiber une
r\'etraction $p : \HC_F(W)\To{} W^\perp$ de l'immersion
$i:W^\perp\injo \HC_F(W)$ dont les fibres sont des polydisques
affino\"ides. En effet, lorsqu'on a une telle r\'etraction et lorsque $\Lambda$ est de torsion (premi\`ere \`a $p$), on
consid\`ere les morphismes canoniques d'adjonction
$$ \Lambda_{W^\perp}\To{can} Rp_*p^*\Lambda_{W^\perp} \To{can}
Rp_*i_*i^*p^*\Lambda_{W^\perp} \simeq \Lambda_{W^\perp}. $$
L'isomorphisme de droite vient de l'\'egalit\'e $p\circ
i=\id_{W^\perp}$ et la compos\'ee de ces deux morphismes est
l'identit\'e du faisceau $\Lambda_{W^\perp}$. D'autre part, d'apr\`es
Berkovich \cite[7.4.2]{Bic2} et notre hypoth\`ese sur les fibres de $p$, la
premi\`ere fl\`eche est un isomorphisme. La deuxi\`eme est donc aussi
un isomorphisme et  en appliquant le foncteur
$R\Gamma({W^\perp}^{ca},.)$ \`a cette deuxi\`eme fl\`eche on obtient le
r\'esultat voulu, \`a savoir que l'application canonique de
restriction $R\Gamma(\HC_F(W)^{ca},\Lambda) \To{} R\Gamma({W^\perp}^{ca},\Lambda)$
est un isomorphisme. Ceci r\`egle le cas de torsion. Comme les espaces $\HC_F(W)$ et $W^\perp$ sont compacts et quasi-alg\'ebriques, on d\'eduit le cas $\lambda$-adique par passage \`a la limite gr\^ace aux isomorphismes $H^q(X^{ca},\Lambda)\simto \limproj \;H^q(X^{ca},\Lambda/\lambda^n)$ valable pour de tels espaces.

Pour construire l'application $p$, on choisit des suppl\'ementaires
$\o{Y_i}$ de $\o{W_i}$ dans $\o{V_i}$ pour $i=0,\cdots,q+1$ (avec les
notations d\'eja introduites plus haut) et on applique le lemme
pr\'ec\'edent pour trouver un sous-espace $Y\subset V^*$ tel que pour
tout $i=0,\cdots, q+1$, l'image de $Y\cap \VC^*_i$ dans $\o{V_i}$ soit
$\o{Y_i}$. On remarque que  $Y^\perp\cap W^\perp=\emptyset$ et
$\dim(Y^\perp)+\dim(W^\perp) = d-2$, de sorte que l'on peut 
d\'efinir  $p :\PM(V)\setminus Y^\perp \To{} W^\perp$ la projection
lin\'eaire sur $W^\perp$ de centre $Y^\perp$. On a
\begin{eqnarray*}
 Y^\perp \cap \HC_F(W) & = & \{x\in \PM(V),\; \wt{x}(Y)=0 \hbox{ et
  }\forall i\in \ZM,\;\;
\wt{x}(W\cap \VC^*_i)\subseteq \wt{x}(\VC^*_{i-1}) \} \\
& \subseteq & 
\{x\in \PM(V),\; \forall i\in \ZM,\;\;
\wt{x}(\VC^*_i)\subseteq \wt{x}(\VC^*_{i-1}) \} \\
& = & \emptyset,
\end{eqnarray*}
la derni\`ere \'egalit\'e venant de la non-nullit\'e de $\wt{x}$
associ\'e \`a $x\in \PM(V)$. Donc $p$ est bien d\'efinie sur
$\HC_F(W)$ et il est clair que c'est une r\'etraction de l'inclusion
$i : W^\perp \injo \HC_F(W)$.

 On veut maintenant \'etudier les fibres
de $p$ et de $p_{|\HC_F(W)}$. Pour cela on fixe un \'el\'ement $x$ de
$W^\perp \subset \PM(V)$. On note $\wh\kappa(x)$ la compl\'etion du
corps r\'esiduel $\kappa(x)$ en $x$ ; la fibre de $p$ au-dessus de $x$
est naturellement munie d'une structure de $\wh\kappa(x)$-espace
analytique {\em affine}. Pour expliciter cette structure, notons
$V_x:=V\otimes_K \wh\kappa(x)$  et $\PM(V_x)$ le
$\wh\kappa(x)$-espace projectif associ\'e. Notons aussi $V^*_x =
Y_x\oplus W_x$ la d\'ecomposition de $V^*_x:=V^*\otimes_K\wh\kappa(x)$
obtenue de la d\'ecomposition $V^*=Y\oplus W$ par extension des
scalaires de $K$ \`a $\wh\kappa(x)$. On notera aussi $Y_x^\perp=\{z\in
\PM(V_x),\; \forall f\in Y_x,\;|f(z)|=0\}$ le sous-espace projectif
``othogonal'' associ\'e \`a $Y_x$ ; on a \'evidemment $Y_x^\perp
\simeq Y^\perp \wh\otimes_K \wh\kappa(x)$.
Posons alors
$$Y(x) =\{f\in Y_x,  |f(x)|=0\} \subseteq Y_x \subset V^*_x$$
(On voit ici $x$ comme un point de l'espace $\PM(V_x)$). La fibre de
$p$ en $x$ munie de sa structure de $\wh\kappa(x)$-espace analytique
est 
$$ p^{-1}(x) \simeq Y(x)^\perp \setminus Y_x^\perp \subset \PM(V_x). $$
C'est un espace affine dont  l'alg\`ebre de
fonctions analytiques est la localisation de l'alg\`ebre gradu\'ee
$S(V_x^*/Y(x))$ par un g\'en\'erateur, disons $f_x$, de
$Y_x/Y(x)$. Ainsi l'application 
$$
\application{\alpha:\;}{S(W_x)}{S(V_x^*/Y(x))_{f_x}}{f}{ff_x^{-1}}$$
induit un isomorphisme d'espaces $\wh\kappa(x)$-affines 
$$ \alpha^* : \; p^{-1}(x) \simto \AM(W_x). $$
On veut maintenant calculer la fibre de $p_{|\HC_F(W)}$ en $x$,
c'est-\`a-dire l'image par l'isomorphisme $\alpha^*$ de l'intersection
$ p^{-1}(x)\cap \HC_F(W)$. 

La suite strictement croissante de
$\OC_K$-r\'eseaux $(\VC^*_i)_{i\in \ZM}$ dans le $K$-espace $V^*$
repr\'esentant la facette  $F$ induit par 
extension des scalaires puis passage au quotient une suite {\em
  croissante} de r\'eseaux $(\XC^*_i)_{i\in \ZM}$ de
$\OC_{\wh\kappa(x)}$-r\'eseaux dans le $\wh\kappa(x)$-espace
$V^*_x/Y(x)$. 
Identifions $W_x$ avec son image dans $V_x^*/Y(x)$ ; notre choix du
$K$-espace $Y$ nous assure que 
$$ \forall i\in \ZM,\;\; \XC^*_i = (\XC^*_i\cap W_x) \oplus
(\XC^*_i\cap (Y_x/Y(x))). $$
Remarquons qu'il existe un unique $i\in \ZM$ tel que $f_x\in
\XC^*_i\cap (Y_x/Y(x)) \setminus \XC^*_{i-1}\cap (Y_x/Y(x))$. Nous
supposerons, quitte \`a effectuer un d\'ecalage 
dans la num\'erotation des r\'eseaux initiaux  $\VC^*_i$,  que $i=-q-1$.
Ceci permet d'\'ecrire pour $i=1,\cdots,q+1$ :
$$\XC_{-i}^* = (W_x\cap \XC^*_{-i}) \oplus \OC_{\wh\kappa(x)}.f_x$$
et pour $i=0$,
$$\XC_{0}^* = (W_x\cap \XC^*_{0}) \oplus
\OC_{\wh\kappa(x)}.\varpi_K^{-1}f_x.$$ 
Maintenant,
 nous
pouvons choisir une $\OC_{\wh\kappa(x)}$-d\'ecomposition $$
W_x\cap\XC^*_0 = \LC^0 \oplus 
\cdots \oplus \LC^q  $$ 
avec \'eventuellement $\LC^i=0$ pour certains $i$ mais telle que 
\begin{eqnarray*}
 W_x\cap \XC^*_{-1} & = & \varpi_K \LC^0 \oplus \LC^1 \cdots \oplus \LC^{q}
 \\
& \vdots & \\
W_x\cap \XC^*_{-q}&  = & \varpi_K \LC^0 \oplus  \cdots \oplus \varpi_K
\LC^{q-1} \oplus \LC^q 
\end{eqnarray*}

Pour un sous-$\OC_{\wh\kappa(x)}$-module $\XC$  de type fini de
$V_x^*/Y(x)$ et un \'el\'ement $y\in p^{-1}(x)$ vu comme une norme sur
l'alg\`ebre $S(V_x^*/Y(x))$, notons
$$ y(\XC):= \sup_{f\in \XC} |f(y)|. $$
Alors la condition pour que $y\in \HC_F(W)$ s'\'ecrit
$$ \left(\forall i\in\{0,\cdots,
  q\},\;\; y(W_x\cap \XC_{-i}) \leq y(\XC_{-i-1})\right)$$  
ou encore
$$
  \left\{ 
  \begin{array}{lcl}
\sup(y(\LC^0),\cdots, y(\LC^q)) & \leq & \sup(q^{-1}_Ky(\LC^0),y(\LC^1)
\cdots, y(\LC^q), |f_x(y)|) \\
\sup(q^{-1}_Ky(\LC^0),\cdots, y(\LC^q))& \leq &
\sup(q^{-1}_Ky(\LC^0),q^{-1}_Ky(\LC^1), 
 y(\LC^2),\cdots , |f_x(y)|)\\ 
& \vdots & \\
\sup(q^{-1}_Ky(\LC^0),\cdots, q^{-1}_K y(\LC^{q-1}), y(\LC^q)) & \leq &
\sup(q^{-1}_Ky(\LC^0), 
\cdots, q^{-1}_Ky(\LC^q), |f_x(y)|) \\ 
  \end{array}\right.
$$
Nous allons montrer que ce syst\`eme d'in\'egalit\'es est \'equivalent
\`a la seule in\'egalit\'e
\ini
\begin{equation}
  \label{ineg2}
  \sup(y(\LC^0),\cdots, y(\LC^q)) \leq |f_x(y)|.
\end{equation}
Tout d'abord, il est clair que si \ref{ineg2} est v\'erifi\'ee alors le
syst\`eme d'in\'egalit\'es ci-dessus l'est aussi. R\'eciproquement,
supposons le syst\`eme d'in\'egalit\'es v\'erifi\'e et notons
$k\in\{0,\cdots,q\}$
le plus grand indice tel que $y(\LC^k)=\sup_{i=0,\cdots,q}(y(\LC^i))$. Alors la
$k+1$-i\`eme ligne du syst\`eme nous dit que 
$y(\LC^k)\leq \sup(q_K^{-1}y(\LC^0),\cdots,
q_K^{-1}y(\LC^k),y(\LC^{k+1}),\cdots, y(\LC^q),|f_x(y)|)$. Mais le
choix de $k$ entraine alors que $y(\LC^k)\leq |f_x(y)|$.

On en d\'eduit que l'isomorphisme $\alpha^*$ envoie $p^{-1}(x)\cap
\HC_F(W)$ sur le polydisque $D_{\XC_0}$ de $\AM(W_x)$ d\'efini par
$$ D_{\XC_0}:=\{z\in \AM(W_x),\;\; \sup_{f\in W_x\cap\XC_0^*} |f(z)| \leq
1\}. $$

\end{proof}

\alin{Le syst\`eme de suites spectrales $\EM_r^{pq}$} \label{systss}
D'apr\`es la proposition pr\'ec\'edente, on a construit un syst\`eme de coefficients $F\mapsto \XG_F$,  $G$-\'equivariant sur $BT$ en arrangements ({\em i.e} \`a valeurs dans la cat\'egorie $\hbox{Arr}(\PM(V))$ d\'efinie en \ref{morarr}). Par fonctorialit\'e, on en d\'eduit plusieurs autres syst\`emes de coefficients $G$-\'equivariants comme 
\begin{itemize}
	\item un syst\`eme en suites spectrales $\EM_r^{pq}:\;\;F\mapsto E^{pq}_r(\XG_F)$ dont le $\EM_1$ est donn\'e par 
	\item des syst\`emes en complexes de $\Lambda W_K$-modules $\EM_1^{p*}:\;\;F\mapsto \CC^*(\XG_F,X\mapsto H^p(X^{ca},\Lambda))$, et le syst\`eme des aboutissements est un gradu\'e 
	\item des syst\`emes $F\mapsto H^{p+q}_c(U_F^{ca},\Lambda)$.
\end{itemize}
Pr\'ecisons seulement que pour construire le syst\`eme en suites spectrales, il faut choisir les r\'esolutions $(\JC^*_n)_n$ de la preuve de \ref{sspec2} de mani\`ere $G$-\'equivariante, {\em i.e.} dans $\FC(\PM(V),\Lambda G_{disc})$.


Introduisons maintenant le sous-ensemble 
$$\XG(p) := \{W\in \XG,\;
H^p(W^\perp,\Lambda)\neq 0 \}.$$
 C'est un sous-ensemble partiellement
ordonn\'e plein de $\XG$ qui est non-vide seulement si $p$ est pair et $\leq 2d-2$. Dans ce cas, l'ensemble simplicial associ\'e est
profini de dimension $d-1-p/2$, et a \'et\'e \'etudi\'e par
Schneider et Stuhler 
\cite[prop 3.6.,proof]{SS1}. Le dictionnaire avec leurs notations est le
suivant :
$ \XG(p)^{(k)}$ s'\'ecrirait avec leurs notations
$\NC\TC_k^{(d-1-p/2)}$. Nous renvoyons \`a {\em loc. cit}, notamment avant
leur lemme 3.3.,  pour la
d\'efinition pr\'ecise de la topologie sur $\XG(p)^{(k)}$ et nous noterons
$C^\infty(\XG(p)^{(k)},\Lambda)$ le $\Lambda$-module des fonctions
localement constantes $\XG(p)^{(k)}\To{} \Lambda$.
Suivant Schneider et Stuhler, on associe \`a $\XG(p)$ le complexe de cocha\^ines
continues $(\CC^*(\XG(p),\Lambda),\partial^*)$ suivant :
$$ \Lambda  \To{\partial^0}
C^\infty(\XG(p)^{(0)},\Lambda) 
 \To{\partial^1} \cdots \To{\partial^{k}} 
C^\infty(\XG(p)^{(k)},\Lambda)  \To{\partial^{k+1}} \cdots
\To{\partial^{d-1-p/2}} C^\infty(\XG(p)^{(d-1-p/2)},\Lambda)
$$
o\`u la diff\'erentielle est d\'efinie comme la somme altern\'ee des
applications de d\'eg\'en\'er\'escence (version continue de
\ref{poset}). Ce complexe est naturellement un complexe de $\Lambda
G$-modules lisses admissibles.

\begin{prop}
Il existe pour tout $p\in\NM$ pair et $\leq 2d-2$ un isomorphisme de syst\`emes de coefficients $G$-\'equivariants en complexes de $\Lambda W_K$-modules :
$$ \EM_1^{p*} \simto \left( F\mapsto \CC^*(\XG(p),\Lambda)^{G_F^+}(-p/2) \right).$$
Pour $p$ impair, $\EM_1^{p*}=0$.
\end{prop}

\begin{proof}
Commen{\c c}ons par souligner les cons\'equences suivantes du point iv) de la proposition \ref{geom}:
\begin{enumerate}
\item La cohomologie de  $\HC_F(W)$ est donn{\'e}e par
  $$
\begin{array}{rcl}
H^p(\HC_F(W)^{ca},\Lambda) & \simeq  & \left\{ 
  \begin{array}{ll}
   \Lambda(-p/2) & \hbox{ si } p\leq 2c(W) \hbox{ est pair.} \\
   0 & \hbox{ sinon}
  \end{array}\right.
\end{array}$$
o\`u $c(W)$ d\'esigne la codimension du sous-espace $W$ de $V^*$. 
\item  Pour toute facette $F$ et tous $W,W'$ tels que $\HC_F(W')\subset \HC_F(W)$,
  l'application  $$H^p(\HC_F(W)^{ca},\Lambda)\To{}
  H^p(\HC_F(W')^{ca},\Lambda)$$ est un isomorphisme pour tout $p\leq 2c(W')$.
 
\item Pour tout $W\in \XG^{(0)}$ et toutes facettes $F'\subset F$,
  l'application canonique $$H^p(\HC_F(W)^{ca},\Lambda)\To{}
  H^p(\HC_{F'}(W)^{ca},\Lambda)$$ est un isomorphisme pour tout $p$.
\end{enumerate}


Fixons maintenant $F$ et introduisons le sous-ensemble $$\XG_F(p) = \{X\in \XG_F,\;
H^p(X^{ca},\Lambda)\neq 0 \}.$$ C'est un sous-ensemble partiellement
ordonn\'e de $\XG_F$ tel que $X\leq X'\in \XG_F \Rightarrow X\in
\XG_F$ (par ii) ci-dessus) et qui est non vide seulement si $p=2i$ est pair et $\leq 2d-2$ (par i)). 
Toujours d'apr\`es i) et ii) ci-dessus, le syst\`eme de coefficients $X\mapsto H^p(X^{ca},\Lambda)$ sur $\XG_F$ est support\'e par $\XG_F(p)$ et y est constant.
On a donc une identit\'e de complexes de cocha\^ines
$$\CC^*(\XG_F(p),\Lambda)(-p/2) \simto \CC^*(\XG_F,X\mapsto H^p(X^{ca},\Lambda))
$$
o\`u $\Lambda$ dans le terme de gauche d\'esigne le syst\`eme de
coefficients sur $\XG_F(p)$ constant $X\mapsto \Lambda$.

De plus, par iii) l'application $\XG_{F'}\To{} \XG_F$ pour $F'\subset F$ envoie
$\XG_{F'}(p)$ dans $\XG_F(p)$ et pour tout $g\in G$,
$\XG_{g(F)}(p)=g(\XG_F(p))$.
En r\'esum\'e, le terme $\EM_1$ est donn\'e en tant que syst\`eme de coefficients $G$-\'equivariant en $\Lambda W_K$-modules par  
$$ \EM_1^{p*} = \left( F\mapsto \CC^*(\XG_F(p),\Lambda)(-p/2) \right). $$

Il r\'esulte maintenant de \ref{geom} ii) que pour tout $k\geq 0$ et toute
facette $F$, l'application canonique $\XG \To{} \XG_F$ induit une
bijection  $\XG(p)^{(k)}/G_F^+ \simto \XG_F(p)^{(k)}$. De cette
bijection on tire un isomorphisme de complexes de $\Lambda $-modules :
$$ \CC^*(\XG_F(p),\Lambda) \simto \CC^*(\XG(p),\Lambda)^{G_F^+}.$$
Par construction les isomorphismes ainsi obtenus sont compatibles 
\`a l'action de $G$ et aux inclusions de facettes : ce sont des  
isomorphismes de syst\`emes de coefficients $G$-\'equivariants sur
$BT$ \`a valeurs dans les complexes de $\Lambda$-modules. On a donc termin\'e la preuve.

\end{proof}

\begin{coro}
Le terme $\EM_2$ du syst\`eme de suites spectrales \ref{systss} est donn\'e par 
  $$
\begin{array}{rcl}
\EM_2^{pq} & \simeq  & \left\{ 
  \begin{array}{ll}
   \gamma_0(\pi_{\{1,\cdots,p/2\}}^\Lambda)(-p/2) & \hbox{ si } p\leq 2d-2 \hbox{ est pair et }\, q=d-1-p/2 \\
   0 & \hbox{ sinon}
  \end{array}\right.
\end{array}$$

\end{coro}
\begin{proof}
D'apr\`es Schneider et Stuhler  \cite[lemma 4.1]{SS1}, la cohomologie du complexe $\CC^*(\XG(p),\Lambda)$ est concentr\'ee en degr\'e $d-1-p/2$ et est isomorphe \`a $\pi_{\{1,\cdots,p/2\}}^\Lambda(-p/2)$.

\end{proof}

\alin{Fin de la preuve de la proposition \ref{cle}}
D'apr\`es le corollaire pr\'ec\'edent, le syst\`eme de suites spectrales (ou si l'on pr\'ef\`ere, la suite spectrale en syst\`emes de coefficients) $\EM_r^{pq}$ d\'eg\'en\`ere en $\EM_2$ et fournit des isomorphismes 
$$ \EM_2^{2i,d-1-i} \simto \left(F \mapsto H^{d-1+i}(U_F^{ca},\Lambda)\right).$$

\section{Appendice 2}

\subsection{Crit{\`e}res de scindage d'un complexe} \label{sectionscin}

Le but de cette section est de donner des crit{\`e}res de scindage d'un
complexe cohomologiquement born{\'e} d'une cat{\'e}gorie
d{\'e}riv{\'e}e. 
Ces crit{\`e}res  
sont tr{\`e}s simples et certainement bien connus des sp{\'e}cialistes
 -- comme l'a signal\'e G. Laumon \`a l'auteur, on trouve
d\'eja dans Deligne \cite{Deldegen} un \'enonc\'e  semblable \`a
\ref{actionphi} ci-dessous, mais un
un peu trop diff\'erent pour \^etre cit\'e tel quel.
 Le cadre naturel, un peu plus g{\'e}n{\'e}ral,
est celui des cat{\'e}gories triangul{\'e}es.

\ali
Soit $\DC$ une cat{\'e}gorie triangul{\'e}e munie d'une $t$-structure
$(\DC^{\leq 0},\DC^{\geq 0})$ non-d{\'e}g{\'e}n{\'e}r{\'e}e, {\em cf} \cite[1.3.7.]{BBD}. 
 On utilise les notations habituelles : $\tau_{\leq q}$ et $\tau_{\geq
   q}$  pour les endofoncteurs de troncation  et 
$\HC^q:=\tau_{\leq q}\tau_{\geq q}[q]=\tau_{\geq q}\tau_{\leq
  q}[q]=\tau_{\geq 0}\tau_{\leq 0}\circ [q]$  les
(endo)foncteurs cohomologiques. On note aussi $\DC^+$, resp. $\DC^-$
et $\DC^b$ les sous-cat{\'e}gories triangul{\'e}es form{\'e}es des objets $X$ tels
que $\HC^q(X)=0$ pour $q << 0$, resp. $q>>0$, resp. $q<<0$ ou $q>>0$
(dans ce dernier cas $X$ est dit ``cohomologiquement born{\'e}'').


Nous dirons qu'un objet $X$ est {\em scindable}, s'il est
cohomologiquement born{\'e} et s'il existe un isomorphisme dans $\DC$
$$ \alpha:\;\; X \simto  \bigoplus_q \HC^q(X)[-q]. $$
Un tel isomorphisme sera appel{\'e} {\em un scindage} s'il induit
l'identit{\'e} en cohomologie, c'est-{\`a}-dire si pour tout $q\in\ZM$, on a
$\HC^q(\alpha)=\id_{\HC^q(X)}$.


\begin{lemme} Soit $X\in D^-$. Fixons $q\in \ZM$ et 
  supposons que $\hom{\HC^l(X)[k]}{\HC^q(X)}{\DC}=0$ 
  pour tout couple d'entiers $(l,k)$ tel que $l+k=q-1$. Alors 
l'application
$$\HC^q:\;\; \hom{X}{\HC^q(X)[-q]}{\DC} \To{}
\hom{\HC^q(X)}{\HC^q(X)}{\DC} $$
est surjective.
\end{lemme}
\begin{proof}
Quitte {\`a} d{\'e}caler, on peut supposer que $q=0$, ce que nous ferons.
  L'application en question s'inscrit dans un diagramme 
$$\xymatrix{ \hom{X}{\HC^0(X)}{\DC} \ar[r]^{\HC^0}
  \ar[d]_{\tau_{\leq 0}} & \hom{\HC^0(X)}{\HC^0(X)}{\DC} \\
\hom{\tau_{\leq 0}(X)}{\HC^0(X)}{\DC}
  \ar[ru]_{\tau_{\geq 0}} & }$$
o{\`u} la fl{\`e}che not{\'e}e $\tau_{\geq 0}$ est bijective et la fl{\`e}che not{\'e}e
$\tau_{\leq 0}$ est la composition avec le morphisme canonique
$\tau_{\leq 0}(X)\To{} X$. Ainsi, l'application de l'{\'e}nonc{\'e} est
surjective \ssi\ pour tout morphisme $\gamma : \tau_{\leq 0}(X) \To{}
\HC^0(X)$, on peut compl{\'e}ter le diagramme suivant
$$ \xymatrix{  \tau_{>0}(X)[-1] \ar[r]^{-1} & \tau_{\leq 0}(X) \ar[r]^{can}
  \ar[d]^\gamma &
  X \ar@{..>}[ld]^\phi \ar[r]^{can} & \\ & \HC^0(X)  &  & }.$$
Une condition n{\'e}cessaire est bien-s{\^u}r que $\gamma\circ -1=0$, et
l'axiome de l'octa{\`e}dre montre que cette condition est suffisante. 
Nous allons en fait montrer par r{\'e}currence sur l'entier $n\in\ZM$ tel que
$X\in \DC^{\leq n}$ que l'hypoth{\`e}se de l'{\'e}nonc{\'e} (c'est-{\`a}-dire :
$\hom{\HC^l(X)[k]}{\HC^0(X)}{\DC}=0$  
  pour tout couple d'entiers $(l,k)$ tel que $l+k=-1$) implique 
$$ \hom{ \tau_{>0}(X)[-1]}{\HC^0(X)}{\DC} =0,$$
ce qui sera suffisant pour prouver le lemme.

Remarquons que si $n\leq 0$, on a $\tau_{>0}(X)=0$ et la propri{\'e}t{\'e}
cherch{\'e}e est imm{\'e}diate. Supposons donc $n>0$ et la propri{\'e}t{\'e} montr{\'e}e
pour $n-1$. Soit $X\in \DC^{\leq n}$. L'objet $\tau_{<n}(X)$ est dans
$\DC^{\leq n-1}$ et v{\'e}rifie l'hypoth{\`e}se du lemme, on peut donc lui
appliquer  l'hypoth{\`e}se de r{\'e}currence. Le triangle distingu{\'e}
$ \tau_{>0}\tau_{<n}(X) \To{} \tau_{>0}(X) \To{} \HC^n(X)[-n]$ fournit
une suite exacte
\begin{eqnarray*}
\hom{\HC^n(X)[-1-n]}{\HC^0(X)}{\DC} & \To{} & \hom{\tau_{>0}(X)[-1]}{\HC^0(X)}{\DC}
 \\ & \To{} & \hom{\tau_{>0}\tau_{<n}(X)[-1]}{\HC^0(X)}{\DC}.
\end{eqnarray*}
Le premier terme est nul par hypoth{\`e}se, ainsi que le dernier  par
hypoth{\`e}se de r{\'e}currence (remarquer que la fl{\`e}che canonique
$\HC^0(\tau_{<n}(X)) \To{} \HC^0(X)$ est un isomorphisme), le terme du
milieu est donc nul aussi. 
\end{proof}

\begin{coro} \label{scindage} Soit $X\in\DC$ cohomologiquement born{\'e}.
  Supposons que $$\hom{\HC^l(X)[k]}{\HC^q(X)}{\DC}=0$$
  pour tout triplet d'entiers $(q,l,k)$ 
  tel que $l+k=q-1$. Alors $X$ est scindable.
\end{coro}
\begin{proof}
D'apr{\`e}s le lemme pr{\'e}c{\'e}dent, pour tout $q\in\ZM$ on peut trouver
$\alpha_q:\; X \To{} \HC^q(X)[-q] $ tel que
$\HC^q(\alpha_q)=\id_{\HC^q(X)}$. Le morphisme somme
$$ \bigoplus_{q} \alpha_q :\;\; X \To{} \bigoplus_{q} \HC^q(X)[-q] $$
induit l'identit{\'e} en cohomologie et par cons{\'e}quent est un
isomorphisme, car la $t$-structure est non-d{\'e}g{\'e}n{\'e}r{\'e}e.
\end{proof}

Dans le lemme suivant on suppose de plus que $\DC$ est $R$-lin{\'e}aire
pour un anneau commutatif $R$ fix{\'e}.

\begin{lemme} \label{actionphi}
  Soit $X$ un objet cohomologiquement born{\'e} de $\DC$ et $\phi \in \endo{\DC}{X}$. On suppose
  donn{\'e}e une famille de polyn{\^o}mes $P_q(T)\in R[T]$,
  $q\in\ZM$, presque tous {\'e}gaux {\`a} $1$, et tels que  pour tout
  $q\in\ZM$, on ait $P_q(\HC^q(\phi))=0$ 
  dans $\endo{\DC}{\HC^q(X)}$. 
  \begin{enumerate}
  \item Posons $P(T):=\prod_{q\in\ZM} P_q(T)$, alors on a $P(\phi)=0$
    dans $\endo{\DC}{X}$.
\item Supposons de plus que pour tous $p\neq q$, on a
  $P_qR[X]+P_pR[X]=R[X]$. Alors il existe un unique scindage  $$\alpha :
   X \To{} \bigoplus_{q} \HC^q(X)[-q]$$ tel que 
 $\alpha \phi \alpha^{-1} = \bigoplus_q \HC^q(\phi)[-q]$.
  \end{enumerate}
\end{lemme}
\begin{proof}
  Pour le point i) on proc{\`e}de par r{\'e}currence sur l'amplitude
  cohomologique de $X$ en utilisant le r{\'e}sultat
  auxiliaire suivant :

{\em Soient $a\in\ZM$ et $P_-(T), P_+(T) \in R[T]$ deux polyn{\^o}mes tels
  que $P_-(\tau_{\leq a}(\phi))=0$ dans $\endo{\DC}{\tau_{\leq a}(X)}$
  et $P_+(\tau_{> a}(\phi))=0$ dans $\endo{\DC}{\tau_{>
      a}(X)}$. Alors $P_-P_+(\phi)=0$ dans $\endo{\DC}{X}$.}

Nous laisserons la r{\'e}currence au lecteur, mais nous allons montrer
l'assertion ci-dessus. On d{\'e}finit un morphisme $\psi_-: \tau_{>a}(X)
\To{} X$ par le diagramme commutatif suivant :
$$\xymatrix{ \tau_{\leq a}(X) \ar[r]^{\iota_a} \ar[d]_{\tau_{\leq a}(P_-(\phi))=0} & X
  \ar[r]^{\rho_a} \ar[d]_{P_-(\phi)}  & \tau_{>a}(X) \ar@{..>}[ld]_{\psi_-}
  \ar[r]^{+1} & \\
\tau_{\leq a}(X) \ar[r] & X \ar[r] & \tau_{>a}(X)
  \ar[r]^{+1} & }$$
L'existence de $\psi_-$ est assur{\'e}e par l'axiome de l'ocat{\`e}dre
appliqu{\'e} {\`a} la composition $P_-(\phi)\circ \iota_a =0$.
De m{\^e}me on d{\'e}finit $\psi_+ : X \To{} \tau_{\leq a}(X)$ par le
diagramme
$$\xymatrix{ \tau_{\leq a}(X) \ar[r]^{\iota_a}  & X
  \ar[r]^{\rho_a} \ar[d]^{P_+(\phi)}  \ar@{..>}[ld]_{\psi_+}  & \tau_{>a}(X)
  \ar[r]^{+1} \ar[d]^{\tau_{> a}(P_+(\phi))=0} & \\
\tau_{\leq a}(X) \ar[r] & X \ar[r] & \tau_{>a}(X)
  \ar[r]^{+1} & }.$$
On calcule alors $(P_-P_+)(\phi)=P_-(\phi)\circ P_+(\phi) =
\psi_-\circ \rho_a\circ \iota_a \circ \psi_+ =0$ car $\rho_a\circ
\iota_a=0$.

Pour l'existence dans le point ii), on proc{\`e}de encore par r{\'e}currence
sur l'amplitude cohomologique de $X$ en utilisant le
r{\'e}sultat auxiliaire suivant :

{\em  En gardant les notations ci-dessus, supposons que $P_-(X)R[X] +
  P_+(X)R[X]=R[X]$. Alors il existe un isomorphisme
$$ \alpha_a:\;\;X \simto \tau_{\leq a}(X) \oplus \tau_{>a}(X) $$
tel que $\tau_{\leq a}(\alpha_a)=\id_{\tau_{\leq a}(X)}$,
$\tau_{>a}(\alpha_a)=\id_{\tau_{>a}(X)}$ et $\alpha_a\phi\alpha_a^{-1}
= \tau_{\leq a}(\phi) \oplus \tau_{>a}(\phi) $.}

Nous laissons {\`a} nouveau la r{\'e}currence au lecteur, mais nous prouvons
cette assertion. Pour cela, remarquons que quitte {\`a} remplacer $P_-$,
resp. $P_+$, par un de ses multiples dans $R[X]$, on peut supposer
$P_-(X)+P_+(X)=1$. Nous allons montrer que les morphismes
$$ \alpha_a :\;\; X \To{\psi_+\oplus \rho_a} \tau_{\leq a}(X)\oplus
\tau_{>a}(X) \;\;\hbox{ et } \;\; \beta_a:\;\; 
\tau_{\leq a}(X)\oplus
\tau_{>a}(X) \To{\iota_a\oplus\psi_-} X $$
sont inverses l'un de l'autre.
On a d'abord $\beta_a\circ \alpha_a= \psi_-\rho_a+\iota_a\psi_+ =
P_-(\phi) +P_+(\phi) = \id_X$.
Par ailleurs $\alpha_a\circ \beta_a= \psi_+\iota_a \oplus
\rho_a\psi_-$. Montrons que $\rho_a\psi_- = \id_{\tau_{>a}(X)}$, le
raisonnement sera identique et omis pour $\psi_+\iota_a$. 
Tout d'abord, puisque $\tau_{>a}(P_+(\phi))=0$, on a
$\tau_{>a}(P_-(\phi)) = \id_{\tau_{>a}(X)}$. 
Il nous suffira donc de montrer que le triangle du bas du diagramme suivant est commutatif
$$\xymatrix{ X
  \ar[r]^{\rho_a} \ar[d]_{P_-(\phi)}  & \tau_{>a}(X)
  \ar@{..>}[ld]_{\psi_-} \ar[d]^{\tau_{>a}(P_-(\phi))=\id} \\
 X \ar[r] & \tau_{>a}(X) }$$
sachant que celui du haut l'est, par d{\'e}finition. La commutativit{\'e} de
celui du haut et du carr{\'e} ext{\'e}rieur donnent l'{\'e}galit{\'e}
$\rho_a\psi_-\rho_a= \tau_{>a}(P_-(\phi)) \circ \rho_a$ dans
$\hom{X}{\tau_{>a}(X)}{\DC}$.
 En appliquant le foncteur
$\tau_{>a}$ et en tenant compte de ce que $\tau_{>a}(\rho_a) =
\id_{\tau_{>a}(X)}$ (en identifiant $\tau_{>a}(\tau_{>a}(X))$ et
$\tau_{>a}(X)$) et de ce que $\tau_{>a}(\rho_a\psi_-) = \rho_a\psi_-$,
on obtient l'{\'e}galit{\'e} cherch{\'e}e $\rho_a\psi=\tau_{>a}(P_-(\phi))$.

L'isomorphisme $\alpha_a$ ainsi construit v{\'e}rifie bien 
 $\tau_{\leq a}(\alpha_a)=\id_{\tau_{\leq a}(X)}$ et
$\tau_{>a}(\alpha_a)=\id_{\tau_{>a}(X)}$. On calcule aussi $\alpha_a
\phi \beta_a = \rho_a\phi\psi_- \oplus \psi_+\phi\iota_a$. Or, 
$\rho_a\phi\psi_- = \tau_{>a}(\phi\psi_-\rho_a)=
\tau_{>a}(\phi)\tau_{>a}(P_-(\phi))= \tau_{>a}(\phi)$ et de m{\^e}me on calcule
$\psi_+\phi\iota_a=\tau_{\leq a}(\phi)$.

Il reste {\`a} voir l'unicit{\'e} de $\alpha$ dans le point ii). Soit $\beta$ un second
isomorphisme v{\'e}rifiant les propri{\'e}t{\'e}s requises par le point
ii). L'automorphisme $\beta\alpha^{-1}$ de $\bigoplus \HC^q(X)[-q]$
commute {\`a} l'endomorphisme $H(\phi):=\bigoplus \HC^q(\phi)[-q]$. Soit $S\subset
\ZM$ le support cohomologique de $X$. Fixons des
polyn{\^o}mes $(Q_q)_{q\in S}$ tels que 
$$Q_q \in \prod_{p\neq q} P_p.R[X] \hbox{
et } \sum_{q\in S} Q_q = 1. $$
Alors l'endomorphime $Q_q(H(\phi))$ est nul sur les $\HC^p(X)[-p]$,
$p\neq q$ et envoie $\HC^q(X)[-q]$ identiquement dans lui-m{\^e}me.
Comme l'automorphisme $\beta\alpha^{-1}$ commute aux $Q_q(H(\phi))$,
il est de la forme $\bigoplus_q \gamma_q[q]$ avec $\gamma_q \in
\endo{\DC}{\HC^q(X)}.$ Mais comme il doit aussi induire l'identit{\'e}
en cohomologie on a $\gamma_q=\id_{\HC^q(X)}$ et par suite
$\beta\alpha^{-1}= \id$.

\end{proof}

\begin{boite}{Remarque}
La preuve du point i) montre plus pr{\'e}cis{\'e}ment que si on se donne une
famille d'endomorphismes $(\phi_q)_{q\in\ZM}$ de $X$, presque tous
{\'e}gaux {\`a} $\id_X$ et tels que pour tout $q\in\ZM$, $\HC^q(\phi_q)=0$,
alors la compos{\'e}e
$$ \cdots\circ \phi_{q-1}\circ \phi_q \circ \phi_{q+1} \circ \cdots
\;\; \in \endo{\DC}{X}$$
est nulle. (Mais pas n{\'e}cessairement celle dans l'autre sens!)
\end{boite}

\def\zp{\ZM[\frac{1}{p}]}

\subsection{D\'ecomposition "par le niveau" de $\Mo{\Lambda}{G}$} \label{decomposition}

Soit $G={\mathbf G}(F)$ le groupe des points rationnels d'un groupe alg\'ebrique r\'eductif ${\mathbf G}$ d\'efini sur un corps local non-archim\'edien $F$ de caract\'eristique r\'esiduelle $p$.
Le but de cette section est d'\'etendre \`a  la cat\'egorie $\Mo{\zp}{G}$ des
  repr\'esentations lisses de $G$ \`a coefficients dans $\ZM[{1\over p}]$ la
  d\'ecomposition "par le niveau", implicite dans les travaux de Moy et
  Prasad en caract\'eristique z\'ero, et \'etendue au cas d'un corps de
  coefficients de caract\'eristique positive $\neq p$ par Vign\'eras dans
  \cite[II.5]{Vig}.
  
Cette d\'ecomposition intervient \`a plusieurs endroits du texte principal (\ref{finit} et \ref{dimcoh}) par son corollaire suivant : si $(\pi,V)$ est une $\zp$-repr\'esentation lisse, alors le complexe de cha\^ines du  syst\`eme de coefficients sur l'immeuble de Bruhat-Tits associ\'e par Schneider-Stuhler ({\em cf} \cite{SS2}, \cite{Vigsheaves}) \`a $(\pi,V)$ est {\em acyclique}. Lorsque ${\mathbf G}=GL(n)$, on a des r\'esultats un peu plus pr\'ecis qui permettent notamment de prouver le fait \ref{hn}.
  
Les arguments reposent {\em in fine} sur les constructions de Moy et Prasad dans \cite{MP2}. Celles-ci fonctionnent bien sous des hypoth\`eses de ramification mod\'er\'ee pour ${\mathbf G}$. Nous n'avons besoin dans ce texte que du cas ${\mathbf G}$  d\'eploy\'e, donc ces hypoth\`eses sont satisfaites.

\alin{D\'ecomposition de cat\'egories ab\'eliennes}
Nous rappelons ici un peu d'{\em abstract nonsense}.
Soit $\CC$ une cat\'egorie ab\'elienne (avec limites inductives exactes).
Pour une famille $(Q_n)_{n\in\NM}$ d'objets de $\CC$ on consid\`ere les propri\'et\'es suivantes :
\begin{itemize}
	\item (PROJ) Chaque $Q_n$ est projectif et de type fini ("compact").
	\item (DISJ) Si $n\neq m$, alors $\hom{Q_n}{Q_m}{\CC}=0$.
  \item (GEN) Pour tout objet $V$ de $\CC$, on a
    $\hom{\bigoplus_nQ_n}{V}{\CC}\neq 0$.
\end{itemize}
Par ailleurs, pour tout objet $V$ de
$\CC$, posons
$$V_n:= \sum_{\phi\in \hom{Q_n}{V}{G}} \im\phi \subseteq V,$$
un sous-objet de $V$. 
Les propri\'et\'es (PROJ) et (GEN) impliquent
que $V = \sum_n V_n$. La propri\'et\'e (DISJ), toujours
avec (PROJ),  assure que la somme est directe, {\em i.e.} 
$V=\bigoplus_n V_n$. On peut paraphraser cela en introduisant la 
sous-cat\'egorie pleine $\CC_n$ de $\CC$ form\'ee des  objets v\'erifiant $\hom{Q_m}{V}{\CC}=0$ pour tout $m\neq n$. On obtient en effet une d\'ecomposition de $\CC$ en une somme directe de sous-cat\'egories "facteurs directs" 
$ \CC \simeq \bigoplus_n \CC_n$.

Plus g\'en\'eralement, pour $I\subset \NM$, notons $\CC_I$ la sous-cat\'egorie pleine de $\CC$
  form\'ee des  objets v\'erifiant $\hom{Q_m}{V}{\CC}=0$ pour $m\notin I$. 
Alors $\CC_I$ est une sous-cat\'egorie "facteur direct" de $\CC$ et
si $P_I$ est un objet projectif et g\'en\'erateur de $\CC_I$ ({\em i.e.} $\hom{P_I}{V}{\CC}\neq 0$ pour tout objet $V$ de $\CC_I$), alors le foncteur
$$ \application{}{\CC}{\Mo{}{\endo{\CC}{P_I}}}{V}{\hom{P_I}{V}{\CC}}$$
induit une \'equivalence de cat\'egories entre $\CC_I$ et la cat\'egorie $\Mo{}{\endo{\CC}{P_I}}$ de tous les modules \`a droite sur l'anneau $\endo{\CC}{P_I}$.
Une \'equivalence inverse est donn\'ee par :
$$\application{}{\Mo{}{\endo{\CC}{P_I}}}{\CC}{M}{M\otimes_{\endo{\CC}{P_I}} P_I},$$
le produit tensoriel \'etant d\'efini comme un conoyau \`a partir d'une pr\'esentation du $\endo{\CC}{P_I}$-module $M$.

\alin{Types non raffin\'es de Moy-Prasad et d\'ecomposition de $\Mo{\zp}{G}$}
Soit $\IC:=\IC({\mathbf G},F)$ l'immeuble de Bruhat-Tits associ\'e \`a $G$. Pour $x\in \IC$, on notera $G_x$ son fixateur  dans $G$. 

Moy et Prasad ont d\'efini (\cite{MP1}, \cite{MP2} et \cite[II.5]{Vig}), une certaine filtration d\'ecroissante de $G_x$ par des
pro-$p$-sous-groupes ouverts $G_{x,r}, r\in \RM_+$. Les sauts de cette filtration sont discrets et on a des relations de commutateurs $(G_{x,r},G_{x,s})\subset G_{x,r+s}$. Si l'on convient de noter $G_{x,r^+}:=\bigcup_{s>r} G_{x,s}$, alors $G_{x,0^+}$ est le pro-$p$-radical de $G_x$, et pour tout $r>0$, le groupe fini $G_{x,r}/G_{x,r^+}$ est naturellement un $\FM_p$-espace vectoriel.
Ils ont ensuite d\'efini certains caract\`eres complexes des gradu\'es $G_{x,r}/G_{x,r^+}$ appel\'es {\em types non raffin\'es minimaux de niveau $r$}, dont nous noterons l'ensemble $NR_{x,r}$. Ces caract\`eres sont donc \`a valeurs dans l'extension $\ZM[\frac{1}{p},\zeta_p]$ si $\zeta_p$ est une racine $p$-i\`eme de l'unit\'e.
Enfin, Moy et Prasad ont d\'efini un ensemble $PO$ de "points optimaux" dans l'immeuble,  fini modulo action de $G$, et nous noterons $(r_n)_{n\in\NM}$ une \'enum\'eration des sauts des filtrations associ\'ees aux points de $PO$.

Posons maintenant $Q_0 : = \bigoplus_{x} \cind{G_{x,0^+}}{G}{\zp}$ o\`u $x$ d\'ecrit un
ensemble (fini) de repr\'esentants des $G$-orbites de sommets de $\IC$.
Pour $r\in \RM_+$, posons (comme dans la remarque de \cite[p. 136]{Vig}) 
$$ P(r) := \bigoplus_{x\in PO,\chi\in NR_{x,r}} \cind{G_{x,r}}{G}{\chi} $$
que l'on voit comme une
repr\'esentation \`a coefficients dans $\ZM[\frac{1}{p},\zeta_p]$. Puisqu'un type non raffin\'e minimal tordu par l'action d'un \'el\'ement de $\gal(\QM(\zeta_p)|\QM)$ est encore un type non raffin\'e minimal, cette repr\'esentation 
se descend  en une $\zp$-repr\'esentation que nous noterons encore $P(r)$.

\begin{lemme}
La famille $Q_n:=P(r_n)$, $n\in\NM$ d'objets de $\Mo{\zp}{G}$ v\'erifie les propri\'et\'es (PROJ), (GEN) et (DISJ).
\end{lemme}
\begin{proof}
D'apr\`es \cite[II.5]{Vig}, pour tout corps alg\'ebriquement clos $R$ de caract\'eristique $\neq p$, la famille de repr\'esentations $(Q_n\otimes_{\zp} R)_{n\in\NM}$ de $\Mo{R}{G}$ v\'erifie les propri\'et\'es (PROJ), (DISJ) \cite[II.5.8]{Vig} et (GEN) \cite[II.5.3]{Vig} de la section pr\'ec\'edente. Nous allons montrer que cela implique formellement qu'il en est de m\^eme de la famille $(Q_n)_{n\in\NM}$ dans $\Mo{\zp}{G}$.

(PROJ) :  En tant que somme d'induites de $\zp$-repr\'esentations de type fini  de pro-$p$-sous-groupes ouverts, $P(r)$ est  
projective et de type fini dans $\Mo{\zp}{G}$. 

(DISJ) :  puisque $Q_m$ est sans torsion,
$\hom{Q_n}{Q_m}{G} \injo \hom{Q_n\otimes \CM}{Q_m\otimes \CM}{G}$. Ce
dernier est nul par \cite[II.5.8]{Vig} appliqu\'e \`a $R=\CM$.

(GEN):  soit $V$ un objet de $\Mo{\zp}{G}$ tel qu'il existe $l\neq
p$ premier tel que $V_l:=\{v\in V, lv=0\} \neq 0$. On peut voir $V_l$
comme une $\FM_l$-repr\'esentation de $G$. 
 On sait alors par
\cite[II.5.3]{Vig} qu'il existe $n\in \NM$ et un morphisme non nul
$\phi :\; Q_n \To{} V_l\otimes \o\FM_l$. Par engendrement fini de $Q_n$,
ce morphisme se factorise par $V_l\otimes \FM_{l^k}$ pour un certain
$k\in \NM$. D'o\`u un morphisme non nul $Q_n \To{} (V_l)^k$
et par suite l'existence d'un morphisme  non nul $Q_n \To{} V_l$
que l'on peut composer avec l'injection $V_l\injo V$.

Si maintenant $V_l=0$ pour tout $l\neq p$, c'est \`a dire si $V$ n'a
pas de torsion, $V$ se plonge dans $V\otimes \QM$. Comme
pr\'ec\'edemment on d\'eduit de \cite[II.5.3]{Vig} l'existence d'un
morphisme non nul $Q_n \To{} V\otimes \QM$. Par engendrement fini de
$Q_n$, on peut multiplier par un ``d\'enominateur commun'' pour
obtenir un morphisme \`a image dans $V$. 

\end{proof} 
 
 Dans \cite{SS2}, Schneider et Stuhler ont d\'efini une autre filtration $(U_x^{(e)})_{e\in\NM}$ du fixateur d'un sommet $x$ de $\IC$. 
Fixons $e\in \NM$ et notons ${\Mo{\zp}{G}}^e$ la sous-cat\'egorie "facteur direct" ({\em cf} paragraphe pr\'ec\'edent) des objets de $\Mo{\zp}{G}$ tels que $\hom{P(r_n)}{V}{G}=0$ pour $r_n> e$.

\begin{lemme}
Les objets de la sous-cat\'egorie pleine $\Mo{\zp}{G}^e$ sont ceux qui sont engendr\'es par leur $U_x^{(e)}$-invariants, $x$ d\'ecrivant l'ensemble des sommets de $\IC$.
\end{lemme}
\begin{proof}
Les arguments de Vign\'eras dans le cas d'un corps alg\'ebriquement clos \cite[Thm 2.2]{Vigsheaves} reposent sur la comparaison des filtrations de Schneider-Stuhler et Moy-Prasad et s'appliquent sans changement ici.
\end{proof}

\begin{coro}
Soit $V$ un objet de $\Mo{\zp}{G}^e$. Le complexe de cha\^ines du  syst\`eme de coefficients $\gamma_e(V)$ sur $\IC$ associ\'e par Schneider-Stuhler ({\em cf} \cite{SS2} et \cite[2.5]{Vigsheaves}) \`a $V$ est {\em acyclique}. 
\end{coro}
\begin{proof}
Comme il est expliqu\'e dans la preuve de \cite[Prop. 2.6]{Vigsheaves}, le lemme pr\'ec\'edent permet de se ramener aux cas $V=\ind{U_x^{(e)}}{G}{\zp}$, pour $x$ sommet de $\IC$, et dans ce dernier cas, l'argument est "g\'eom\'etrique" (le complexe associ\'e s'identifie au complexe de cocha\^ines du syst\`eme de coefficients constant d'un certain sous-espace simplicial de l'immeuble dont on montre la contractibilit\'e).
\end{proof}

\begin{coro} \label{corofinit}
Soit $Z\subset Z(G)$ un sous-groupe cocompact du centre de $G$ et $\Lambda$ un anneau o\`u $p$ est inversible. Alors 
\begin{enumerate}
	\item tout $V\in \Mo{\Lambda}{G/Z}^e$ admet une r\'esolution par des objets projectifs de type fini de $\Mo{\Lambda}{G/Z}^e$.
	\item si $\Lambda$ est banal pour $G$, alors il existe une telle r\'esolution de longueur $\hbox{rang s/s}(G)+1$.
\end{enumerate}
\end{coro}
\begin{proof}
L'argument est le m\^eme que celui de \cite[2.11]{Vigsheaves}.
Les termes du complexe de cha\^ines du syst\`eme de coefficients $\gamma_e(V)$ du corollaire pr\'ec\'edent sont des sommes d'objets de la forme $\cind{P_F^+}{G_d}{V^{U_F^{(e)}}}$ o{\`u} $P_F^+$ est
le stabilisateur d'une facette $F$ et $U_F^{(e)}$ est le
sous-groupe de conguences de $P_F^+$ de niveau $e$ d{\'e}fini dans
\cite[ch. 1]{SS2}. Comme l'induction \`a supports compacts respecte le
caract\`ere projectif,  ces objets sont projectifs lorsque $\Lambda$ est banal, d'o\`u le point ii). Lorsque $\Lambda$ n'est pas banal
on choisit pour chaque tel facteur une
r\'esolution de $V^{U_F^{(e)}}$ par des $\Lambda (P_F^+/ZU_F^{(e)})$-modules
projectifs et finis  et on en d\'eduit la r\'esolution cherch\'ee.
\end{proof}

\begin{rema} Dans le cas ${\mathbf G}=GL(n)$, on a une seule classe de conjugaison de sommets dans $\IC$ et l'on peut prendre pour prog\'en\'erateur de $\Mo{\zp}{G}^e$, la repr\'esentation $\ind{U_x^{(e)}}{G}{\zp}$. Par d\'efinition,  $U_x^{(e)}$ est conjugu\'e \`a $1+\varpi_F^e\MC_n(F)$ pour $e>0$.
\end{rema}

 \subsection{Variations sur une construction de Berkovich}

Dans cette section, on reprend le contexte et les notations de la partie \ref{berkovich}. Nous expliquons les arguments de Berkovich pour prouver la proposition \ref{torsion0}. Comme on l'a d\'eja dit dans le texte, nous adoptons un langage plus \'el\'ementaire que celui de \cite{Bicunp}, et les erreurs eventuelles sont exclusivement d\^ues \`a l'auteur de ces lignes.
Puis nous prolongeons un peu ces arguments pour prouver le lemme \ref{ladic0} et la proposition \ref{ladic2}, en y ajoutant des techniques de Jannsen.

\alin{Un crit{\`e}re d'acyclicit{\'e}} \label{critacy}
Soit  $X$  un espace analytique, que l'on suppose  paracompact. Soit $\pi^X :
X_{et}\To{} |X|$ le morphisme de sites du site {\'e}tale de
$X$ vers le site topologique de $X$.
Nous dirons qu'un faisceau {\'e}tale $\FC$ sur $X$ est ``{\'e}tale-c-mou'' si 
\begin{enumerate}
\item Pour tout $x\in X$, $\FC_x$ est un
  $\gal(\HC(x)^a/\HC(x))$-module acyclique.
\item Le faisceau ${\pi^X_*}(\FC)$ est c-mou.
\end{enumerate}
D'apr{\`e}s \cite{Bic2}, tout faisceau {\'e}tale-c-mou  est
$\Gamma_c$-acyclique. On peut aussi formuler i) de mani{\`e}re plus
canonique (sans choix d'une cl{\^o}ture s{\'e}parable) : en consid{\'e}rant $\FC_x$ comme un
faisceau sur le site {\'e}tale de $\HC(x)$, la
condition i) est {\'e}quivalente {\`a} demander que pour toute extension
Galoisienne finie $K$ de $\HC(x)$, le groupe ab{\'e}lien
$\FC_x(K)$ soit $\gal(K/\HC(x))$-acyclique.

\ali
Soit $X$ un espace $K$-analytique localement paracompact muni d'une
action continue d'un groupe topologique $G$. On d{\'e}finit
$$X_{disc} := \bigsqcup_{x\in X} \MC(\HC(x)) \;\; \hbox{ et } \;\;
 \application{\nu:\;}{X_{disc}}{X}{x}{x}.$$
Avec les d{\'e}finitions de \cite{Bic2}, $X_{disc}$ n'est pas un espace
analytique sur $K$. On peut n{\'e}anmoins 
lui associer un site {\'e}tale
(et c'est tout ce qui nous
importe) d{\'e}fini  comme le produit
des sites {\'e}tales des $\MC(\HC(x))$. Il est muni d'une action de $G$ et
 $\nu$ induit un morphisme de sites $G$-{\'e}quivariant.
Un faisceau {\'e}tale ab{\'e}lien sur $X_{disc}$ est donc simplement une famille
$(\FC_x)_{x\in X}$ de
$\gal(\HC(x)^a/\HC(x))$-modules lisses, pour $x\in X$. Tout point $x$
est ouvert dans $X_{disc}$ et les sections au-dessus de $x$ sont
donn{\'e}es par $\FC(x)=\FC_x^{\gal(\HC(x)^a/\HC(x))}$. 

Rappelons aussi qu'on a not{\'e} $\pi^X : X_{et}\To{} |X|$ le morphisme du site {\'e}tale
de $X$ vers le site topologique de $X$.

\begin{lemme} \label{lemmemou}
  Soit $\FC \in \FC(X_{disc},\Lambda G_{disc})$,
  \begin{enumerate}
  \item $\pi^X_* \nu_*(\FC)$ est un faisceau flasque sur $|X|$.
\item $\pi^X_* (\nu_*(\FC)^\infty)$ est un faisceau c-mou sur $|X|$.
  \end{enumerate}
\end{lemme}
\begin{proof}
Pour tout ouvert $U\subset X$, on a 
$$ \pi^X_*\nu_*(\FC)(U) = \prod_{x\in U} \FC(x),$$
expression qui montre clairement l'assertion i).

Pour le ii) on suit un argument de Berkovich. Soit $\Sigma$ un compact
de $X$ et $f \in \Gamma(\Sigma,\pi^X_* (\nu_*(\FC)^\infty))$. Il
existe un voisinage ouvert $U \supset \Sigma$ et une section $f \in
\pi^X_* (\nu_*(\FC)^\infty)(U)$ qui induit $f$ sur $\Sigma$. Soit
alors $U\supset \o{V} \supset \Sigma$ un domaine analytique compact de
$X$. On d{\'e}finit 
$$ f'\in \Gamma_c(X,\pi^X_* (\nu_*(\FC))) \subset \prod_{x\in X} \FC(x)
$$
par $f'(x) =0 $ si $x\notin \o{V}$ et $f'(x)=f(x)$ pour $x\in
\o{V}$. On a {\'e}videmment $f'_{|\Sigma}=f_{|\Sigma}$ ; il nous suffira
donc de montrer que $f'\in \Gamma_c(X,\pi^X_*
(\nu_*(\FC))^\infty)$, c'est-{\`a}-dire que $f'$ est stabilis{\'e}e par un
sous-groupe ouvert compact de $G$. Par d{\'e}finition de $f$ et par
compacit{\'e} de $\o{V}$, il existe des ouverts distingu{\'e}s $(U_i)_{i\in
  I}$, avec $I$ fini tels que $\o{V}=\cup_{i\in I} (\o{V}\cap U_i)$ et
$f_{|U_i}$ est $J_i$-invariante pour un certain sous-groupe ouvert compact
$J_i$ de $G$. Comme le stabilisateur $G_{\o{V}}$ de $\o{V}$ dans $G$ est ouvert,
il s'ensuit que $f'_{|U_i}$ est invariante sous le sous-groupe ouvert
compact $J_i\cap G_{\o{V}}$, et par cons{\'e}quent que $f$ est invariante sous
le sous-groupe ouvert compact $\cap_{i\in I} (J_i\cap G_{\o{V}})$.
\end{proof}

Pour pr{\'e}parer le point ii) du lemme suivant, nous introduisons la
notation $G_x$ pour le stabilisateur de $x$ dans $G$ et pour toute
$\HC(x)$-alg{\`e}bre {\'e}tale $K$, nous notons $(G_x)_K$ un sous-groupe
ouvert compact de $G_x$ suffisamment petit pour que l'action de $(G_x)_K$ sur
$\HC(x)$ s'{\'e}tende canoniquement {\`a} une action sur $K$ (cas particulier
de \cite[Key Lemma 7.2]{Bic3}).

\begin{lemme} \label{lemmegalacy}
   Soit $\FC \in \FC(X_{disc},\Lambda G_{disc})$.
   \begin{enumerate}
   \item Si pour tout $x\in X$, le $\gal(\HC(x)^a/\HC(x))$-module lisse
     $\FC_x$ est acyclique (pour la cohomologie ``continue''), alors
     il en est de m{\^e}me pour les fibres $\nu_*(\FC)_x$ de $\nu_*(\FC)$ en
     tout $x$.
   \item
    Si pour tout $x\in X$, toute extension Galoisienne finie $K$ de
    $\HC(x)$ et tout  sous-groupe ouvert compact $H$
     de $(G_x)_K$, le  $\gal(K/\HC(x))$-module
     $\FC_x(K)^{H}$ est acyclique 
 alors  les fibres
     $(\nu_*(\FC)^\infty)_x$ de $\nu_*(\FC)^\infty$ sont acycliques en tout $x$.

   \end{enumerate}
\end{lemme}

\begin{proof} i) Ici encore, tous les arguments sont d{\^u}s {\`a} Berkovich (et
  les erreurs {\`a} l'auteur).
Soit $K$ une extension  Galoisienne finie de
$\HC(x)$. 
On a la
description suivante
de $\nu_*(\FC)_x(K)$ :
fixons un
morphisme {\'e}tale $X'\To{f} X$ tel que $f^{-1}(x)=\{x'\}$ et
$\HC(x')\simeq K$ (on fixe alors un tel isomorphisme).
Alors
$$ (\nu_*(\FC))_x(K) = \limi{x\in V\subset X} \Gamma(f^{-1}(V),\nu_*(\FC)),$$
la limite inductive portant sur les voisinages ouverts (pour la
topologie analytique) de $x$ dans
$X$. Par un argument de cofinalit{\'e}, on peut se restreindre aux $V$
tels que $f^{-1}(V) \To{} V$ est {\'e}tale Galoisien de groupe
$\gal(K/\HC(x))$. Par commutation de la cohomologie aux limites inductives, il
nous suffira de montrer que pour tout tel $V$, le $\gal(K/\HC(x))$-module
$\Gamma(f^{-1}(V),\nu^*(\FC))$ est cohomologiquement trivial.
On a
$$ \Gamma(f^{-1}(V),\nu_*(\FC)) = \prod_{y\in V} \prod_{y'\in f^{-1}(y)}
\FC_y(\HC(y')). $$ 
  Comme l'action de $\gal(K/\HC(x))$ sur $f^{-1}(V)$ respecte les fibres de $f$,
  il nous suffit de voir que pour tout $y\in V$, le $\gal(K/\HC(x))$-module 
\ini\begin{equation}\label{galmod}
 \prod_{y'\in f^{-1}(y)}
\FC_y(\HC(y')) 
\end{equation}  est cohomologiquement trivial. Fixons $y'\in
f^{-1}(y)$ et notons $\gal_{y'}$ son stabilisateur dans $\gal(K/\HC(x))$. L'action de $\gal_{y'}$
sur le $\HC(y)$-espace analytique $\MC(\HC(y'))$ induit un isomorphisme $\gal_{y'} \simto
\gal(\HC(y')/\HC(y))$. On voit alors que le $\gal(K/\HC(x))$-module \ref{galmod} est induit du $\gal_{y'}$-module $\FC_y(\HC(y'))$ qui
par hypoth{\`e}se est cohomologiquement trivial. Le lemme de Shapiro
montre donc la $\gal(K/\HC(x))$-acyclicit{\'e} de ce module.

Preuve de ii). Reprenons l'argument pr{\'e}c{\'e}dent avec maintenant le faisceau
$\nu_*(\FC)^\infty$. On a cette fois
$$ (\nu_*(\FC)^\infty)_x(K) = \limi{x\in V\subset X}\limi{H\subset
  G_{f^{-1}(V)}} \Gamma(f^{-1}(V),\nu_*(\FC))^H$$
o{\`u} la premi{\`e}re limite porte sur les voisinages ouverts distingu{\'e}s de
$x$ tels que le morphisme $f^{-1}(V)\To{} V$ soit Galoisien de groupe
$\gal(K/\HC(x))$, et la deuxi{\`e}me limite porte sur les pro-$p$-sous-groupes ouverts
 $H$ de
$G_{f^{-1}(V)}$.
Comme pr{\'e}c{\'e}demment il nous suffira de montrer que pour tout tel couple
$(V,H)$ le $\gal(K/\HC(x))$-module 
$$ \Gamma(f^{-1}(V), \nu_*(\FC))^H = \left(\prod_{y\in V} \prod_{y'\in f^{-1}(y)}
\FC_y(\HC(y')) \right)^H $$
est cohomologiquement trivial.
Soit $[V/H]$ un ensemble de repr{\'e}sentants des orbites de $H$ dans $V$
et pour $y\in [V/H]$, notons $H_y$ son fixateur dans $H$. Alors 
$$\left(\prod_{y\in V} \prod_{y'\in f^{-1}(y)}
\FC_y(\HC(y')) \right)^H = \prod_{y\in [V/H]} \left(  \prod_{y'\in
  f^{-1}(y)} \FC_y(\HC(y'))\right)^{H_y}.$$
Fixons alors $y\in [V/H]$ ; tous les $y'\in f^{-1}(y)$ ont le m{\^e}me
fixateur $H_y$ que $y$ dans $H$ (car $\gal(K/\HC(x))$ les permute transitivement).
 Ainsi le $\gal(K/\HC(x))$-module
$$  \left(  \prod_{y'\in
  f^{-1}(y)} \FC_y(\HC(y'))\right)^{H_{y}}$$
est induit du  $\gal_{y'}$-module $\FC_y(\HC(y'))^{H_{y}}$. Par
notre hypoth{\`e}se et le lemme de Shapiro, il s'ensuit que ce
$\gal(K/\HC(x))$-module est  acyclique. 

\end{proof}

Consid{\'e}rons  maintenant  la r{\'e}union disjointe d'espaces analytiques
sur $K$
$$GX_{disc}:=\bigsqcup_{(g,x)\in G\times X}
\MC(\HC(x))$$
 sur laquelle $G$ agit par $h(g,x):=(hg,x)$. On a un
morphisme  $G$-{\'e}quivariant 
$$ \application{\mu:\;}{GX_{disc}}{X_{disc}}{(g,x)}{gx} $$
qui induit un morphisme $G$-{\'e}quivariant de sites {\'e}tales.

\begin{lemme} \label{lemmetech}
Soit $\FC \in \FC(GX_{disc},\Lambda G_{disc})$. Si pour
  tout $(g,x)\in G\times X$, le $\gal(\HC(x)^a/\HC(x))$-module lisse
  $\FC_{g,x}$ est cohomologiquement trivial, alors
 le faisceau $\mu_*(\FC)$
satisfait l'hypoth{\`e}se du point ii) du lemme pr{\'e}c{\'e}dent.
\end{lemme}
\begin{proof}
Soit $K$ une $\HC(x)$-alg{\`e}bre {\'e}tale. Un calcul rapide de la fibre en $x\in X$ du faisceau
$\mu_*(\FC)$ montre que :
$$ \mu_*(\FC)_x(K) = \prod_{g\in G}
\FC_{g,g^{-1}x}(K^g) $$
o{\`u} $K^g$ est la $\HC(g^{-1}x)$-alg{\`e}bre {\'e}tale  $K\otimes_{\HC(x)}
\HC(g^{-1}x)$, produit tensoriel pris pour l'isomorphisme
$g^*:\;\HC(x)\simto \HC(g^{-1}x)$ donn{\'e} par l'action de $g$.
Soit $H$ un sous-groupe de $(G_x)_K$,
il agit par permutation $\FC_{g,g^{-1}x}(K^g) \To{\tau_\FC(h)} \FC_{hg,g^{-1}x}(K^g)
= \FC_{hg,(hg)^{-1}x}(K^{hg})$. Il s'ensuit que 
$$\mu_*(\FC)_x(K)^{H} \simeq  \prod_{g\in H\ba G} \FC_{g,g^{-1}x}(K^g).$$
Supposons $K$ extension Galoisienne de $\HC(x)$. Par notre hypoth{\`e}se,
les $\gal(K/\HC(x))$-modules 
$\FC_{g,g^{-1}x}(K^g)^{H^g}$ sont acycliques et par cons{\'e}quent leur produit
aussi.

\end{proof}

Remarquons ici que l'astuce consistant {\`a} introduire l'espace
$GX_{disc}$ est encore une fois d{\^u}e {\`a} Berkovich.

\begin{prop} \label{alaBerk}
  Tout faisceau {\'e}tale $G$-{\'e}quivariant discret/lisse sur $X$ se plonge
  dans un objet injectif de $\FC(X,\Lambda G)$ dont le faisceau
  sous-jacent est {\'e}tale-c-mou.
\end{prop}
\begin{proof}
Soit $\FC \in \FC(X,\Lambda G)$.
Choisissons un plongement $\mu^*\nu^*(\FC)\injo \IC$ dans un objet injectif
$\IC$ de $\FC(GX_{disc},\Lambda G_{disc})$ dont le faisceau sous-jacent est un objet injectif dans
$\FC(GX_{disc},\Lambda)$ (il en existe toujours, {\em cf} lemme \ref{induction}).
 Alors $\nu_*\mu_*(\IC)$ est un
objet injectif de $\FC(X,\Lambda G_{disc})$ car 
le foncteur $\nu_*\mu_*$ admet un
adjoint {\`a} gauche exact, $\mu^*\nu^*$.
Comme on le v{\'e}rifie immediatement sur les fibres, on  a un plongement
$$ \FC \injo (\nu\mu)_*(\nu\mu)^*(\FC) \injo \nu_*\mu_*(\IC). $$
 Appliquant le foncteur de lissification (exact \`a gauche), on trouve alors un
plongement
$$ \FC \injo (\nu_*\mu_*(\IC))^\infty$$
et le faisceau $(\nu_*\mu_*(\IC))^\infty$ est un objet injectif de
$\FC(X,\Lambda G)$ car le foncteur de
lissification est adjoint {\`a} droite du plongement canonique.

Maintenant le faisceau $\IC$ {\'e}tant injectif, les
$\gal(\HC(x)^a/\HC(x))$-modules lisses $\IC_{g,x}$ sont acycliques. Donc
par le lemme pr{\'e}c{\'e}dent, le faisceau $\mu_*(\IC_{g,x})$
v{\'e}rifie l'hypoth{\`e}se du point ii) du
  lemme \ref{lemmegalacy}. Par ce point ii) et le point ii)
du lemme \ref{lemmemou},
le faisceau $\nu_*\mu_*(\IC)^\infty$ 
est {\'e}tale-c-mou.

\end{proof}

\alin{Preuve de la proposition \ref{ladic2}} \label{preuveladic2} Il
s'agit de montrer que le foncteur 
 $\limproj^\infty$ envoie assez d'objets injectifs sur des objets
$\Gamma_!$-acycliques. Soit $(\FC_n)_n$ un objet de
$\FC(X,\Lambda_\bullet G)$. Comme dans la preuve de la proposition
\ref{alaBerk}, pour chaque $n$ on peut 
trouver un plongement $\FC_n\To{i_n} \nu_*\mu_*(\IC_n)^\infty$ avec
$\IC_n$ un objet injectif de $\FC(GX_{disc},\Lambda/\mG^n)$.
Rappelons que 
les faisceaux
$\nu_*\mu_*(\IC_n)^\infty$ sont injectifs dans $\FC(X,\Lambda_n G)$.
Posons alors pour tout $n$
$$\JC_n := \prod_{k=1}^n \nu_*\mu_*(\IC_k)^\infty.$$  Les morphismes de
projections sur les premiers termes  $\JC_n \To{} \JC_m$ pour $m<n$
d{\'e}finissent un objet de $\FC(X,\Lambda_\bullet G)$ qui est injectif.
De plus les morphismes 
$$ \prod_{k=1}^n i_k \circ f_{n,k} :\;\; \FC_n \To{} \JC_n$$ 
(o{\`u} $f_{n,k}
:\FC_n \To{} \FC_k$ est le morphisme de transition) d{\'e}finissent un
plongement du syst{\`e}me $(\FC_n)_n$ dans $(\JC_n)_n$.
Par ailleurs on a 
$$ \limproj^\infty( (\JC_n)_n) = \left(\prod_{n}
\mu_*\nu_*(\IC_n)^\infty\right)^\infty \simto \left(\prod_{n}
\mu_*\nu_*(\IC_n)\right)^\infty \simeq \left(
\mu_*\nu_*(\prod_n \IC_n)\right)^\infty. $$
L'isomorphisme du milieu est un fait g{\'e}n{\'e}ral : le morphisme canonique
qui va du lissifi{\'e} d'un
produit de lissifi{\'e}s vers le  lissifi{\'e} du produit est un isomorphisme.
Le dernier isomorphisme vient de la commutation des images directes et
des produits.
Maintenant la "fibre" du produit des $\IC_n$ en $(g,x)\in GX_{disc}$
s'{\'e}crit :
$$ \left(\prod_n \IC_n \right)_x \simeq \left(\prod_n (\IC_n)_x
\right)^\infty $$
o{\`u} le signe $\infty$ d{\'e}signe le $\gal(\HC(x)^a/\HC(x))$-lissifi{\'e} du
module produit ({\`a} ne pas confondre avec le signe $\infty$ au-dessus qui
d{\'e}signait le $G$-lissifi{\'e} d'un faisceau). Comme les $(\IC_n)_x$ sont
$\gal(\HC(x)^a/\HC(x))$-acycliques, leur produit lisse l'est aussi.
On peut alors appliquer les lemmes \ref{lemmetech}, \ref{lemmegalacy}
ii) et \ref{lemmemou} ii) pour en d{\'e}duire que le faisceau
$\limproj^\infty((\JC_n)_n)$ est {\'e}tale-c-mou et par cons{\'e}quent
$\Gamma_!$-acyclique.

\alin{Preuve du lemme \ref{ladic0}}
 Il s'agit de voir que $\iota_{disc}$ envoie suffisamment
d'objets de $\FC(X,\Lambda_\bullet G)$ sur des  objets
$\Gamma_c^{disc}$-acycliques.  Soit $\UM(X)$ l'ensemble des ouverts
distingu{\'e}s de $X$. D'apr{\`e}s \cite[Prop. 4.1.8]{Fargues}, on a pour tout
$p>0$ une suite exacte
$$ 0\To{} \limi{U\in\UM(X)} R^1\limp{n} H^{p-1}_c(U,\FC_n) \To{}
H^p_c(X,(\FC_n)_{n}) \To{} \limi{U\in \UM(X)} \limp{n} H^p_c(U,\FC_n)
\To{} 0 $$
de laquelle on tire que si $(\FC_n)$ est un syst{\`e}me projectif de
faisceaux qui pour tout $U\in \UM(X)$ sont $\Gamma_!(U,.)$-acycliques 
et tel que pour tout $U\in \UM(X)$, le syst{\`e}me
projectif de $\Lambda$-modules $(\Gamma_c(U,\FC_n))_n$ satisfait la
condition de Mittag-Leffler, alors le syst{\`e}me $(\FC_n)_n$ est
$\Gamma_c$-acyclique. 
En particulier, le syst{\`e}me $(\JC_n)_n$ construit au paragraphe
\ref{preuveladic2} est $\Gamma_c$-acyclique (et donc
$\Gamma_c^{disc}$-acyclique).
D'apr{\`e}s \ref{preuveladic2}, il y a assez de tels objets.

\alin{Induction} \label{induction}
 Soit $G$ un groupe agissant contin{\^u}ment sur $X$ et
$H\subset G$ un sous-groupe ferm{\'e}. Alors le foncteur de restriction
${\Res{G}{H}}:\; {\FC(X,\Lambda G)}\To{}{\FC(X,\Lambda H)}$ admet un
adjoint {\`a} droite $\Ind{H}{G}$. 

Pour construire $\Ind{H}{G}$, introduisons la petite cat{\'e}gorie $[G/H]$ dont
l'ensemble des objets est $G$ et les morphismes sont donn{\'e}s par 
$$ \cas{\hom{g_1}{g_2}{[G/H]}}{g_1\To{h} g_2}{g_1=g_2h \hbox{ pour }
  h\in H}{\emptyset}{g_1H\neq g_2 H}$$
avec une loi de composition {\'e}vidente.
Soit $\FC$ un faisceau $H$-{\'e}quivariant dont on note simplement $\tau$
la structure $H$-{\'e}quivariante. On d{\'e}finit un foncteur $[G/H]\To{}
\FC(X,\Lambda)$ en associant {\`a} tout $g\in G$ le faisceau $g_*\FC$ et {\`a}
toute fl{\`e}che $g_1\To{h} g_2$ le morphisme ${g_2}_*(\tau(h)) :\;
{g_1}_*\FC \To{} {g_2}_*\FC$. On pose alors 
$$ \Ind{H}{G_{disc}}(\FC):= \limp{[G/H]} g_*\FC $$ que l'on munit de
la structure $G$-{\'e}quivariante induite par les isomorphismes de
``permutation'' 
$$\sigma(\gamma) : \gamma_*\limp{[G/H]} g_*\FC  = \limp{[G/H]} (\gamma
  g)_*\FC \simto \limp{[G/H]} g_*\FC .$$
Enfin on pose $\Ind{H}{G}(\FC):=\Ind{H}{G_{disc}}(\FC)^\infty$. Si
$(\FC,\tau)$ est un faisceau $G$-{\'e}quivariant, on a un morphisme
canonique $\FC \To{\tau} \limp{[G/H]}(g_*\FC)$ associ{\'e} {\`a} la
collection de morphismes $\FC \To{g_*(\tau(g^{-1}))} g_*\FC$. D'autre
part si $\GC$ est un faisceau $H$-{\'e}quivariant, on a la projection
$\limp{[G/H]}(g_*\GC) \To{} {1_G}_*\FC=\FC$. On v{\'e}rifie alors
simplement que les $2$ morphismes de foncteurs $\id \To{} \Ind{H}{G}\circ
\Res{G}{H}$ et $\Res{G}{H}\circ \Ind{H}{G} \To{} \id$  ainsi obtenus
font de $(\Res{G}{H},\Ind{H}{G})$ une paire de foncteurs adjoints.

\begin{boite}{Remarque :}
Ici aussi le formalisme de \cite{Bicunp} est beaucoup plus efficace :
avec les notations de ce paragraphe et de la remarque du paragraphe \ref{fgeq}, il y
a un morphisme canonique de 
sites $X(H)_{et}\To{\alpha} X(G)_{et}$. Via les isomorphismes canoniques $\FC(X,\Lambda G)
\simeq   \FC(X(G),\Lambda)$ et respectivement pour $H$, le foncteur de
restriction $\Res{G}{H}$ n'est autre que $\alpha^*$ et par cons{\'e}quent
l'induction $\Ind{H}{G}$ s'identifie {\`a} $\alpha_*$.
\end{boite}

\alin{La r{\'e}solution fonctorielle de Godement-Berkovich}

Ce paragraphe n'est finalement pas utilis\'e dans ce texte. Mais il l'est dans d'autres articles, comme \cite{HaCusp} ou \cite{HaKyoto}. Nous esp\'erons qu'il pourra servir de r\'ef\'erence en attendant la publication vivement souhait\'ee du manuscrit de Berkovich \cite{Bicunp}.
Introduisons,
toujours suivant Berkovich, l'espace
$$ \wt{GX}_{disc} :=\bigsqcup_{(g,x)\in G\times X} \MC(\wh{\HC(x)^a}) $$
o{\`u} $\wh{\HC(x)^a}$ est la compl{\'e}tion d'une cl{\^o}ture alg{\'e}brique de
$\HC(x)$. Le choix des cl{\^o}tures est fait de sorte que pour tout
$(g,x)$ on a $\HC(gx)^a=\HC(x)^a\otimes_{\HC(x),g_*} \HC(gx)$. 
 On a un morphisme analytique $G$-{\'e}quivariant tautologique
$\rho:\;\;\wt{GX}_{disc} \To{} GX_{disc}$. 
 \begin{lemme}
   Pour tout faisceau $\FC\in\FC(\wt{GX}_{disc},\Lambda G_{disc})$, le
   $\gal(\HC(x)^a/\HC(x))$-module lisse 
  $\rho_*(\FC)_{g,x}$ est cohomologiquement trivial.
 \end{lemme}
 \begin{proof} 
   En effet, on v{\'e}rifie imm{\'e}diatement que pour tout $(g,x)\in G\times
   X$, on a
$$ \rho_*(\FC)_{g,x} = \CC^\infty(\gal(\HC(x)^a/\HC(x)),\Lambda)\otimes
\FC_{g,x}.$$

 \end{proof}

Le r{\'e}sultat suivant est un des points d'orgue de \cite{Bicunp}. La
formulation dans {\em loc. cit} est diff{\'e}rente et plus
intelligente.
 Celle ci-dessous est
peut-{\^e}tre plus {\'e}l{\'e}mentaire. 
\begin{prop} (Berkovich \cite[Prop. 1.5.1]{Bicunp}) Soit $B$
  l'endofoncteur de $\FC(X,\Lambda G)$ d{\'e}fini par 
  $B:= \infty\circ(\nu\mu\rho)_*(\nu\mu\rho)^*$. 
Alors 
\begin{enumerate}
\item $B$ est exact.
\item Le morphisme d'adjonction de foncteurs
  $\id_{\FC(X,\Lambda G)} \To{} B$ est un monomorphisme et pour tout $\FC \in
  \FC(X,\Lambda G)$, le faisceau {\'e}tale sous-jacent {\`a} $B(\FC)$ est
  {\'e}tale-c-mou.
\item Si $H$ est un autre groupe agissant sur $X$ et commutant {\`a}
  l'action de $G$, alors $B$ d{\'e}finit aussi un endofoncteur de
  $\FC(X,\Lambda G\times H_{disc})$.
\end{enumerate}
\end{prop}
\begin{proof}
Le point ii) est une cons{\'e}quence du lemme pr{\'e}c{\'e}dent et des lemmes \ref{lemmemou},
  \ref{lemmegalacy} et \ref{lemmetech}.

Pour prouver le point i), remarquons d'abord qu'il suffit de montrer
que le foncteur $\infty\circ(\nu\mu\rho)_*$ est exact, ce que l'on va
v{\'e}rifier sur les fibres.
Soit $\FC\in \FC(\wt{GX}_{disc},\Lambda G_{disc})$ et $x\in X$.
Reprenons les notations
de la preuve du lemme \ref{lemmegalacy}, c'est-{\`a}-dire une extension
finie $K$ de $\HC(x)$, un morphisme $X'\To{f} X$ etc... On a 
$$ ((\nu\mu\rho)_*(\FC)^\infty)_x(K) = \limi{x\in V\subset X}\limi{H\subset
  G_{f^{-1}(V)}} \Gamma(f^{-1}(V),(\nu\mu\rho)_*(\FC))^H$$
avec
\begin{eqnarray*}
\Gamma(f^{-1}(V),(\nu\mu\rho)_*(\FC))^H & = & \left(\prod_{y\in V} \prod_{y'\in f^{-1}(y)}
\mu_*\rho_*(\FC)_y(\HC(y'))\right)^H \\
& = & \prod_{y\in [V/H]} \left(  \prod_{y'\in
  f^{-1}(y)} \mu_*\rho_*(\FC)_y(\HC(y'))\right)^{H_y}. 
\end{eqnarray*}
Par ailleurs, les preuves du lemme pr{\'e}c{\'e}dent et de \ref{lemmetech}
montrent que
\begin{eqnarray*}
 \mu_*\rho_*(\FC)_y(\HC(y')) & = &  \prod_{g\in G}
 \rho_*(\FC)_{g,g^{-1}y}(\HC(y')^g) \\
& = &   \prod_{g\in G}
\Lambda[\gal_{y'}]\otimes_\Lambda \FC_{g,g^{-1}y}
\end{eqnarray*}
en posant $\gal_{y'}=\gal(\HC(y')/\HC(y))\simeq
\gal(\HC(y')^g/\HC(g^{-1}y))$, de sorte que
$$ 
 \mu_*\rho_*(\FC)_y(\HC(y'))^{H_y}= \prod_{g\in [G/H_y]}
\Lambda[\gal_{y'}]\otimes_\Lambda \FC_{g,g^{-1}y}.$$
Sous cette forme, l'exactitude du foncteur $\FC\mapsto
\Gamma(f^{-1}(V),(\nu\mu\rho)_*(\FC))^H $ 
appara{\^\i}t clairement, et par suite celle du foncteur
$\infty\circ(\nu\mu\rho)_*$.

Enfin le point iii) est {\'e}vident : il suffit de remplacer les
cat{\'e}gories $\FC( ?, \Lambda G_{disc})$ par $\FC(?,\Lambda
G_{disc}\times H_{disc})$  pour $?=X_{disc}$, $GX_{disc}$ ou
$\wt{GX}_{disc}$. L'action de $H$ sur ces deux derniers espaces est
donn{\'e}e par $h(g,x)=(g,hx)$.
\end{proof}

Renommons l'endofoncteur $B_0:=B$ et d{\'e}finissons $B_1:=B\circ(B/\id)$
puis par r{\'e}currence pour $i\in \NM$, $B_i:=B\circ(B_{i-1}/B_{i-2})$.

\begin{coro} \label{coroBerk}
  La suite $\id\To{}B_0\To{} \cdots B_i\To{} \cdots$ est une
  r{\'e}solution du foncteur identit{\'e} dans la cat{\'e}gorie des endofoncteurs
  de $\FC(X,\Lambda G)$ (resp. de $\FC(X,\Lambda G\times H_{disc}$)
  par des foncteurs exacts tels que pour tout $\FC\in\FC(X,\Lambda)$,
  le faisceau sous-jacent {\`a} $B_i(\FC)$ est {\'e}tale-c-mou. 
\end{coro}

\bibliography{artbiblio}

\end{document}